\documentclass[11pt]{amsart}

\usepackage{latexsym,amsfonts,amssymb,exscale,enumerate,comment}
\usepackage{amsmath,amsthm,amscd}

\usepackage{fullpage}
\usepackage{appendix}

\usepackage{url}
\usepackage[bookmarks=true,%
    colorlinks=true,%
    linkcolor=blue,%
    citecolor=blue,%
    filecolor=blue,%
    menucolor=blue,%
    urlcolor=blue,%
    breaklinks=true]{hyperref}

\input xy
\usepackage[all]{xy}
\SelectTips{cm}{}
\usepackage{tikz}
\usepackage{tikz-cd} %% added by Benjamin
\usetikzlibrary{shapes,snakes}
\usetikzlibrary{decorations.markings}
\usetikzlibrary{decorations.pathreplacing}
\tikzstyle directed=[postaction={decorate,decoration={markings,
    mark=at position #1 with {\arrow{>}}}}]

\newcommand{\hackcenter}[1]{
 \xy (0,0)*{#1}; \endxy}

\usepackage{graphicx}
\usepackage{color}
\newcommand{\onel}{\1_{\lambda}}

\newcommand{\onen}{\1_{n}}

\theoremstyle{plain}
\newtheorem{theorem}{Theorem}

\newtheorem{corollary}[theorem]{Corollary}
\newtheorem{proposition}[theorem]{Proposition}
\newtheorem{lemma}[theorem]{Lemma}

\newtheorem{notation}[theorem]{Notation}

\theoremstyle{definition}

\newtheorem{definition}[theorem]{Definition}

\theoremstyle{definition}
\newtheorem{remark}[theorem]{Remark}
\numberwithin{equation}{section}
\numberwithin{theorem}{section}
\newcommand{\maps}{\colon}
\newcommand{\und}[1]{\underline{#1}}

\newcommand{\refequal}[1]{\xy {\ar@{=}^{#1}
(-1,0)*{};(1,0)*{}};
\endxy}
\newcommand{\cat}[1]{\ensuremath{\mbox{\bfseries {\upshape {#1}}}}}

\newcommand{\To}{\Rightarrow}

\newcommand{\Hom}{{\rm Hom}}

\renewcommand{\to}{\rightarrow}

\def\Id{\mathrm{Id}}

\def\mf{\mathfrak}

\newcommand{\oUcat}{\mf{U}}

\numberwithin{equation}{section}
\def\ME#1{\textcolor[rgb]{0.40,0.00,0.90}{[ME: #1]}}%

\def\b{$\blacktriangleright$}

\let\epsilon=\varepsilon

\usepackage{bbm}

\def\N{{\mathbbm N}}

\def\Z{{\mathbbm Z}}
\def\cal#1{\mathcal{#1}}%
\def\1{\mathbbm{1}}%
\def\nn{\notag}

\def\la{\langle}
\def\ra{\rangle}

\renewcommand{\l}{\lambda}

\usepackage{bbm}

\def\cal#1{\mathcal{#1}}

\newcommand\nc{\newcommand}
\nc\rnc{\renewcommand}
\nc\Kar{\operatorname{Kar}}
\nc\End{\operatorname{End}}

\newcommand{\scs}{\scriptstyle}
\nc\Sym{\operatorname{Sym}}

\input{macros.tex}

\allowdisplaybreaks

\title { Super rewriting theory and nondegeneracy of odd categorified $\mf{sl}_2$}

\begin{document}
\setcounter{tocdepth}{1}

\author{Benjamin Dupont}
\email{bdupont@math.univ-lyon1.fr}
\address{Institut Camille Jordan, University of Lyon, France}

\author{Mark Ebert}
\email{markeber@usc.edu}
\address{Department of Mathematics\\ University of Southern California \\ Los Angeles, CA}

\author{Aaron D. Lauda}
\email{lauda@usc.edu}
\address{Department of Mathematics\\ University of Southern California \\ Los Angeles, CA}
\thanks{Research was sponsored by the Army Research Office and was
accomplished under Grant Number W911NF-20-1-0075. The views and conclusions contained in this
document are those of the authors and should not be interpreted as representing the official policies, either
expressed or implied, of the Army Research Office or the U.S. Government. The U.S. Government is authorized to reproduce and distribute reprints for Government purposes notwithstanding any copyright
notation herein.}

\maketitle

\begin{abstract}
We develop the rewriting theory for monoidal supercategories and 2-supercategories.  This extends the theory of higher-dimensional rewriting established for (linear) 2-categories to the super setting, providing a suite of tools for constructing bases and normal forms for 2-supercategories given by generators and relations.  We then employ this newly developed theory to prove the non-degeneracy conjecture for the odd categorification of quantum sl(2) from \cite{Lau-odd,BE2}.  As a corollary, this gives a classification of dg-structures on the odd 2-category conjectured in \cite{Odd-diff}.
\end{abstract}
\setcounter{tocdepth}{2}
%\tableofcontents

% #################################
\section{Introduction}
% #################################

Higher representation theory studies the higher categorical structure present when an associative algebra $A$ acts on an additive/abelian category $\cal{V}$, with algebra generators acting by additive or exact functors and algebra relations lifting to explicit natural isomorphisms of functors.  In its most refined form, this involves a categorification of an algebra $A$ itself, lifting $A$ to a monoidal category $\cal{A}$.  The algebra $A$ is categorified in the sense that there is an isomorphism from the (additive or abelian) Grothendieck group $K(\cal{A})$ to $A$.  The monoidal structure equips $K(\cal{A})$ with the structure of an algebra, where the $[X\otimes Y] = [X]\cdot[Y]$ and the class $[\1]$ of the unit in the  monoidal category becomes the unit element for algebra.

If the algebra $A$ is equipped with a system of mutually orthogonal idempotents, the most natural setting for categorification is to lift $A$ to an additive linear 2-category.  Since any monoidal category can be regarded as a 2-category with one object, the 2-categorical setting is often the most natural.  In particular, the diagrammatic calculus of 2-categorical string diagrams often appear in categorification, where the 2-categories $\cal{A}$ are defined diagrammatically via generating 2-morphisms modulo certain diagrammatic relations.  Then the categorification isomorphisms $K(\cal{A}) \cong A$ often requires significant effort to demonstrate that the diagrammatic presentation does not   collapse.  In particular, finding a basis for the spaces of 2-morphisms in $\cal{A}$ becomes a fundamental problem.  This can be viewed as the higher representation theoretic analog of studying PBW bases and related bases for enveloping algebras.  In the same way that those more traditional bases are a basic tool in the study of these algebras, the analogous bases for the spaces of 2-morphisms are equally relevant in higher representation theory.

Higher-dimensional rewriting theory applies the tools of rewriting theory in higher categorical settings.  It provides a set of tools for determining when a presentation of a 2-category will be coherent and allows for a determination of a normal form for a given 2-morphism within a given rewriting class,  constructively providing bases from a specific presentation of a 2-category.    The techniques of higher-dimensional rewriting have been effectively applied in a number of important examples in higher-representation theory~\cite{AL16,AL-cat,DUP19,DUP19bis} including cases where a determination of these bases have eluded experts for some time~\cite{DUP19bis}.
%This has the advantage  constructively providing bases from a specific presentation of a 2-category.

More recently, the field of higher representation theory has taken on the categorification of superalgebras $A$.  Superalgebras no longer lift to monoidal categories or 2-categories. Rather, they lift to so-called monoidal supercategories or 2-supercategories where the familiar interchange law is replaced by a super interchange law that depends on an additional $\Z_2$-grading on 2-morphisms~\cite{BE1,BE2}.
Monoidal supercategories and 2-supercategories are becoming increasingly common place in modern representation theory with examples ranging from categorification (Heisenberg categories~\cite{CautisLicata,BSW-found}, super 2-Kac-Moody algebras~\cite{EKL,Lau-odd,BE2, KKT, KKO, KKO2}, affine oriented Brauer-Clifford supercategory \cite{BCK},  Frobenius nilHecke \cite{Frob-nil}),  descriptions of the representation category of Lie superalgebras of Type $Q$~\cite{brown2020quantum}, Deligne categories for periplectic superalgebras~\cite{EAS}, and super analogs of modular/fusion tensor categories~\cite{ALK, Lacabanne,usher2016fermionic}.

Here we extend the theory of higher-dimensional rewriting to the super setting, allowing for these techniques to be applied to  monoidal supercategories and 2-supercategories.  This allows for a constructive approach to constructing bases in 2-supercategories and provides a suite of techniques for identifying Grothendieck groups needed for categorification.
As an application, we prove the non-degeneracy conjecture for the odd categorification of quantum $\mf{sl}_2$.

The \emph{odd} categorification of quantum $\mf{sl}_2$ arose as an attempt to provide a higher representation theoretic explanation for a phenomena discovered in link homology theory.   Ozsvath, Rassmusen, and Szabo showed that Khovanov's categorification of the Jones polynomial was not unique~\cite{ORS}. They defined what they called odd Khovanov homology that was similar in many ways to ordinary Khovanov homology (the theories agree when coefficients are reduced modulo two), but rather than being based on 2D TQFT, this theory was based on a strange type of 2D TQFT where signs appear when heights of handles are interchanged~\cite{Putyra}.  These theories are inequivalent in the sense that each can distinguish knots the other cannot~\cite{Shum}.  Since Khovanov homology has a higher representation theoretic interpretation coming from the categorification of quantum $\mf{sl}_2$~\cite{LQR,Web}, Ellis, Khovanov, and the third author initiated a program~\cite{EKL} to define odd analogs of quantum $\mf{sl}_2$ and related structures.  The result was the discovery of odd, noncommutative, analogs of many of the structures appearing in connection with $\mf{sl}_2$ categorification including odd analogs of the Hopf algebra of symmetric functions~\cite{EK,EKL}, cohomologies of Grassmannians~\cite{EKL} and Springer varieties~\cite{LauR}.  Subsequent work has shown these odd categorifications extend to arc algebras and constructions of odd Khovanov homology for tangles~\cite{NV,NK20,OddArc}.

%A closely related spin Hecke algebra associated to the affine Hecke-Clifford superalgebra appeared in earlier work of Wang~\cite{Wang} and many of the essential features of the odd nilHecke algebra including skew-polynomials appears much earlier in this and related works on spin symmetric groups~\cite{KW1,KW2,KW4}.

These investigations into odd categorification turned out to be closely connected with parallel investigations into super Kac-Moody algebra categorifications~\cite{KKT, KKO, KKO2}, with the odd categorification of $\mf{sl}_2$ lifting the rank one super Kac-Moody algebra.  These odd categorifications are also closely connected with the theory of covering Kac-Moody algebras~\cite{HillWang,ClarkWang, CHW,CHW2}.  Covering algebras $U_{q,\pi}(\mf{g})$ generalize quantum enveloping algebras, depending on an additional parameter $\pi$ with $\pi^2 = 1$.  When $\pi=1$, it reduces to the usual quantum enveloping algebra $U_{q}(\mf{g})$, while the $\pi = -1$ specialization recovers the quantum group of a super Kac-Moody algebra.  Covering algebras, and the novel introduction of the parameter $\pi$, allow for the first construction of canonical bases for Lie superalgebras~\cite{CHW2,ClarkWang}.

In the rank one case, the $\pi=1$ specialization is $U_{q}(\mf{sl}_2)$, while for $\pi=-1$ it gives the quantum group $U_{q}(\mf{osp}(1|2))$ associated with the super algebra $\mf{osp}(1|2)$.  Following a categorification of the positive parts of these algebra in \cite{HillWang}, Ellis and the third author categorified the full rank one covering algebra proving a conjecture from \cite{ClarkWang}.  In doing so, a 2-supercategory $\mf{U}:=\mf{U}(\mf{sl}_2)$ was defined~\cite{Lau-odd} for the rank one covering algebra whose Grothendieck group recovers $U_{q,\pi}(\mf{sl}_2)$.  This categorification was later greatly simplified in \cite{BE2}, where the 2-supercategory formalism was better developed, building off of the work~\cite{BE1}.  This covering formalism and the connection with $\mf{osp}(1|2)$ also informs the realization of odd Khovanov homology in theoretical physics~\cite{Witten-odd}.

Despite being able to establish the categorification isomorphism for $\mf{U}(\mf{sl}_2)$, a basis for the space of 2-morphisms was not achieved in \cite{Lau-odd}.  A spanning set was given in \cite{Lau-odd} and conjectured to form a basis --  the \emph{ non-degeneracy conjecture for odd categorified $\mf{sl}_2$}.  The need for a basis result was highlighted in \cite{Odd-diff} where dg-structures were defined on $\mf{U}$ extending differentials on the positive part from \cite{EllisQi}.  These differentials make the dg-Grothendieck group of its compact derived category isomorphic to the small quantum group $\dot{u}_{\sqrt{-1}}(\mf{sl}_2)$ that plays a role in quantum approaches to the Alexander polynomial.  Such dg-structures were conjecturally classified on $\mf{U}$ assuming the non-degeneracy conjecture~\cite[Proposition 7.1]{Odd-diff}.  As a corollary of the basis results achieved here, we prove this conjectured classification is complete.

\bigskip

This paper is organized as follows.  In Section~\ref{sec:super-rewrite} we adapt the theory of rewriting in linear 2-categories to the context of super 2-categories.   In Section~\ref{sec:SIso} we give a convergent presentation of the 2-supercategory we call odd isotopies.  This is analogous to the polygraph of isotopies from   \cite{GM09,DUP19},  but adapted to the context of 2-super Kac Moody algebras.  Section~\ref{sec:ONil} presents the 2-supercategory associated to the odd nilHecke algebra;  the resulting normal form is shown to recover the basis of the odd nilHecke algebra from \cite{EKL}.  Section~\ref{sec:odd2Cat} gives a presentation of the odd 2-category $\mf{U}(\mf{sl}_2)$.  Finally, in Section~\ref{sec:basis} we show that the (3,2)-superpolygraph presenting $\mf{U}(\mf{sl}_2)$ is quasi-terminating and confluent modulo.  The resulting quasi-normal form proves the non-degeneracy conjecture for  $\mf{U}(\mf{sl}_2)$.

%Rewriting theory
%\begin{itemize}
%  \item finite \cite{GM-finite}
%  \item Alleaume~\cite{AL16},
%\end{itemize}
%Polygraphs~\cite{Bur93} or computad~\cite{Street}

\subsection*{Acknowledgements}     A.D.L. and M.E were partially supported by NSF grant DMS-1902092 and Army Research Office W911NF2010075. B.D. would like to thank his PhD advisors Philippe Malbos, St\'ephane Gaussent and Alistair Savage for their help and support.

%\BD{Notations: $x,y$ for 0-cells $p,q$ for 1-cells, $u,v$ for 2-cells, $f,g$ for 3-cells and $\alpha, \beta$ for generating $3$-cells (that is elements of $P_3$)}

% #################################
\section{Super rewriting theory} \label{sec:super-rewrite}
% #################################

% ---------------------------------
\subsection{2-supercategories}
%----------------------------------
Here we review Brundan and Ellis~\cite{BE1,BE2} notion of a $2$-supercategory.

% - - - - - - - - - - - - - - - -
\subsubsection{Super vector spaces} \label{sec:superspaces}
% - - - - - - - - - - - - - - - -

Let $\Bbbk$ be a field with characteristic not equal to 2.   A \emph{superspace} is a $\Z_2$-graded vector space
$
V = V_{ \bar{0}} \oplus V_{ \bar{1}}.
$
For a homogeneous element $v \in V$, write $|v|$ for the parity of~$v$.

Let \cat{SVect} denote the category of superspaces and all linear maps. Note that  $\Hom_{\cat{SVect}}(V,W)$ has the structure of a superspace since  a linear map $f\maps V \to W$ between superspaces decomposes uniquely into an even and odd map.
The usual tensor product of $\Bbbk$-vector spaces is again a superspace with
$
(V\otimes W)_{ \bar{0}} = V_{ \bar{0}} \otimes W_{ \bar{0}}\oplus V_{ \bar{1}} \otimes W_{ \bar{1}}
$
and
$
(V\otimes W)_{ \bar{1}} = V_{ \bar{0}} \otimes W_{ \bar{1}}\oplus V_{ \bar{1}} \otimes W_{ \bar{0}}.
$
Likewise, the tensor product $f\otimes g$ of two linear maps between superspaces is defined by
\begin{equation}
  (f\otimes g) (v \otimes w) := (-1)^{|g||v|} f(v) \otimes g(w).
\end{equation}
Note that this tensor product does \emph{not} define a tensor product on \cat{SVect}, as the usual interchange law between tensor product and composition has a sign in the presence of odd maps
\begin{equation}
  (f\otimes g) \circ (h \otimes k) = (-1)^{|g||h|} (f \circ h) \otimes (g\circ k).
\end{equation}
This failure of the interchange law depending on parity is the primary structure differentiating monoidal supercategories from their non-super analogs.

If we set $\underline{\cat{SVect}}$ to be the subcategory consisting of only even maps, then the tensor product equips $\underline{\cat{SVect}}$ with a monoidal structure.  The map $u \otimes v \mapsto (-1)^{|u||v|} v \otimes u$ makes  $\underline{\cat{SVect}}$ into a symmetric monoidal category.

% - - - - - - - - - - - - - - - -
\subsubsection{Supercategories}
% - - - - - - - - - - - - - - - -
Supercategories, superfunctors, and supernatural transformations are defined~\cite{BE1} via the theory of enriched categories by enriching over the symmetric monoidal category $\underline{\cat{SVect}}$.   See \cite{kel1} for a review of the enriched category theory.  Unpacking this definition we have the following.

\begin{definition} [Supercategories] \hfill \newline
A supercategory $\cal{C}$ is a category enriched in the monoidal category $\underline{\cat{SVect}}$.
This consists of the data of a set $C_0$ called the objects, or $0$-cells, of $\cal{C}$ and
\begin{itemize}
 \item For each $x,y \in C_0$ a superspace of 1-cells $\cal{C}(x,y)$,
 \item For each $x\in  C_0$ an identity assigning map $i_x \maps \cal{I} \to \cal{C}(x,x)$ where $\cal{I}$ is the superspace $\Bbbk$ concentrated in degree zero.
 \item For each $x,y,z \in C_0$,   composition is given by a even linear map
 $$\star_0^{xyz} \maps \cal{C}(x,y) \otimes \cal{C}(y,z) \to \cal{C}(x,z)$$

\end{itemize}
such that composition is associative and unital with respect to identities.
\end{definition}

Superfunctors are functors between supercategories that give even linear maps on hom spaces.  For more details see \cite[Definition 1.1]{BE1}.

% - - - - - - - - - - - - - - - -
\subsubsection{2-supercategories} \label{sec:2supercat}
% - - - - - - - - - - - - - - - -

\begin{definition} [2-supercategories] \label{def:super2cat} \hfill \newline
A \emph{$2$-supercategory} $\cal{C}$ is a category enriched in the monoidal category of supercategories $\cat{SCat}$.
Namely, a 2-supercategory $\cal{C}$ is the data of a set $C_0$ called the objects of $\cal{C}$ and
\begin{itemize}
 \item For each $x,y \in C_0$ a supercategory $\cal{C}(x,y)$,
 \item For each $x\in  C_0$ an identity assigning superfunctor $i_x \maps \cal{I} \to \cal{C}(x,x)$ where $\cal{I}$ is the supercategory with
 \begin{itemize}
  \item one object $I$;
  \item $Hom(I,I)=\Bbbk$ where everything is even;
  \item Composition is the linear map $\circ: \Bbbk\otimes \Bbbk\to \Bbbk$ sending $c\otimes d \to cd$.
  \end{itemize}
 \item For each $x,y,z \in C_0$, a composition superfunctor $\star_0^{xyz} \maps \cal{C}(x,y) \otimes \cal{C}(y,z) \to \cal{C}(x,z)$
\end{itemize}
such that
\begin{itemize}
  \item $\star_0^{xzw} \circ (\star_0^{xyz} \otimes id_{\cal{C}(z,w)}) =\star_0^{xyw} \circ (id_{\cal{C}(x,y)} \otimes \star_0^{yzw})$ (associativity of composition)
  \item $\star_0^{xxy} \circ (i_x \times id_{\cal{C}(x,y)}) \circ is_l=id_{\cal{C}(x,y)}=\star_0^{xyy} \circ (id_{\cal{C}(x,y)} \times i_b) \circ is_r$ where $is_l$ and $is_r$ are the canonical isomorphisms $ C(a,b) \to I \otimes C(a,b)$ and $C(a,b) \to C(a,b) \otimes I$ (unitors).
\end{itemize}
The objects of the hom supercategories $\cal{C}(x,y)$ taken over all $x$ and $y$ define the set $C_1$ of 1-cells of $\cal{C}$ and the 1-cells in $\cal{C}(x,y)$ form the set $C_2$ of 2-cells in $\cal{C}$.  We use $\star_1$ to denote the composition operation in the supercategory $\cal{C}(x,y)$ and call this vertical composition of 2-cells.

For $p$ an object of the supercategory $\cal{C}(x,y)$ we define the \emph{0-source} of $p$ as $s_0(p)=x$ and \emph{0-target} of $p$ as $t_0(p) = y$.  The source and target maps in $\cal{C}(x,y)$ give \emph{1-source} and \emph{1-target} maps $s_1, t_1 \maps C_1 \to C_0$.
\end{definition}

The fact that composition is given by a monoidal superfunctor implies that the usual interchange axiom of a 2-category must be replaced by the superinterchange law.  That is,
given 2-cells $u \maps p \To q \maps x \to y$, $u' \maps p'\To q' \maps y \to z$,  $v \maps q \To r \maps x \to y$, $v' \maps q'\To r' \maps y \to z$,  then the \emph{superinterchange equation}
%\begin{equation} \label{eq:superinterchange}
%    (v\star_0 v') \star_1(u \star_0 u')
%    \;\; = \;\; (-1)^{|v'| |u|}
%    (v \star_1 u)\star_0 (v'\star_1 u')
%  %  =(\Id_{q'}v) \star_0 (u \Id_{p})
%\end{equation}
\begin{equation} \label{eq:superinterchange}
    (u \star_0 u')\star_1(v\star_0 v')
    \;\; = \;\; (-1)^{|u'| |v|}
    (u \star_1 v)\star_0 (u'\star_1 v')
  %  =(\Id_{q'}v) \star_0 (u \Id_{p})
\end{equation}
holds in a 2-supercategory $\cal{C}$.

\begin{definition}
A 2-supercategory $\cal{C}$ with one object is a monoidal supercategory. The tensor product operation is given by the $\star_0$-composition and composition of morphisms by $\star_1$.  The unit for the monoidal structure is given by the identity morphism of the unique object.  For more details see \cite[Definition 1.4]{BE1}.
\end{definition}

\begin{definition} \label{def:hom-basis}
A \emph{hom-basis} for a 2-supercategory $\cal{C}$ is a family of sets $(B_{p,q})_{p,q \in C_1}$ such that $B_{p,q}$ is a linear basis of the $\Bbbk$-superspace $C_2(p,q)$.
\end{definition}

The standard 2-categorical string diagrams can be adapted to the super setting.  The primary difference is that  the interchange law is replaced by the  superinterchange.  Since odd parity 2-morphisms now skew commute with each other, this means that for 2-supercategories one must be careful with the heights of 2-morphisms.  In particular, the superinterchange axiom \eqref{eq:superinterchange} implies that given 2-cells $u \maps p \To q \maps x \to y$ and $v \maps p'\To q' \maps y \to z$  then

\begin{equation}
% (V \star_1 Id_p') \star_0 ( \Id_{q} \star_1 u)
%    =
    (\Id_{p}\star_0 v) \star_1( u\star_0\Id_{q'})
    = (-1)^{|u||v|}(u \star_0 v)
    = (-1)^{|u||v|}
    ( u\star_0 \Id_{p'})\star_1 (  \Id_{q}\star_0 v)
\end{equation}
 \[
  \hackcenter{
\begin{tikzpicture}
\draw[very thick, black] (0,-1) to (0,1);
\node at (0,1.3) {$q'$};
\node at (0,-1.3) {$p'$};
\node[draw, thick, fill=black!20,rounded corners=4pt,inner sep=4pt] (p') at (0,-.4) {$v$};
\draw[very thick, black] (1,-1) to (1,1);
\node at (1,1.3) {$q$};
%\node at (1.3,.5) {$q$};
\node at (1,-1.3) {$p$};
\node[draw, thick, fill=black!20,rounded corners=4pt,inner sep=4pt] (p') at (1,.4) {$u$};
\node at (-.5,.75) {$z$};
\node at (.5,.75) {$y$};
\node at (1.5,.75) {$x$};
\end{tikzpicture} }
\quad = \quad
(-1)^{|u||v|}\;
  \hackcenter{
\begin{tikzpicture}
\draw[very thick, black] (0,-1) to (0,1);
\draw[very thick, black] (0,-1) to (0,1);
\node at (0,1.3) {$q'$};
\node at (0,-1.3) {$p'$};
\node[draw, thick, fill=black!20,rounded corners=4pt,inner sep=4pt] (p') at (0,0) {$v$};
\draw[very thick, black] (1,-1) to (1,1);
\node at (1,1.3) {$q$};
\node at (1,-1.3) {$p$};
\node[draw, thick, fill=black!20,rounded corners=4pt,inner sep=4pt] (p') at (1,0) {$u$};
%\node[draw, thick, fill=black!20,rounded corners=4pt,inner sep=4pt] (p') at (1,-.25) {$f'$};
\node at (-.5,.75) {$z$};
\node at (.5,.75) {$y$};
\node at (1.5,.75) {$x$};
\end{tikzpicture} }
\quad = \quad
(-1)^{|u||v|} \;
\hackcenter{
\begin{tikzpicture}
\draw[very thick, black] (0,-1) to (0,1);
\node at (0,1.3) {$q'$};
\node at (0,-1.3) {$p'$};
\node[draw, thick, fill=black!20,rounded corners=4pt,inner sep=4pt] (p') at (0,.4) {$v$};
\draw[very thick, black] (1,-1) to (1,1);
\node at (1,1.3) {$q$};
\node at (1,-1.3) {$p$};
\node[draw, thick, fill=black!20,rounded corners=4pt,inner sep=4pt] (p') at (1,-.4) {$u$};
\node at (-.5,.75) {$z$};
\node at (.5,.75) {$y$};
\node at (1.5,.75) {$u$};
\end{tikzpicture}  }
\]

\begin{remark}
Throughout this paper, we read our compositions cells as is common in higher category theory, just as the first author does in \cite{DUP19, DUP19bis}. This composition is read backwards from the more prevalent way of reading composition used by Brundan and Ellis \cite[Definition 2.1]{BE1}.
That is, $f\star_i g$ in this paper translates to $g\star_i f$ in \cite{BE1}.  So for example, we would have
\[
\hackcenter{\begin{tikzpicture}
  \draw[thick, black, ->] (0,0) to (0,1);
  \draw[thick, black, ->] (1,0) to (1,1);
  \node at (0,.5) {$\bullet$};
  \node at (1.5,.75) {$\scriptstyle{\lambda}$};
  \node at (0.5,.75) {$\scriptstyle{\lambda+2}$};
  \node at (-0.5,.75) {$\scriptstyle{\lambda+4}$};
\end{tikzpicture} }
\;\;=\;\;
\hackcenter{\begin{tikzpicture}
  \draw[thick, black, ->] (0,0) to (0,1);
  \node at (0.5,.75) {$\scriptstyle{\lambda}$};
  \node at (-0.5,.75) {$\scriptstyle{\lambda+2}$};
\end{tikzpicture}}
\star_0
\hackcenter{\begin{tikzpicture}
  \draw[thick, black, ->] (0,0) to (0,1);
  \node at (0,.5) {$\bullet$};
  \node at (0.5,.75) {$\scriptstyle{\lambda+2}$};
  \node at (-0.5,.75) {$\scriptstyle{\lambda+4}$};
\end{tikzpicture}}
\]
\end{remark}

% - - - - - - - - - - - - - - - -
\subsubsection{2-superpolygraphs and free  2-supercategories}
% - - - - - - - - - - - - - - - -

We now introduce the notion of superpolygraphs extending the notion of linear polygraphs developed in \cite{AL16}.  The theory of linear polygraphs is quite general, providing presentations of linear $(n,p)$-categories; these are defined using a combination of globular $n$-category objects and $p$-fold iterative enrichment (see \cite[Definition 2.2.1 \& 2.2.2]{AL16}) so that a linear $(n+1,p+1)$-category is a category enriched in $(n,p)$-categories, with the base case of linear $(n,0)$-category corresponding to an internal $n$-category in \cat{Vect}.  This means that a linear (1,1)-category is a linear category, a linear (1,0)-category is a category object in vector spaces, and a linear (2,2)-category is a linear 2-category.    Within the higher dimensional rewriting framework, a linear $(n,p)$-category is presented by a linear~$(n+1,p)$-polygraph.

Here we will need to extend several instances of the general linear $(n,p)$-category framework to the super setting.  This is because a $(2,2)$-supercategory is just a 2-supercategory as defined in Definition~\ref{def:super2cat} and these will be presented by $(3,2)$-superpolygraphs.  It is not hard to generalize Alleaume's theory of linear $(n,p)$-polygraphs to the super setting more generally, but as we do not have interesting examples of these structures in higher dimensions, we focus on unpacking the general inductive definitions in the cases of interest. To ease the exposition in this article, we make use of the definitions and notation of linear $(n,p)$-polygraphs from \cite[Section 3.2]{AL16}. We start with $(2,2)$-superpolygraphs which will be used to form the free 2-supercategory on a given set of generating cells.

Following \cite{AL16}, we will denote by $P_n^{\ast}$ the free strict $n$-category on a globular set
\[
\xymatrix{
P_n \ar@<1ex>[r]^-{s_{n-1}}
& \;\; \dots\;\;\ar@<1ex> @{<-}[l]^-{t_{n-1}}
P_{p} \ar@<1ex>[r]^{s_{p-1}}
& P_{p-1}  \ar@<1ex> @{<-}[l]^{t_{p-1}}  \ar @<1ex> [r]^{s_{p-2}}
& \;\; \dots\;\;\ar@<1ex> @{<-}[l]^{t_{p-2}}
\ar @<1ex> [r]^{s_0}
& P_0\ar@<1ex> @{<-}[l]^{t_0}
}
\]

\begin{definition}
  A \emph{$(2,2)$-superpolygraph} is a collection $P=(P_0, P_1, P_2)$ of sets equipped with set maps $s_k,t_k:P_{k+1}\to P_k^*$ for $k<2$,  such that:
   \begin{itemize}
     \item $(P_0, P_1)$ with $s_j,t_j$ for $j<1$ is a $1$-polygraph  as defined in \cite[Section 3]{AL16};
     \item $P_2$ is a \emph{super globular extension} of the free $1$-category $P_{1}^\ast$ on $(P_0,P_1)$, that is a $\Z_2$-graded set equipped with source and target maps
        $s_{1}, t_{1}:P_2 \to P_{1}^*$ satisfying globular relations $s_{0}\circ s_{1}=s_{0} \circ t_{1}$ and $t_{0}\circ s_{1}=t_{0} \circ t_{1}$.
   \end{itemize}
\end{definition}
We sometimes refer to $(2,2)$-superpolygraphs as $2$-superpolygraphs for convenience.

\begin{definition}
  A pasting diagram on $(2,2)$-superpolygraph $P=(P_0,P_1,P_2)$ is a formal composite of elements of $P_2':=P_2\cup \{ \1_x\maps x \To x \mid x \in P_1^{\ast} \}$ of the form:
  \begin{itemize}
    \item $u$ for $\alpha \in P_2'$,
    \item $u\star_1 v$ for $u,v$ pasting diagrams on $P$ with $t_1(u)=s_1(v)$,
    \item $u\star_0 v$ for $u,v$ pasting diagrams with $t_0s_1(u)=s_0s_1(v)$.
  \end{itemize}
 Such a composite inherits a $\Z_2$-grading determined by the parity of elements in $P_2$ as follows: $|\1_u|=0$, and $|u \star_k v|=|u|+|v|$ for $k=0,1$.
  We define a source $s_1(D)$ and target $t_1(D)$ of a composition $D$ iteratively by
  \begin{itemize}
    \item $s_1(u)$ and $t_1(u)$ are the normal 1-source and 1-target for $u\in P_2'$,
    \item $s_1(u \star_1 v)=s_1(u)$, \; $t_1(u \star_1 v)=t_1(v)$,
    \item $s_1(u \star_0 v)=s_1(u)\star_0 s_1(v)$, \; $t_1(u \star_0 v)=t_1(u)\star_0 t_1(v)$.
  \end{itemize}
  Then pasting diagrams on $P$ are such formal compositions quotiented by associativity of $\star_0$ and $\star_1$:
  \[
  u \star_k (v \star_k w)=(u\star_k v) \star_k w \quad \text{for}\; k=0,1
  \]
\end{definition}

We can now define the free $(2,2)$-supercategory on a $(2,2)$-superpolygraph by adapting the definition \cite[Def. 2.4.3]{GM-finite}.  A (2,2)-supercategory is the same thing as a 2-supercategory, so we will interchange freely between these two terminologies.

\begin{definition} \label{def:free-2-supercat}
Let~$P$ be a $2$-superpolygraph. The \emph{free $(2,2)$-supercategory over~$P$}, denoted by $P_2^s$, is defined as follows:
\begin{itemize}
\item the $0$-cells of $P_2^s$ are the $0$-cells of $P_0$,
\item for all $0$-cells~$x$ and~$y$ of~$P$, $P_2^s(x,y)$ is the supercategory whose
\begin{itemize}
\item $0$-cells are the $1$-cells $f \in P_1^\ast (x,y)$, where $P_1^\ast$ is the free $1$-category generated by the $1$-polygraph $(P_0,P_1)$,
\item  set of $1$-cells is the disjoint union of superspaces $P_2^s(p,q):=\text{Past} (p,q)$  where $\text{Past} (p,q)$ is the free superspace on the set of pasting diagrams with $1$-source $p$ and $1$-target $q$ for any $p,q \in P_1^\ast(x,y)$,
\end{itemize}
and quotiented by the congruence generated by the cellular extensions made of all the possible
\[
(u \star_0 v)\star_1(u' \star_0 v')=-1^{|v||u'|}(u \star_1 u')\star_0( v \star_1 v'), \qquad \1_{s_1(u)}\star_1 u=u=u \star_1 \1_{t_1(u)} \]
for all pasting diagrams $u,v,u',v'$ composable in this way.
The $0$-cells (resp. $1$-cells) of the hom supercategories $P_2^s(x,y)$ will be the $1$-cells (resp. $2$-cells) of $P_2^s$.
For any $0$-cells $p,q$ and $r$ in $P_2^s(x,y)$, there is an even linear map $\star_1 \maps P_2^s(p,q) \otimes P_2^s(q,r) \to P_2^s(p,r)$ given by gluing two $2$-cells $u : p \Rightarrow q$ and $v : q\Rightarrow r$ in $P_2^s$ along their common $1$-cell $q$.
For any $0$-cells $x,y,z\in P_0$, there is a composition map $\star_0:P_1^*(x,y)\otimes P_1^*(y,z) \to P_1^*(x,z)$ defined as the composition map on $P_1^*$.
Let $p,q$ and $r,s$ be any $0$-cells in the supercategories $P_2^s(x,y)$ and $P_2^s(y,z)$ respectively.
Then there is an even linear map $\star_0 : P_2^s(p,q) \otimes P_2^s(r,s)\to P_2^s(p\star_0 r,q\star_0 s)$ given by gluing two $2$-cells $u: p \Rightarrow q$ and $v:r\Rightarrow s$ in $P_2^s$ along their common $0$-cell $y$.
The $\star_0$ maps above give the data of a composition superfunctor $\star_0: P_2^s(x,y)\otimes P_2^s(y,z)\to P_2^s(x,z)$.
For any $1$-cells $u_1$, $\dots$, $u_m$ in $P_2^s(x,y)$ and $v_1$, $\dots$,$v_n$ in $P_2^s(y,z)$, these compositions satisfy
\begin{align*}
&\big( u_1 \star_1 \cdots \star_1 u_m  \big)
	\:\star_0\: \big( v_1  \star_1 \cdots \star_1 v_n  \big)
\:   \\
& \quad =
\: \left( u_1 \star_0 s(v_1) \right) \star_1 \cdots \star_1 \left( u_m \star_0 s(v_1) \right) {\star_1}
	 \left( t(u_m)  \star_0 v_1 \right) \star_1 \cdots \star_1 \left( t(u_m) \star_0 v_n \right)
\end{align*}
\end{itemize}
\end{definition}

\begin{remark}
If the $\Z_2$-grading of $P_2$ in a $(2,2)$-superpolygraph $P$ is all concentrated in even parity, then a $(2,2)$-superpolygraph is just a linear $(2,2)$-polygraph \cite[Definition 3.2.3]{AL16}, and the free $(2,2)$-supercategory $P_2^s$ generated by $P$ will be a linear $(2,2)$-category $P_2^{\ell}$ defined as in \cite[Definition 3.2.4]{AL16}.
\end{remark}

\begin{notation}
Let $P$ be a $2$-superpolygraph. Consider a subset $Q_2$  of the set $P_2$ of generating $2$-cells. For a given $2$-cell $u$ of $P_2^s$, denote by $||u||_{Q_2}$ the number of generating cells of $Q_2$ appearing in $u$. When $Q_2 = \{ w \}$ is a singleton, $||u||_{Q_2}$ counts the number of occurrences of the generating $2$-cell $w$ in $u$.
\end{notation}

The notion of monomial in a free 2-supercategory is defined by disregarding the $\Z_2$-grading and utilizing the definition of monomial for free linear 2-categories from \cite[Definition 4.1.4]{AL16}.

\begin{definition} \label{def:MonomialsOfSuperpolygraphs}
  Let $P=(P_0,P_1,P_2)$ be a $2$-superpolygraph and let $U(P)$ be the linear $(2,2)$-polygraph obtained by forgetting the parity of the elements $P_2$.
  Then a \emph{monomial} of the free 2-supercategory $P_2^s$ is a monomial of the free linear $(2,2)$-category $U(P)_2^\ell$ equipped with a parity determined by $P_2$.
\end{definition}

  The set of monomials of $U(P)_2^\ell$ is the set of 2-cells of the free 2-category $U(P)_2^*$, so equipping each element in the set of 2-cells of $U(P)_2^*$ with the parity determined by $P_2$ gives the monomials of $P_2^s$.

%\begin{lemma}
%The monomial decomposition in the free 2-supercategory $P_2^s$ generated by a (3,2)-superpolygraph is unique.
%\end{lemma}
%
%\begin{proof}
%Suppose $a$ and $b$ are pasting diagrams in $P$ such that $a=\pm b$ as 2-cells in $P_2^s$.
%If $a$ and $b$ are both monomials of $P_2^s$, then there exist two monomials $a',b'$ of $U(P)_2^l$ such that $a'=b'$.
%However, we know from \cite[Remark 4.1.3]{AL16} that every 2-cell in $U_2^\ell$ has a unique decomposition as a linear combination of monomials, which implies that the set of monomials for $U_2^\ell$ must be linearly independent.
%Thus, only one of $a', b'$ can be a monomial of $U^\ell$ and hence only one of $a,b$ is a monomial of $P_2^s$.
%This tells us that the monomials of $P_2^s$ are linearly independent, so if $u$ can be written as a linear combination of monomials, then this monomial decomposition must be unique.
%We also know by the definition of the free $2$-supercategory that every 2-cell in $P_2^s$ is a linear combination of pasting diagrams modulo superinterchange and identity, so every 2-cell has a monomial decomposition.
%Thus, every 2-cell in $P_2^s$ has a unique decomposition as a linear combination of monomials.
%\end{proof}

\begin{remark}
\label{rmk:set-of-monomials}
  For a $2$-superpolygraph $P$, let $A$ be a set of 2-cells of $U(P)_2^*$ containing one element from each exchange equivalence class of pasting diagrams of $U(P)$, where $2$-cells $u,v \in U(P)_2^*$ are in the same exchange equivalence class if $u=v$ via the exchange and identity relations.
  Then every 2-cell of $U(P)_2^*$ is equal to a unique element in $A$ by exchange and identity relations, so $A$ is the set of monomials of $U(P)_2^\ell$ and a linear combination of elements in $A$ is a monomial decomposition.
  We then obtain a set $B$ of $2$-cells of $P_2^s$ by assigning to each element of $A$ the parity determined by $P_2$.
  Then $B$ is the set of monomials of $P_2^s$ and so a linear combination of elements in $B$ is a monomial decomposition.
\end{remark}
It is known from  \cite[Definition 4.1.4]{AL16} that every 2-cell of $U(P)_2^\ell$ has a unique monomial decomposition.
This is true because there are no relations in $U_2^\ell$ other than the exchange and identity relations and no two elements of $A$ are related via these relations.
We now prove a lemma that gives this result for $2$-supercategories using similar principles.

\begin{lemma}
  Every $2$-cell in the free $2$-supercategory $P_2^s$ generated by a $2$-superpolygraph $P$ admits a unique monomial decomposition.
\end{lemma}

\begin{proof}
  Let $A$ and $B$ be the sets of monomials of $U(P)_2^\ell$ and $P_2^s$ described in Remark~\ref{rmk:set-of-monomials}.
  If $u=\pm v$ in $B$ by superinterchange, then there are corresponding elements $u'$ and $v'$ in $A$ that satisfy $u'=v'$ by exchange and identity relations.
  But we know that no two elements of $A$ are equal, so there are no two elements of $B$ that are scalar multiplies of each other by superinterchange and identity relations.
  Hence, $B$ is a linearly independent set of 2-cells because $P_2^s$ has no other relations other than the superinterchange and identity relations.
  Furthermore, every pasting diagram of $P$ is equal as a 2-cell by the superinterchange and identity relations to an element in $B$ up to a sign, so every $2$-cell admits a decomposition as a linear combination of elements of $B$ by construction.
  Hence, every 2-cell of $P_2^s$ admits a unique decomposition into a linear combination of elements of $B$.
\end{proof}

Given a 2-cell $u$ of the free 2-supercategory $P_2^s$  expressed as a  linear combination of monomials $u=\sum \lambda_i u_i$, we set
\[
\text{Supp}(u):=\{u_i \mid u_i \: \text{appears in the monomial decomposition of $u$}\}.
\]
%\begin{definition} Given a 2-superpolygraph $P$, every
%  2-cell $u$ of the free 2-supercategory $P_2^s$  has a unique decomposition into a linear combination of monomials $u=\sum \lambda_i u_i$.
%  Define
%\[
%\text{Supp}(u):=\{u_i \mid u_i \: \text{appears in the monomial decomposition of $u$}\}.
%\]
%\end{definition}

\begin{definition}
 Let $C$ be a $2$-(super)category. For a $k$-cell $f$ in $C$, with $k=1,2$,  
  define the \emph{boundary of $f$} as the ordered pair of $(k-1)$-cells defined by $\partial f := (s_{k-1}(f),t_{k-1}(f))$.
   A $k$-sphere of $C$ is a pair of $k$-cells $(f,g)$ such that $\partial f = \partial g$.  That is, $s_{k-1}(f)=s_{k-1}(g)$ and $t_{k-1}(f)=t_{k-1}(g)$.
  % $s_{k-1}(f)=s_{k-1}(g)$ and $t_{k-1}(f)=t_{k-1}(g)$ ($\partial f = \partial g$).
\end{definition}

Let us recall some key definitions needed to prove termination using the derivation method from \cite{GM09}: that of a context of a 2-category.

\begin{definition} \label{label:context-nonlin}
 A \emph{context} of a 2-category $C$ is a pair $(S,c)$ where $S$ is a 1-sphere of $C$ and $c$  is a 2-cell in the 2-category $C[S]$, defined as $C$ extended by a formal $2$-cell tiling the sphere $S$ as in \cite[Section 1.3]{GM09} such that this $2$-cell occurs exactly once in $c$. In other words, it is a $2$-cell $c$ that contains one `hole' with boundary the sphere $S$.

When $C$ is a $2$-category freely generated by a $2$-polygraph, a context of $C$ has the form \linebreak $c=m_1\star_1(m_2\star_0 S \star_0 m_3)\star_1 m_4$ where $m_i$ are monomials of $C$. For a 2-cell $u$ in $C_2$ such that $\partial u=S$, we denote by $c[u]$ the $2$-cell $m_1\star_1(m_2\star_0 u \star_0 m_3)\star_1 m_4$ in $C_2$.
\end{definition}

\begin{definition} \label{def:cat-contexts}
  Let $C$ be a 2-category. Then define the category of contexts $\cat{C}(C)$ as the category with
  \begin{itemize}
    \item Objects: 2-cells in $C$
    \item Morphisms: $\Hom(u,v)$ is the set of contexts $(\partial u, c)$ of $C$ such that $c[u]=v$.
    \item composition: If $x=(\partial u, c) \in \Hom(p,q)$ and $y=(\partial v, c') \in \Hom(q,r)$, then $x\circ y:= (\partial u, c' \circ c)\in \Hom(p,r)$ where $(c'\circ c)[w]:=c'[c[w]]$. %\ME{$(c'\circ c)[f]:=c'[c[f]]=c'[g]=h$}
    \item For any object $u$, there is an identity morphism $1_u:=(S=\partial u, c=S\in C[S])$. For $w\in C_2$ with $\partial w=\partial u$, $\partial u [w]=w$, so $(c\circ \partial u) [w]=c[w]$.
  \end{itemize}
%  Composition is associative and there is an identity object.
\end{definition}

In order to define rewriting steps of $(3,2)$-superpolygraph  we need to extend Definition~\ref{label:context-nonlin} to the case of contexts of 2-supercategories.
\begin{definition}
 A \emph{context} of a 2-supercategory $C$ is a pair $(S,c)$ where $S:=(p,q)$ is a 1-sphere of $C$ and $c$ is a 2-cell in the 2-supercategory $C[S]$, defined as the 2-supercategory $C$ extended with additional even $2$-cells $\lambda w$, for $\lambda \in \Bbbk$, tiling the sphere $S$ such that one of these $2$-cells appears exactly once in $c$.

In the case where $C$ is freely generated by a $2$-superpolygraph, that is $C = P_2^s$, a context of $P_2^s$ has the form $c=\lambda m_1\star_1(m_2\star_0 S \star_0 m_3)\star_1 m_4 +u$ for some scalar $\lambda$, monomials $m_i$ in $P_2^s$ and a $2$-cell $u$ in $P_2^s$. %\emph{such that the monomial $m_1\star_1(m_2\star_0 g \star_0 m_3)\star_1 m_4 \notin \text{Supp}(u)$.}
For a 2-cell $v$ of $P_2^s$ with $\partial v=(p,q)$, denote by $c[v]$ the $2$-cell $\lambda m_1\star_1(m_2\star_0 v \star_0 m_3)\star_1 m_4 +u$ in $P_2^s$.
\end{definition}

% ---------------------------------
\subsection{  (3,2)-superpolygraphs}
%----------------------------------
We now define (3,2)-superpolygraphs as a means of presenting $(2,2)$-supercategories.  This extends linear (3,2)-polygraphs from \cite[Definition 3.2.4]{AL16}.

\begin{definition}
  A $(3,2)$-superpolygraph is the data of $P=(P_0,P_1,P_2,P_3)$ where $(P_0,P_1,P_2)$ is a 2-superpolygraph and $P_3$ is a super globular extension of the free 2-supercategory $P_2^s$ on $(P_0,P_1,P_2)$, that is $P_3$ is a $\Z_2$-graded set equipped with even set maps $s_2, t_2:P_3\to P_2^s$ such that $s_{1}\circ s_2 = s_{1}\circ t_2$ and $t_{1}\circ s_2= t_{1} \circ t_2$ where $s_1,t_1$ are the 1-source and 1-target maps of $P_2^s$.
\end{definition}

The evenness of the set maps $s_2$ and $t_2$ in the definition of a $(3,2)$-superpolygraph implies they preserve the $\Z_2$ parity, so that the elements in $P_3$ with even parity have even sources and targets, while the elements in $P_3$ with odd parity have odd source and target.
% ---------------------------------
\subsection{$(3,2)$-supercategory}
%----------------------------------

\begin{definition}
A \emph{$(1,0)$-supercategory} is a category object in \underline{\cat{SVect}}.
A (2,1)-supercategory is a category enriched in $(1,0)$-supercategories.
A (3,2)-supercategory is a category enriched in $(2,1)$-supercategories.
\end{definition}
We will unpack these definitions in the cases of interest below.

%----------------------------------
\subsubsection{Free $(3,2)$-supercategory}
%----------------------------------

\begin{definition}
  A \emph{pasting diagram} on a $(3,2)$-superpolygraph $P=(P_0,P_1,P_2,P_3)$ is a formal composite of elements of the form
  \begin{itemize}
    \item $\alpha$ for $\alpha \in P_3':=P_3\cup \{\1_u:u \Rrightarrow u \mid u\in P_2^s\}$
    \item $f \star_2 g$ for pasting diagrams $f,g$ with $t_2(f)=s_2(g)$
    \item $f \star_1 g$ for pasting diagrams $f,g$ with $t_1t_2(f)=s_1s_2(g)$
    \item $f \star_0 g$ for pasting diagrams $f,g$ with $t_0t_1t_2(f)=s_0s_1s_2(g)$
  \end{itemize}
  quotiented by associativity relations for $\star_0$, $\star_1$ and $\star_2$:
   \[
  f \star_k (g\star_k h)=(f\star_k g) \star_k h \quad \text{for}\; 0 \leq k \leq 2.
  \]
  The source $s_2(f)$ and target $t_2(f)$ of a such a composition is defined by
  \begin{itemize}
    \item $s_2(f \star_2 g)=s_2(f)$,  \;  $t_2(f \star_2 g)=t_2(g)$
    \item $s_2(f \star_i g)=s_2(f)\star_i s_2(g)$ for $i\in \{ 0,1 \}$
    \item $t_2(f \star_i g)=t_2(f)\star_i t_2(g)$ $i\in \{ 0,1 \}$
  \end{itemize}
  The parity of such a composition is defined by $|f\star_k g|=|f|+|g|$  for $k=0,1$ and $|f\star_2 g|=|f|=|s_2(f)|$.
\end{definition}

\begin{definition}
Let $P = (P_0, P_1, P_2 ,P_3)$ be a $(3,2)$-superpolygraph. The \emph{free $(3,2)$-supercategory} generated by $P$, denoted by $P_3^s$, is defined as follows: its $0$-cells are the $0$-cells of $P_0$. For any $0$-cells $x$ and $y$ of $P$, we define the Hom $(2,1)$-supercategory $P_3^s(x,y)$ as follows:
\begin{itemize}
\item its $0$-cells are the $1$-cells $p \in P_1^\ast(x,y)$, where $P_1^\ast$ is the free $1$-category generated by the $1$-polygraph $(P_0,P_1)$,
\item for any $0$-cells $p$ and $q$ in $P_3^s(x,y)$, let us define the 2Hom $(1,0)$-supercategory $P_3^s(p,q)$ as follows:
\begin{itemize}
\item its set of $0$-cells  $P_2(p,q)$  is given by the superspace  $P_2^s(p,q)$  of $2$-cells of the free $(2,2)$-supercategory $P_2^s$ with $1$-source $p$ and $1$-target $q$,
\item its set of $1$-cells  $P_3(p,q)$  is the superspace given by the free superspace on $(3,2)$-pasting diagrams with $1$-source $p$ and $1$-target $q$ quotiented by relations
\[ (f \star_i g)\star_j( h \star_i k)=(-1)^{|g||h|}(f\star_j h)\star_i(g\star_j k), \qquad \1_{s_1(f)}\star_2 f=f=f \star_2 \1_{t_1(f)} \] for any $0 \leq i < j \leq 2$ and for all pasting diagrams $f,g,k,h$ composable in this way.
\end{itemize}
%The $0$-cells of the Hom $(2,1)$-supercategories $P_3^s(x,y)$ (respectively the $0$-cells, $1$-cells of the 2Hom $(1,0)$-categories $P_3^s(f,g)$) will be the $1$-cells (resp. $2$-cells, $3$-cells) of $P_3^s$.
The $\star_0$-composition for $1$-cells and $\star_0$, $\star_1$-composition for $2$-cells of $P_3^s$ are defined as in the free $(2,2)$-supercategory $P_2^s$.
For any $0$-cells $p$, $q$ in $P_3^s(x,y)$ and $r$,$s$ in $ P_3^s(y,z)$, there is an even linear map
$ \star_0 \maps P_3(p,q) \otimes P_3(r,s) \to P_3(p \star_0 r, q \star_0 s)$ given by gluing two $3$-cells along their common $0$-cell $y$. For any $0$-cells $p$, $q$, $r$ in $P_3^s(x,y)$, there is an even linear map
$ \star_1 \maps P_3(p,q) \otimes P_3(q,r) \to P_3(p, r)$ given by gluing two $3$-cells along their common $1$-cell $q$.
For any $0$-cells $p$, $q$ in $P_3^s(x,y)$, there is an even linear map
$ \star_2 \maps P_3(p,q) \times_{P_2(p, q)} P_3(p,q) \to P_3(p,q)$ given by gluing two $1$-cells %$3$-cells 
$f: u \Rrightarrow v$ and $g: v \Rrightarrow w$ of the 2Hom $(1,0)$-supercategory $P_3^s(p, q)$ along their common $0$-cell  %$2$-cell 
$v \in P_2(p, q)$.
For any $2$-cells $f_1$, $\dots$, $f_n$, $g_1$, $\dots$,$g_n$ in $P_3^s(x,y)$, these compositions satisfy
\begin{align*}
&\big( f_1 \star_2 \cdots \star_2 f_m  \big)
	\:\star_1\: \big( g_1  \star_2 \cdots \star_1 g_n  \big)
\: = \\
&\quad \: \left( f_1 \star_1 s(g_1) \right) \star_2 \cdots \star_2 \left( f_m \star_1 s(g_1) \right) {\star_2}
	 \left( t(f_m)  \star_1 g_1 \right) \star_2 \cdots \star_2 \left( t(f_m) \star_1 g_n \right) .
\end{align*}
\end{itemize}
\end{definition}

\begin{remark}
When the $\Z_2$-grading on the sets $P_2$ and $P_3$ are concentrated in degree zero, then a $(3,2)$-superpolygraph and (3,2)-supercategory reduce to a linear (3,2)-polygraphs and linear (3,2)-categories from \cite{AL16}.
\end{remark}

% ---------------------------------
\subsection{Presenting 2-supercategories by $(3,2)$-superpolygraphs}
%----------------------------------
\begin{definition}
  Let $P$ be a $(3,2)$-superpolygraph, and let $P_3^s$ be the free $(3,2)$-supercategory on $P$.
  Define an equivalence relation $\equiv$ on $P_2^s$ by
  \[
  u \equiv v \quad \text{if there is a 3-cell} \; f \in P_3^s \; \text{such that} \; s_2(f)=u \; \text{and} \; t_2(f)=v.
  \]
  We say that a 2-supercategory $C$ is presented by the $(3,2)$-superpolygraph $P$ if  $C$ is isomorphic to the quotient 2-supercategory $P_2^s / \equiv $.
\end{definition}

\begin{definition}
  A \emph{rewriting step} of a $(3,2)$-superpolygraph $P$ is a 3-cell $c[\alpha]\in P_3^s$ of the form
  \[
  c[\alpha]:c[s_2(\alpha)] \to c[t_2(\alpha)]
  \]
  where $\alpha \in P_3$ is a generating 3-cell, and $c=\lambda m_1\star_1(m_2\star_0 S \star_0 m_3)\star_1 m_4 +u$ is a context of $P_2^s$ such that the monomial $m_1\star_1(m_2\star_0 s_2(\alpha) \star_0 m_3)\star_1 m_4$ does not appear in the monomial decomposition of $u$.
  A {\em rewriting sequence} is a sequence of rewriting steps. A $3$-cell $f$ of $P_3^s$ is called \emph{positive} if it is an identity $3$-cell or a $\star_2$-composition $f= f_1 \star_2 \dots \star_2 f_n$ of rewriting steps of $P$. The \emph{length} of a positive $3$-cell $f$ in $P_3^s$, denoted by $\ell(f)$, is the number of rewriting steps of $P$ needed to write $f$ as a $\star_2$-composition of these rewriting steps. As a consequence, the terminologies rewriting path of $P$ (resp. rewriting step of $P$) and positive $3$-cell of $P_3^s$ (resp. positive $3$-cell of $P_3^s$ of length $1$) can both be used to represent the same notion.
\end{definition}

% ---------------------------------
\subsection{Termination and confluence}
%----------------------------------
A \emph{branching} (resp. \emph{local branching}) of a $(3,2)$-superpolygraph $P$ is a pair of rewriting sequences (resp. rewriting steps) of $P$ which have the same $2$-cell as $2$-source.
Such a branching (resp. local branching) is \emph{confluent} if it can be completed by rewriting sequences $f'$ and $g'$ of $P$ as follows:
\vskip-6pt
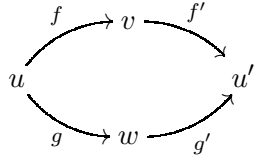
\[ \xymatrix @C=2.6em@R=0.8em{
& v
	\ar @/^1.5ex/ [dr] ^-{f'}
\\
u
	\ar @/^1.5ex/ [ur] ^-{f}
	\ar @/_1.5ex/ [dr] _-{g}
&& u'
\\
& w
	\ar @/_1.5ex/ [ur] _-{g'}
}	\]

\vskip-6pt
\noindent A~$(3,2)$-superpolygraph $P$ is said to be:
\begin{enumerate}[{\bf i)}]
\item \emph{left-monomial} if for any $\alpha$ in $P_3$, $s_2(\alpha)$ is a monomial of $P_2^s$.
\item \emph{terminating} if there is no infinite rewriting sequences in $P$.
\item \emph{quasi-terminating} if for each sequence $(u_n)_{n \in \N}$ of $2$-cells such that there is a rewriting step from $u_n$ to $u_{n+1}$ for each $n$ in $\N$, the sequence $(u_{n})_{n \in \N}$ contains an infinite number of occurrences of the same $2$-cell.
\item \emph{confluent} (resp. \emph{locally confluent}) if all the branchings (resp. local branchings) of $P$ are confluent.
\item \emph{convergent} if it is both terminating and confluent.
%\item \emph{exponentiation free} if for any $2$-cell $u$, there does not exist a $3$-cell $\alpha$ in $P_3^\ell$ such that $$
%\text{$u \overset{\alpha}{\fl} \lambda u + h$ with $\lambda \in \Bbbk \backslash \{ 0 \}$ and $h \ne 0$}. $$
\end{enumerate}

From now on, we will only consider left-monomial $(3,2)$-superpolygraphs.
Let us fix a (3,2)-superpolygraph $P$. A \emph{normal form} of $P$ is a $2$-cell $u$ that cannot be rewritten by any rewriting step of $P$. When $P$ is terminating, any $2$-cell admits at least one normal form, and exactly one when it is also confluent. A \emph{quasi-normal form} is a $2$-cell $u$ such that for any rewriting step from $u$ to another $2$-cell $v$, there exists a rewriting sequence from $v$ to $u$.

If $P$ is a terminating $(3,2)$-superpolygraph, Newman's lemma  
states that its confluence is equivalent to its local confluence.   
Following \cite[Section 4]{AL16}, branchings of a $(3,2)$-superpolygraph may be divided into four families: aspherical branchings, additive branchings, Peiffer branchings and overlapping branchings. A \emph{critical branching} of $P$ is an overlapping local branching that is minimal for the order $\sqsubseteq$ on monomials of $P_2^s$ defined by $f \sqsubseteq g$ if there exists a context $c$ of the free  $2$-category $U(P)_2^\ast$ generated by $P$ such that $g = c[f]$.

Following \cite{AL16}, we prove that a terminating $(3,2)$-superpolygraph is locally confluent if and only if its critical branchings are confluent. Indeed, the proofs of \cite[Lemma 4.2.12 $\&$ Theorem 4.2.13]{AL16} would remain the same: first proving that additive branchings are confluent and then prove that confluence of critical branchings implies confluence of all the overlapping branchings using implicit rewriting modulo superinterchange instead of the usual interchange.
Moreover,  with the definition of monomials from Definition~\ref{def:MonomialsOfSuperpolygraphs}, we obtain that if $P$ is a convergent $(3,2)$-polygraph presenting a $(2,2)$-supercategory $C$, then the set of monomials in normal form with respect to $P$ gives a hom-basis of $C$ in the sense of Definition~\ref{def:hom-basis}. Indeed, the same linear algebra argument as in the proof of \cite[Prop. 4.2.15]{AL16} would apply in this context since monomials of $P_2^s$ are defined in such a way that a $2$-cell of $P_2^s$ admits a unique monomial decomposition.

\subsubsection{Termination by Derivation} \label{subsec:termination-derivation}
Recall from \cite{GM09} a method to prove termination of a 3-polygraph using derivations of a 2-category.
The first author extended this method to the setting of linear 2-categories in \cite{DUP19, DUP19bis}, giving a method to prove termination of a $(3,2)$-linear polygraph using derivations of a 2-category.
Inspired by this extension to the linear setting, we describe a method to prove termination of a $(3,2)$-superpolygraph using derivations of a 2-category in this subsection.

The linear extension of proving termination by derivation from \cite{DUP19, DUP19bis} utilizes monomials of a linear $(2,2)$-category, which up to parity, are the same as  monomials of supercategories. This suggests the following definition.

\begin{definition} \label{def:UP}
  Let $P=(P_0,P_1,P_2,P_3)$ be a $(3,2)$-superpolygraph. Then define $U(P)$ as the linear $(3,2)$-polygraph with
  \begin{enumerate}
    \item $U(P)_i=P_i$ except that we forget the parity of elements
    \item The same source and target maps as in $P$ (forgetting parity of elements sends map $s_2,t_2: P_3 \to P_2^s$ to maps $s_2,t_2: U(P)_3 \to U(P)_2^l$)
  \end{enumerate}
\end{definition}

\begin{definition}
  Let $C$ be a 2-category.
  A $C$-module is a functor $M: \cat{C}(C) \to \cat{Ab}$ where $\cat{C}(C)$ is the category of context from Definition~\ref{def:cat-contexts} and $\cat{Ab}$ is the category of abelian groups.
\end{definition}

Let $\cat{Ord}$  denote the category of partially ordered sets and monotone maps. This is a monoidal category under the cartesian product.  As in \cite{GM09}, thinking of $\cat{Ord}$   as a 2-category with one object, we build examples of $C$-modules as follows:
\begin{definition}
    Let $C$ be a 2-category, $G$ be an internal abelian group in $\cat{Ord}$, and $X:C\to \cat{Ord}$ and $Y: C^{{\rm op}} \to \cat{Ord}$ be 2-functors. Then we can a define $C$-module $M:=M_{X,Y,G}$ as follows.
    \begin{itemize}
      \item Every 2-cell $u:p\Rightarrow q$ in $C$ is sent to the abelian group of morphisms $M(u)=\Hom_{\cat{Ord}}(X(p)\times Y(q),G)$

      \item If $p,q$ are 1-cells of $C$ and $c=p' \star_0 S \star_0 q'$ is a context from $u:p \Rightarrow q$ to $p' \star_0 u \star_0 q'$, then $M(c)$ sends a morphism $a:X(p)\times Y(q) \to G$ in $\cat{Ord}$ to the morphism $X(p')\times X(p)\times X(q')\times Y(p')\times Y(q)\times Y(q')\to G$ in $\cat{Ord}$ sending $(x',x,x'',y',y,y'')\to a(x,y)$.

      \item If $u:p'\to p$, $w:q\to q'$, are 2-cells and $c=u\star_1 x \star_1 w$ is a context from a $2$-cell $v: p \Rightarrow q$ to $u \star_1 v \star_1 w$, then $M(c)$ sends morphism $a:X(p)\times Y(q) \to G$ in $\cat{Ord}$ to the morphism to $a\circ (X\times Y)$, which is the map $X(p')\times Y(q') \to G$ sending $(x,y)\to a(X(g)(x),Y(h)(y))$.
    \end{itemize}
    When $C=U(P)_2^*$ is freely generated by a 2-polygraph $U(P)_{\leq 2}$, then such a $C$-module is uniquely determined by $X(p)$ and $Y(p)$ for $p\in P_1$ and the morphisms $X(u):X(p)\to X(q)$ and $Y(u): Y(q) \to Y(p)$ for every generating $2$-cell $u: p \Rightarrow q$ in $U(P)_2$.
   \end{definition}

We also recall the notion of a derivation of a 2-category:
\begin{definition}
  A derivation of a 2-category $C$ into a $C$-module $M$ is a map sending every 2-cell $u$ in $C$ to an element $d(u)\in M(f)$ such that
  \[
    d(u\star_i v)=u \star_i d(v) + d(u) \star_i v
  \]
  where $u \star_i d(v)=M(u \star_i x)(d(v))$ and $d(u)\star_i v=M(x \star_i v)(d(u))$.
\end{definition}

Then following \cite{DUP19}, we get the following result:

\begin{theorem} \label{T:TerminationDerivationTheorem}
  Let $P$ be a $(3,2)$-superpolygraph and $U(P)$ be the linear $(3,2)$-polygraph defined in Definition~\ref{def:UP}.
  Then if there exist
  \begin{enumerate}
    \item Two 2-functors $X:U(P)_2^*\to \cat{Ord}$ and $Y: (U(P)_2^*)^{{\rm op}}\to \cat{Ord}$ such that for every 1-cell $p$ in $P_1$, the sets $X(p)$ and $Y(p)$ are non-empty and for every generating 3-cell $\alpha$ in $P_3$, the inequalities $X(s_2(\alpha))\ge X(h)$ and $Y(s_2(\alpha))\ge Y(h)$ hold for every $h\in \text{Supp}(t_2(\alpha))$.
    \item An abelian group G in $\cat{Ord}$ whose addition is strictly monotone in both arguments and such that every decreasing sequence of non-negative elements of G is stationary.
        %We usually take $G=\Z$
    \item A derivation of $U(P)_2^*$ into the $U(P)_2^*$-module $M_{X,Y,G}$ such that for every 2-cell of $u\in U(P)_2^*$, we have $d(u)\ge 0$, and for every generating 3-cell $\alpha$ in $P_3$, $d(s_2(\alpha))> d(h)$ for every $h\in \text{Supp}(t_2(\alpha))$.
  \end{enumerate}
  Then the $(3,2)$-superpolygraph $P$ terminates.
\end{theorem}

\begin{proof}
Like in the linear setting, the proof works in a similar manner to the proof of termination by derivation given in \cite[Theorem 4.2.1]{GM09}.
%  Let $P$ be a $(3,2)$-superpolygraph such that there exists such $X,Y,d$ as described in the theorem above.
%  Now, assume for contradiction that $f_1\to f_2 \to \dots \to f_n \to \dots $ is an infinite rewriting sequence in $P$.
%  Let $h_k=d(f_k):=\text{max}{d(h) \mid h \in \text{Supp}(f_k)}$.
%  Then $(d(f_k))_{k\in \N}$ is a non-increasing sequence of natural numbers by properties of $X$, $Y$, and $d$ \cite{GM09}.
%  Furthermore, since $d(s_2(x))> d(t_2(x))$ for every $x\in P_3$ and there are no 2-cells whose support contains an infinite number of monomials, there is an infinite subsequence $(h_{k_i})_{i\in \N} \subset (h_k)_{k\in \N}$ of natural numbers that is strictly decreasing.
%  This is a contradiction, as an infinite sequence of natural numbers that is strictly decreasing cannot exist.
%  Therefore we cannot have an infinite rewriting sequence in $P$, which is to say that $P$ is terminating.
\end{proof}

\begin{remark}
 Usually we take internal abelian group $G=\Z$ and consider derivations with values into a $C$-module of the form $M_{X,Y,\Z}$.  We often consider $C$-module where $X$ or $Y$ are the trivial 2-functor and write $M_{X,\ast, \Z}$ or $M_{ \ast,Y, \Z}$.
\end{remark}

\subsubsection{Termination by context stable maps}\label{sec:context-stable}
Derivations were introduced in order to define termination orders by requiring some inequalities on sources and targets of generating $3$-cells;  the properties of derivations make this order stable by context of $2$-categories. Instead of a derivation, we can equivalently use maps $d \: : \: \mathcal{C}_2 \to \N$ that are stable under context, that is $d(a) \geq d(b)$ implies $d(c[a]) \geq d(c[b])$ for any context $c$ of $\mathcal{C}$.

\subsubsection{Derivation by steps}\label{sec:derivationbysteps}
The process of proving termination can be achieved in steps, proving termination for subsets of generating 3-cells at a time.

\begin{lemma} \label{lemma:steps}
 Let $P=(P_0,P_1,P_2,P_3)$ be a superpolygraph with $P_3=A \sqcup B$
 %.  Consider $X$ and $Y$ 2-functors as in Theorem~\ref{T:TerminationDerivationTheorem}
 and let $d \maps \mathcal{C}_2 \to \N$ be a context stable map satisfying the inequalities
 %$X,Y,d$ satisfy
%\[
%X(s_2(f))\geq X(t_2(f)) , \quad  $Y(s_2(f))\geq Y(t_2(f))\]
%for all $f \in P_3$
\[
d(s_2(f))>d(t_2(f)) \quad  \text{for $f\in A$, } \quad   d(s_2(g))\geq d(t_2(g)) \quad  \text{for $g\in B$}.
\]
  Then $P$ terminates if $P'=(P_0,P_1,P_2,B)$ terminates.
\end{lemma}

\begin{proof}
Suppose $P'$ terminates and \[v_1\overset{f_1}{\rightarrow} v_2 \overset{f_2}{\rightarrow} v_3 \to \dots \] is an infinite rewriting sequence in $P$.
Define   $d(u):=\text{max}\{d(u')\mid u'\in \text{Supp}(u)\}$.
Then since $P'$ terminates, there are an infinite number of $f_i$ that are in $A$.
Then consider the non-increasing infinite sequence $(d(v_n))_{n\in \N}$ of natural numbers.
The inequality is strict for $f_n\in A$ and $\text{Supp}(u)$ is a finite set, so $d(v_n)$ must decrease after a finite number of rewriting steps from $A$.
Hence, there is an infinite subsequence $(d(v_{n_k}))_{k\in \N}$ of natural numbers that is strictly decreasing giving a contradiction.
%However, no infinite sequence of natural numbers can be strictly decreasing, so we can't have an infinite rewriting sequence $v_1\to v_2 \to \dots$ in $P$.
\end{proof}

Lemma~\ref{lemma:steps} allows us to prove termination, progressively eliminating 3-cells.  When one of these steps is constructed from a context stable map arising from a derivation, we will need the conditions
\[
X(s_2(f))\geq X(t_2(f)) , \quad  Y(s_2(f))\geq Y(t_2(f)) \quad \text{for all $f \in P_3$}
\]
to hold at each step for the 2-functors used in defining the derivations.

One can view the process of proving derivations in steps
%in order to progressively eliminate $3$-cells used in this paper amounts to
as defining a termination lexicographic order.  If we denote the context stable map used at step $j$ by $d_j$, then a $k$ step procedure amounts to considering one large context stable map $d=(d_1, \: d_2, \dots, d_k)$ satisying
\[ \left( d_1(s_2(\alpha)), d_2(s_2(\alpha)), \dots, d_k(s_2(\alpha)) \right) >_{\text{lex}} \left( d_1(t_2(\alpha)), d_2(t_2(\alpha)), \dots, d_k (t_2(\alpha)) \right)
\]
for any generating $3$-cell $\alpha$ of the $(3,2)$-superpolygraph $P$, where $>_{\text{lex}}$ denotes the lexicographic order on $\N^k$. Each of these components being stable by context, we thus obtain that if there is an infinite rewriting sequence
\[ u_1 \fl u_2 \fl \dots \]
with respect to $P$, this yields an infinite strictly decreasing sequence
\[ \left( d_1 (u_1), d_2 (u_1), \dots, d_k (u_1) \right) >_{\text{lex}} \left( d_1 (u_2), d_2 (u_2), \dots, d_k (u_2) \right) >_{\text{lex}} \dots
\]
for the lexicographic order on $\N^k$, which is impossible since this order is well-founded.

% ---------------------------------
\subsection{  (3,2)-superpolygraphs modulo}
%----------------------------------

In this section we introduce the notion of rewriting modulo in 2-supercategories extending the work of the second author~\cite{DUP19,DUP19bis}.  This is tool for breaking termination and confluence arguments into incremental steps.  We utilize this to first prove that `odd isotopies' have a convergent presentation.  We then study presentations of the odd 2-category $\mf{U}$ modulo these odd isotopies.

A \emph{(3,2)-superpolygraph modulo} is a data $(R,E,S)$ made of two (3,2)-super polygraphs $R$ and $E$ such that $R_{\leq 1} = E_{\leq 1}$ and $E_2 \subseteq R_2$, and a cellular extension $S$ of the free 2-supercategory generated by $R_{\leq 2}$ satisfying $R \subseteq S \subseteq \ERE$, where the cellular extension $\ERE$ is made of elements of triples of the form $(e, f, e')$ for $3$-cells $e$, $e'$ in $E_3^s$ and a rewriting step $f$ of $R$ such that $t_2(e)= s_2(f)$ and $t_2(f) = s_2(e')$ as follows:
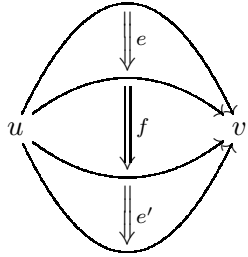
\[ \xymatrix@C=3em@R=1em{
u \ar @/^10ex/ [rr] ^-{} ^-{}="1"
 \ar @/^4ex/  [rr] ^{}   ^-{}="2"
 \ar @/_4ex/ [rr] _-{}  _-{}="3"
 \ar @/_10ex/ [rr] _-{} _-{}="4"
  & &  v
 \ar@2 "1"!<0pt,-2pt>;"2"!<0pt,2pt> ^-{e}
  \ar@2 "2"!<0pt,-2pt>;"3"!<0pt,2pt> ^-{f}
   \ar@2 "3"!<0pt,-2pt>;"4"!<0pt,2pt> ^-{e'}
  } \]
The rewriting sequences with respect to $\ERE$ thus correspond to application of rewriting sequences of $R$ by allowing sources and targets of $3$-cells to be transformed by a zig-zag sequence of rewriting steps of $E$. We refer to \cite{DM18} for a detailed definition of higher-dimensional polygraphs modulo. Given a $(3,2)$-superpolygraph modulo $(R,E,S)$, the data of $R_{\leq 2}$ and $S$ gives a $(3,2)$-superpolygraph, that we denote by $S$ in the sequel.

% - - - - - - - - - - - - - - - -
\subsubsection{Branchings and confluence modulo}
% - - - - - - - - - - - - - - - -
A \emph{branching modulo $E$} of a (3,2)-superpolygraph $(R,E,S)$ is a triple~$(f,e,g)$ where $f$ and $g$ are rewriting sequences of $S$, with $f$ non-identity, and $e$ is a $3$-cell in $E_3^\ell$ such that $s_2 (f) = s_2(e)$ and $s_2(g) = t_2(e)$. Such a branching modulo is \emph{confluent modulo $E$} if there exist rewriting sequences $f'$ and $g'$ of $S$, and a $3$-cell $e'$ in $E_3^s$ as in the following diagram:
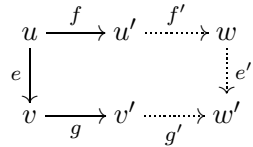
\[
\raisebox{0.55cm}{
\xymatrix @R=1.5em @C=2em {
u
  \ar[r] ^-{f}
  \ar[d] _-{e}
&
u'
  \ar@{.>}[r] ^-{f'}
&
w
  \ar@{.>}[d] ^-{e'}
\\
v
  \ar [r] _-{g}
&
v'
  \ar@{.>}[r] _-{g'}
&
w'
}}
\]
We then say that the triple $(f',e',g')$ is a confluence modulo $E$ of the branching $(f,e,g)$ modulo $E$. The $(3,2)$-superpolygraph $S$ is \emph{confluent modulo $E$} if all its branchings modulo $E$ are confluent modulo~$E$. A branching $(f,e,g)$ modulo $E$ is \emph{local} if $f$ is a rewriting step of $S$, $g$ is a positive $3$-cell of $S_3^s$ and $e$ is a $3$-cell of $E_3^s$ such that $\ell(g) + \ell(e) = 1$. Following \cite[Section 2.2.6]{DUP19}, local branchings are classified in the following families: local aspherical, local Peiffer, local additive, local Peiffer modulo, local additive modulo and overlappings modulo which are all the remaining local branchings modulo. A \emph{critical branching modulo $E$} is an overlapping branching modulo which is minimal for the order $\sqsubseteq$ defined by $(f,e,g) \subseteq (c[f],c[e],c[g])$ for any context $c$ of the 2-supercategory $R_2^s$.

% - - - - - - - - - - - - - - - -
\subsubsection{(Quasi)-Normal forms modulo}
% - - - - - - - - - - - - - - - -
Let us consider a $(3,2)$-superpolygraph modulo $(R,E,S)$ such that $S$ is confluent modulo $E$.  If $S$ is terminating (resp. quasi-terminating), each $2$-cell $u$ of $R_2^s$ admits at least one normal form (resp. quasi-normal form) with respect to $S$, and all these normal forms (resp. quasi-normal forms) are congruent modulo $E$ by confluence of $S$ modulo $E$. We fix such a normal form (resp. quasi-normal form), that we denote by $\widehat{u}$. By convergence of $E$, any $2$-cell $u$ of $R_2^s$ admits a unique normal form with respect to $E$, that we denote by $\widetilde{u}$. Note that when $S$ is confluent modulo $E$, the element $\tilda{u}$ does not depend on the chosen normal form $\widehat{u}$ for $u$ with respect to $S$, since two normal forms of $u$ being equivalent with respect to $E$, they have the same $E$-normal form. A \emph{normal form for $(R,E,S)$} (resp. \emph{quasi-normal form for $(R,E,S)$}) of a $2$-cell $u$ in $R_2^s$ is a $2$-cell $v$ such that $v$ appears in the monomial decomposition of $\widetilde{w}$, where $w$ is a monomial in the support of $\widehat{u}$. Such a set is obtained by reducing a $2$-cell $u$ in $R_2^s$ into its chosen normal form (resp. quasi-normal form) with respect to $S$, then taking all the monomials appearing in the $E$-normal form of each element in $\text{Supp}(\widehat{u})$.

% - - - - - - - - - - - - - - - -
\subsubsection{Decreasingness modulo}
% - - - - - - - - - - - - - - - -
The property of decreasingness modulo has been introduced in \cite{DUP19} following Van Oostrom's abstract decreasingness property for a rewriting system to give confluence criteria with respect to a well-founded labelling on the rewriting steps of a linear~$(3,2)$-polygraph modulo. When this polygraph is quasi-terminating, one may consider the quasi-normal form labelling, given by measuring the distance between a $2$-cell and a fixed quasi-normal form. It is proven in \cite{DUP19} that if a linear~$(3,2)$-polygraph is decreasing with respect to this labelling, which can be proved by proving the confluence of its critical branchings, it is confluent modulo. Note that this extends to the case of $(3,2)$-superpolygraphs since it is an abstract property. Another proof of the critical branching lemma modulo in the quasi-terminating setting may be found in \cite{CDM20}, based on induction on the distance to the quasi-normal form.

% ---------------------------------
\subsection{Linear bases from confluence modulo}
%----------------------------------
 Given a $(3,2)$-superpolygraph $P$, we define a \emph{splitting} of $P$ as a pair $(E,R)$ of $(3,2)$-superpolygraphs such that:
\begin{enumerate}[{\bf i)}]
\item $E$ is a sub-superpolygraph of $P$ such that $E_{\leq 1} = P_{\leq 1}$ and $E_2 \subseteq P_2$,
\item $R$ is a $(3,2)$-superpolygraph such that $R_{\leq 2} = P_{\leq 2}$ and $P_3 = R_3 \coprod E_3$.
\end{enumerate}
Such a splitting is called \emph{convergent} if we require that $E$ is convergent. The data of a splitting of a $(3,2)$-superpolygraph $P$ gives two distinct $(3,2)$-superpolygraphs $E$ and $R$ from which we can construct $(3,2)$-superpolygraphs modulo. Then, since the definition of monomials imply that every $2$-cell $u$ of $P_2^s$ admits a unique monomial decomposition, we prove in the same fashion as in the non-super setting the following statement:
\begin{theorem}
\label{T:BasisFromConfluenceModulo}
Let $P$ be a $(3,2)$-superpolygraph presenting a $(2,2)$-supercategory $C$, $(E,R)$ a convergent splitting of $P$ and $(R,E,S)$ a $(3,2)$-superpolygraph modulo such that
\begin{enumerate}[{\bf i)}]
\item $S$ is terminating (resp. quasi-terminating),
\item $S$ is confluent modulo $E$,
\end{enumerate}
then the set of all normal forms (resp. of all quasi-normal forms) for $(R,E,S)$ is a hom-basis of $\mathcal{C}$.
\end{theorem}

\begin{remark}
Note that we require $E$ to be convergent to ensure that any quasi-normal form with respect to the polygraph modulo $S$ admits a unique normal form with respect to $E$. However, even if we will still require $E$ to be terminating, the whole confluence assumption can be weakened. In particular, when $E$ is convergent with a set of $2$-cells that does not contain all the generating $2$-cells of $P$, the generating $2$-cells of $P_2 - E_2$ could create new indexed critical branchings, and thus obstructions to confluence. But confluence outside of these indexed critical branchings might be enough provided that these obstructions can be removed using the $3$-cells of $S$, so that any (quasi-)normal form with respect to $S$ still admit a unique normal form with respect to $E$. This is the case for the $(3,2)$-superpolygraph $\mathbf{Osl(2)}$ in which the $(3,2)$-superpolygraph will be confluent outside of crossing indexations as in \eqref{eq:newindexations}, but the polygraph modulo ${}_E R$ admits $3$-cells allowing the removal of self-intersections, as explained in Section \ref{sec:splittingofOsl2}.
\end{remark}

% #################################
\section{A convergent presentation of the super isotopy category}  \label{sec:SIso}
% #################################

% ---------------------------------
\subsection{Definition of supercategory of super isotopies}
%----------------------------------

Let $I$ be a possibly infinite index set equipped with a parity function
\begin{align}
  I &\to \Z/2, \qquad
 i \mapsto |i| \nn
\end{align}
We say that $i \in I$ is \emph{odd} if $|i|=\bar{1}$ and \emph{even} if $|i|=\bar{0}$.

 Let $(-d_{ij})_{i,j\in I}$ be a generalized Cartan matrix with $d_{ii}=-2$, $d_{ij}\geq 0$ for $i \neq j$, and $d_{ij}=0$ if and only if $d_{ji}=0$.  Under the additional assumption that $d_{ij}$ is even whenever $i$ is odd, Brundan and Ellis define a super 2-Kac-Moody algebra as a certain 2-supercategory $\mf{U}(\mf{g})$ associated to the Kac-Moody algebra $\mf{g}$ determined by the generalized Cartan matrix $(-d_{ij})_{i,j\in I}$.    In particular, associated to this Cartan matrix pick one can choose a complex vector space $\mf{h}$ and linearly independent subsets $\{\alpha_i \mid i \in I\}\subset \mf{h}^{\ast}$, $\{h_i \mid i \in I\} \subset \mf{h}$, such that the natural pairing $\mf{h}^* \times \mf{h} \to \Z$ is given by $\la h_i,\alpha_j \ra = -d_{ij}$ for all $i,j \in I$.
We denote the \emph{weight lattice} of $\mf{g}$ by $X=\{ \lambda \in \mf{h}^{\ast} \mid \la h_i, \lambda \ra \in \Z \text{ for all $i\in I$}\}$ and the \emph{root lattice} by $Y = \bigoplus_{i\in I} \Z\alpha_i$.  We sometimes write $\l_i := \la h_i,\lambda \ra$.

In what follows we consider a certain sub super 2-category of the super 2-Kac-Moody category $\mf{U}(\mf{g})$ defined by Brundan and Ellis~\cite[Definition 1.5]{BE2}.  This can be thought of as a super analog of the 2-category of pearls from \cite{GM09}.

\begin{definition} \label{def:SIso-cat}
Define the 2-supercategory of $\mf{g}$-valued isotopies $\mf{SIso}(\mf{g})$ to have
\begin{enumerate}[i)]
  \item objects consisting of  weights $\lambda \in X$ of the Kac-Moody algebra $\mf{g}$;

  \item the 1-morphisms  are generated by
\[
\1_{\lambda} \maps \lambda \to \lambda, \quad
\cal{E}_i\1_{\lambda} \maps \lambda \to \lambda + \alpha_i, \quad
\cal{F}_i\1_{\lambda} \maps \lambda \to \lambda - \alpha_i
\]

\item 2-morphisms generated by
\begin{alignat}{2}
&\hackcenter{\udott{i}} \maps   \cal{E}_i\1_{\lambda} \to \cal{E}_i\1_{\lambda}
     \qquad
&&\hackcenter{\ddott{i}}\maps \cal{F}_i\1_{\lambda} \to \cal{F}_i\1_{\lambda}
\qquad
%&& \bigotimes_i^{\l} \maps \1_{\lambda} \to \1_{\lambda}
    \\ \nn
& \quad \text{(parity $|i|$)} \qquad
&&\quad \text{(parity $|i|$)}
%&& \text{$i$ odd, (parity $\bar{1}$)}
\end{alignat}
\begin{alignat*}{4}
\capl{i} &\maps  \cal{F}_i\cal{E}_i\1_{\lambda} \to \1_{\lambda}
\qquad
&&\cupl{i}   \maps \cal{E}_i\cal{F}_i\1_{\lambda} \to \1_{\lambda}
\qquad
&&\capr{i}   \maps\1_{\lambda} \to \cal{F}_i\cal{E}_i\1_{\lambda}
\qquad
&&\cupr{i}   \maps\1_{\lambda} \to \cal{E}_i\cal{F}_i\1_{\lambda}
\\ \nn
&  \text{(parity $|i,\l|$)} \qquad
&&\quad\text{(parity $|i,\l|$)} \qquad
&&\quad\text{(parity $\bar{0}$)}\qquad
&&\quad\text{(parity $\bar{0}$)}
\end{alignat*}
 where
\[
|i,\lambda| := |i|(\la h_i, \l\ra +1) =  |i|(\l_i +1).
\]
\end{enumerate}
These 2-morphisms are required to satisfy the following axioms.

\begin{enumerate}[a)]
  \item Super zig-zag identities
  \begin{align*}
\label{E:IsotopyCells}
& \mathord{
\begin{tikzpicture}[baseline = 0]
  \draw[->,thick,black] (0.3,0) to (0.3,.4);
	\draw[-,thick,black] (0.3,0) to[out=-90, in=0] (0.1,-0.4);
	\draw[-,thick,black] (0.1,-0.4) to[out = 180, in = -90] (-0.1,0);
	\draw[-,thick,black] (-0.1,0) to[out=90, in=0] (-0.3,0.4);
	\draw[-,thick,black] (-0.3,0.4) to[out = 180, in =90] (-0.5,0);
  \draw[-,thick,black] (-0.5,0) to (-0.5,-.4);
   \node at (-0.5,-.6) {$\scriptstyle{i}$};
   \node at (0.5,0) {$\scriptstyle{\lambda}$};
\end{tikzpicture}
} =
\mathord{\begin{tikzpicture}[baseline=0]
  \draw[->,thick,black] (0,-0.4) to (0,.4);
   \node at (0,-.6) {$\scriptstyle{i}$};
   \node at (0.2,0) {$\scriptstyle{\lambda}$};
\end{tikzpicture}
}, \qquad
\mathord{
\begin{tikzpicture}[baseline = 0]
  \draw[->,thick,black] (0.3,0) to (0.3,-.4);
	\draw[-,thick,black] (0.3,0) to[out=90, in=0] (0.1,0.4);
	\draw[-,thick,black] (0.1,0.4) to[out = 180, in = 90] (-0.1,0);
	\draw[-,thick,black] (-0.1,0) to[out=-90, in=0] (-0.3,-0.4);
	\draw[-,thick,black] (-0.3,-0.4) to[out = 180, in =-90] (-0.5,0);
  \draw[-,thick,black] (-0.5,0) to (-0.5,.4);
   \node at (-0.5,.6) {$\scriptstyle{i}$};
   \node at (0.5,0) {$\scriptstyle{\lambda}$};
\end{tikzpicture}
}
=
\mathord{\begin{tikzpicture}[baseline=0]
  \draw[<-,thick,black] (0,-0.4) to (0,.4);
   \node at (0,.6) {$\scriptstyle{i}$};
   \node at (0.2,0) {$\scriptstyle{\lambda}$};
\end{tikzpicture}
}, \qquad
\mathord{
\begin{tikzpicture}[baseline = 0]
  \draw[-,thick,black] (0.3,0) to (0.3,-.4);
	\draw[-,thick,black] (0.3,0) to[out=90, in=0] (0.1,0.4);
	\draw[-,thick,black] (0.1,0.4) to[out = 180, in = 90] (-0.1,0);
	\draw[-,thick,black] (-0.1,0) to[out=-90, in=0] (-0.3,-0.4);
	\draw[-,thick,black] (-0.3,-0.4) to[out = 180, in =-90] (-0.5,0);
  \draw[->,thick,black] (-0.5,0) to (-0.5,.4);
   \node at (0.3,-.6) {$\scriptstyle{i}$};
   \node at (0.5,0) {$\scriptstyle{\lambda}$};
\end{tikzpicture}
}
=
(-1)^{|i , \lambda |} \mathord{\begin{tikzpicture}[baseline=0]
  \draw[->,thick,black] (0,-0.4) to (0,.4);
   \node at (0,-.6) {$\scriptstyle{i}$};
   \node at (0.2,0) {$\scriptstyle{\lambda}$};
\end{tikzpicture}
}, \qquad
\mathord{
\begin{tikzpicture}[baseline = 0]
  \draw[-,thick,black] (0.3,0) to (0.3,.4);
	\draw[-,thick,black] (0.3,0) to[out=-90, in=0] (0.1,-0.4);
	\draw[-,thick,black] (0.1,-0.4) to[out = 180, in = -90] (-0.1,0);
	\draw[-,thick,black] (-0.1,0) to[out=90, in=0] (-0.3,0.4);
	\draw[-,thick,black] (-0.3,0.4) to[out = 180, in =90] (-0.5,0);
  \draw[->,thick,black] (-0.5,0) to (-0.5,-.4);
   \node at (0.3,.6) {$\scriptstyle{i}$};
   \node at (0.5,0) {$\scriptstyle{\lambda}$};
\end{tikzpicture}
}
=
\mathord{\begin{tikzpicture}[baseline=0]
  \draw[<-,thick,black] (0,-0.4) to (0,.4);
   \node at (0,.6) {$\scriptstyle{i}$};
   \node at (0.2,0) {$\scriptstyle{\lambda}$};
\end{tikzpicture}
},
\end{align*}

\item  For $i\in I$ of parity $1$, define the odd bubble by
    \begin{equation} \label{eq:def-odd-bubble}
\raisebox{-4mm}{$\oddbubble{i}$}
\; :=\; \left\{
    \begin{array}{cl}
       (-1)^{\lfloor \frac{\l}{2}\rfloor }
       \hackcenter{\begin{tikzpicture}
\begin{scope} [ x = 10pt, y = 10pt, join = round, cap = round, thick, scale=1.8]
  \draw[<-,thick,black,scale=1.3] (0.1,1.2) to[out=180,in=90] (-.3,0.8);
  \draw[-,thick,black,scale=1.3] (0.5,0.8) to[out=90,in=0] (0.1,1.2);
 \draw[-,thick,black,scale=1.3] (-.3,0.8) to[out=-90,in=180] (0.1,0.4);
  \draw[-,thick,black,scale=1.3] (0.1,0.4) to[out=0,in=-90] (0.5,0.8);
 \node at (-0.5,0.8) {$\scriptstyle{i}$};
 %\node at (-.3,1.3) {$\bullet$};
 %\node at (-.6,1.4) {$\scriptstyle{\l}$};
 \node at (.5,1.3) {$\bullet$};
 \node at (.7,1.4) {$\scriptstyle{\l}$};
 \node at (1.2,0.6) {$\scriptstyle{\lambda}$};
 \end{scope}
\end{tikzpicture}  }
 & \text{if $\l \geq 0$} \\
        \raisebox{-4mm}{$ \begin{tikzpicture}
\begin{scope} [ x = 10pt, y = 10pt, join = round, cap = round, thick, scale=1.8]
  \draw[->,thick,black,scale=1.3] (0.1,1.2) to[out=180,in=90] (-.3,0.8);
  \draw[-,thick,black,scale=1.3] (0.5,0.8) to[out=90,in=0] (0.1,1.2);
 \draw[-,thick,black,scale=1.3] (-.3,0.8) to[out=-90,in=180] (0.1,0.4);
  \draw[-,thick,black,scale=1.3] (0.1,0.4) to[out=0,in=-90] (0.5,0.8);
 \node at (-0.5,0.8) {$\scriptstyle{i}$};
 \node at (0.55,1.3) {$\bullet$};
 \node at (1.0,1.5) {$\scriptstyle{-\l}$};
 \node at (1.2,0.4) {$\scriptstyle{\lambda}$};
 \end{scope}
\end{tikzpicture}$}  & \text{if $\l\leq 0$}
    \end{array}
\right.
\end{equation}
%\begin{equation} \label{eq:def-odd-bubble}
%\raisebox{-4mm}{$\oddbubble{i}$}
%\; :=\; \left\{
%    \begin{array}{cl}
%       %(-1)^{[\frac{h}{2}]}
%       \hackcenter{\begin{tikzpicture}
%\begin{scope} [ x = 10pt, y = 10pt, join = round, cap = round, thick, scale=1.8]
%  \draw[<-,thick,black,scale=1.3] (0.1,1.2) to[out=180,in=90] (-.3,0.8);
%  \draw[-,thick,black,scale=1.3] (0.5,0.8) to[out=90,in=0] (0.1,1.2);
% \draw[-,thick,black,scale=1.3] (-.3,0.8) to[out=-90,in=180] (0.1,0.4);
%  \draw[-,thick,black,scale=1.3] (0.1,0.4) to[out=0,in=-90] (0.5,0.8);
% \node at (-0.5,0.8) {$\scriptstyle{i}$};
% \node at (0.55,1.3) {$\bullet$};
% \node at (1.0,1.4) {$\scriptstyle{h}$};
% \node at (1.2,0.6) {$\scriptstyle{\lambda}$};
% \end{scope}
%\end{tikzpicture}  }
% & \text{if $h>0$} \\
%        \raisebox{-4mm}{$ \begin{tikzpicture}
%\begin{scope} [ x = 10pt, y = 10pt, join = round, cap = round, thick, scale=1.8]
%  \draw[->,thick,black,scale=1.3] (0.1,1.2) to[out=180,in=90] (-.3,0.8);
%  \draw[-,thick,black,scale=1.3] (0.5,0.8) to[out=90,in=0] (0.1,1.2);
% \draw[-,thick,black,scale=1.3] (-.3,0.8) to[out=-90,in=180] (0.1,0.4);
%  \draw[-,thick,black,scale=1.3] (0.1,0.4) to[out=0,in=-90] (0.5,0.8);
% \node at (-0.5,0.8) {$\scriptstyle{i}$};
% \node at (0.55,1.3) {$\bullet$};
% \node at (1.05,1.4) {$\scriptstyle{-h}$};
% \node at (1.2,0.4) {$\scriptstyle{\lambda}$};
% \end{scope}
%\end{tikzpicture}$}  & \text{if $h\leq 0$}
%    \end{array}
%\right.
%\end{equation}
Then the odd `cyclicity' relations
\begin{align*}
&\surdd{i} = \upd{i},
\qquad
\suld{i} = \left\{
    \begin{array}{ll}
        \dpd{i} & \text{if $i$ is even} \\
        2 \raisebox{-6mm}{$\downoddbubbleright{i}{i}$} - \dpd{i} & \text{if $i$ is odd}
    \end{array}
\right.
\\
&\sdrd{i} = \dpd{i}, \qquad
 \sdld{i}  \quad = \left\{
    \begin{array}{ll}
        \upd{i} & \text{if $i$ is even} \\
        2 \: \raisebox{-6mm}{$\upoddbubbleleft{i}{i}$} - (-1)^{\l_i+1} \upd{i} & \text{if $i$ is odd}
    \end{array}
\right.
\end{align*}
hold.
 %\begin{align*}
% &
%\cupldl{i}{} = \left\{
%    \begin{array}{ll}
%       \cupldr{i}{} & \text{if $i$ is even} \\
%        (-1)^h \cupldr{i}{} + 2 \raisebox{-4mm}{$\cupoddbubble{i}{i}$} & \text{if $i$ is odd}
%    \end{array}
%\right. , \quad \capldl{i}{} \overset{i_{\lambda}^4}{\Rrightarrow} \left\{
%    \begin{array}{ll}
%       \capldr{i}{} & \text{if $i$ is even} \\
%        (-1)^h \capldr{i}{} + 2 \: \: \raisebox{-4mm}{$\capoddbubblenest{i}{i}$} & \text{if $i$ is odd}
%    \end{array}
%\right.
%\end{align*}
\end{enumerate}
\end{definition}

% ---------------------------------
\subsection{The super (3,2)-polygraph SIso} \label{sec:OIso}
%----------------------------------

In this section we define a $(3,2)$-super polygraph presenting the super 2-category $\mf{SIso}(\mf{g})$ of $\mf{g}$-valued isotopies. Let $\mathbf{SIso}(\mathfrak{g})$ be the super~$(3,2)$-polygraph defined by:
\begin{enumerate}[{\bf i)}]
\item the elements of $\mathbf{SIso}(\mathfrak{g})_0$ are the weights $\lambda \in X$ of the Kac-Moody algebra;
\item the elements of $\mathbf{SIso}(\mathfrak{g})_1$ are given by
\[
 1_{\lambda '} \mathcal{E}_{\varepsilon_1 i_1} \dots \mathcal{E}_{\varepsilon_m i_m} 1_{\lambda}
\]
for any signed sequence of vertices with $\epsilon_\ell \in \{ \pm\}$ and $i_\ell\in I$ for all $1 \leq \ell \leq m$. Here we identify $\cal{E}_{+} := \cal{E}$ and $\cal{E}_{-} :=\cal{F}$.    Such a $1$-cell has for $0$-source $\lambda$ and $0$-target $\lambda ' = \l + \sum_{\ell} \epsilon_l \alpha_{i_{\ell}}$, and
\[
1_{\lambda '} \mathcal{E}_{\varepsilon_1 i_1} \dots \mathcal{E}_{\varepsilon_m i_m} 1_{\lambda} \star_0 1_{\lambda ''} \mathcal{E}_{\varepsilon'_1 j_1} \dots \mathcal{E}_{\varepsilon_l j_l} 1_{\lambda'}  = 1_{\lambda ''} \mathcal{E}_{\varepsilon'_1 j_1} \dots  \mathcal{E}_{\varepsilon'_l j_l} \mathcal{E}_{\varepsilon_1 i_1}\dots
\mathcal{E}_{\varepsilon_m i_m} 1_{\lambda} \]
\item the elements of $\mathbf{SIso}(\mathfrak{g})_2$ are the following generating $2$-cells: for any $i$ in $I$ and $\lambda '$ in $X$,
\begin{align*}
\udott{i} \qquad \ddott{i} \qquad
\capl{i} \qquad  \cupl{i} \qquad \capr{i} \qquad \cupr{i}
\end{align*}
with respective parity $|i|$, $|i|$, $|i,\l|$, $|i,\l|$, $0$, $0$.

\item $\mathbf{SIso}(\mathfrak{g})_3$ consists of the following $3$-cells:
\begin{align*}
\label{E:IsotopyCells}
& \mathord{
\begin{tikzpicture}[baseline = 0]
  \draw[->,thick,black] (0.3,0) to (0.3,.4);
	\draw[-,thick,black] (0.3,0) to[out=-90, in=0] (0.1,-0.4);
	\draw[-,thick,black] (0.1,-0.4) to[out = 180, in = -90] (-0.1,0);
	\draw[-,thick,black] (-0.1,0) to[out=90, in=0] (-0.3,0.4);
	\draw[-,thick,black] (-0.3,0.4) to[out = 180, in =90] (-0.5,0);
  \draw[-,thick,black] (-0.5,0) to (-0.5,-.4);
   \node at (-0.5,-.6) {$\scriptstyle{i}$};
   \node at (0.5,0) {$\scriptstyle{\lambda}$};
\end{tikzpicture}
} \overset{u_{\lambda,0}}{\Rrightarrow}
\mathord{\begin{tikzpicture}[baseline=0]
  \draw[->,thick,black] (0,-0.4) to (0,.4);
   \node at (0,-.6) {$\scriptstyle{i}$};
   \node at (0.2,0) {$\scriptstyle{\lambda}$};
\end{tikzpicture}
}, \qquad
\mathord{
\begin{tikzpicture}[baseline = 0]
  \draw[->,thick,black] (0.3,0) to (0.3,-.4);
	\draw[-,thick,black] (0.3,0) to[out=90, in=0] (0.1,0.4);
	\draw[-,thick,black] (0.1,0.4) to[out = 180, in = 90] (-0.1,0);
	\draw[-,thick,black] (-0.1,0) to[out=-90, in=0] (-0.3,-0.4);
	\draw[-,thick,black] (-0.3,-0.4) to[out = 180, in =-90] (-0.5,0);
  \draw[-,thick,black] (-0.5,0) to (-0.5,.4);
   \node at (-0.5,.6) {$\scriptstyle{i}$};
   \node at (0.5,0) {$\scriptstyle{\lambda}$};
\end{tikzpicture}
}
\overset{d_{\lambda,0}}{\Rrightarrow}
\mathord{\begin{tikzpicture}[baseline=0]
  \draw[<-,thick,black] (0,-0.4) to (0,.4);
   \node at (0,.6) {$\scriptstyle{i}$};
   \node at (0.2,0) {$\scriptstyle{\lambda}$};
\end{tikzpicture}
}, \qquad
\mathord{
\begin{tikzpicture}[baseline = 0]
  \draw[-,thick,black] (0.3,0) to (0.3,-.4);
	\draw[-,thick,black] (0.3,0) to[out=90, in=0] (0.1,0.4);
	\draw[-,thick,black] (0.1,0.4) to[out = 180, in = 90] (-0.1,0);
	\draw[-,thick,black] (-0.1,0) to[out=-90, in=0] (-0.3,-0.4);
	\draw[-,thick,black] (-0.3,-0.4) to[out = 180, in =-90] (-0.5,0);
  \draw[->,thick,black] (-0.5,0) to (-0.5,.4);
   \node at (0.3,-.6) {$\scriptstyle{i}$};
   \node at (0.5,0) {$\scriptstyle{\lambda}$};
\end{tikzpicture}
}
\overset{u'_{\lambda,0}}{\Rrightarrow}
(-1)^{|i , \lambda |} \mathord{\begin{tikzpicture}[baseline=0]
  \draw[->,thick,black] (0,-0.4) to (0,.4);
   \node at (0,-.6) {$\scriptstyle{i}$};
   \node at (0.2,0) {$\scriptstyle{\lambda}$};
\end{tikzpicture}
}, \qquad
\mathord{
\begin{tikzpicture}[baseline = 0]
  \draw[-,thick,black] (0.3,0) to (0.3,.4);
	\draw[-,thick,black] (0.3,0) to[out=-90, in=0] (0.1,-0.4);
	\draw[-,thick,black] (0.1,-0.4) to[out = 180, in = -90] (-0.1,0);
	\draw[-,thick,black] (-0.1,0) to[out=90, in=0] (-0.3,0.4);
	\draw[-,thick,black] (-0.3,0.4) to[out = 180, in =90] (-0.5,0);
  \draw[->,thick,black] (-0.5,0) to (-0.5,-.4);
   \node at (0.3,.6) {$\scriptstyle{i}$};
   \node at (0.5,0) {$\scriptstyle{\lambda}$};
\end{tikzpicture}
}
\overset{d'_{\lambda,0}}{\Rrightarrow}
\mathord{\begin{tikzpicture}[baseline=0]
  \draw[<-,thick,black] (0,-0.4) to (0,.4);
   \node at (0,.6) {$\scriptstyle{i}$};
   \node at (0.2,0) {$\scriptstyle{\lambda}$};
\end{tikzpicture}
},
\end{align*}
%
%
%
%
%Let $\l_i := \la h_i,\l \ra$ then
\begin{align*}
 &\cuprdl{i}{}  \overset{i_{\lambda}^1}{\Rrightarrow} \cuprdr{i}{}, \qquad \quad
\cupldl{i}{} \overset{i_{\lambda}^3}{\Rrightarrow} \left\{
    \begin{array}{ll}
       \cupldr{i}{} & \text{if $i$ is even} \\
        (-1)^{\lambda_i} \cupldr{i}{} + 2 \raisebox{-4mm}{$\cupoddbubble{i}{i}$} & \text{if $i$ is odd}
    \end{array}
\right.
\\
&\caprdl{i}{} \overset{i_{\lambda}^2}{\Rrightarrow} \caprdr{i}{}, \qquad \quad
 \capldl{i}{} \overset{i_{\lambda}^4}{\Rrightarrow} \left\{
    \begin{array}{ll}
       \capldr{i}{} & \text{if $i$ is even} \\
        (-1)^{\lambda_i} \capldr{i}{} + 2 \: \: \raisebox{-7mm}{$\capoddbubblenest{i}{i}$} & \text{if $i$ is odd}
    \end{array}
\right.
\end{align*}
For the definition of an odd bubble in weight spaces $\lambda_i=\la h_i, \l \ra = 0$ with $|i| = \overline{1}$, we also add $3$-cells
\begin{align*}
  \raisebox{-4mm}{\begin{tikzpicture}
\begin{scope} [ x = 10pt, y = 10pt, join = round, cap = round, thick, scale=1.8]
  \draw[<-,thick,black,scale=1.3] (0.1,1.2) to[out=180,in=90] (-.3,0.8);
  \draw[-,thick,black,scale=1.3] (0.5,0.8) to[out=90,in=0] (0.1,1.2);
 \draw[-,thick,black,scale=1.3] (-.3,0.8) to[out=-90,in=180] (0.1,0.4);
  \draw[-,thick,black,scale=1.3] (0.1,0.4) to[out=0,in=-90] (0.5,0.8);
 \node at (0.85,1.3) {$\scriptstyle{0}$};
 \node at (-0.3,0.3) {$\scriptstyle{i}$};
 \end{scope}
\end{tikzpicture}} \:
\overset{I_0}{\Rrightarrow} \:
  \hackcenter{\begin{tikzpicture}
\begin{scope} [ x = 10pt, y = 10pt, join = round, cap = round, thick, scale=1.8]
  \draw[->,thick,black,scale=1.3] (0.1,1.2) to[out=180,in=90] (-.3,0.8);
  \draw[-,thick,black,scale=1.3] (0.5,0.8) to[out=90,in=0] (0.1,1.2);
 \draw[-,thick,black,scale=1.3] (-.3,0.8) to[out=-90,in=180] (0.1,0.4);
  \draw[-,thick,black,scale=1.3] (0.1,0.4) to[out=0,in=-90] (0.5,0.8);
 \node at (0.85,1.3) {$\scriptstyle{0}$};
  \node at (-0.3,0.3) {$\scriptstyle{i}$};
 \end{scope}
\end{tikzpicture}}
\end{align*}
and for any $i \in I$ of parity $1$ and any endomorphism $2$-cell $k$ of the identity $\1_{\l}$  in normal form with respect to the set of $3$-cells above:

\begin{align*} \raisebox{-9mm}{$\sourcealphamk{i}{i}{m}{k}$} \overset{\alpha_{m,k}}{\Rrightarrow} \left\{
    \begin{array}{ll}
       (-1)^{m + |k|} \raisebox{-3mm}{$\targetalphamk{i}{i}{m+1}{k}$} & \text{if $m+\l_i+1 $ is even} \\
        0 & \text{if $m+\l_i+1$ is odd}
    \end{array}
\right.
\end{align*}
\begin{align*} \raisebox{-9mm}{$\sourcebetamk{i}{i}{m}{k}$} \overset{\beta_{m,k}}{\Rrightarrow} \left\{
    \begin{array}{ll}
      \raisebox{-3mm}{$\targetbetamk{i}{i}{m+1}{k}$}  & \text{if $m+\l_i+1$ is even} \\
        0 & \text{if $m+\l_i+1$ is odd}
    \end{array}
\right.
\end{align*}
where the odd bubble $\raisebox{-5mm}{$\oddbubble{i}$}$ is the 2-cell  defined as in \eqref{eq:def-odd-bubble}.
\end{enumerate}

\begin{lemma}
One can define 3-cells
\begin{align*}
&\surdd{i} \overset{u_{\lambda,1}}{\Rrightarrow} \upd{i}, \qquad
 \suld{i} \overset{d'_{\lambda,1}}{\Rrightarrow} \left\{
    \begin{array}{ll}
        \dpd{i} & \text{if $i$ is even} \\
        2 \raisebox{-6mm}{$\downoddbubbleright{i}{i}$} - \dpd{i} & \text{if $i$ is odd}
    \end{array}
\right. ,
\\
&\sdrd{i} \overset{d_{\lambda,1}}{\Rrightarrow} \dpd{i}, \qquad \sdld{i} \overset{u'_{\lambda,1}}{\Rrightarrow} \left\{
    \begin{array}{ll}
        \upd{i} & \text{if $i$ is even} \\
        2 \: \raisebox{-6mm}{$\upoddbubbleleft{i}{i}$} - (-1)^{\l_i+1} \upd{i} & \text{if $i$ is odd}
    \end{array}
\right.
\end{align*}
from the generating 3-cells of $\mathbf{SIso}(\mathfrak{g})$
\end{lemma}

\begin{remark}
Note that in every $3$-cell except for $\alpha_{m,k}$ and $\beta_{m,k}$, every strand of the source and target are labeled by the same $i\in I$.
In $\alpha_{m,k}$ and $\beta_{m,k}$, the strands of $k$ can be labeled with any $j\in I$.
However, we cannot rewrite $k$ using any rewriting step since it is in normal form by definition.
Knowing this, we can write $\mathbf{SIso}_{3}=\bigsqcup\limits_{i\in I}\mathbf{SIso}_{3}^{i}$, where $\mathbf{SIso}_{3}^{i}$ is the set of $3$-cells where all of the strands of the source and target are labeled by $i$ along with the $3$-cells $\alpha_{m,k}$ and $\beta_{m,k}$ where the strands of the odd bubble and the bubble surrounding $k$ are labeled with $i$.
Then there can be no critical branchings between $3$-cells in $\mathbf{SIso}_{3}^{i}$ and $\mathbf{SIso}_{3}^{j}$ unless $i=j$.

We can prove that the $(3,2)$-superpolygraph $(\mathbf{SIso}_{0},\mathbf{SIso}_{1},\mathbf{SIso}_{2},\mathbf{SIso}_{3}^i)$ is convergent for any $i\in I$ of parity $|i|=0$ by using an argument similar to the proof that the polygraph of pearls from \cite[Section 5.5]{GM09} is convergent.
Thus, if we prove that the $(3,2)$-superpolygraph $(\mathbf{SIso}_{0},\mathbf{SIso}_{1},\mathbf{SIso}_{2},\mathbf{SIso}_{3}^i)$ is convergent for an arbitrary $i\in I$ of parity $|i|=1$, then we will have proved that the entire $(3,2)$-superpolygraph $\mathbf{SIso}(\mathfrak{g})$ is convergent.
\end{remark}

% - - - - - - - - - - - - - - - -
\subsubsection{Termination}
% - - - - - - - - - - - - - - - -

We now prove the termination of the $(3,2)$-superpolygraph $\mathbf{SIso}(\mathfrak{g})$ using the derivation method  from Section~\ref{subsec:termination-derivation}.

\begin{lemma} \label{lem:step1}
  Let $U(\mathbf{SIso}(\mathfrak{g}))$  be the linear $(3,2)$-polygraph given by $U(\mathbf{SIso}(\mathfrak{g}))_i=\mathbf{SIso}(\mathfrak{g})_i$ forgetting the parity of elements in $\mathbf{SIso}(\mathfrak{g})$ as in Definition~\ref{def:UP}.  Then   the map $d \maps U(\mathbf{SIso}(\mathfrak{g}))_2^* \to \N $ given by
  \[ d(u) = ||u||_{ \{ \: \cuplAs{i}, \: \cuprAs{i}, \: \raisebox{-3mm}{$\caplAs{i}$}, \: \raisebox{-3mm}{$\caprAs{i}$} \: \} }-2\text{ times the number of odd bubbles} \]
is stable under contexts as described in \ref{sec:context-stable}.
\end{lemma}

\begin{proof}
   For $f\in \{u_{\l,0}, d_{\l,0}, u_{\l,0}', d_{\l,0}'\}$, we have $d(s_2(f))=2>0=d(t_2(f))$.
  Furthermore, for any context $c$ of $U(\mathbf{SIso}(\mathfrak{g}))_2^*$ such that $c[f]$ is defined, we have that
  \[||c[s_2(f)]||_{ \{ \: \cuplAs{i}, \: \cuprAs{i}, \: \raisebox{-3mm}{$\caplAs{i}$}, \: \raisebox{-3mm}{$\caprAs{i}$} \: \}}=||c[t_2(f)]||_{ \{ \: \cuplAs{i}, \: \cuprAs{i}, \: \raisebox{-3mm}{$\caplAs{i}$}, \: \raisebox{-3mm}{$\caprAs{i}$} \: \}}+2\] and $c[t_2(f)]$ must have at least as many odd bubbles as $c[s_2(f)]$.
  Thus, $d(c[s_2(f)])\geq d(c[t_2(f)])+2>d(c[t_2(f)])$ for $f\in \{u_{\l,0}, d_{\l,0}, u_{\l,0}', d_{\l,0}'\}$.

For $f\in \{i_1^\l, i_2^\l \}$, and any context $c$ for which $c[f]$ is defined, we have that $d(c[s_2(f)])\geq d(c[t_2(f)])$ because $c[s_2(f)]$ and $c[t_2(f)]$ have the same number of caps and cups and $c[s_2(f)]$ cannot have more odd bubbles than $c[t_2(f)]$ by the definition of the odd bubble in \ref{eq:def-odd-bubble}.
%  Thus, $d(c[s_2(f)])\geq d(t_2(f))$ for $f\in \{i_1^\l, i_2^\l \}$.
%
For $f\in \{i_\l^3, i_\l^4\}$, the context $c[t_2(f)]$ has two terms.
We have $d(c[s_2(f)])\geq d(c[h])$ for all $h\in \text{Supp}(t_2(f))$ using a similar argument for the first term of the target and observing that in the second term  the target has exactly two more caps and cups and at least one more odd bubble than the source.
%  For $f\in \{i_3^\l, i_4^\l \}$, a similar argument gives the desired inequality between the source and the first term of the target.
%  The inequality holds for the second term of the target as well because it has exactly two more caps and cups and at least one more odd bubble than the source.

The remaining 3-cells are endomorphism 2-cells of the identity $1_\l$ and it is straightforward to verify the desired inequality.
%
%  Since the source and target of the rest of the 3-cells are endomorphism 2-cells of the identity $1_\l$, it is enough to show that $d(s_2(f))\geq d(t_2(f))$.
%  We omit this check since it is nearly trivial.
\end{proof}

\begin{proposition}
The $(3,2)$-superpolygraph $\mathbf{SIso}(\mathfrak{g})$ terminates.
\end{proposition}

\begin{proof}
We prove the termination of $\mathbf{SIso}(\mathfrak{g})$ in steps as described in Section~\ref{sec:derivationbysteps}.

\noindent {\bf Step 1.} %First, consider
%  While not a derivation, this map is stable under contexts, so it can be used like a derivation as described in \ref{sec:context-stable}.
Using the context stable map from Lemma~\ref{lem:step1}, we have that $d(c[s_2(f)])>d(c[t_2(f)])$ for $f\in \{u_{\l,0}, d_{\l,0}, u_{\l,0}', d_{\l,0}'\}$ and $d(c[s_2(f)])\geq d[c[t_2(f)]]$ for the remaining 3-cells.  Hence, the map $d$ allows us to
 reduce termination of $\mathbf{SIso}(\mathfrak{g})$ to termination of \[\mathbf{SIso}(\mathfrak{g})':=(\mathbf{SIso}(\mathfrak{g})_0,\mathbf{SIso}(\mathfrak{g})_1,\mathbf{SIso}(\mathfrak{g})_2,\mathbf{SIso}(\mathfrak{g})_3-\{u_{\l,0}, d_{\l,0}, u_{\l,0}', d_{\l,0}'\}). \]
\smallskip

\noindent {\bf Step 2.}
Define 2-functors $X\maps U(\mathbf{SIso}(\mathfrak{g})')_2^* \to \cat{Ord}$ and $Y\maps (U(\mathbf{SIso}(\mathfrak{g})')_2^*)^{{\rm co}} \to \cat{Ord}$
%\begin{align}
%  X &\maps U(\mathbf{SIso}(\mathfrak{g})')_2^* \to \cat{Ord} \nn \\ \nn
%  Y &\maps (U(\mathbf{SIso}(\mathfrak{g})')_2^*)^{{\rm co}} \to \cat{Ord}
%\end{align}
%$X \maps U(\mathbf{SIso}(\mathfrak{g}))_2^* \to \cat{Ord}$, $Y \maps (U(\mathbf{SIso}(\mathfrak{g}))_2^*)^{{\rm co}} \to \cat{Ord}$
 whose nonempty values are given on generators by
\begin{align*}
  & X\left( \hackcenter{\begin{tikzpicture}
  \draw[-,thick,black] (0,-0.4) to (0,.4);
   \node at (0,-.6) {$\scriptstyle{i}$};
   \node at (0.2,0) {$\scriptstyle{\lambda}$};
\end{tikzpicture}} \right )
=Y
\left( \hackcenter{ \begin{tikzpicture}
  \draw[-,thick,black] (0,-0.4) to (0,.4);
   \node at (0,-.6) {$\scriptstyle{i}$};
   \node at (0.2,0) {$\scriptstyle{\lambda}$};
\end{tikzpicture} } \right)=\N,
\quad
X\left( \hackcenter{ \cuplA{i}{} } \right)=
  X\left( \hackcenter{ \cuprA{i}{}  }\right)=(0,0),
\quad
  Y\left( \hackcenter{\caplA{i}{} } \right)=
Y\left( \hackcenter{ \caprA{i}{} } \right)=(0,0),
\\
 &
 X\left(\hackcenter{ \updA{i} }\right)(n)
 =Y\left(\hackcenter{ \updA{i}}\right)(n)
 =X\left( \hackcenter{ \dpdA{i} }\right)(n)
=
Y\left( \hackcenter{ \dpdA{i} }\right)(n)=n+1
\end{align*}
Then a derivation $d\maps U(\mathbf{SIso}(\mathfrak{g})')_2^* \to M_{X,Y,\Z}$ is defined from
\begin{align*}
 & d\left( \hackcenter{ \upd{i} } \right)(n,m)=0 \qquad
  d\left( \hackcenter{ \dpd{i} } \right)(n,m)=0
 \\
  &d\left( \hackcenter{\cupl{i}{} } \right)(n,m) =
   d\left( \hackcenter{\capl{i}{}} \right)(n,m)=m
  \qquad
  d\left( \hackcenter{\capr{i}{} } \right)(n,m)=
d\left( \hackcenter{\cupr{i}{}} \right)(n,m)=0.
\end{align*}
It is immediate from this definition that
\begin{alignat*}{4}
  d \left( \cupldl{i}{a} \right) (n,m) &=m+a \qquad
    &&d\left(\capldl{i}{a}\right)(n,m) =m+a \qquad
    &&      d \left( \hackcenter{\begin{tikzpicture}
\begin{scope} [ x = 10pt, y = 10pt, join = round, cap = round, thick, scale=1.8]
  \draw[<-,thick,black,scale=1.3] (0.1,1.2) to[out=180,in=90] (-.3,0.8);
  \draw[-,thick,black,scale=1.3] (0.5,0.8) to[out=90,in=0] (0.1,1.2);
 \draw[-,thick,black,scale=1.3] (-.3,0.8) to[out=-90,in=180] (0.1,0.4);
  \draw[-,thick,black,scale=1.3] (0.1,0.4) to[out=0,in=-90] (0.5,0.8);
 \node at (-0.45,0.8) {$\scriptstyle{i}$};
 \node at (0.55,1.3) {$\bullet$};
 \node at (0.85,1.3) {$\scriptstyle{a}$};
 \end{scope}
\end{tikzpicture}} \right) = 0 \qquad
&& d \left(\hackcenter{\oddbubble{i}} \right) = 0
  \\
  d\left(\cupldr{i}{a}\right)(n,m)&=m \qquad
    &&d\left(\capldr{i}{a}\right)(n,m)=m  \qquad
    &&d \left( \hackcenter{\begin{tikzpicture}
\begin{scope} [ x = 10pt, y = 10pt, join = round, cap = round, thick, scale=1.8]
  \draw[->,thick,black,scale=1.3] (0.1,1.2) to[out=180,in=90] (-.3,0.8);
  \draw[-,thick,black,scale=1.3] (0.5,0.8) to[out=90,in=0] (0.1,1.2);
 \draw[-,thick,black,scale=1.3] (-.3,0.8) to[out=-90,in=180] (0.1,0.4);
  \draw[-,thick,black,scale=1.3] (0.1,0.4) to[out=0,in=-90] (0.5,0.8);
 \node at (-0.45,0.8) {$\scriptstyle{i}$};
 \node at (0.55,1.3) {$\bullet$};
 \node at (0.85,1.3) {$\scriptstyle{a}$};
 \end{scope}
\end{tikzpicture}} \right) = 0
&&
\end{alignat*}
Then for every generating 3-cell $x \in \mathbf{SIso}(\mathfrak{g})_3'$, the inequalities $X(s_2(x))\geq X(h)$, $Y(s_2(x))\geq Y(h)$, and $d(s_2(x))\geq d(h)$  hold for every $h\in \text{Supp}(t_2(x))$.
%
%
%For convenience, denote $\text{max}\{d(h) \mid h\in \text{Supp}(t_2(x))\}$ by $d(t_2(x))$
%\begin{alignat*}{2}
%  & d(s_2(i_\lambda^1))(n,m)=0=d(t_2(i_\lambda^1))(n,m)
%    \qquad && d(s_2(i_\lambda^2))(n,m)=0=d(t_2(i_\lambda^2))(n,m)
%    \\
%  & d(s_2(i_\lambda^3))(n,m)=m+1>m=d(t_2(i_\lambda^3))(n,m)
%    \qquad && d(s_2(i_\lambda^4))(n,m)=m+1>m=d(t_2(i_\lambda^4))(n,m)
%    \\
% & d(s_2(\alpha_{m,k}))=d(k)=d(t_2(\alpha_{m,k}))
%    \qquad && d(s_2(\beta_{m,k}))=d(k)=d(t_2(\beta_{m,k}))
%    \\
%   & d(s_2(I_0))=0=d(t_2(I_0)) &&
%\end{alignat*}
%
Furthermore, for $f\in \{i_\lambda ^3, i_\lambda^4\}$ we have the inequalities $X(s_2(x))\geq X(h)$, $Y(s_2(x))\geq Y(h)$ and a strict inequality $d(s_2(x))>d(h)$ for every $h\in \text{Supp}(t_2(x))$.
This reduces the termination of $\mathbf{SIso}(\mathfrak{g})'$ to the termination of the $(3,2)$-superpolygraph $R$ with $R_{\leq 2}:=\mathbf{SIso}(\mathfrak{g})_{\leq 2}'$ and $R_3:=\mathbf{SIso}(\mathfrak{g})_3' - \{ i_\lambda^3, i_\lambda^4\}$.

%Since the inequality is strict for $x=i_\lambda^3$ and $x=i_\lambda^4$ , termination of $\mathbf{SIso}(\mathfrak{g})$ terminates is reduced to termination of the superpolygraph $R$ with $R_i:=\mathbf{SIso}(\mathfrak{g})_i$ for $i=0,1,2$ and $R_3:=\mathbf{SIso}(\mathfrak{g})_3 - \{ i_\lambda^3, i_\lambda^4\}$.
\smallskip

\noindent {\bf Step 3.} To prove termination of $R$, consider the derivation $d$ into the trivial $U(R)_2^*$-module $M_{*,*,\Z}$ counting the number caps and cups, that is
\[ d(u) = ||u||_{ \{ \: \cuplAs{i}, \: \cuprAs{i}, \: \raisebox{-3mm}{$\caplAs{i}$}, \: \raisebox{-3mm}{$\caprAs{i}$} \: \} } \]
for any $2$-cell $u$ of $R_2^s$. For every generating 3-cell in $\alpha \in R_3$, we have the inequality $d(s_2(\alpha))\geq d(h)$ for every $h\in \text{Supp}(t_2(\alpha))$.
\begin{alignat*}{2}
 & d(s_2(i_\lambda^1))=1=d(t_2(i_\lambda^1))
    \qquad
    && d(s_2(i_\lambda^2))=1=d(t_2(i_\lambda^2))
    \\
 & d(s_2(\alpha_{m,k}))=d(k)+4>d(k)+2=d(t_2(\alpha_{m,k}))
    \qquad
    && d(s_2(\beta_{m,k}))=d(k)+4>d(k)+2=d(t_2(\beta_{m,k}))
  \\
  & d(s_2(I_0))=2=d(t_2(I_0))  &&
\end{alignat*}

Furthermore, for $\alpha \notin \{i_\lambda^1, i_\lambda^2,I_0 \}$ we have strict inequalities $d(s_2(\alpha))> d(h)$ for every $h\in \text{Supp}(t_2(\alpha))$.
This reduces the termination of $R$ to the termination of the $(3,2)$-superpolygraph $R'$ with $R'_{\leq 2}:=R_{\leq 2}$ and $R'_3=\{ i_\lambda^1, i_\lambda^2,I_0\}$.

%\[d(f)=\text{the number of leftwards facing caps and cups in f}\]
%\begin{enumerate}
 % \item $d(s_2 (u_{\lambda, 0}))=0=d(t_2(u_{\lambda, 0}))$
 % \item $d(s_2(d_{\lambda, 0}))=0=d(t_2(d_{\lambda,0}))$
  %\item $d(s_2 (u_{\lambda, 1}))=0=d(t_2(u_{\lambda, 1}))$
  %\item $d(s_2(d_{\lambda, 1}))=0=d(t_2(d_{\lambda,1}))$
 % \item $d(s_2 (u'_{\lambda, 0}))=2 >0=d(t_2(u'_{\lambda, 0}))$
 % \item $d(s_2 (d'_{\lambda, 0})) =2> 0=d(t_2(d'_{\lambda, 0}))$
  %\item $d(s_2 (u'_{\lambda, 1})) =2>1= d(t_2(u'_{\lambda, 1}))$
  %\item $d(s_2 (d'_{\lambda, 1})) =2>1= d(t_2(d'_{\lambda, 1}))$
 % \item $d(s_2(i_\lambda^1))=0=d(t_2(i_\lambda^1))$
  %\item $d(s_2(i_\lambda^2))=0=d(t_2(i_\lambda^2))$
  %\item $d(s_2(\alpha_{m,k}))=d(k)+2>=d(k)+1=d(t_2(\alpha_{m,k}))$
  %\item $d(s_2(\beta_{m,k}))=d(k)+2>d(k)+1=d(t_2(\beta_{m,k}))$
%\end{enumerate}

%Since the inequality is strict for $u'_{\lambda, 0}, d'_{\lambda, 0}, u'_{\lambda, 1}, d'_{\lambda, 1}, \alpha_{m,k}, \beta_{m,k}$, termination of $R$ is reduced to termination of $R'=(R_0,R_1,R_2, R_3':=R_3-\{  u'_{\lambda, 0}, d'_{\lambda, 0}, u'_{\lambda, 1}, d'_{\lambda, 1}, \alpha_{m,k}, \beta_{m,k} \})$.

%Thus, $\mathbf{SIso}(\mathfrak{g})$ terminates if $R'=(\mathbf{SIso}(\mathfrak{g})_0, \mathbf{SIso}(\mathfrak{g})_1,\mathbf{SIso}(\mathfrak{g})_2,R_3')$ terminates where
%$R_3'=\{u_{\lambda, 0}, d_{\lambda, 0}, u_{\lambda, 1}, d_{\lambda, 0}, i_\lambda^1, i_\lambda^2 \}$

\smallskip

\noindent {\bf Step 4.}  Now consider 2-functors $X:U(R')_2^* \to \cat{Ord}$,$Y:(U(R')_2^*)^{co} \to \cat{Ord}$
\begin{align*}
 & X\left( \hackcenter{\begin{tikzpicture}
  \draw[-,thick,black] (0,-0.4) to (0,.4);
   \node at (0,-.6) {$\scriptstyle{i}$};
   \node at (0.2,0) {$\scriptstyle{\lambda}$};
\end{tikzpicture}} \right )
=Y
\left( \hackcenter{ \begin{tikzpicture}
  \draw[-,thick,black] (0,-0.4) to (0,.4);
   \node at (0,-.6) {$\scriptstyle{i}$};
   \node at (0.2,0) {$\scriptstyle{\lambda}$};
\end{tikzpicture} } \right)=\N,
\quad \;
X\left( \cuplA{i}{} \right)=X\left( \cuprA{i}{} \right)=(0,0),
\quad \;
 Y\left( \caplA{i}{} \right)= Y\left( \caprA{i}{} \right)=(0,0)
\\
& X\left(\updA{i}\right)(n)=Y\left(\updA{i}\right)(n)=X\left( \dpdA{i} \right)(n)=Y\left(\dpdA{i}\right)(n)=n+1
\end{align*}
and derivation $d:U(R')_2^* \to M_{X,Y,\Z}$ given by
\begin{align*}
&  d\left( \updA{i} \right)(n,m)=0, \qquad  d\left( \dpdA{i} \right)(n,m)=0,
\\
&  d\left( \cuplA{i}{} \right)(n,m) =m, \quad\; d\left( \caplA{i}{} \right)(n,m)=m, \quad\; d\left( \caprA{i}{} \right)(n,m)=m, \quad\; d\left( \cuprA{i}{} \right)(n,m)=m.
\end{align*}
Then we have the desired inequalities $X(s_2(\alpha))\geq X(h)$, $Y(s_2(\alpha))\geq Y(h)$, and $d(s_2(\alpha))\geq d(h)$ for every $h \in \text{Supp}(t_2(\alpha))$ for every generating 3-cell $\alpha$ of $R'$, with strict inequalities for $\alpha \in \{i_\lambda^1,i_\lambda^2 \}$.
%\begin{enumerate}
%  \item $d(s_2(i_\lambda^1))(n,m)=m+1>m=d(t_2(i_\lambda^1))(n,m)$
%  \item $d(s_2(i_\lambda^2))(n,m)=m+1>m=d(t_2(i_\lambda^2))(n,m)$
%  \item $d(s_2(I_0))=0=d(t_2(I_0))$
%\end{enumerate}
%
So termination of $R'$ reduces to termination of $R'':=(R'_0,R'_1,R'_2,\{I_0\})$.
\smallskip

\noindent {\bf Step 5. }Consider the derivation  $d$ into the trivial module $M_{*,*,\Z}$ defined by
\[d(u)= ||u||_{ \cuplAs{i} } .\]
Then we have that $d(s_2(I_0))=1>0=d(t_2(I_0))$, so $R''$ terminates.
Hence, $R'$ terminates.  Therefore $\mathbf{SIso}(\mathfrak{g})$ terminates.
\end{proof}

% - - - - - - - - - - - - - - - -
\subsubsection{Convergence of $\mathbf{SIso}(\mathfrak{g})$}
% - - - - - - - - - - - - - - - -

\begin{proposition}
The $(3,2)$-superpolygraph $\mathbf{SIso}(\mathfrak{g})$ defined in Section~\ref{sec:OIso} is convergent.
\end{proposition}

\begin{proof}
Since $\mathbf{SIso}(\mathfrak{g})$ is terminating, following \cite[Theorem 4.2.13]{AL16} its confluence is equivalent to the confluence of its critical branchings, that are all proved confluent in Appendix \ref{appendix:criticalbranchingsSIso}.
\end{proof}

% #################################
\section{A convergent presentation of the odd nilHecke algebra} \label{sec:ONil}
% #################################

% ---------------------------------
\subsection{Definition of odd nilHecke 2-supercategory}
%----------------------------------

Here we recall the odd nilHecke algebra and its associated 2-supercategory.  This algebra appeared independently in~\cite{EKL,KKO} and is closely related to the spin Hecke algebra associated to the affine Hecke-Clifford superalgebra appearing in earlier work of Wang~\cite{Wang}. %and many of the essential features of the odd nilHecke algebra including skew-polynomials appears much earlier in this and related works on spin symmetric groups~\cite{KW1,KW2,KW4}.

\begin{definition} \label{def:oddNil}
Define the odd nilHecke 2-supercategory to have
\begin{enumerate}[i)]
  \item one object $\ast$,
  \item 1-morphisms $n \in \N$,
  \item 2-morphisms generated by
$\identdot{}: 1 \to 1\vspace{-0.15in}$ and $\raisebox{-7mm}{$\crossing{}{} \vspace{-0.15in}$}: 2 \to 2$
both of parity $\bar{1}$
\end{enumerate}
modulo the relations
\[
 \raisebox{-8mm}{$\dbcrossing{}{}$} \:= 0,
 \quad
 \raisebox{-9mm}{$\ybleft{}{}{}$} \: = \:
 \raisebox{-8mm}{$\ybright{}{}{}$},
\quad
\hackcenter{\begin{tikzpicture}
\begin{scope} [ x = 10pt, y = 10pt, join = round, cap = round, thick, scale = 2 ] \draw[-] (0.00,0.75)--(0.00,0.50) ;
\draw[-] (1.00,0.75)--(1.00,0.50) ; \draw (0.00,0.50)--(1.00,0.00) (1.00,0.50)--(0.00,0.00) ; \draw (0.00,0.00)--(0.00,-0.25) (1.00,0.00)--(1.00,-0.25) ;
\node at (1,0.50) {$\bullet$} ;
\end{scope}
\end{tikzpicture}  }
\: +
\hackcenter{\begin{tikzpicture}
\begin{scope} [ x = 10pt, y = 10pt, join = round, cap = round, thick, scale = 2 ] \draw[-] (0.00,0.75)--(0.00,0.50) ;
\draw[-] (1.00,0.75)--(1.00,0.50) ; \draw (0.00,0.50)--(1.00,0.00) (1.00,0.50)--(0.00,0.00) ; \draw (0.00,0.00)--(0.00,-0.25) (1.00,0.00)--(1.00,-0.25) ;
\node at (0,0.00) {$\bullet$} ;
\end{scope}
\end{tikzpicture}}
\;=\;  \raisebox{-7mm}{$\did{}{}$},
\quad
\hackcenter{\begin{tikzpicture}
\begin{scope} [ x = 10pt, y = 10pt, join = round, cap = round, thick, scale = 2 ] \draw[-] (0.00,0.75)--(0.00,0.50) ;
\draw[-] (1.00,0.75)--(1.00,0.50) ; \draw (0.00,0.50)--(1.00,0.00) (1.00,0.50)--(0.00,0.00) ; \draw (0.00,0.00)--(0.00,-0.25) (1.00,0.00)--(1.00,-0.25) ;
\node at (0,0.50) {$\bullet$} ;
\end{scope}
\end{tikzpicture}}
\; + \hackcenter{\begin{tikzpicture}
\begin{scope} [ x = 10pt, y = 10pt, join = round, cap = round, thick, scale = 2 ] \draw[-] (0.00,0.75)--(0.00,0.50) ;
\draw[-] (1.00,0.75)--(1.00,0.50) ; \draw (0.00,0.50)--(1.00,0.00) (1.00,0.50)--(0.00,0.00) ; \draw (0.00,0.00)--(0.00,-0.25) (1.00,0.00)--(1.00,-0.25) ;
\node at (1,0.00) {$\bullet$} ;
\end{scope}
\end{tikzpicture}} = \raisebox{-7mm}{$\did{}{}$} \: . \vspace{-0.15in}\]
\end{definition}

% ---------------------------------
\subsection{The super (3,2)-polygraph ONH}
%----------------------------------
\label{subsec:ConvergentPresentationONH}

% - - - - - - - - - - - - - - - -
\subsubsection{Definition}
% - - - - - - - - - - - - - - - -
In this section we define a $(3,2)$-superpolygraph presenting the odd Nilhecke 2-supercategory.
Let $\mathbf{ONH}$ be the $(3,2)$-superpolygraph defined by
\begin{enumerate}
  \item one object denoted by $\lambda$,
  \item one generating 1-cell denoted $1$, with $n$ denoting the $\star_0$-composition of $1$ with itself $n$ times. Since there are only one $0$-cell and one generating $1$-cell, we omit them in the string diagrams below.
  \item generating 2-cells $\identdot{}$ and $\raisebox{-7mm}{$\crossing{}{}$}$ both of parity $\bar{1}$,
  \item generating 3-cells
  \[ \hspace{-0.5cm} \raisebox{-6mm}{$\dbcrossing{}{}$} \: \overset{dc}{\Rrightarrow} 0,
\qquad
\raisebox{-8mm}{$\ybleft{}{}{}$} \: \overset{yb}{\Rrightarrow} \: \raisebox{-8mm}{$\ybright{}{}{}$},
\qquad
\begin{array}{l}
   \hackcenter{\begin{tikzpicture}
    \begin{scope} [ x = 10pt, y = 10pt, join = round, cap = round, thick, scale = 2 ] \draw[-] (0.00,0.75)--(0.00,0.50) ;
    \draw[-] (1.00,0.75)--(1.00,0.50) ; \draw (0.00,0.50)--(1.00,0.00) (1.00,0.50)--(0.00,0.00) ; \draw (0.00,0.00)--(0.00,-0.25) (1.00,0.00)--(1.00,-0.25) ;
    \node at (0,0.50) {$\bullet$} ;
    \end{scope}
    \end{tikzpicture}} \:
    \overset{on_1}{\Rrightarrow} -
        \hackcenter{\begin{tikzpicture}
        \begin{scope} [ x = 10pt, y = 10pt, join = round, cap = round, thick, scale = 2 ] \draw[-] (0.00,0.75)--(0.00,0.50) ;
        \draw[-] (1.00,0.75)--(1.00,0.50) ; \draw (0.00,0.50)--(1.00,0.00) (1.00,0.50)--(0.00,0.00) ; \draw (0.00,0.00)--(0.00,-0.25) (1.00,0.00)--(1.00,-0.25) ;
        \node at (1,0.00) {$\bullet$} ;
        \end{scope}
        \end{tikzpicture}} + \raisebox{-7mm}{$\did{}{}$} ,
\\
         \hackcenter{\begin{tikzpicture}
        \begin{scope} [ x = 10pt, y = 10pt, join = round, cap = round, thick, scale = 2 ] \draw[-] (0.00,0.75)--(0.00,0.50) ;
        \draw[-] (1.00,0.75)--(1.00,0.50) ; \draw (0.00,0.50)--(1.00,0.00) (1.00,0.50)--(0.00,0.00) ; \draw (0.00,0.00)--(0.00,-0.25) (1.00,0.00)--(1.00,-0.25) ;
        \node at (1,0.50) {$\bullet$} ;
        \end{scope}
        \end{tikzpicture}} \: \overset{on_2}{\Rrightarrow} -
        \hackcenter{\begin{tikzpicture}
        \begin{scope} [ x = 10pt, y = 10pt, join = round, cap = round, thick, scale = 2 ] \draw[-] (0.00,0.75)--(0.00,0.50) ;
        \draw[-] (1.00,0.75)--(1.00,0.50) ; \draw (0.00,0.50)--(1.00,0.00) (1.00,0.50)--(0.00,0.00) ; \draw (0.00,0.00)--(0.00,-0.25) (1.00,0.00)--(1.00,-0.25) ;
        \node at (0,0.00) {$\bullet$} ;
        \end{scope}
        \end{tikzpicture}} + \raisebox{-7mm}{$\did{}{}$}.
        \vspace{-0.15in}
\end{array}
\]
\end{enumerate}

\subsubsection{Termination}
% - - - - - - - - - - - - - - - -

We closely follow \cite[Section 2.3.3]{DUP19bis} to prove the termination of the $(3,2)$-superpolygraph $\mathbf{ONH}$ in two steps.
\smallskip

\begin{proposition} \label{prop:ONH-terminating}
The (3,2)-superpolygraph $\mathbf{ONH}$ terminates.
\end{proposition}

\begin{proof}
\noindent {\bf Step 1.}  Define a $2$-functor $X: U(\mathbf{ONH})_2^\ast \to \mathbf{Ord}$ by setting $X(i) = \N$, so that $X (i \star_0 i) = \N \times \N$, and on generating $2$-cells of $\mathbf{ONH}$ by:
\[ X \big( \ident{} \big) (n)=n \hspace{1cm} X \big( \identdot{} \big) (n)= n \hspace{1cm} X \big( \raisebox{-7mm}{$\scalebox{0.9}{\crossing{}{}}$} \big) (n,m) = (m,n+1 )\vspace{-0.15in} \]
for all $n,m \in \N$.
%We consider the $U(\mathbf{ONH})_2^\ast$-module $M_{X,*,\Z}$.
Define a derivation $d \maps U(\mathbf{ONH})_2^* \to M_{X,*,\Z}$ on the generating $2$-cells of $\mathbf{ONH}$ by \[ d \big( \ident{} \big) (n) = 0, \qquad d \big( \raisebox{-7mm}{$\scalebox{0.9}{\crossing{}{}}$} \big) (n,m) = m, \qquad d \big( \identdot{} \big)(n) = 0  \vspace{-0.15in} \]
for any $n,m \in \N$.
Then by the same calculation in \cite[Section 2.3.3]{DUP19bis} for the even nilHecke algebra, we obtain the inequalities $X(s_2(f)) \geq X(t_2(f))$ and $d(s_2(f)) \geq d(t_2(f))$ for all 3-cells $f$ and $d(s_2(\alpha)) > d(t_2(\alpha))$ for $\alpha \in \{yb, dc\}$.
Thus termination of $\mathbf{ONH}$ is reduced to termination of $\mathbf{ONH}':=(\mathbf{ONH}_0,\mathbf{ONH}_1,\mathbf{ONH}_2,\{on_1, on_2\})$
\smallskip

\noindent {\bf Step 2.}
%To prove that the $(3,2)$-superpolygraph $\mathbf{ONH}'$ terminates,
%we use
Define a $2$-functor $X: U(\mathbf{ONH})_2^\ast \to \mathbf{Ord}$ on the generating $2$-cells of $\mathbf{ONH}$ by:
\[ X \big( \ident{} \big) (n)=n \hspace{1cm} X \big( \identdot{} \big) (n)= n \hspace{1cm} X \big( \raisebox{-7mm}{$\scalebox{0.9}{\crossing{}{}}$} \big) (n,m) = (m+2,n+1 ) \vspace{-0.15in}\]
for all $n,m \in \N$, and a derivation $d:U(\mathbf{ONH}')_2^* \to M_{X,*,\Z}$ given by \[ d \big( \ident{} \big) (n) = 0, \qquad d \big( \raisebox{-7mm}{$\scalebox{0.9}{\crossing{}{}}$} \big) (n,m) = n \qquad d \big( \identdot{} \big)(n) = n \vspace{-0.15in}\]
for any $n,m \in \N$.
Then by \cite[Section 2.3.3]{DUP19bis}, we obtain the desired inequalities $X(s_2(\alpha)) \geq X(t_2(\alpha))$ and $d(s_2(\alpha)) > d(t_2(\alpha))$ for $\alpha \in \{on_1, on_2\}$, so that
%
%Therefore, the $2$-functor $X$ and the derivation $d$ satisfy the conditions of
Theorem \ref{T:TerminationDerivationTheorem} implies that $\mathbf{ONH}'$ is terminating, and thus $\mathbf{ONH}$ is terminating.
\end{proof}

Moreover, we now prove the following result:
\begin{proposition}
The $(3,2)$-superpolygraph $\mathbf{ONH}$ is convergent.
\end{proposition}

\begin{proof}
Since $\mathbf{ONH}$ is terminating by Proposition~\ref{prop:ONH-terminating}, following \cite[Theorem 4.2.13]{AL16}, its confluence is equivalent to the confluence of its critical branchings, whose classification follows from \cite{DUP19bis}, and are all proved confluent in Appendix \ref{appendix:criticalbranchingsONH}.
\end{proof}

% - - - - - - - - - - - - - - - -
\subsubsection{Bases of ONH}
% - - - - - - - - - - - - - - - -

\begin{definition}
Define the \emph{normal form basis} for the 2-supercategory $ONH$ to be the basis obtained from the convergent (3,2)-superpolygraph $\mathbf{ONH}$. This basis is obtained by choosing a fixed representative from each equivalence class of normal forms modulo superinterchange.
\end{definition}

In practice, to rewrite a 2-cell in $\mathbf{ONH}_2^s$, one checks if there is a representative in its equivalence class modulo superinterchange that is reducible by a 3-cell. If there is more than one representative where a 3-cell can be applied, the convergence of the superpolygraph ensures that it does not matter which representative is chosen to apply a 3-cell.   Then a 2-cell is in its normal form if and only if for any representative modulo superinterchange, this representative is irreducible using the set of 3-cells in the (3,2)-superpolygraph $\mathbf{ONH}$.

In the case of the odd nilHecke algebra, we can further specify the resulting normal form basis by making a preferred choice of representative of the superinterchange class for the order of dots; for example, choosing that dots will decrease in height going from left to right.   With our fixed choice or ordering of dots, we can represent these dot sequences as $\und{x}^{\alpha} =x_1^{\alpha_1} \dots x_n^{\alpha_n}$ with $\alpha_1$ dots appearing on the first strand, $\alpha_2$ dots below these on the second, and so on.

 The 3-cells in $\mathbf{ONH}$ ensure that all dots appearing in a given normal form 2-cell appear below any crossings. Then for each reduced expression of $w=s_{i_1} \dots s_{i_k}$ of a permutation in the symmetric group $S_n$, there is a corresponding crossing diagram $\partial_w = \partial_{i_1} \dots \partial_{i_k}$ in the odd nilHecke algebra, where $\partial_i$ is the crossing of the $i$th and $(i+1)$st lines.  The crossings appearing at the top of a normal form diagram will have reduced expressions $\partial_w$ where no equivalence class under superinterchange admits a reduction $\partial_i \partial_{i+1} \partial_i \Rrightarrow \partial_{i+1} \partial_{i}\partial_{i+1}$.  The superinterchange equivalence class may still be undetermined if the reduced expression contains a subsequence of the form $\partial_i \partial_j = - \partial_j \partial_i$ with $|i-j| >1$.  We can then uniquely specify a representative by choosing the ordering $\partial_i \partial_j$ where $i\leq j$.  An example is given below with the reduced expression $s_2s_1s_3s_2$, rather than $s_2s_1s_3s_2$, illustrating this choice of ordering.
\[
\partial_2\partial_1\partial_3\partial_2x_1^{\alpha_1}x_2^{\alpha_2}x_3^{\alpha_3}x_4^{\alpha_4} \;\; := \;\;
\hackcenter{
\begin{tikzpicture} [scale =0.65, thick]
    \draw (1,-1) -- (1,.25) (2,-1) -- (2,.25) (3,-1) -- (3,.25) (4,-1) -- (4,.25) ;
    \draw (1,.75) -- (1,.25)  (2,.25) -- (3,.75) (3,.25) -- (2,.75) (4,.25) -- (4,.75) ;
    \draw (1,.75) -- (1,1) (2,.75) -- (2,1) (3,.75) -- (3,1) (4,.75) -- (4,1) ;
    \draw (1,1) -- (1,1.5) (2,1) -- (2,1.5) (3,1) -- (4,1.5) (4,1) -- (3,1.5) ;
    \draw (1,1.5) to (1,1.75) -- (2,2.25)  (2,1.5) to (2,1.75) -- (1,2.25) (3,1.5) -- (3,2.25) (4,1.5) -- (4,2.25) ;
    \draw (1,2.25) -- (1,2.5) (2,2.25) -- (2,2.5) (3,2.25) -- (3,2.5) (4,2.25) -- (4,2.5) ;
    \draw (1,2.5) -- (1,3)  (2,2.5) -- (3,3) (3,2.5) -- (2,3) (4,2.5) -- (4,3) ;
    \draw (1,3) -- (1,3.25) (2,3) -- (2,3.25) (3,3) -- (3,3.25) (4,3) -- (4,3.25) ;
      \node at (1,0) {$\bullet$};
      \node at (2,-.25) {$\bullet$};
      \node at (3,-.5) {$\bullet$};
      \node at (4,-.75) {$\bullet$};
     \node at (.65, .1) {$\scriptstyle \alpha_1$};
     \node at (1.65, -.15) {$\scriptstyle \alpha_2$};
     \node at (2.65, -.4) {$\scriptstyle \alpha_3$};
     \node at (3.65, -.65) {$\scriptstyle \alpha_4$};
\end{tikzpicture}}
\]

In \cite[Proposition 2.11]{EKL}, bases for the odd nilHecke algebra are defined by making a choice of a reduced expression for each element $w \in S_n$ and considering elements $\{\partial_{w} \und{x}^{\alpha} \}$ or $\{\und{x}^{\alpha}\partial_{w}\}$ where $\partial_w$.

\begin{proposition}
The superpolygraph $\mathbf{ONH}$ presents the odd nilHecke 2-supercategory.  The resulting normal form basis recovers the basis  $\{\partial_{w} \und{x}^{\alpha} \}$ from \cite[Proposition 2.11]{EKL} where the choice of reduced expressions cannot be simplified further by any application of the identity $\partial_i \partial_{i+1} \partial_i =\partial_{i+1} \partial_{i}\partial_{i+1}$ for any representative of the superinterchange equivalence class of $\partial_{w} \und{x}^{\alpha}$.
\end{proposition}

% #################################
\section{Rewriting modulo in the odd 2-category}  \label{sec:odd2Cat}
% #################################

% ---------------------------------
\subsection{Definition of the odd 2-category}
%----------------------------------
\label{subsec:definitionodd2cat}

Ellis and Brundan give a description of the odd 2-category $\mf{U}(\mf{sl}_2)$ involving a minimal number of relations by requiring the invertibility of certain maps lifting the $\mf{sl}_2$-relations.  They show that the invertibility of these maps imply the relations given below.  In the definition that follows we do not attempt to provide a minimal set of relations.  In section~\ref{subsec:polygraphOsl2} we will explain how to reduce the number of generating 2-morphisms and defining relations in a way that will be helpful for presenting this super 2-category by a $(3,2)$-superpolygraph.

\begin{definition} \label{def:oddU}
The odd $2$-supercategory $\oUcat=\oUcat(\mf{sl}_2)$ is the $2$-supercategory
consisting of
\begin{itemize}
\item objects $\l$ for $\l \in \Z$,
\item
for a signed sequence $\underline{\epsilon}=(\epsilon_1,\epsilon_2, \dots, \epsilon_m)$, with $\epsilon_1,\dots,\epsilon_m \in \{+,-\}$, define
\[\cal{E}_{\underline{\epsilon}}:= \cal{E}_{\epsilon_1} \cal{E}_{\epsilon_2} \dots \cal{E}_{\epsilon_m}
\]
where $\cal{E}_{+}:=\cal{E}$ and $\cal{E}_{-}:= \cal{F}$.  A
$1$-morphisms from $\l$ to $\l'$ is a formal finite direct sum of strings
\[
 \cal{E}_{\underline{\epsilon}}\onel   =
  \1_{\l'}\cal{E}_{\underline{\epsilon}}
\]
for any   signed sequence $\underline{\epsilon}$ such that $\l'= \l+ 2\sum_{j=1}^m\epsilon_j1$.

\item $2$-morphisms are generated by the $\Z\times \Z_2$-graded generating 2-morphisms
\begin{align}
\hackcenter{\begin{tikzpicture}[scale=0.8]
    \draw[thick, ->] (0,0) -- (0,1.25)
        node[pos=.5, shape=coordinate](DOT){};
    \filldraw  (DOT) circle (2.5pt);
    \node at (-.85,.85) {\tiny $\lambda +2$};
    \node at (.5,.85) {\tiny $\lambda$};
    %\node at (-.2,.1) {\tiny $i$};
\end{tikzpicture}}
    &\maps \cal{E}\onel \to \cal{E}\onel %\la 2\ra
& \quad &
\hackcenter{\begin{tikzpicture}[scale=0.8]
    \draw[thick, <-] (0,0) -- (0,1.25)
        node[pos=.5, shape=coordinate](DOT){};
    \filldraw  (DOT) circle (2.5pt);
    \node at (-.85,.85) {\tiny $\lambda -2$};
    \node at (.5,.85) {\tiny $\lambda$};
    %\node at (-.2,1.4) {\tiny $i$};
\end{tikzpicture}}
\maps \cal{F}\onel \to \cal{F}\onel %\la i\cdot i \ra  \nn
\nn \\
 & \text{degree } (2, \bar{1}) & & \qquad \text{degree } (2, \bar{1}) \nn
\\
  \hackcenter{\begin{tikzpicture}[scale=0.8]
    \draw[thick, ->] (0,0) .. controls (0,.5) and (.75,.5) .. (.75,1.0);
    \draw[thick, ->] (.75,0) .. controls (.75,.5) and (0,.5) .. (0,1.0);
    \node at (1.1,.55) {\tiny $\lambda$};
    %\node at (-.2,.1) {\tiny $i$};
%    \node at (.95,.1) {\tiny $j$};
\end{tikzpicture}} \;\; &\maps \cal{E}\cal{E}\onel  \to \cal{E}\cal{E}\onel %\la -2 \ra
 &  \quad &
 \quad
  \hackcenter{\begin{tikzpicture}[scale=0.8]
    \draw[thick, <-] (0,0) .. controls (0,.5) and (.75,.5) .. (.75,1.0);
    \draw[thick, <-] (.75,0) .. controls (.75,.5) and (0,.5) .. (0,1.0);
    \node at (1.1,.55) {\tiny $\lambda$};
%    \node at (-.25,.1) {\tiny $i$};
%    \node at (1,.1) {\tiny $j$};
\end{tikzpicture}}\;\;
\maps \cal{F}\cal{F}\onel  \to \cal{F}\cal{F}\onel\la -2 \ra
 \nn
 \\
 & \text{degree } (-2, \bar{1}) & & \qquad \text{degree } (-2, \bar{1}) \nn
 \\
 %% THIRD ROW
\hackcenter{\begin{tikzpicture}[scale=0.8]
    \draw[thick, <-] (.75,2) .. controls ++(0,-.75) and ++(0,-.75) .. (0,2);
    \node at (.4,1.2) {\tiny $\lambda$};
    %\node at (-.2,1.9) {\tiny $i$};
\end{tikzpicture}} \;\; &\maps \onel  \to \cal{F}\cal{E}\onel\la  1 + {\lambda}  \ra
&
&
\hackcenter{\begin{tikzpicture}[scale=0.8]
    \draw[thick, ->] (.75,2) .. controls ++(0,-.75) and ++(0,-.75) .. (0,2);
    \node at (.4,1.2) {\tiny $\lambda$};
    %\node at (.95,1.9) {\tiny $i$};
\end{tikzpicture}} \;\; \maps \onel  \to\cal{E}\cal{F}\onel\la  1 - {\lambda}  \ra   \nn
 \\
 &\text{degree } (1+\l,\bar{0} )
  & &
 \quad \text{degree } (1-\l, \overline{\l+1}) ) \nn
 \\
\hackcenter{\begin{tikzpicture}[scale=0.8]
    \draw[thick, ->] (.75,-2) .. controls ++(0,.75) and ++(0,.75) .. (0,-2);
    \node at (.4,-1.2) {\tiny $\lambda$};
   % \node at (.95,-1.9) {\tiny $i$};
\end{tikzpicture}} \;\; & \maps \cal{F}\cal{E}\onel \to\onel\la  1 + {\lambda} \ra   &
    &
\hackcenter{\begin{tikzpicture}[scale=0.8]
    \draw[thick, <-] (.75,-2) .. controls ++(0,.75) and ++(0,.75) .. (0,-2);
    \node at (.4,-1.2) {\tiny $\lambda$};
    %\node at (-.2,-1.9) {\tiny $i$};
\end{tikzpicture}} \;\;\maps\cal{E}\cal{F}\onel  \to\onel\la  1 - {\lambda}  \ra  \nn
\\
 &\text{degree } (1+\l, \overline{\l+1}))
  & &
 \quad \text{degree } (1-\l,\bar{0} ) \nn
\nn
\end{align}
where we have indicated a $Q$-grading and parity as an ordered tuple $(x,\bar{y})$.
% Note  that the $\Z_2$ degree of the right pointing cap and cup are not the mod 2 reductions of the $\Z$-degree.
\end{itemize}
The identity $2$-morphism of the $1$-morphism $\cal{E} \onen$ is represented by an upward oriented line (likewise, the identity $2$-morphism of $\cal{F} \onen$ is represented by a downward oriented line).

Horizontal and vertical composites of the above diagrams are interpreted using the conventions for supercategories explained in Section~\ref{sec:2supercat}.
 The rightmost region in our diagrams is usually colored by $\l$.  The fact that we are defining a 2-supercategory means that diagrams with odd parity skew commute.
The $2$-morphisms satisfy the following relations
 (see \cite{BE2} for more details).
\begin{enumerate}

\item {\bf (Odd nilHecke)} The odd nilHecke relations from Definition~\ref{def:oddNil} are satisfied for upward oriented strands and any $\l\in \Z$.

\item {\bf (Odd isotopies)} The odd isotopy relations from Definition~\ref{def:SIso-cat} for a Cartan data with a single odd $i\in I$.

\item {\bf (Bubble relations)}
\noindent Dotted bubbles of negative degree are zero, so that
for all $m \geq 0$
\begin{equation}
 \hackcenter{ \begin{tikzpicture} [scale=.65]
 %\draw (-.15,.35) node { $\scs i$};
 \draw[ ]  (0,0) arc (180:360:0.5cm) [thick];
 \draw[<- ](1,0) arc (0:180:0.5cm) [thick];
\filldraw  [black] (.1,-.25) circle (2.5pt);
 \node at (-.2,-.5) {\tiny $m$};
 \node at (1.15,.8) { $\lambda  $};
\end{tikzpicture} } \;  = 0\quad \text{if $m < \l -1$}, \qquad \quad
\;
\hackcenter{ \begin{tikzpicture} [scale=.65]
 %\draw (-.15,.35) node { $\scs i$};
 \draw  (0,0) arc (180:360:0.5cm) [thick];
 \draw[->](1,0) arc (0:180:0.5cm) [thick];
\filldraw  [black] (.9,-.25) circle (2.5pt);
 \node at (1,-.5) {\tiny $m$};
 \node at (1.15,.8) { $\lambda $};
\end{tikzpicture} } \;  = 0 \quad  \text{if $m < -\l -1$}.
\end{equation}
Dotted bubbles of degree $0$ are equal to the identity $2$-morphism:
\begin{equation} \label{eq:degreezero}
 \hackcenter{ \begin{tikzpicture} [scale=.65]
 %\draw (-.15,.35) node { $\scs i$};
 \draw  (0,0) arc (180:360:0.5cm) [thick];
 \draw[,<-](1,0) arc (0:180:0.5cm) [thick];
\filldraw  [black] (.1,-.25) circle (2.5pt);
 \node at (-.5,-.5) {\tiny $\l -1$};
 \node at (1.15,1) { $\lambda  $};
\end{tikzpicture} }
\;\; =
\Id_{\1_{\l}}
\quad \text{for $  \l \geq 1$}, \qquad \quad
\;
\hackcenter{ \begin{tikzpicture} [scale=.65]
 %\draw (-.15,.35) node { $\scs i$};
 \draw  (0,0) arc (180:360:0.5cm) [thick];
 \draw[->](1,0) arc (0:180:0.5cm) [thick];
\filldraw  [black] (.9,-.25) circle (2.5pt);
 \node at (1.35,-.5) {\tiny $-\l -1$};
 \node at (1.15,1) { $\lambda $};
\end{tikzpicture} }
\;\; =
\Id_{\1_{\l}} \quad \text{if $   \l \leq -1$}.
\end{equation}
We will sometimes make use of the shorthand notation
\begin{equation}
   \vcenter{\hbox{$\posbubd{}{ n+\ast\quad }$}} :=
     \vcenter{\hbox{$\posbubd{}{\l-1+n\qquad} $}}
\qquad
   \vcenter{\hbox{$\negbubd{}{\; \;n+\ast }$}} :=
     \vcenter{\hbox{$\negbubd{}{\qquad \;-\l-1+n\;} $}}
\end{equation}
 The degree two bubble is given a special notation as in \eqref{eq:def-odd-bubble} and squares to zero by the superinterchange law.

We call a clockwise (resp. counterclockwise) bubble  fake
if $m+\l-1<0$ and (resp. if $m-\l-1<0$).   These correspond to positive degree bubbles that are labeled by a negative number of dots. These are to be interpreted as formal symbols recursively defined by
 the odd infinite Grassmannian relations
\begin{align} \label{eq:fake-def}
\vcenter{\hbox{$\posbubd{}{2n + \ast \quad \;}$}}
&:=  - \sum\limits_{l=1}^{ n }
\raisebox{2.5mm}{\hbox{$\posbubdfffsl{}{ 2(n-\ell) + \ast \qquad }$}}
\raisebox{-2.5mm}{\hbox{$\negbubdfffsl{}{\; 2l+ \ast }$}}
\qquad \text{for $0 \leq 2n <-\l$,}
\\ \nn
\vcenter{\hbox{$\negbubd{}{\;\;\;2n + \ast }$}}
&:=  - \sum\limits_{l=1}^{ n }
\raisebox{2.5mm}{\hbox{$\posbubdfffsl{}{ 2l+ \ast \;\;}$}}
\raisebox{-2.5mm}{\hbox{$\negbubdfffsl{}{\qquad\;  2(n-\ell) + \ast }$}}
\qquad \text{for $0 \leq 2n <\l$,}
\\ \nn
\vcenter{\hbox{$\posbubd{}{2n+1 + \ast \quad \;}$}}
&:=
\stackrel{\vcenter{\hbox{$\posbubdfffsl{}{ 2n+ \ast \quad }
                    $}}}{\xy (-10,0)*{\;}; (0,0)*{\bigotimes}; \endxy}
\qquad \text{for $0 \leq 2n <-\l$,}
\\ \nn
\vcenter{\hbox{$\negbubd{}{\qquad   2n+1 + \ast }$}}
&:= \quad
\stackrel{\vcenter{\hbox{$\negbubdfffsl{}{ 2n+ \ast  }
                    $}}}{\xy (-4,0)*{\;}; (0,0)*{\bigotimes}; \endxy}
\qquad \text{for $0 \leq 2n +1 <\l$,}
\end{align}

\item {\bf (Centrality of odd bubbles) } Odd bubbles are central
\begin{equation} \label{rel: centrality of odd bubbles}
\xy 0;/r.17pc/:
 (0,8);(0,-8); **\dir{-} ?(1)*\dir{>}+(2.3,0)*{\scriptstyle{}}
 ?(.1)*\dir{ }+(2,0)*{\scs };
 (-5,-2)*{\txt\large{$\bigotimes$}};
 (-8,6)*{ \l};
% (8,5)*{ \l+2};
 (-10,0)*{};(10,0)*{};(-2,-8)*{\scs };
 \endxy
\;\; = \;\;
\xy 0;/r.17pc/:
 (0,8);(0,-8); **\dir{-} ?(1)*\dir{>}+(2.3,0)*{\scriptstyle{}}
 ?(.1)*\dir{ }+(2,0)*{\scs };
 (5,-2)*{\txt\large{$\bigotimes$}};
 (-8,6)*{ \l};
% (8,5)*{ \l+2};
 (-10,0)*{};(10,0)*{};(-2,-8)*{\scs };
 \endxy
 \qquad \qquad
\xy 0;/r.17pc/:
 (0,8);(0,-8); **\dir{-} ?(0)*\dir{<}+(2.3,0)*{\scriptstyle{}}
 ?(.1)*\dir{ }+(2,0)*{\scs };
 (-5,-2)*{\txt\large{$\bigotimes$}};
 (-8,6)*{ \l};
% (8,5)*{ \l+2};
 (-10,0)*{};(10,0)*{};(-2,-8)*{\scs };
 \endxy
\;\; = \;\;
\xy 0;/r.17pc/:
 (0,8);(0,-8); **\dir{-} ?(0)*\dir{<}+(2.3,0)*{\scriptstyle{}}
 ?(.1)*\dir{ }+(2,0)*{\scs };
 (5,-2)*{\txt\large{$\bigotimes$}};
 (-8,6)*{ \l};
% (8,5)*{ \l+2};
 (-10,0)*{};(10,0)*{};(-2,-8)*{\scs };
 \endxy
\end{equation}

\item \label{item_cycbiadjoint} {\bf (Odd crossing cyclicity)}
%
%\begin{equation} \label{eq:cyclic_dot}
%\hackcenter{\begin{tikzpicture}[scale=0.8]
%    \draw[thick, <-]  (0,-1) to (0,1);
%    %\node at (.8,-.4) { $\lambda+\alpha_i$};
%    \node at (-.5,-.4) { $\lambda$};
%    \filldraw  (0,.2) circle (2.5pt);
%   % \node at (-.2,.8) {\tiny $i$};
%\end{tikzpicture}}
% \; := \;
%\hackcenter{\begin{tikzpicture}[scale=0.8]
%    \draw[thick, ->]  (0,.4) .. controls ++(0,.6) and ++(0,.6) .. (.75,.4) to (.75,-1);
%    \draw[thick, <-](0,.4) to (0,-.4) .. controls ++(0,-.6) and ++(0,-.6) .. (-.75,-.4) to (-.75,1);
%    \filldraw  (0,-.2) circle (2.5pt);
%    %\node at (1,-.9) { $\lambda+\alpha_i$};
%    %\node at (1.3,.9) { $\lambda + \alpha_i$};
%   % \node at (-.95,.8) {\tiny $i$};
%\end{tikzpicture}}
%     \quad = \quad 2 \;
%\hackcenter{\begin{tikzpicture}[scale=0.8]
%    \draw[thick, <-]  (0,-1) to (0,1);
%   % \node at (.8,-.4) { $\lambda+\alpha_i$};
%    \node at (-.5,-.4) { $\lambda$};
%    \filldraw  (0,.2) circle (2.5pt);
%    %\node at (-.2,.8) {\tiny $i$};
%    \node at (-1,0) {\large{$\bigotimes$}};
%\end{tikzpicture}}
% \;\; - \;\;
%\hackcenter{\begin{tikzpicture}[scale=0.8]
%    \draw[thick, ->]  (0,.4) .. controls ++(0,.6) and ++(0,.6) .. (-.75,.4) to (-.75,-1);
%    \draw[thick, <-](0,.4) to (0,-.4) .. controls ++(0,-.6) and ++(0,-.6) .. (.75,-.4) to (.75,1);
%    \filldraw  (0,-.2) circle (2.5pt);
%    %\node at (1,-.9) { $\lambda+\alpha_i$};
%    \node at (-1,.9) { $\lambda$};
%    \node at (.95,.8) {\tiny $i$};
%\end{tikzpicture}}
%\end{equation}
%
The cyclic relations for crossings\footnote{Equation~\ref{eq:cyclic} differs by a sign from \cite[Equation (1.28)]{BE2}, but is consistent with the original formulation of the odd 2-category from \cite{Lau-odd}. } are given by
\begin{equation} \label{eq:cyclic}
\hackcenter{
\begin{tikzpicture}[scale=0.8]
    \draw[thick, <-] (0,0) .. controls (0,.5) and (.75,.5) .. (.75,1.0);
    \draw[thick, <-] (.75,0) .. controls (.75,.5) and (0,.5) .. (0,1.0);
    \node at (1.1,.65) { $\lambda$};
%    \node at (-.2,.1) {\tiny $i$};
%    \node at (.95,.1) {\tiny $j$};
\end{tikzpicture}}
\;\; := \;\;
\hackcenter{\begin{tikzpicture}[scale=0.7]
    \draw[thick, ->] (0,0) .. controls (0,.5) and (.75,.5) .. (.75,1.0);
    \draw[thick, ->] (.75,0) .. controls (.75,.5) and (0,.5) .. (0,1.0);
    \draw[thick] (0,0) .. controls ++(0,-.4) and ++(0,-.4) .. (-.75,0) to (-.75,2);
    \draw[thick] (.75,0) .. controls ++(0,-1.2) and ++(0,-1.2) .. (-1.5,0) to (-1.55,2);
    \draw[thick, ->] (.75,1.0) .. controls ++(0,.4) and ++(0,.4) .. (1.5,1.0) to (1.5,-1);
    \draw[thick, ->] (0,1.0) .. controls ++(0,1.2) and ++(0,1.2) .. (2.25,1.0) to (2.25,-1);
    \node at (-.35,.75) {  $\lambda$};
%    \node at (1.3,-.7) {\tiny $i$};
%    \node at (2.05,-.7) {\tiny $j$};
%    \node at (-.9,1.7) {\tiny $i$};
%    \node at (-1.7,1.7) {\tiny $j$};
\end{tikzpicture}}
\quad = \quad -
\hackcenter{\begin{tikzpicture}[xscale=-1.0, scale=0.7]
    \draw[thick, ->] (0,0) .. controls (0,.5) and (.75,.5) .. (.75,1.0);
    \draw[thick, ->] (.75,0) .. controls (.75,.5) and (0,.5) .. (0,1.0);
    \draw[thick] (0,0) .. controls ++(0,-.4) and ++(0,-.4) .. (-.75,0) to (-.75,2);
    \draw[thick] (.75,0) .. controls ++(0,-1.2) and ++(0,-1.2) .. (-1.5,0) to (-1.55,2);
    \draw[thick, ->] (.75,1.0) .. controls ++(0,.4) and ++(0,.4) .. (1.5,1.0) to (1.5,-1);
    \draw[thick, ->] (0,1.0) .. controls ++(0,1.2) and ++(0,1.2) .. (2.25,1.0) to (2.25,-1);
    \node at (1.2,.75) {  $\lambda$};
%    \node at (1.3,-.7) {\tiny $j$};
%    \node at (2.05,-.7) {\tiny $i$};
%    \node at (-.9,1.7) {\tiny $j$};
%    \node at (-1.7,1.7) {\tiny $i$};
\end{tikzpicture}} .
\end{equation}

Sideways crossings satisfy the following identities:
\begin{align} \label{eq:crossl-gen-cyc}
\hackcenter{
\begin{tikzpicture}[scale=0.8]
    \draw[thick, ->] (0,0) .. controls (0,.5) and (.75,.5) .. (.75,1.0);
    \draw[thick, <-] (.75,0) .. controls (.75,.5) and (0,.5) .. (0,1.0);
    \node at (1.1,.65) { $\lambda$};
%    \node at (-.2,.1) {\tiny $i$};
%    \node at (.95,.1) {\tiny $j$};
\end{tikzpicture}}
  := \;\;
\hackcenter{\begin{tikzpicture}[scale=0.7]
    \draw[thick, ->] (0,0) .. controls (0,.5) and (.75,.5) .. (.75,1.0);
    \draw[thick, ->] (.75,-.5) to (.75,0) .. controls (.75,.5) and (0,.5) .. (0,1.0) to (0,1.5);
    \draw[thick] (0,0) .. controls ++(0,-.4) and ++(0,-.4) .. (-.75,0) to (-.75,1.5);
    %\draw[thick] (.75,0) .. controls ++(0,-1.2) and ++(0,-1.2) .. (-1.5,0) to (-1.55,2);
    \draw[thick, ->] (.75,1.0) .. controls ++(0,.4) and ++(0,.4) .. (1.5,1.0) to (1.5,-.5);
    %\draw[thick, ->] (0,1.0) .. controls ++(0,1.2) and ++(0,1.2) .. (2.25,1.0) to (2.25,-1);
    \node at (1.85,.55) {  $\lambda$};
%    \node at (1.75,-.2) {\tiny $j$};
%    \node at (.55,-.2) {\tiny $i$};
%    \node at (-.9,1.2) {\tiny $j$};
%    \node at (.25,1.2) {\tiny $i$};
\end{tikzpicture}}
   = \;\;
\hackcenter{\begin{tikzpicture}[xscale=-1.0, scale=0.7]
    \draw[thick, <-] (0,0) .. controls (0,.5) and (.75,.5) .. (.75,1.0);
    \draw[thick, <-] (.75,-.5) to (.75,0) .. controls (.75,.5) and (0,.5) .. (0,1.0) to (0,1.5);
    \draw[thick] (0,0) .. controls ++(0,-.4) and ++(0,-.4) .. (-.75,0) to (-.75,1.5);
    %\draw[thick] (.75,0) .. controls ++(0,-1.2) and ++(0,-1.2) .. (-1.5,0) to (-1.55,2);
    \draw[thick, ] (.75,1.0) .. controls ++(0,.4) and ++(0,.4) .. (1.5,1.0) to (1.5,-.5);
    %\draw[thick, ->] (0,1.0) .. controls ++(0,1.2) and ++(0,1.2) .. (2.25,1.0) to (2.25,-1);
    \node at (-1.1,.55) {  $\lambda$};
%    \node at (1.75,-.2) {\tiny $i$};
%    \node at (1,-.2) {\tiny $j$};
%    \node at (-.9,1.2) {\tiny $i$};
%    \node at (.25,1.2) {\tiny $j$};
\end{tikzpicture}}
\qquad
 \hackcenter{
\begin{tikzpicture}[scale=0.8]
    \draw[thick, <-] (0,0) .. controls (0,.5) and (.75,.5) .. (.75,1.0);
    \draw[thick, ->] (.75,0) .. controls (.75,.5) and (0,.5) .. (0,1.0);
    \node at (1.1,.65) { $\lambda$};
%    \node at (-.2,.1) {\tiny $i$};
%    \node at (.95,.1) {\tiny $j$};
\end{tikzpicture}}
     := (-1)^{\l+1}
\hackcenter{\begin{tikzpicture}[xscale=-1.0, scale=0.7]
    \draw[thick, ->] (0,0) .. controls (0,.5) and (.75,.5) .. (.75,1.0);
    \draw[thick, ->] (.75,-.5) to (.75,0) .. controls (.75,.5) and (0,.5) .. (0,1.0) to (0,1.5);
    \draw[thick] (0,0) .. controls ++(0,-.4) and ++(0,-.4) .. (-.75,0) to (-.75,1.5);
    %\draw[thick] (.75,0) .. controls ++(0,-1.2) and ++(0,-1.2) .. (-1.5,0) to (-1.55,2);
    \draw[thick, ->] (.75,1.0) .. controls ++(0,.4) and ++(0,.4) .. (1.5,1.0) to (1.5,-.5);
    %\draw[thick, ->] (0,1.0) .. controls ++(0,1.2) and ++(0,1.2) .. (2.25,1.0) to (2.25,-1);
    \node at (-1.1,.55) {  $\lambda$};
%    \node at (1.75,-.2) {\tiny $i$};
%    \node at (1,-.2) {\tiny $j$};
%    \node at (-.9,1.2) {\tiny $i$};
%    \node at (.25,1.2) {\tiny $j$};
\end{tikzpicture}}   - \;\;
\hackcenter{\begin{tikzpicture}[scale=0.7]
    \draw[thick, <-] (0,0) .. controls (0,.5) and (.75,.5) .. (.75,1.0);
    \draw[thick, <-] (.75,-.5) to (.75,0) .. controls (.75,.5) and (0,.5) .. (0,1.0) to (0,1.5);
    \draw[thick] (0,0) .. controls ++(0,-.4) and ++(0,-.4) .. (-.75,0) to (-.75,1.5);
    %\draw[thick] (.75,0) .. controls ++(0,-1.2) and ++(0,-1.2) .. (-1.5,0) to (-1.55,2);
    \draw[thick, ] (.75,1.0) .. controls ++(0,.4) and ++(0,.4) .. (1.5,1.0) to (1.5,-.5);
    %\draw[thick, ->] (0,1.0) .. controls ++(0,1.2) and ++(0,1.2) .. (2.25,1.0) to (2.25,-1);
    \node at (1.85,.55) {  $\lambda$};
%    \node at (1.75,-.2) {\tiny $j$};
%    \node at (.55,-.2) {\tiny $i$};
%    \node at (-.9,1.2) {\tiny $j$};
%    \node at (.25,1.2) {\tiny $i$};
\end{tikzpicture}}
\end{align}

\item {\bf (Odd sl(2) relations)}
\begin{equation}
\begin{split}
 \hackcenter{\begin{tikzpicture}[scale=0.8]
    \draw[thick] (0,0) .. controls ++(0,.5) and ++(0,-.5) .. (.75,1);
    \draw[thick,<-] (.75,0) .. controls ++(0,.5) and ++(0,-.5) .. (0,1);
    \draw[thick] (0,1 ) .. controls ++(0,.5) and ++(0,-.5) .. (.75,2);
    \draw[thick, ->] (.75,1) .. controls ++(0,.5) and ++(0,-.5) .. (0,2);
    %    \node at (-.2,.15) {\tiny $i$};
 %   \node at (.95,.15) {\tiny $i$};
     \node at (1.1,1.44) { $\lambda $};
\end{tikzpicture}}
\;\; + \;\;
\hackcenter{\begin{tikzpicture}[scale=0.8]
    \draw[thick, ->] (0,0) -- (0,2);
    \draw[thick, <-] (.75,0) -- (.75,2);
  %   \node at (-.2,.2) {\tiny $i$};
  %  \node at (.95,.2) {\tiny $i$};
     \node at (1.1,1.44) { $\lambda $};
\end{tikzpicture}}
\;\; = \;\;
\sum_{\overset{f_1+f_2+f_3}{=\l-1}}(-1)^{f_2}\; \hackcenter{
 \begin{tikzpicture}[scale=0.8]
 \draw[thick,->] (0,-1.0) .. controls ++(0,.5) and ++ (0,.5) .. (.8,-1.0) node[pos=.75, shape=coordinate](DOT1){};
  \draw[thick,<-] (0,1.0) .. controls ++(0,-.5) and ++ (0,-.5) .. (.8,1.0) node[pos=.75, shape=coordinate](DOT3){};
 \draw[thick,->] (0,0) .. controls ++(0,-.45) and ++ (0,-.45) .. (.8,0)node[pos=.25, shape=coordinate](DOT2){};
 \draw[thick] (0,0) .. controls ++(0,.45) and ++ (0,.45) .. (.8,0);
 %\draw (-.15,.7) node { $\scs i$};
%\draw (1.05,0) node { $\scs i$};
%\draw (-.15,-.7) node { $\scs i$};
% \draw (.3,.125) node {};
% \draw  (0,0) arc (180:360:0.5cm) [thick];
% \draw[,<-](1,0) arc (0:180:0.5cm) [thick];
%\filldraw  [black] (.1,-.25) circle (2.5pt);
 \node at (.95,.65) {\tiny $f_3$};
 \node at (-.55,-.05) {\tiny $\overset{-\l-1}{+f_2}$};
  \node at (.95,-.65) {\tiny $f_1$};
 \node at (1.65,.3) { $\lambda $};
 \filldraw[thick]  (DOT3) circle (2.5pt);
  \filldraw[thick]  (DOT2) circle (2.5pt);
  \filldraw[thick]  (DOT1) circle (2.5pt);
\end{tikzpicture} }
\\
 \hackcenter{\begin{tikzpicture}[scale=0.8]
    \draw[thick,<-] (0,0) .. controls ++(0,.5) and ++(0,-.5) .. (.75,1);
    \draw[thick] (.75,0) .. controls ++(0,.5) and ++(0,-.5) .. (0,1);
    \draw[thick, ->] (0,1 ) .. controls ++(0,.5) and ++(0,-.5) .. (.75,2);
    \draw[thick] (.75,1) .. controls ++(0,.5) and ++(0,-.5) .. (0,2);
    %    \node at (-.2,.15) {\tiny $i$};
%    \node at (.95,.15) {\tiny $i$};
     \node at (1.1,1.44) { $\lambda $};
\end{tikzpicture}}
\;\; + \;\;
\hackcenter{\begin{tikzpicture}[scale=0.8]
    \draw[thick, <-] (0,0) -- (0,2);
    \draw[thick, ->] (.75,0) -- (.75,2);
    % \node at (-.2,.2) {\tiny $i$};
%    \node at (.95,.2) {\tiny $i$};
     \node at (1.1,1.44) { $\lambda $};
\end{tikzpicture}}
\;\; = \;\;
\sum_{\overset{f_1+f_2+f_3}{=-\l-1}}(-1)^{f_2}\; \hackcenter{
 \begin{tikzpicture}[scale=0.8]
 \draw[thick,<-] (0,-1.0) .. controls ++(0,.5) and ++ (0,.5) .. (.8,-1.0) node[pos=.75, shape=coordinate](DOT1){};
  \draw[thick,->] (0,1.0) .. controls ++(0,-.5) and ++ (0,-.5) .. (.8,1.0) node[pos=.75, shape=coordinate](DOT3){};
 \draw[thick ] (0,0) .. controls ++(0,-.45) and ++ (0,-.45) .. (.8,0)node[pos=.25, shape=coordinate](DOT2){};
 \draw[thick, ->] (0,0) .. controls ++(0,.45) and ++ (0,.45) .. (.8,0);
% \draw (-.15,.7) node { $\scs i$};
%\draw (1.05,0) node { $\scs i$};
%\draw (-.15,-.7) node { $\scs i$};
% \draw (.3,.125) node {};
% \draw  (0,0) arc (180:360:0.5cm) [thick];
% \draw[,<-](1,0) arc (0:180:0.5cm) [thick];
%\filldraw  [black] (.1,-.25) circle (2.5pt);
 \node at (.95,.65) {\tiny $f_3$};
 \node at (-.55,-.05) {\tiny $\overset{\l-1}{+f_2}$};
  \node at (.95,-.65) {\tiny $f_1$};
 \node at (1.65,.3) { $\lambda $};
 \filldraw[thick]  (DOT3) circle (2.5pt);
  \filldraw[thick]  (DOT2) circle (2.5pt);
  \filldraw[thick]  (DOT1) circle (2.5pt);
\end{tikzpicture} } . \label{eq:sl2}
\end{split}
\end{equation}
\end{enumerate}
\end{definition}

\begin{remark} \label{rem:Symd}
Let ${\sf Sym}$ denote the algebra of symmetric functions over $\Bbbk$.  This algebra is generated by elementary symmetric functions ${\sf e}_r$ for $r\geq 0$ and by the complete symmetric functions ${\sf h}_s$ with $s\geq 0$.  By convention ${\sf e}_0={\sf h}_0=1$.   These generators are related by the equations
\[
\sum_{r+s=n} (-1)^s {\sf e}_r {\sf h}_s =0 \quad \text{for all $n\geq 0$. }
\]
Let ${\sf Sym}[d]$ be the supercommutative superalgebra obtained by placing ${\rm Sym}$ in even degree and adjoining an odd generator ${\sf d}$ with ${\sf d}^2=0$.  Then consider the unique surjective homomorphism
\begin{align}
  \beta_{\l} \maps {\sf Sym} \longrightarrow \End_{\mf{U}}(\1_{\l})
\end{align}
such that
\begin{alignat}{2}
&  {\sf e}_n \mapsto \;\; \vcenter{\hbox{$\posbubd{}{\l-1+2n \qquad \;\;}$}} \quad  \text{if $n> -\frac{h}{2}$},
\quad &&
{\sf h}_n \mapsto   (-1)^n\l \vcenter{\hbox{$\negbubd{}{\quad \quad \; -\l-1+2n }$}} \;\;
 \text{if $n> \frac{h}{2}$},
\\
&  {\sf de}_n \mapsto \vcenter{\hbox{$\posbubd{}{\l-1+2n+1 \qquad \quad }$}} \;  \text{if $n> -\frac{h}{2}$},
\quad &&
{\sf  dh}_n \mapsto (-1)^n\;   \vcenter{\hbox{$\negbubd{}{\qquad \quad  \;-\l-1+2n+1 }$}} \;\;
 \text{if $n> \frac{h}{2}$}.
\end{alignat}
The relations in $\mf{U}$ imply that this is a homomorphism and that the relations \eqref{eq:fake-def} defining the fake bubbles hold for all values of $\l$ and for all $n\geq 0$, see \cite[Proposition 5.1]{BE2}.
\end{remark}

% ---------------------------------
\subsection{The super (3,2)-polygraph $Osl(2)$}
%----------------------------------
\label{subsec:polygraphOsl2}

%% - - - - - - - - - - - - - - - -
%\subsubsection{Definition}
%% - - - - - - - - - - - - - - - -

\begin{definition}
Let $\mathbf{Osl(2)}$ be the linear~$(3,2)$-polygraph defined by:
\begin{enumerate}[{\bf i)}]
\item the elements of $\mathbf{Osl(2)}_0$ are the weights $\lambda \in \mathbb{Z}$ of $\mathfrak{sl}_2$,
\item the elements of $\mathbf{Osl(2)}_1$ are given by
\[ 1_{\lambda'} \mathcal{E}_{\varepsilon_1} \dots \mathcal{E}_{\varepsilon_m} 1_{\lambda} \] for any sequence of signs $(\varepsilon_1, \dots, \varepsilon_m)$ and $\lambda$,$\lambda'$ in $\Z$. Such a $1$-cell has for $0$-source $\lambda$ and $0$-target $\lambda '$, and
%\[ 1_{n''} \mathcal{E}_{\varepsilon'_1} \dots \mathcal{E}_{\varepsilon_l} 1_{n'} \star_0 1_{n'} \mathcal{E}_{\varepsilon_1} \dots \mathcal{E}_{\varepsilon_m} 1_{n} = 1_{n''} \mathcal{E}_{\varepsilon'_1} \dots \mathcal{E}_{\varepsilon_m} 1_{n} \]
\[1_{\lambda '} \mathcal{E}_{\varepsilon_1} \dots \mathcal{E}_{\varepsilon_m } 1_{\lambda} \star_0 1_{\lambda ''} \mathcal{E}_{\varepsilon'_1} \dots \mathcal{E}_{\varepsilon_l} 1_{\lambda'}  = 1_{\lambda''} \mathcal{E}_{\varepsilon'_1 } \dots \mathcal{E}_{\varepsilon_m} 1_{\lambda} \]
\item the elements of $\mathbf{Osl(2)}_2$ are the following generating $2$-cells: for $\lambda \in \Z$:
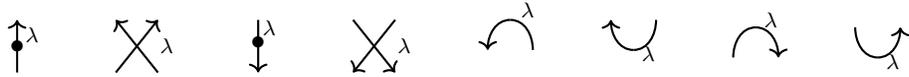
\begin{align*}
\udott{} \qquad \crossup{}{} \qquad \ddott{} \qquad \crossdn{}{} \qquad
\capl{} \qquad  \cupl{} \qquad \capr{} \qquad \cupr{}
\end{align*}
with respective parity $1$, $1$, $1$, $1$, $\l+1$, $\l+1$, $0$, $0$.
\item $\mathbf{Osl(2)}_3$ consists of the following $3$-cells:
\begin{itemize}
\item[1)] The odd nilHecke $3$-cells, given by
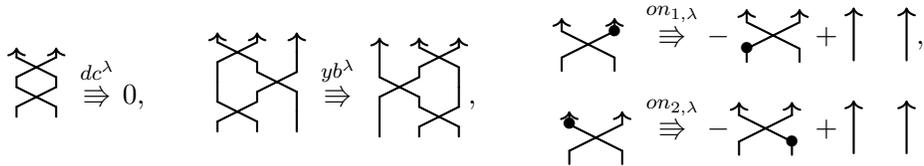
\[
\hspace{-0.5cm} \vcenter{\hbox{$\dbcrossingup{}{}$}} \: \overset{dc^\lambda}{\Rrightarrow} 0,
\qquad
\vcenter{\hbox{$\ybleftup{}{}{}$}} \: \overset{yb^\lambda}{\Rrightarrow} \: \vcenter{\hbox{$\ybrightup{}{}{}$}},
\qquad
\begin{array}{l}
\hackcenter{\begin{tikzpicture}
\begin{scope} [ x = 10pt, y = 10pt, join = round, cap = round, thick, scale = 2 ] \draw[<-] (0.00,0.75)--(0.00,0.50) ;
\draw[<-] (1.00,0.75)--(1.00,0.50) ; \draw (0.00,0.50)--(1.00,0.00) (1.00,0.50)--(0.00,0.00) ; \draw (0.00,0.00)--(0.00,-0.25) (1.00,0.00)--(1.00,-0.25) ;
\node at (1,0.50) {$\bullet$} ;
\end{scope}
\end{tikzpicture} }\: \overset{on_{1,\lambda}}{\Rrightarrow} - \hackcenter{\begin{tikzpicture}
\begin{scope} [ x = 10pt, y = 10pt, join = round, cap = round, thick, scale = 2 ] \draw[<-] (0.00,0.75)--(0.00,0.50) ;
\draw[<-] (1.00,0.75)--(1.00,0.50) ; \draw (0.00,0.50)--(1.00,0.00) (1.00,0.50)--(0.00,0.00) ; \draw (0.00,0.00)--(0.00,-0.25) (1.00,0.00)--(1.00,-0.25) ;
\node at (0,0.00) {$\bullet$} ;
\end{scope}
\end{tikzpicture}} + \raisebox{-7mm}{$\didup{}{}$},
\\
\hackcenter{\begin{tikzpicture}
\begin{scope} [ x = 10pt, y = 10pt, join = round, cap = round, thick, scale = 2 ] \draw[<-] (0.00,0.75)--(0.00,0.50) ;
\draw[<-] (1.00,0.75)--(1.00,0.50) ; \draw (0.00,0.50)--(1.00,0.00) (1.00,0.50)--(0.00,0.00) ; \draw (0.00,0.00)--(0.00,-0.25) (1.00,0.00)--(1.00,-0.25) ;
\node at (0,0.50) {$\bullet$} ;
\end{scope}
\end{tikzpicture} }
    \: \overset{on_{2,\lambda}}{\Rrightarrow} -
\hackcenter{\begin{tikzpicture}
\begin{scope} [ x = 10pt, y = 10pt, join = round, cap = round, thick, scale = 2 ] \draw[<-] (0.00,0.75)--(0.00,0.50) ;
\draw[<-] (1.00,0.75)--(1.00,0.50) ; \draw (0.00,0.50)--(1.00,0.00) (1.00,0.50)--(0.00,0.00) ; \draw (0.00,0.00)--(0.00,-0.25) (1.00,0.00)--(1.00,-0.25) ;
\node at (1,0.00) {$\bullet$} ;
\end{scope}
\end{tikzpicture} }
+ \raisebox{-7mm}{$\didup{}{}$}
\end{array}
\]
with the rightmost region of the diagram being labeled $\lambda$. When no confusion is likely to arise we often drop the $\l$ subscript from this notation.

\item[2)] The super isotopy $3$-cells of $\mathbf{SIso}_3$,

\item[3)] The cyclicity $3$-cell for the definition of the downward crossing:
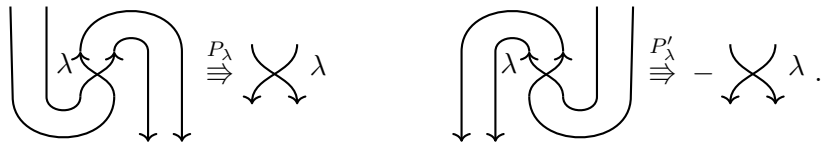
\[
\hackcenter{\begin{tikzpicture}[scale=0.6]
    \draw[thick, ->] (0,0) .. controls (0,.5) and (.75,.5) .. (.75,1.0);
    \draw[thick, ->] (.75,0) .. controls (.75,.5) and (0,.5) .. (0,1.0);
    \draw[thick] (0,0) .. controls ++(0,-.4) and ++(0,-.4) .. (-.75,0) to (-.75,2);
    \draw[thick] (.75,0) .. controls ++(0,-1.2) and ++(0,-1.2) .. (-1.5,0) to (-1.55,2);
    \draw[thick, ->] (.75,1.0) .. controls ++(0,.4) and ++(0,.4) .. (1.5,1.0) to (1.5,-1);
    \draw[thick, ->] (0,1.0) .. controls ++(0,1.2) and ++(0,1.2) .. (2.25,1.0) to (2.25,-1);
    \node at (-.35,.75) {  $\lambda$};
\end{tikzpicture}}
\: \overset{P_\lambda}{\Rrightarrow} \hackcenter{
\begin{tikzpicture}[scale=0.8]
    \draw[thick, <-] (0,0) .. controls (0,.5) and (.75,.5) .. (.75,1.0);
    \draw[thick, <-] (.75,0) .. controls (.75,.5) and (0,.5) .. (0,1.0);
    \node at (1.1,.65) { $\lambda$};
\end{tikzpicture}} \qquad \qquad
\hackcenter{\begin{tikzpicture}[xscale=-1.0, scale=0.6]
    \draw[thick, ->] (0,0) .. controls (0,.5) and (.75,.5) .. (.75,1.0);
    \draw[thick, ->] (.75,0) .. controls (.75,.5) and (0,.5) .. (0,1.0);
    \draw[thick] (0,0) .. controls ++(0,-.4) and ++(0,-.4) .. (-.75,0) to (-.75,2);
    \draw[thick] (.75,0) .. controls ++(0,-1.2) and ++(0,-1.2) .. (-1.5,0) to (-1.55,2);
    \draw[thick, ->] (.75,1.0) .. controls ++(0,.4) and ++(0,.4) .. (1.5,1.0) to (1.5,-1);
    \draw[thick, ->] (0,1.0) .. controls ++(0,1.2) and ++(0,1.2) .. (2.25,1.0) to (2.25,-1);
    \node at (1.2,.75) {  $\lambda$};
\end{tikzpicture}} \: \overset{P'_\lambda}{\Rrightarrow} \: - \: \hackcenter{
\begin{tikzpicture}[scale=0.8]
    \draw[thick, <-] (0,0) .. controls (0,.5) and (.75,.5) .. (.75,1.0);
    \draw[thick, <-] (.75,0) .. controls (.75,.5) and (0,.5) .. (0,1.0);
    \node at (1.1,.65) { $\lambda$};
\end{tikzpicture}} .
\]
together with their respective images $Q_\lambda$ and $Q'_\lambda$ through the Chevalley involution $\omega$ defined in \cite[Proposition 3.5]{BE2} giving the same cyclicity condition for the upward crossing in terms of the downward crossing.

%\item[3)] The $3$-cells coming from the leftward cap and cup generators:
%\begin{equation}
%\label{E:NewGenerators}
%\vcenter{\hbox{$\tfishdlp{}{- \lambda }$}} \overset{D_{\lambda}^{-}}{\Rrightarrow} \cupl{} \quad \text{for $\lambda \leq 0$}, \hspace{0.5cm}
%\vcenter{\hbox{$\tfishulpfME{}{ \lambda }$}} \overset{B_{\lambda}^{+}}{\Rrightarrow}  \capl{} \quad \text{for $\lambda \geq 0$}
%\end{equation}

\item[4)] The $3$-cells for the degree conditions on bubbles: for every $\lambda \in \Z$ and $n\in \N$
\begin{equation}
\label{E:ReductionPosBub}
\vcenter{\hbox{$\posbubd{}{n}$}} \underset{b_{\lambda}^{0,n}}{\overset{b_{\lambda}^1}{\Rrightarrow}}  \left\{
    \begin{array}{ll}
        1_{1_\lambda} & \mbox{if } n = \lambda - 1 \\
        0 & \mbox{if } n < \lambda - 1
    \end{array}
\right.
\end{equation}
\begin{equation}
\label{E:ReductionNegBub}
\vcenter{\hbox{$\negbubd{}{n}$}} \underset{c_{\lambda}^{0,n}}{\overset{c_{\lambda}^1}{\Rrightarrow}} \left\{
    \begin{array}{ll}
        1_{1_\lambda} & \mbox{if } n = - \lambda - 1 \\
        0 & \mbox{if } n < - \lambda - 1
    \end{array}
\right.
\end{equation}
\item[5)] The Infinite-Grassmannian $3$-cells: for any $\lambda \in \Z$ and $n \geq 1$ such that $2n+\lambda-1\geq 0$
\begin{align*}
\vcenter{\hbox{$\posbubdff{}{2n + \ast}$}} \overset{\text{ig}_{2n, \lambda}}{\Rrightarrow} - \sum\limits_{l=1}^{ n }
\raisebox{2.5mm}{\hbox{$\posbubdfffsl{}{ 2(n-\ell) + \ast \qquad }$}}
\raisebox{-2.5mm}{\hbox{$\negbubdfffsl{}{\; 2l+ \ast }$}}
\end{align*}
\item[6)] Bubble Slide $3$-cells
 \begin{equation} \label{eq:ccbubslideLR}
\hspace{-1.5cm}  \raisebox{-2mm}{$\posbubdfffsl{}{n+ \ast}$} \identu{}
  \overset{s_{\lambda,n}^+}{\Rrightarrow} \sum\limits_{r \geq 0} (2r+1) \identdotsufsl{}{2r} \posbubdffr{}{n-2r+\ast}
\end{equation}
\begin{equation*}
\raisebox{-2mm}{$\negbubdfffsl{}{n+ \ast}$} \identu{} \: \overset{s_{\lambda,n}^-}{\Rrightarrow} \: \identusl{} \negbubdff{}{n+ \ast} -3   \identdotsusl{}{2} \negbubdfff{}{n-2+ \ast} + 4 \sum\limits_{r \geq 2} (-1)^r \identdotsusl{}{2r} \negbubdfff{}{n-2r+ \ast}
\end{equation*}
and their reflections across the horizontal axis $r_{\lambda,n}^+$ and $r_{\lambda,n}^-$, which allow a bubble to go through a downwards strand.  The reflections correspond to the images of these relations via the Chevalley involution $\omega$ defined in \cite[Proposition 3.5]{BE2}.  By \eqref{eq:def-odd-bubble} and the definition of fake bubbles \eqref{eq:fake-def},   we simplify notation and write $s_{\l,1} = s_{\l,1}^+ = s_{\l,1}^-$ and $r_{\l,1} = r_{\l,1}^+ = r_{\l,1}^-$.
 These are added to the presentation to reach confluence modulo.
\item[7)] The invertibility $3$-cells:
\begin{align*}
\vcenter{\hbox{$\tdcrosslr{}{}$}}
%\raisebox{4mm}{$\tdcrosslr{}{}$}
& \overset{F_\lambda}{\Rrightarrow} -(-1)^{\lambda+1} \mathord{
\begin{tikzpicture}[baseline = 0]
	\draw[<-,thick,black] (0.08,-.3) to (0.08,.4);
	\draw[->,thick,black] (-0.28,-.3) to (-0.28,.4);
   \node at (-0.28,-.4) {$\scriptstyle{i}$};
   \node at (0.08,.5) {$\scriptstyle{i}$};
   \node at (.3,.05) {$\scriptstyle{\lambda}$};
\end{tikzpicture}
} +
\!\!\sum_{n=0}^{\lambda-1}
\sum_{r \geq 0} (-1)^{n+r}
\stdb{}{r}{-n-r-2}{n}, \\
\vcenter{\hbox{$\tdcrossrl{}{}$}}
&\overset{E_\lambda}{\Rrightarrow} - (-1)^{\lambda+1} \mathord{
\begin{tikzpicture}[baseline = 0]
	\draw[->,thick,black] (0.08,-.3) to (0.08,.4);
	\draw[<-,thick,black] (-0.28,-.3) to (-0.28,.4);
   \node at (-0.28,.5) {$\scriptstyle{i}$};
   \node at (0.08,-.4) {$\scriptstyle{i}$};
   \node at (.3,.05) {$\scriptstyle{\lambda}$};
\end{tikzpicture}
} +
\!\!\sum_{n=0}^{-\lambda-1}
\sum_{r \geq 0} (-1)^{n+r}
\stda{}{n}{-n-r-2}{r}.
\end{align*}
\item[8)] Remaining $3$-cells:
\begin{align*}
\tfishur{}
&\overset{C_{\lambda}}{\Rrightarrow}
\sum_{n=0}^{\lambda}
(-1)^n \rulec{}{n}; \qquad
  \tfishdr{}
\overset{A_{\lambda}}{\Rrightarrow}
\sum_{n=0}^{-\lambda}
(-1)^n \rulea{}{n}; \\
\tfishul{} & \overset{B_{\lambda}}{\Rrightarrow} \sum_{n=0}^{- \lambda}
(-1)^n \ruleb{}{n}; \qquad
\tfishdl{}  \overset{D_{\lambda}}{\Rrightarrow} \sum_{n=0}^{\lambda}
(-1)^n \ruled{}{n}.
\end{align*}
\begin{align*}
\raisebox{0.8mm}{$\begin{tikzpicture}[baseline = 0,scale=0.55]
	\draw[->,thick,black] (0.3,-.5) to (-0.3,.5);
	\draw[-,thick,black] (-0.2,-.2) to (0.2,.3);
        \draw[-,thick,black] (0.2,.3) to[out=50,in=180] (0.5,.5);
        \draw[->,thick,black] (0.5,.5) to[out=0,in=90] (0.8,-0.5);
        \draw[-,thick,black] (-0.2,-.2) to[out=230,in=0] (-0.5,-.4);
        \draw[-,thick,black] (-0.5,-.4) to[out=180,in=-90] (-0.8,.5);
        \draw[-,thick,black] (0.3,-1.5) to (-0.3,-2.6);
        \draw[-,thick,black] (0.5,-2.6) to[out=0,in=270] (0.8,-1.8);
        \draw[-,thick,black] (0.2,-2.3) to[out=130,in=180] (0.5,-2.6);
        \draw[-,thick,black] (-0.2,-1.8) to (0.2,-2.3);
        \draw[-,thick,black] (-0.2,-1.8) to[out=130,in=0] (-0.5,-1.6);
        \draw[->,thick,black] (-0.5,-1.6) to[out=180,in=-270] (-0.8,-2.6);
        \node at (1,0.5) {$\scriptstyle{\lambda}$};
        %% LEFT CROSSING PART
        \draw[-,thick,black] (0.3,-1.5) to (-1.5,-0.5);
        \draw[->,thick,black] (-1.5,-0.5) to (-1.5,0.5);
        \draw[-,thick,black]  (0.3,-0.5) to (-1.5,-1.5);
        \draw[-,thick,black] (-1.5,-1.5) to (-1.5,-2.6);
        \draw (0.8,-0.5) to (0.8,-2);
\end{tikzpicture}$}
\:
\raisebox{-4mm}{$\overset{\Gamma_\lambda}{\Rrightarrow}$} \:
\raisebox{0.8mm}{$\begin{tikzpicture}[baseline = 0,scale=0.55]
	\draw[<-,thick,black] (0.3,.5) to (-0.3,-.5);
	\draw[-,thick,black] (-0.2,.2) to (0.2,-.3);
        \draw[-,thick,black] (0.2,-.3) to[out=130,in=180] (0.5,-.4);
        \draw[-,thick,black] (0.5,-.4) to[out=0,in=270] (0.8,.5);
        \draw[-,thick,black] (-0.2,.2) to[out=130,in=0] (-0.5,.5);
        \draw[-,thick,black] (-0.5,.5) to[out=180,in=-270] (-0.8,-.5);
      \draw[-,thick,black] (-0.3,-1.5) to (0.3,-2.5);
         \draw[-,thick,black] (-0.5,-2.5) to[out=180,in=-90] (-0.8,-1.5);
         \draw[-,thick,black] (-0.2,-2.2) to[out=230,in=0] (-0.5,-2.5);
         \draw[-,thick,black] (-0.2,-2.2) to (0.2,-1.7);
         \draw[-,thick,black] (0.2,-1.7) to[out=50,in=180] (0.5,-1.6);
          \draw[->,thick,black] (0.5,-1.6) to[out=0,in=90] (0.7,-2.5);
    \node at (1.8,0.5) {$\scriptstyle{\lambda}$};
    %% RIGHT CROSSING PART
        \draw[-,thick,black] (-0.3,-1.5) to (1.5,-0.5);
        \draw[->,thick,black] (1.5,-0.5) to (1.5,0.5);
        \draw[-,thick,black] (-0.3,-0.5) to (1.5,-1.5);
        \draw[-,thick,black] (1.5,-1.5) to (1.5,-2.5);
        \draw[-,thick,black]  (-0.8,-0.5) to (-0.8,-1.5);
\end{tikzpicture}$} \: \raisebox{-4mm}{$- \:
\sum\limits_{r,s,t \geq 0}
(-1)^{r+s}
\begin{tikzpicture}[baseline = 0]
	\draw[-,thick,black] (0.3,0.7) to[out=-90, in=0] (0,0.3);
	\draw[->,thick,black] (0,0.3) to[out = 180, in = -90] (-0.3,0.7);
    \node at (1.4,-0.32) {$\scriptstyle{\lambda}$};
  \draw[->,thick,black] (0.2,0) to[out=90,in=0] (0,0.2);
  \draw[-,thick,black] (0,0.2) to[out=180,in=90] (-.2,0);
\draw[-,thick,black] (-.2,0) to[out=-90,in=180] (0,-0.2);
  \draw[-,thick,black] (0,-0.2) to[out=0,in=-90] (0.2,0);
   \node at (0.2,0) {$\color{black}\bullet$};
   \node at (.6,0) {$\color{black}\scriptstyle{\substack{-r-s\\-t-3}}$};
   \node at (-0.23,0.43) {$\color{black}\bullet$};
   \node at (-0.43,0.43) {$\color{black}\scriptstyle{r}$};
	\draw[<-,thick,black] (0.3,-.7) to[out=90, in=0] (0,-0.3);
	\draw[-,thick,black] (0,-0.3) to[out = 180, in = 90] (-0.3,-.7);
   \node at (-0.25,-0.5) {$\color{black}\bullet$};
   \node at (-.4,-.5) {$\color{black}\scriptstyle{s}$};
	\draw[->,thick,black] (.98,-0.7) to (.98,0.7);
   \node at (0.98,0.53) {$\color{black}\bullet$};
   \node at (1.15,0.53) {$\color{black}\scriptstyle{t}$};
\end{tikzpicture} + \sum\limits_{r,s,t \geq 0}
(-1)^{r+s}
\begin{tikzpicture}[baseline = 0]
	\draw[<-,thick,black] (0.3,0.7) to[out=-90, in=0] (0,0.3);
	\draw[-,thick,black] (0,0.3) to[out = 180, in = -90] (-0.3,0.7);
    \node at (.6,-0.32) {$\scriptstyle{\lambda}$};
  \draw[-,thick,black] (0.2,0) to[out=90,in=0] (0,0.2);
  \draw[<-,thick,black] (0,0.2) to[out=180,in=90] (-.2,0);
\draw[-,thick,black] (-.2,0) to[out=-90,in=180] (0,-0.2);
  \draw[-,thick,black] (0,-0.2) to[out=0,in=-90] (0.2,0);
   \node at (-0.2,0) {$\color{black}\bullet$};
   \node at (-0.6,0) {$\color{black}\scriptstyle{\substack{-r-s\\-t-3}}$};
   \node at (0.23,0.43) {$\color{black}\bullet$};
   \node at (0.43,0.43) {$\color{black}\scriptstyle{s}$};
	\draw[-,thick,black] (0.3,-.7) to[out=90, in=0] (0,-0.3);
	\draw[->,thick,black] (0,-0.3) to[out = 180, in = 90] (-0.3,-.7);
   \node at (0.25,-0.5) {$\color{black}\bullet$};
   \node at (.4,-.5) {$\color{black}\scriptstyle{r}$};
	\draw[->,thick,black] (-.98,-0.7) to (-.98,0.7);
   \node at (-0.98,0.53) {$\color{black}\bullet$};
   \node at (-1.15,0.53) {$\color{black}\scriptstyle{t}$};
\end{tikzpicture}$}
\end{align*}
\end{itemize}
\end{enumerate}
\end{definition}

Note that the last $3$-cell is added to the presentation to recover the Yang-Baxter relation for sideways crossing\footnote{The 3-cell $\Gamma_{\l}$ corrects a minor typo from \cite[Equation (7.20)]{BE2}.}, see \cite[Equation (7.20)]{BE2}, and is needed to reach confluence modulo and to fix a preferred choice of representative for all   possible orientations of the Yang-Baxter equations.

\begin{remark}
  By the definition of fake bubbles \eqref{eq:fake-def} in terms of positively dotted bubbles from $\oUcat$, we can use $ig_{2n,\lambda}$ for all $n\geq 1$  by using it as an equality for $2n+\lambda-1 <0$ and as an oriented $3$-cell for $2n+\lambda-1 \geq 0$, see also Remark~\ref{rem:Symd}.  Likewise, we can use $b_\lambda, c_\lambda$ for all $n\in \Z$ by using it as an equality for $n<0$ and as an oriented $3$-cell for $n\geq 0$.
\end{remark}

\begin{remark}
 The summations appearing in the targets of $r_{\l,n}^{\pm}$, $s_{\l,n}^{\pm}$, $E_\l$, $F_\l$, and $\Gamma_\l$ are assumed to be restricted so that no negative degree bubbles appear.  For example, the target of $F_{\l}$ has the summation with $r$ ranging from 0 to $\l-1-n$ and $E_{\l}$ the $r$ summation runs from 0 to $-\l-1-n$.
%
%  When we have a sum $\sum\limits_{r\geq 0}$ as in the targets of $r_{\l,n}^{\pm}$, $s_{\l,n}^{\pm}$, $E_\l$, $F_\l$, and $\Gamma_\l$ this actually denotes $\sum\limits_{r=0}^k$ where $k$ is the largest value of $r$ such that the $r=k$ term doesn't have a negative degree bubble.
%  For example,
The first sum in $t_2(\Gamma)$ implicitly has the restriction $-r-s-t+\l \geq 0$ since the degree of the bubble in that summand is $-r-s-t+\l$.
\end{remark}

%------------------
\subsection{Splitting of $\mathbf{Osl(2)}$}
%-----------------------
\label{sec:splittingofOsl2}
Let us split the $(3,2)$-superpolygraph $\mathbf{Osl(2)}$ into two parts: consider the $(3,2)$-superpolygraph $E$ defined by $E_i = \mathbf{Osl(2)}_i$ for $0 \leq i \leq 1$, $E_2 = \mathbf{Osl(2)} - \{ \hackcenter{
\begin{tikzpicture}[scale=0.45]
    \draw[thick, <-] (0,0) .. controls (0,.5) and (.75,.5) .. (.75,1.0);
    \draw[thick, <-] (.75,0) .. controls (.75,.5) and (0,.5) .. (0,1.0);
    %\node at (1.1,.65) { $\lambda$};
%    \node at (-.2,.1) {\tiny $i$};
%    \node at (.95,.1) {\tiny $j$};
\end{tikzpicture}} \} = \mathbf{SIso}_2 \cup \{ \hackcenter{
\begin{tikzpicture}[scale=0.45]
    \draw[thick, ->] (0,0) .. controls (0,.5) and (.75,.5) .. (.75,1.0);
    \draw[thick, ->] (.75,0) .. controls (.75,.5) and (0,.5) .. (0,1.0);
    %\node at (1.1,.65) { $\lambda$};
%    \node at (-.2,.1) {\tiny $i$};
%    \node at (.95,.1) {\tiny $j$};
\end{tikzpicture}} \}$ and $E_3 = \mathbf{SIso}_3 \cup \{ yb, \: dc \}$. Let  $R$  be the $(3,2)$-superpolygraph such that $R_i = \mathbf{Osl(2)}_i$ for $0 \leq i \leq 2$  and containing all the remaining $3$-cells.

\begin{proposition}
The (3,2)-superpolygraph $E$ is terminating.
\end{proposition}

\begin{proof}
The proof goes in steps as explained in Section~\ref{sec:derivationbysteps}:
\smallskip

\noindent {\bf Step 1.}  Eliminate the zigzag $3$-cells using the first step of the proof of termination of $\mathbf{SIso}$.
\smallskip

\noindent {\bf Step 2.} Eliminate $yb$ and $dc$ using the first step of the proof of termination of $\mathbf{ONH}$, extending values of $X$ and $d$ by
\vspace{-0.15in}
\[   X \big( \cupl{}{} \big)= X\big( \cupr{}{} \big)=(0,0), \quad d \big( \cupl{}{} \big)= d\big( \cupr{}{} \big)=0 \hspace{.5cm} d \big( \capl{}{} \big)(n,m)= d\big( \capr{}{} \big)(n,m)=0 \]
so that the  inequalities $d(s_2(\alpha)) \geq d(t_2(\alpha))$ hold for any $\alpha \in \mathbf{SIso} - \{ u_{\lambda,0}, \: u'_{\lambda,0}, d_{\lambda,0}, \: d'_{\lambda,0} \}$.

\noindent {\bf Step 3.} Finish the proof by eliminating the $3$-cells in the same order as in the proof of termination of $\mathbf{SIso}$.
\end{proof}

Since $E_2 = \mathbf{SIso}_2 \cup \{ \hackcenter{
\begin{tikzpicture}[scale=0.45]
    \draw[thick, ->] (0,0) .. controls (0,.5) and (.75,.5) .. (.75,1.0);
    \draw[thick, ->] (.75,0) .. controls (.75,.5) and (0,.5) .. (0,1.0);
    %\node at (1.1,.65) { $\lambda$};
%    \node at (-.2,.1) {\tiny $i$};
%    \node at (.95,.1) {\tiny $j$};
\end{tikzpicture}} \}$,
additional indexed critical branchings appear in $E$ between $i_\lambda^1$ and $i_\lambda^4$ of the form
   \begin{equation}
   \label{eq:newindexations}
    \hackcenter{
\begin{tikzpicture} [scale =0.5, thick]
\node at (1.5,1.5) {$\bullet$};
\draw  (3.5,-.25) to (3.5,.75) .. controls ++(0,.3) and  ++(0,-.3) .. (2.5,1.5) ;
\draw[->]   (2.5,.75) .. controls ++(0,.3) and  ++(0,-.3) .. (3.5,1.5) to (3.5,2.5) ;
\draw[->] (2.5,1.5) .. controls ++(0,.65) and  ++(0,.65) ..(1.5,1.5) to (1.5,.75);
\draw (1.5,.75) .. controls ++(0,-.65) and  ++(0,-.65) ..(2.5,.75);
    \node at (4,1.4) {$\scs \l$};
\end{tikzpicture}}
\end{equation}
that are not confluent. However, we still have the following.

\begin{lemma}
\label{lem:confluenceoutsideofselfintersection}
Any $2$-cell $u$ that does not contain a strand that self-intersects admits a unique decomposition into monomials in normal form with respect to $E$.
\end{lemma}

\begin{proof}
Let $u$ be a $2$-cell that does not contain a self-intersecting strand, that is up to application of $yb$ that does not contain any element of the form \eqref{eq:newindexations}. Since $E$ is terminating and left-monomial, $u$ admits at least a linear decomposition into monomials in normal form with respect to $E$. If two such decompositions exist, then the two reductions leading to these results give a branching, that is either a Peiffer branching or come from a critical branching in a context. However, since $u$ does not contain a self-intersection, this critical branching is not given by a crossing indexation as in \eqref{eq:newindexations}, and thus from confluence of critical branchings of $\mathbf{SIso}$ and $\{ yb, dc \}$, there exists a confluence of that branching, so that these two decompositions are equal.
\end{proof}

Lemma~\ref{lem:confluenceoutsideofselfintersection} is enough to get the hom-basis of $\oUcat$ since the $3$-cells $A_\lambda$, $B_\lambda$, $C_\lambda$ and $D_\lambda$ in ${}_E R^s$ can be used to remove all self-intersections, so that any quasi-normal form with respect to ${}_E R$ will admit a unique normal form with respect to $E$.

% ---------------------------------
\subsection{Quasi termination of ${}_E R$}
%----------------------------------
\label{subsec:quasi-termination}
In this subsection, we will prove that the $(3,2)$-superpolygraph $R$ is terminating without bubble slide and cyclicity $3$-cells, and quasi-terminating with these $3$-cells. We also give a procedure showing that ${}_E R$ is quasi-terminating with rewriting cycles being induced by bubble slide cycles  as in \cite{AL16},   isotopy cycles created by dots moving on cups and caps, and cyclicity for crossings.

% - - - - - - - - - - - - - - - -
\subsubsection{Termination without bubble slide and cyclicity 3-cells}
% - - - - - - - - - - - - - - - -

\begin{lemma} \label{lem:R'-terminates}
The (3,2)-superpolygraph $R':=R-\{s_{\lambda}^{+},s_{\lambda}^{-},r_{\l}^+, r_{\l}^-, P_\lambda, P'_\lambda, Q_\lambda, Q'_\lambda \}$ terminates.
\end{lemma}
 
\begin{notation}
  For a 3-cell $\alpha$, define $d(t_2(\alpha)):=max\{d(h)\mid h\in\text{Supp}(t_2(\alpha))\}$ and similarly $X(t_2(\alpha)):=max\{X(h)\mid h\in\text{Supp}(t_2(\alpha))\}$.
\end{notation}

\begin{proof}
We prove termination in steps.

\noindent {\bf Step 1.}  First, consider the derivation $d$ into the trivial $U(R')_2^*$-module $M_{*,*,\Z}$ given by
  \[d(u )=||u||_{ \begin{tikzpicture}
    \begin{scope} [ x = 10pt, y = 10pt, join = round, cap = round, thick, scale =1.1] \draw[<-] (0.00,0.75)--(0.00,0.50) ;
    \draw[<-] (1.00,0.75)--(1.00,0.50) ; \draw (0.00,0.50)--(1.00,0.00) (1.00,0.50)--(0.00,0.00) ; \draw (0.00,0.00)--(0.00,-0.25) (1.00,0.00)--(1.00,-0.25) ;
    \end{scope}
    \end{tikzpicture}}\]
    for any $2$-cell $u$ of $R_2^s$.
  Then $d(s_2(\alpha)) > d(t_2(\alpha))$ for $\alpha \in \{ A_\lambda, B_\lambda, C_\lambda, D_\lambda, E_\lambda, F_\lambda\}$ and $d(s_2(\alpha)) \geq d(t_2(\alpha))$ for all other $\alpha$ in $R'_3$. Thus, termination of $R'$ is reduced to termination of
%$R''$, where $R''$ is the $(3,2)$-superpolygraph
\[ R'':=(R'_0,R'_1,R'_2,R'_3 - \{A_\lambda, B_\lambda, C_\lambda, D_\lambda, E_\lambda, F_\lambda\})= (R'_0,R'_1,R'_2, R''_3=\{ on_1, on_2, \Gamma, b_\lambda^{n,0}, b_\lambda^1, c_\lambda^{n,0}, c_\lambda^1, ig_{2n}\}). \]

\noindent {\bf Step 2.} Consider the 2-functor $X:U(R'')_2^* \to \cat{Ord}$ and derivation $d: U(R'')_2^* \to \Z$ defined by extending the second derivation used for $\mathbf{ONH}$ as follows:
  \[ X \big( \identu{} \big) (n)= X \big( \identd{} \big) (n) = n,
   \quad X \big( \identdotu{} \big) (n)= n,
   \quad  X \big( \raisebox{-1mm}{$\scalebox{0.9}{ \begin{tikzpicture}
    \begin{scope} [ x = 10pt, y = 10pt, join = round, cap = round, thick, scale =1.1] \draw[<-] (0.00,0.75)--(0.00,0.50) ;
    \draw[<-] (1.00,0.75)--(1.00,0.50) ; \draw (0.00,0.50)--(1.00,0.00) (1.00,0.50)--(0.00,0.00) ; \draw (0.00,0.00)--(0.00,-0.25) (1.00,0.00)--(1.00,-0.25) ;
    \end{scope}
    \end{tikzpicture}}$} \big) (n,m) = (m+2,n+1), \vspace{-0.15in}
    \]
    \[
    X \big( \ddott{} \big) (n)= n+1,
    \quad X \big( \raisebox{-1mm}{$\scalebox{0.9}{ \begin{tikzpicture}
    \begin{scope} [ x = 10pt, y = 10pt, join = round, cap = round, thick, scale =1.1]
    \draw (0.00,0.75)--(0.00,0.50) (1.00,0.75)--(1.00,0.50) ;
    \draw (0.00,0.50)--(1.00,0.00) (1.00,0.50)--(0.00,0.00) ;
    \draw (0.00,0.00)--(0.00,-0.25) (1.00,0.00)--(1.00,-0.25) ;
    \draw[->] (0.00,-0.25)--(0.00,-0.50) ;
    \draw[->] (1.00,-0.25)--(1.00,-0.50) ;
    \end{scope}
    \end{tikzpicture}}$} \big) (n,m) = (m,n),
     \quad X \big( \cupl{}{} \big)= X\big( \cupr{}{} \big)=(0,0) \vspace{-0.05in}\]
  \[ d \big( \ident{} \big) (n) = 0, \hspace{.25cm}d\big( \identdotu{} \big)(n) = n,  \hspace{.25cm}
  d \big( \raisebox{-1mm}{$\scalebox{0.9}{ \begin{tikzpicture}
    \begin{scope} [ x = 10pt, y = 10pt, join = round, cap = round, thick, scale =1.1] \draw[<-] (0.00,0.75)--(0.00,0.50) ;
    \draw[<-] (1.00,0.75)--(1.00,0.50) ; \draw (0.00,0.50)--(1.00,0.00) (1.00,0.50)--(0.00,0.00) ; \draw (0.00,0.00)--(0.00,-0.25) (1.00,0.00)--(1.00,-0.25) ;
    \end{scope}
    \end{tikzpicture}}$} \big) (n,m) = n, \hspace{.25cm} d \big( \ddott{} \big) (n)= n, \hspace{.25cm}
     d \big( \raisebox{-1mm}{$\scalebox{0.9}{ \begin{tikzpicture}
    \begin{scope} [ x = 10pt, y = 10pt, join = round, cap = round, thick, scale =1.1]
    \draw (0.00,0.75)--(0.00,0.50) (1.00,0.75)--(1.00,0.50) ;
    \draw (0.00,0.50)--(1.00,0.00) (1.00,0.50)--(0.00,0.00) ;
    \draw (0.00,0.00)--(0.00,-0.25) (1.00,0.00)--(1.00,-0.25) ;
    \draw[->] (0.00,-0.25)--(0.00,-0.50) ;
    \draw[->] (1.00,-0.25)--(1.00,-0.50) ;
    \end{scope}
    \end{tikzpicture}}$} \big) (n,m) = n+m \vspace{-0.15in} \]
  \[ d \big( \cupl{}{} \big)= d\big( \cupr{}{} \big)=0 \hspace{.5cm} d \big( \capl{}{} \big)(n,m)= d\big( \cupr{}{} \big)(n,m)=0. \]
  Then we have $X(s_2(\alpha)) \geq X(t_2(\alpha))$ and $d(s_2(\alpha)) \geq d(t_2(\alpha))$ for all 3-cells $\alpha \in R_3''$.
  Furthermore, $d(s_2(\alpha)) > d(t_2(\alpha))$ for $\alpha \in \{on_1, on_2, \Gamma\}$.
%  \begin{align*}
%    d(s_2(on_1))(n,m)=2n+1>2n=d(t_2(on_1))(n,m) \\
%    d(s_2(on_2))(n,m)=n+m+2>n+m=d(t_2(on_2))(n,m) \\
%    d(s_2(\Gamma_\lambda))(n,m,k)=n+k+4>n+k=d(t_2(\Gamma_\lambda))(n,m,k)
%  \end{align*}
  Thus, termination of $R''$ is reduced to termination of the $(3,2)$-superpolygraph
  \[\check{R}:=(R'_0,R'_1,R'_2, \{ b_\lambda^{n,0}, b_\lambda^1, c_\lambda^{n,0}, c_\lambda^1, ig_{2n}\})\]

\noindent {\bf Step 3.} To prove termination of $\check{R}$, we use a context stable map as in Section~\ref{sec:context-stable}.
  For any $u\in U(\check{R})_2^*$, define a map $d': U(\check{R})_2^* \to \N$ by
\[d'(u):=\text{number of bubbles in $u$}+\sum_{ \text{$\pi$ clockwise bubble in $u$}} |deg(\pi)|
\]
where the sum is over all clockwise bubbles appearing in $u$ and $|deg(\pi)|$ denotes the absolute value of the $\mathbb{Z}$-grading   defined in Definition~\ref{def:oddU}.

  %
%  Then $d'$ is stable under contexts and satisfies $d'(s_2(f))>d'(t_2(f))$ for all $f\in \check{R}$.

 % For $P_\lambda$ and any context $c$ of $U(\check{R})_2^*$ such that $c[P_\lambda]$ is defined, we can see that $c[s_2(P_\lambda)]$ and $c[t_2(P_\lambda)]$ have the same bubble data, so we have that $d'(c[s_2(P_\l)])=d'(c[t_2(P_\l)])$.
%
  For $\alpha \in \check{R}_3$, we have $d'(s_2(\alpha))>d'(t_2(\alpha))$.
  Let $c$ be any context of $U(\check{R})_2^*$ such that $c[\alpha]$ is defined.
  Then we have $d'(c[s_2(\alpha)])=d'(s_2(\alpha))+d'(c[1_{1_\l}])>d'(t_2(\alpha))+d'(c[Id_{\1_\lambda}])=d'(c[t_2(\alpha)])$ since both $s_2(\alpha)$ and $t_2(\alpha)$ are endomorphism 2-cells on the identity 1-cell $\1_\lambda$.
%  Therefore, we can use this map to reduce termination of $\check{R}$ to termination of $R''':=(R'_0,R'_1,R'_2,\{P_\l \})$l
%
 % \item  To prove termination of $R'''$, consider the derivation $d$ into the trivial $U(R''')_2^\ast$-module $M_{\ast, \ast, \Z}$ given by
%  \[ d(u) = ||u||_{\cuplAs{}} \]
%  Then we have $d(s_2(P_\lambda)) = 2 > 0 = d(t_2(P_\lambda))$.
%  This tells us that $R'''$ terminates, and since we have shown that $R'$ terminates if $R'''$ terminates, we conclude that $R'$ must be terminating.
Therefore, $\check{R}$ terminates, implying $R'$ also terminates.
\end{proof}

%-----------------------------
\subsubsection{Indexed cycles}
%-----------------------------
The super~$(3,2)$-polygraph $R' = R-\{s_\lambda^{\pm},r_\lambda^{\pm}, P_\lambda, P'_\lambda, Q_\lambda, Q'_\lambda \}$  terminates by Lemma~\ref{lem:R'-terminates}. However,  ${}_E R'$, and thus ${}_E R$ do  not. Closing off crossing diagrams with caps and cups can create cycles where a dot slides around a closed strand and arrives back where it started as in the configurations:
%Indeed, the isotopy relations making dots move on caps and cups create rewriting cycles that come from closing a diagram involving crossings with caps and cups. Indeed, one might be able to make a dot move through all the crossings of a strand, make it move on top of the diagram using non-oriented dot move isotopy relations, and make it move through all crossings of the new strand again, creating some cycle. For instance, one checks that there are rewriting cycles of the form
\begin{equation}
\label{E:IndexedCyclesIso}
\isocbaka{}{} \quad \isocbakb{}{}
\end{equation}
for $k > 0$  even and $l \geq 1$  odd, where the label $n$ stands for a $\star_1$-composition of $n$ crossings.
By successive application of $on_1$ and $on_2$, these give a rewriting cycle in ${}_E R$. However, for $k$ being even and $l \ne 1$ they do not have to be taken into account since the whole diagram will become $0$ when taking the normal form with respect to $E$. The case $l=1$ gives a rewriting cycle as follows:
 \begin{align}
\label{eq:AA1}
 \hackcenter{\begin{tikzpicture}[scale=0.65]
    \draw[thick,->] (0,0) .. controls ++(0,.5) and ++(0,-.5) .. (.75,1);
    \draw[thick,->] (.75,0) .. controls ++(0,.5) and ++(0,-.5) .. (0,1);
    % LEFT LOOP
    \draw[thick] (0,1) .. controls ++(0,.5) and ++(0,.5) .. (-.75,1) to (-.75,0);
    \draw[thick] (-.75,0) .. controls ++(0,-.5) and ++(0,-.5) .. (0,0);
    %% RIGHT LOOP
    \draw[thick] (.75,1) to (.75,1.25) .. controls ++(0,.5) and ++(0,.5) .. (1.5,1.35) to (1.5,-.35)
                               .. controls ++(0,-.5) and ++(0,-.5) .. (.75,-.35) to (.75,0);
                               %% DOT
    \filldraw[thick]  (-.75,1) circle (2.25pt);
     \node at (.15,1.75) { $\lambda $};
\end{tikzpicture}}
& \; \xy (-3.5,0)*{\;}="1"; (3.5,0)*{\;}="2";  {\ar@3^{({i_\lambda^4.s_{\lambda,1}})_\ast}  "1"; "2"  }; \endxy
(-1)^{\lambda}\;
\hackcenter{\begin{tikzpicture}[scale=0.65]
    \draw[thick,->] (0,0) .. controls ++(0,.5) and ++(0,-.5) .. (.75,1);
    \draw[thick] (.75,0) .. controls ++(0,.5) and ++(0,-.5) .. (0,1);
    % LEFT LOOP
    \draw[thick,directed=0.6] (0,1) .. controls ++(0,.5) and ++(0,.5) .. (-.75,1) to (-.75,0);
    \draw[thick] (-.75,0) .. controls ++(0,-.5) and ++(0,-.5) .. (0,0);
    %% RIGHT LOOP
    \draw[thick] (.75,1) to (.75,1.25) .. controls ++(0,.5) and ++(0,.5) .. (1.5,1.35) to (1.5,-.35)
                               .. controls ++(0,-.5) and ++(0,-.5) .. (.75,-.35) to (.75,0);
                               %% DOT
    \filldraw[thick]  (0,1) circle (2.25pt);
     \node at (.15,1.75) { $\lambda $};
\end{tikzpicture}}
\; - \; 2 (-1)^{\l}\;
\hackcenter{\begin{tikzpicture}[scale=0.65]
    \draw[thick,->] (0,0) .. controls ++(0,.5) and ++(0,-.5) .. (.75,1);
    \draw[thick,->] (.75,0) .. controls ++(0,.5) and ++(0,-.5) .. (0,1);
    % LEFT LOOP
    \draw[thick,] (0,1) .. controls ++(0,.5) and ++(0,.5) .. (-.75,1) to (-.75,0);
    \draw[thick] (-.75,0) .. controls ++(0,-.5) and ++(0,-.5) .. (0,0);
    %% RIGHT LOOP
    \draw[thick] (.75,1) to (.75,1.25) .. controls ++(0,.5) and ++(0,.5) .. (1.5,1.35) to (1.5,-.35)
                               .. controls ++(0,-.5) and ++(0,-.5) .. (.75,-.35) to (.75,0);
                               %% DOT
   % \filldraw[thick]  (0,1) circle (2.25pt);
     %\node at (2.15,1.75) { $\lambda $};
     \node at (-.35,2) { $\bigotimes$};
\end{tikzpicture}}
\\
\nn
&
\; \xy (-3.5,0)*{\;}="1"; (3.5,0)*{\;}="2";  {\ar@3^{on_1}  "1"; "2"  }; \endxy \;
(-1)^{\lambda + 1} \left(
\hackcenter{\begin{tikzpicture}[scale=0.65]
    \draw[thick,->] (0,0) .. controls ++(0,.5) and ++(0,-.5) .. (.75,1);
    \draw[thick,->] (.75,0) .. controls ++(0,.5) and ++(0,-.5) .. (0,1);
    % LEFT LOOP
    \draw[thick] (0,1) .. controls ++(0,.5) and ++(0,.5) .. (-.75,1) to (-.75,0);
    \draw[thick] (-.75,0) .. controls ++(0,-.5) and ++(0,-.5) .. (0,0);
    %% RIGHT LOOP
    \draw[thick] (.75,1) to (.75,1.25) .. controls ++(0,.5) and ++(0,.5) .. (1.5,1.35) to (1.5,-.35)
                               .. controls ++(0,-.5) and ++(0,-.5) .. (.75,-.35) to (.75,0);
                               %% DOT
    \filldraw[thick]  (.75,0) circle (2.25pt);
     \node at (.15,1.75) { $\lambda $};
\end{tikzpicture}}
\; - \;   \;
 \hackcenter{\begin{tikzpicture}[scale=0.65]
    \draw[thick,->] (0,0) to (0,1);
    \draw[thick,->] (.75,0) to (.75,1);
    % LEFT LOOP
    \draw[thick] (0,1) .. controls ++(0,.5) and ++(0,.5) .. (-.75,1) to (-.75,0);
    \draw[thick] (-.75,0) .. controls ++(0,-.5) and ++(0,-.5) .. (0,0);
    %% RIGHT LOOP
    \draw[thick] (.75,1) to (.75,1.25) .. controls ++(0,.5) and ++(0,.5) .. (1.5,1.35);
     \draw[thick]  (1.5,1.35) to (1.5,-.35).. controls ++(0,-.5) and ++(0,-.5) .. (.75,-.35) to (.75,0);
                               %% DOT
    %\filldraw[thick]  (.75,.95) circle (2.25pt);
     \node at (.15,1.75) { $\lambda $};
\end{tikzpicture}}
\; + \; 2  \;
\hackcenter{\begin{tikzpicture}[scale=0.65]
    \draw[thick,->] (0,0) .. controls ++(0,.5) and ++(0,-.5) .. (.75,1);
    \draw[thick,->] (.75,0) .. controls ++(0,.5) and ++(0,-.5) .. (0,1);
    % LEFT LOOP
    \draw[thick,] (0,1) .. controls ++(0,.5) and ++(0,.5) .. (-.75,1) to (-.75,0);
    \draw[thick] (-.75,0) .. controls ++(0,-.5) and ++(0,-.5) .. (0,0);
    %% RIGHT LOOP
    \draw[thick] (.75,1) to (.75,1.25) .. controls ++(0,.5) and ++(0,.5) .. (1.5,1.35) to (1.5,-.35)
                               .. controls ++(0,-.5) and ++(0,-.5) .. (.75,-.35) to (.75,0);
                               %% DOT
   % \filldraw[thick]  (0,1) circle (2.25pt);
    % \node at (2.15,1.75) { $\lambda $};
     \node at (-.35,2) { $\bigotimes$};
\end{tikzpicture}}
\right)
\end{align}
further sliding the dot term produces
\begin{align}
&(-1)^{\lambda + 1}
\hackcenter{\begin{tikzpicture}[scale=0.65]
    \draw[thick,->] (0,0) .. controls ++(0,.5) and ++(0,-.5) .. (.75,1);
    \draw[thick,->] (.75,0) .. controls ++(0,.5) and ++(0,-.5) .. (0,1);
    % LEFT LOOP
    \draw[thick] (0,1) .. controls ++(0,.5) and ++(0,.5) .. (-.75,1) to (-.75,0);
    \draw[thick] (-.75,0) .. controls ++(0,-.5) and ++(0,-.5) .. (0,0);
    %% RIGHT LOOP
    \draw[thick] (.75,1) to (.75,1.25) .. controls ++(0,.5) and ++(0,.5) .. (1.5,1.35) to (1.5,-.35)
                               .. controls ++(0,-.5) and ++(0,-.5) .. (.75,-.35) to (.75,0);
                               %% DOT
    \filldraw[thick]  (.75,0) circle (2.25pt);
     \node at (.15,1.75) { $\lambda $};
\end{tikzpicture}}
\; \xy (-3.5,0)*{\;}="1"; (3.5,0)*{\;}="2";  {\ar@{=}^{SInt}  "1"; "2"  }; \endxy \;
(-1)^{\lambda + 1}\;
\hackcenter{\begin{tikzpicture}[scale=0.65]
    \draw[thick,->] (0,0) .. controls ++(0,.5) and ++(0,-.5) .. (.75,1);
    \draw[thick,->] (.75,0) .. controls ++(0,.5) and ++(0,-.5) .. (0,1);
    % LEFT LOOP
    \draw[thick] (0,1) .. controls ++(0,.5) and ++(0,.5) .. (-.75,1) to (-.75,0);
    \draw[thick] (-.75,0) .. controls ++(0,-.5) and ++(0,-.5) .. (0,0);
    %% RIGHT LOOP
    \draw[thick] (.75,1) to (.75,1.25) .. controls ++(0,.5) and ++(0,.5) .. (1.5,1.35) to (1.5,-.35)
                               .. controls ++(0,-.5) and ++(0,-.5) .. (.75,-.35) to (.75,0);
                               %% DOT
    \filldraw[thick]  (.75,-.4) circle (2.25pt);
     \node at (.15,1.75) { $\lambda $};
\end{tikzpicture}}
\; \xy (-3.5,0)*{\;}="1"; (3.5,0)*{\;}="2";  {\ar@3^{i_\lambda^3}  "1"; "2"  }; \endxy \;
-\;
\hackcenter{\begin{tikzpicture}[scale=0.65]
    \draw[thick,->] (0,0) .. controls ++(0,.5) and ++(0,-.5) .. (.75,1);
    \draw[thick,->] (.75,0) .. controls ++(0,.5) and ++(0,-.5) .. (0,1);
    % LEFT LOOP
    \draw[thick] (0,1) .. controls ++(0,.5) and ++(0,.5) .. (-.75,1) to (-.75,0);
    \draw[thick] (-.75,0) .. controls ++(0,-.5) and ++(0,-.5) .. (0,0);
    %% RIGHT LOOP
    \draw[thick] (.75,1) to (.75,1.25) .. controls ++(0,.5) and ++(0,.5) .. (1.5,1.35) to (1.5,-.35)
                               .. controls ++(0,-.5) and ++(0,-.5) .. (.75,-.35) to (.75,0);
                               %% DOT
    \filldraw[thick]  (1.5,-.4) circle (2.25pt);
     \node at (.15,1.75) { $\lambda $};
\end{tikzpicture}}
\; + \; 2(-1)^{\l+1}\;
 \hackcenter{\begin{tikzpicture}[scale=0.7]
    \draw[thick,->] (0,0) .. controls ++(0,.5) and ++(0,-.5) .. (.75,1);
    \draw[thick,->] (.75,0) .. controls ++(0,.5) and ++(0,-.5) .. (0,1);
    % LEFT LOOP
    \draw[thick] (0,1) .. controls ++(0,.5) and ++(0,.5) .. (-.75,1) to (-.75,0);
    \draw[thick] (-.75,0) .. controls ++(0,-.5) and ++(0,-.5) .. (0,0);
    %% RIGHT LOOP
    \draw[thick] (.75,1) to (.75,1.25) .. controls ++(0,.5) and ++(0,.5) .. (1.5,1.35);
     \draw[thick ]  (1.5,1.35) to (1.5,-.35);
    \draw[thick] (1.5,-.35).. controls ++(0,-.5) and ++(0,-.5) .. (.75,-.35) to (.75,0);
                               %% DOT
    %\filldraw[thick]  (.75,-.4) circle (2.25pt);
    \node at (1.2,-1.25) { $\bigotimes$};
     \node at (.15,1.75) { $\lambda $};
\end{tikzpicture}}
\nn
\\
& \quad
\; \xy (-3.5,0)*{\;}="1"; (3.5,0)*{\;}="2";  {\ar@{=}^{SInt}  "1"; "2"  }; \endxy \;
(-1)^{\l+1}\;
\hackcenter{\begin{tikzpicture}[scale=0.65]
    \draw[thick,->] (0,0) .. controls ++(0,.5) and ++(0,-.5) .. (.75,1);
    \draw[thick,->] (.75,0) .. controls ++(0,.5) and ++(0,-.5) .. (0,1);
    % LEFT LOOP
    \draw[thick] (0,1) .. controls ++(0,.5) and ++(0,.5) .. (-.75,1) to (-.75,0);
    \draw[thick] (-.75,0) .. controls ++(0,-.5) and ++(0,-.5) .. (0,0);
    %% RIGHT LOOP
    \draw[thick] (.75,1) to (.75,1.25) .. controls ++(0,.5) and ++(0,.5) .. (1.5,1.35) to (1.5,-.35)
                               .. controls ++(0,-.5) and ++(0,-.5) .. (.75,-.35) to (.75,0);
                               %% DOT
    \filldraw[thick]  (1.5,1.4) circle (2.25pt);
     \node at (.15,1.75) { $\lambda $};
\end{tikzpicture}}
\; + \; 2(-1)^{\l}\;
 \hackcenter{\begin{tikzpicture}[scale=0.7]
    \draw[thick,->] (0,0) .. controls ++(0,.5) and ++(0,-.5) .. (.75,1);
    \draw[thick,->] (.75,0) .. controls ++(0,.5) and ++(0,-.5) .. (0,1);
    % LEFT LOOP
    \draw[thick] (0,1) .. controls ++(0,.5) and ++(0,.5) .. (-.75,1) to (-.75,0);
    \draw[thick] (-.75,0) .. controls ++(0,-.5) and ++(0,-.5) .. (0,0);
    %% RIGHT LOOP
    \draw[thick] (.75,1) to (.75,1.25) .. controls ++(0,.5) and ++(0,.5) .. (1.5,1.35);
     \draw[thick ]  (1.5,1.35) to (1.5,-.35);
    \draw[thick] (1.5,-.35).. controls ++(0,-.5) and ++(0,-.5) .. (.75,-.35) to (.75,0);
                               %% DOT
    %\filldraw[thick]  (.75,-.4) circle (2.25pt);
    \node at (2,1.85) { $\bigotimes$};
     \node at (.15,1.75) { $\lambda $};
\end{tikzpicture}}
\; \xy (-3.5,0)*{\;}="1"; (3.5,0)*{\;}="2";  {\ar@3^{(i_\lambda^2)^-}  "1"; "2"  }; \endxy \;
(-1)^{\l+1}\;
\hackcenter{\begin{tikzpicture}[scale=0.65]
    \draw[thick,->] (0,0) .. controls ++(0,.5) and ++(0,-.5) .. (.75,1);
    \draw[thick,->] (.75,0) .. controls ++(0,.5) and ++(0,-.5) .. (0,1);
    % LEFT LOOP
    \draw[thick] (0,1) .. controls ++(0,.5) and ++(0,.5) .. (-.75,1) to (-.75,0);
    \draw[thick] (-.75,0) .. controls ++(0,-.5) and ++(0,-.5) .. (0,0);
    %% RIGHT LOOP
    \draw[thick] (.75,1) to (.75,1.25) .. controls ++(0,.5) and ++(0,.5) .. (1.5,1.35) to (1.5,-.35)
                               .. controls ++(0,-.5) and ++(0,-.5) .. (.75,-.35) to (.75,0);
                               %% DOT
    \filldraw[thick]  (.75,1.4) circle (2.25pt);
     \node at (.15,1.75) { $\lambda $};
\end{tikzpicture}}
\; + \; 2(-1)^{\l}\;
 \hackcenter{\begin{tikzpicture}[scale=0.7]
    \draw[thick,->] (0,0) .. controls ++(0,.5) and ++(0,-.5) .. (.75,1);
    \draw[thick,->] (.75,0) .. controls ++(0,.5) and ++(0,-.5) .. (0,1);
    % LEFT LOOP
    \draw[thick] (0,1) .. controls ++(0,.5) and ++(0,.5) .. (-.75,1) to (-.75,0);
    \draw[thick] (-.75,0) .. controls ++(0,-.5) and ++(0,-.5) .. (0,0);
    %% RIGHT LOOP
    \draw[thick] (.75,1) to (.75,1.25) .. controls ++(0,.5) and ++(0,.5) .. (1.5,1.35);
     \draw[thick ]  (1.5,1.35) to (1.5,-.35);
    \draw[thick] (1.5,-.35).. controls ++(0,-.5) and ++(0,-.5) .. (.75,-.35) to (.75,0);
                               %% DOT
    %\filldraw[thick]  (.75,-.4) circle (2.25pt);
    \node at (-.35,1.85) { $\bigotimes$};
     \node at (.35,1.75) { $\lambda $};
\end{tikzpicture}} \nn
\end{align}
The term with the odd bubble cancels with the corresponding term in \eqref{eq:AA1}. Continuing with the dot term we have
\begin{align}
&
(-1)^{\l+1}\;
\hackcenter{\begin{tikzpicture}[scale=0.65]
    \draw[thick,->] (0,0) .. controls ++(0,.5) and ++(0,-.5) .. (.75,1);
    \draw[thick,->] (.75,0) .. controls ++(0,.5) and ++(0,-.5) .. (0,1);
    % LEFT LOOP
    \draw[thick] (0,1) .. controls ++(0,.5) and ++(0,.5) .. (-.75,1) to (-.75,0);
    \draw[thick] (-.75,0) .. controls ++(0,-.5) and ++(0,-.5) .. (0,0);
    %% RIGHT LOOP
    \draw[thick] (.75,1) to (.75,1.25) .. controls ++(0,.5) and ++(0,.5) .. (1.5,1.35) to (1.5,-.35)
                               .. controls ++(0,-.5) and ++(0,-.5) .. (.75,-.35) to (.75,0);
                               %% DOT
    \filldraw[thick]  (.75,1.4) circle (2.25pt);
     \node at (.15,1.75) { $\lambda $};
\end{tikzpicture}}
 \; \xy (-3.5,0)*{\;}="1"; (3.5,0)*{\;}="2";  {\ar@{=}^{SInt}  "1"; "2"  }; \endxy \;
 \hackcenter{\begin{tikzpicture}[scale=0.7]
    \draw[thick,] (0,0) .. controls ++(0,.5) and ++(0,-.5) .. (.75,1);
    \draw[thick,->] (.75,0) .. controls ++(0,.5) and ++(0,-.5) .. (0,1);
    % LEFT LOOP
    \draw[thick] (0,1) .. controls ++(0,.5) and ++(0,.5) .. (-.75,1) to (-.75,0);
    \draw[thick] (-.75,0) .. controls ++(0,-.5) and ++(0,-.5) .. (0,0);
    %% RIGHT LOOP
    \draw[thick] (.75,1) to (.75,1.25) .. controls ++(0,.5) and ++(0,.5) .. (1.5,1.35);
     \draw[thick,->]  (1.5,1.35) to (1.5,-.35);
    \draw[thick] (1.5,-.35).. controls ++(0,-.5) and ++(0,-.5) .. (.75,-.35) to (.75,0);
                               %% DOT
    \filldraw[thick]  (.75,.95) circle (2.25pt);
     \node at (.15,1.75) { $\lambda $};
\end{tikzpicture}}
 \; \xy (-3.5,0)*{\;}="1"; (3.5,0)*{\;}="2";  {\ar@3^{on_2}  "1"; "2"  }; \endxy \;
- \;
\hackcenter{\begin{tikzpicture}[scale=0.65]
    \draw[thick,->] (0,0) .. controls ++(0,.5) and ++(0,-.5) .. (.75,1);
    \draw[thick,->] (.75,0) .. controls ++(0,.5) and ++(0,-.5) .. (0,1);
    % LEFT LOOP
    \draw[thick] (0,1) .. controls ++(0,.5) and ++(0,.5) .. (-.75,1) to (-.75,0);
    \draw[thick] (-.75,0) .. controls ++(0,-.5) and ++(0,-.5) .. (0,0);
    %% RIGHT LOOP
    \draw[thick] (.75,1) to (.75,1.25) .. controls ++(0,.5) and ++(0,.5) .. (1.5,1.35) to (1.5,-.35)
                               .. controls ++(0,-.5) and ++(0,-.5) .. (.75,-.35) to (.75,0);
                               %% DOT
    \filldraw[thick]  (0,0) circle (2.25pt);
     \node at (.15,1.75) { $\lambda $};
\end{tikzpicture}}
\; + \;
 \hackcenter{\begin{tikzpicture}[scale=0.7]
    \draw[thick,->] (0,0) to (0,1);
    \draw[thick,->] (.75,0) to (.75,1);
    % LEFT LOOP
    \draw[thick] (0,1) .. controls ++(0,.5) and ++(0,.5) .. (-.75,1) to (-.75,0);
    \draw[thick] (-.75,0) .. controls ++(0,-.5) and ++(0,-.5) .. (0,0);
    %% RIGHT LOOP
    \draw[thick] (.75,1) to (.75,1.25) .. controls ++(0,.5) and ++(0,.5) .. (1.5,1.35);
     \draw[thick]  (1.5,1.35) to (1.5,-.35).. controls ++(0,-.5) and ++(0,-.5) .. (.75,-.35) to (.75,0);
                               %% DOT
    %\filldraw[thick]  (.75,.95) circle (2.25pt);
     \node at (.15,1.75) { $\lambda $};
\end{tikzpicture}}
\end{align}
The double bubble term combines with the corresponding term in~\eqref{eq:AA1} with coefficient $(1+(-1)^{\lambda})$, so for $\lambda$ odd these cancel. But since negative degree  bubbles vanish, this diagram is only nonzero if $\l=0$ in which the two bubbles are both multiples of the odd bubbles that squares to zero.  Sliding the remaining dot term completes the cycle.
\begin{align}
-\;
\hackcenter{\begin{tikzpicture}[scale=0.65]
    \draw[thick,->] (0,0) .. controls ++(0,.5) and ++(0,-.5) .. (.75,1);
    \draw[thick,->] (.75,0) .. controls ++(0,.5) and ++(0,-.5) .. (0,1);
    % LEFT LOOP
    \draw[thick] (0,1) .. controls ++(0,.5) and ++(0,.5) .. (-.75,1) to (-.75,0);
    \draw[thick] (-.75,0) .. controls ++(0,-.5) and ++(0,-.5) .. (0,0);
    %% RIGHT LOOP
    \draw[thick] (.75,1) to (.75,1.25) .. controls ++(0,.5) and ++(0,.5) .. (1.5,1.35) to (1.5,-.35)
                               .. controls ++(0,-.5) and ++(0,-.5) .. (.75,-.35) to (.75,0);
                               %% DOT
    \filldraw[thick]  (0,0) circle (2.25pt);
     \node at (.15,1.75) { $\lambda $};
\end{tikzpicture}}
 \; \xy (-3.5,0)*{\;}="1"; (3.5,0)*{\;}="2";  {\ar@3^{(i_\lambda^1)^-}  "1"; "2"  }; \endxy \;
-\;
\hackcenter{\begin{tikzpicture}[scale=0.65]
    \draw[thick,->] (0,0) .. controls ++(0,.5) and ++(0,-.5) .. (.75,1);
    \draw[thick,->] (.75,0) .. controls ++(0,.5) and ++(0,-.5) .. (0,1);
    % LEFT LOOP
    \draw[thick] (0,1) .. controls ++(0,.5) and ++(0,.5) .. (-.75,1) to (-.75,0);
    \draw[thick] (-.75,0) .. controls ++(0,-.5) and ++(0,-.5) .. (0,0);
    %% RIGHT LOOP
    \draw[thick] (.75,1) to (.75,1.25) .. controls ++(0,.5) and ++(0,.5) .. (1.5,1.35) to (1.5,-.35)
                               .. controls ++(0,-.5) and ++(0,-.5) .. (.75,-.35) to (.75,0);
                               %% DOT
    \filldraw[thick]  (-.75,0) circle (2.25pt);
     \node at (.15,1.75) { $\lambda $};
\end{tikzpicture}}
 \; \xy (-3.5,0)*{\;}="1"; (3.5,0)*{\;}="2";  {\ar@{=}^{SInt}  "1"; "2"  }; \endxy \;
\hackcenter{\begin{tikzpicture}[scale=0.65]
    \draw[thick,->] (0,0) .. controls ++(0,.5) and ++(0,-.5) .. (.75,1);
    \draw[thick,->] (.75,0) .. controls ++(0,.5) and ++(0,-.5) .. (0,1);
    % LEFT LOOP
    \draw[thick] (0,1) .. controls ++(0,.5) and ++(0,.5) .. (-.75,1) to (-.75,0);
    \draw[thick] (-.75,0) .. controls ++(0,-.5) and ++(0,-.5) .. (0,0);
    %% RIGHT LOOP
    \draw[thick] (.75,1) to (.75,1.25) .. controls ++(0,.5) and ++(0,.5) .. (1.5,1.35) to (1.5,-.35)
                               .. controls ++(0,-.5) and ++(0,-.5) .. (.75,-.35) to (.75,0);
                               %% DOT
    \filldraw[thick]  (-.75,1) circle (2.25pt);
     \node at (.15,1.75) { $\lambda $};
\end{tikzpicture}}
\end{align}
This may seem like a special coincidence that the cycle completed, however the diagram that we started with vanishes unless $\l=0,-1$ using $3$-cells $C_\l$ and $B_\l$, so that an element of the form \eqref{eq:AA1} will never appear in a quasi-normal form with respect to ${}_E R$.  In general, if there are more dots inside the figure \eqref{eq:AA1} the cycle can be shown to complete more generally. In fact, simplifying a diagram of this form with additional dots leads directly to the odd infinite Grassmannian equation.
The cycles built in this way are called \emph{indexed cycles}, and are rewriting cycles proper to the context of rewriting modulo.

% - - - - - - - - - - - - - - - -
\subsubsection{Quasi reduced monomials}
% - - - - - - - - - - - - - - - -
Alleaume showed in \cite{AL16} that linear 2-categories with bubble slide relations cannot be presented by terminating polygraphs, but rather by quasi-terminating polygraphs.
  For the same reason, 2-supercategories with bubble slide relations can't be presented with terminating superpolygraphs, but rather quasi-terminating superpolygraphs. Furthermore, rewriting modulo isotopies with the existence of cyclicity $3$-cells for crossings
%the definition of downward crossings in terms of upwards crossings
imply the existence of cycles of the form:
  \begin{equation} \label{eq:cyclePlambda1}
  \hackcenter{
\begin{tikzpicture}[scale=0.8]
    \draw[thick, ->] (0,0) .. controls (0,.5) and (.75,.5) .. (.75,1.0);
    \draw[thick, ->] (.75,0) .. controls (.75,.5) and (0,.5) .. (0,1.0);
    \node at (1.1,.65) { $\lambda$};
\end{tikzpicture}}
\;\; \equiv \;\;
\hackcenter{\begin{tikzpicture}[scale=0.55]
    \draw[thick, ->] (0,0) .. controls (0,.5) and (.75,.5) .. (.75,1.0);
    \draw[thick, ->] (.75,0) .. controls (.75,.5) and (0,.5) .. (0,1.0);
    \draw[thick] (0,0) .. controls ++(0,-.4) and ++(0,-.4) .. (-.75,0) to (-.75,2);
    \draw[thick] (.75,0) .. controls ++(0,-1.2) and ++(0,-1.2) .. (-1.5,0) to (-1.55,2);
    \draw[thick, -] (.75,1.0) .. controls ++(0,.4) and ++(0,.4) .. (1.5,1.0) to (1.5,-1);
    \draw[thick, -] (0,1.0) .. controls ++(0,1.2) and ++(0,1.2) .. (2.25,1.0) to (2.25,-1);
    \draw[thick, -] (2.25,-1) .. controls ++(0,-0.4) and ++(0,-0.4) .. (3,-1) to (3,2.5);
    \draw[thick,-] (1.5,-1) .. controls ++(0,-1.2) and ++(0,-1.2) .. (3.75,-1) to (3.75,2.5);
    \draw[thick,-] (-1.55,2) .. controls ++(0,0.4) and ++(0,0.4) .. (-2.3,2) to (-2.3,-1.5);
    \draw[thick,-] (-0.75,2) .. controls ++(0,1.2) and ++(0,1.2) .. (-3,2) to (-3,-1.5);
    \node at (-.35,.75) { $\lambda$};
\end{tikzpicture}} \: \overset{P_\lambda}{\Rrightarrow} \hackcenter{\begin{tikzpicture}[xscale=-1.0, scale=0.6]
    \draw[thick, <-] (0,0) .. controls (0,.5) and (.75,.5) .. (.75,1.0);
    \draw[thick, <-] (.75,0) .. controls (.75,.5) and (0,.5) .. (0,1.0);
    \draw[thick,->] (0,0) .. controls ++(0,-.4) and ++(0,-.4) .. (-.75,0) to (-.75,2);
    \draw[thick,->] (.75,0) .. controls ++(0,-1.2) and ++(0,-1.2) .. (-1.5,0) to (-1.55,2);
    \draw[thick, -] (.75,1.0) .. controls ++(0,.4) and ++(0,.4) .. (1.5,1.0) to (1.5,-1);
    \draw[thick, -] (0,1.0) .. controls ++(0,1.2) and ++(0,1.2) .. (2.25,1.0) to (2.25,-1);
    \node at (1.2,.75) {  $\lambda$};
\end{tikzpicture}} \: \overset{Q'_\lambda}{\Rrightarrow}  \: \hackcenter{
\begin{tikzpicture}[scale=0.8]
    \draw[thick, ->] (0,0) .. controls (0,.5) and (.75,.5) .. (.75,1.0);
    \draw[thick, ->] (.75,0) .. controls (.75,.5) and (0,.5) .. (0,1.0);
    \node at (1.1,.65) { $\lambda$};
\end{tikzpicture}}
  \end{equation}
   \begin{equation} \label{eq:cyclePlambda2}
  \hackcenter{
\begin{tikzpicture}[scale=0.8]
    \draw[thick, ->] (0,0) .. controls (0,.5) and (.75,.5) .. (.75,1.0);
    \draw[thick, ->] (.75,0) .. controls (.75,.5) and (0,.5) .. (0,1.0);
    \node at (1.1,.65) { $\lambda$};
\end{tikzpicture}}
\;\; \equiv \;\; \hackcenter{\begin{tikzpicture}[xscale=-1.0, scale=0.55]
    \draw[thick, ->] (0,0) .. controls (0,.5) and (.75,.5) .. (.75,1.0);
    \draw[thick, ->] (.75,0) .. controls (.75,.5) and (0,.5) .. (0,1.0);
    \draw[thick] (0,0) .. controls ++(0,-.4) and ++(0,-.4) .. (-.75,0) to (-.75,2);
    \draw[thick] (.75,0) .. controls ++(0,-1.2) and ++(0,-1.2) .. (-1.5,0) to (-1.55,2);
    \draw[thick, ->] (.75,1.0) .. controls ++(0,.4) and ++(0,.4) .. (1.5,1.0) to (1.5,-1);
    \draw[thick, ->] (0,1.0) .. controls ++(0,1.2) and ++(0,1.2) .. (2.25,1.0) to (2.25,-1);
    \node at (1.2,.75) {$\lambda$};
    \draw[thick,-] (-1.55,2) .. controls ++(0,0.4) and ++(0,0.4) .. (-2.25,2) to (-2.25,-1.5);
    \draw[thick,-] (-0.75,2) .. controls ++(0,1.2) and ++(0,1.2) .. (-3,2) to (-3,-1.5);
    \draw[-,thick] (2.25,-1) .. controls ++(0,-0.4) and ++(0,-0.4) .. (3,-1) to (3,2.5);
    \draw[-,thick] (1.5,-1) .. controls ++(0,-1.2) and ++(0,-1.2) .. (3.75,-1) to (3.75,2.5);
\end{tikzpicture}} \: \overset{P'_\lambda}{\Rrightarrow} \: - \:  \hackcenter{\begin{tikzpicture}[scale=0.6]
    \draw[thick, <-] (0,0) .. controls (0,.5) and (.75,.5) .. (.75,1.0);
    \draw[thick, <-] (.75,0) .. controls (.75,.5) and (0,.5) .. (0,1.0);
    \draw[thick] (0,0) .. controls ++(0,-.4) and ++(0,-.4) .. (-.75,0) to (-.75,2);
    \draw[thick] (.75,0) .. controls ++(0,-1.2) and ++(0,-1.2) .. (-1.5,0) to (-1.55,2);
    \draw[thick, -] (.75,1.0) .. controls ++(0,.4) and ++(0,.4) .. (1.5,1.0) to (1.5,-1);
    \draw[thick, -] (0,1.0) .. controls ++(0,1.2) and ++(0,1.2) .. (2.25,1.0) to (2.25,-1);
    \node at (-.35,.75) {  $\lambda$};
\end{tikzpicture}} \: \overset{Q_\lambda}{\Rrightarrow} \:  \hackcenter{
\begin{tikzpicture}[scale=0.8]
    \draw[thick, ->] (0,0) .. controls (0,.5) and (.75,.5) .. (.75,1.0);
    \draw[thick, ->] (.75,0) .. controls (.75,.5) and (0,.5) .. (0,1.0);
    \node at (1.1,.65) { $\lambda$};
\end{tikzpicture}}
\end{equation}
The image of these cycles through the Chevalley involution $\omega$ give rise to similar cycles for the downward crossings. If we consider sideways crossings as defined in \eqref{eq:crossl-gen-cyc} in terms of upward crossings, we can derive their definition using downward crossing using $P_\lambda$, $P'_\lambda$, and come back to the upward version using $Q_\lambda$,  $Q'_\lambda$. As a consequence, the cyclicity $3$-cells provide cycles from any kind of crossing to itself. A monomial in $\mathcal{P}$ is \emph{quasi-reduced} if it is not $E$-equivalent to $0$ and, up to indexed cycles, it can be rewritten only using rewriting cycles generated by \eqref{eq:cyclePlambda1} and \eqref{eq:cyclePlambda2} and cycles that slide a bubble through a cap or cup.
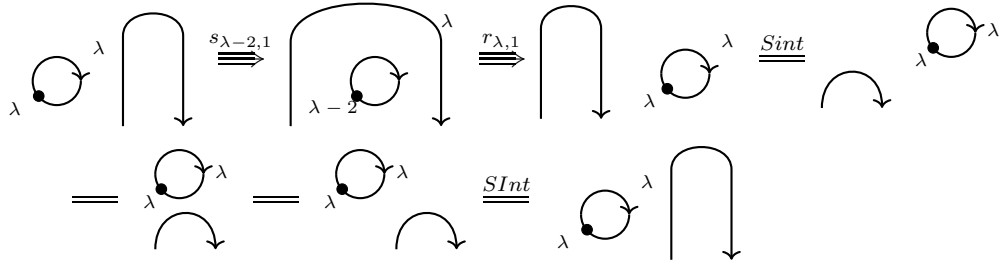
\begin{align*}
&
\hackcenter{\begin{tikzpicture}[scale=0.8]
%% BUBBLE
\draw  (0,0) arc (180:360:0.4cm) [thick];
 \draw[,<-](.8,0) arc (0:180:0.4cm) [thick];
\filldraw  [black] (.1,-.25) circle (2.5pt);
 \node at (-.3,-.45) {\tiny $\l$};
 %% CAP
 \draw[thick, ->] (1.5,-.75) to (1.5,.75).. controls ++(0,.5) and ++(0,.5) .. (2.5,.75) to (2.5,-.75);
 %%Cross
 %\draw[thick, black] (.1,-.25) to (.7,.35);
 %% LABEL
 \node at (1.1,.55) {\tiny $\lambda$};
\end{tikzpicture}}
% TRIPLE ARROW
\; \xy (-3.5,0)*{\;}="1"; (3.5,0)*{\;}="2";  {\ar@3^{s_{\l-2,1}}  "1"; "2"  }; \endxy \;
%%
%% 2nd diagram
%%
\hackcenter{\begin{tikzpicture}[scale=0.8]
%% BUBBLE
\draw  (0,0) arc (180:360:0.4cm) [thick];
 \draw[,<-](.8,0) arc (0:180:0.4cm) [thick];
\filldraw  [black] (.1,-.25) circle (2.5pt);
 \node at (-.3,-.45) {\tiny $\l -2$};
 %% CAP
 \draw[thick, <-] (1.5,-.75) to (1.5,.65).. controls ++(0,.85) and ++(0,.85) .. (-1,.65) to (-1,-.75);
 %% LABEL
 \node at (1.6,1) {\tiny $\lambda$};
\end{tikzpicture}}
% TRIPLE ARROW
\; \xy (-3.5,0)*{\;}="1"; (3.5,0)*{\;}="2";  {\ar@3^{r_{\l,1}}  "1"; "2"  }; \endxy \;
%%
%% 3rd diagram
%%
\hackcenter{\begin{tikzpicture}[scale=0.8]
%% BUBBLE
\draw  (0,0) arc (180:360:0.4cm) [thick];
 \draw[,<-](.8,0) arc (0:180:0.4cm) [thick];
\filldraw  [black] (.1,-.25) circle (2.5pt);
 \node at (-.2,-.45) {\tiny $\l$};
 %% CAP
 \draw[thick, <-] (-1,-.75) to (-1,.75).. controls ++(0,.5) and ++(0,.5) .. (-2,.75) to (-2,-.75);
 %% LABEL
 \node at (1.1,.55) {\tiny $\lambda$};
 \end{tikzpicture}}
 %
 % TRIPLE ARROW
\; \xy (-3.5,0)*{\;}="1"; (3.5,0)*{\;}="2";  {\ar@3{=}^{Sint}  "1"; "2"  }; \endxy \;
%%
%% 4th diagram
%%
\hackcenter{\begin{tikzpicture}[scale=0.8]
%% BUBBLE
\draw  (0,.5) arc (180:360:0.4cm) [thick];
 \draw[,<-](.8,.5) arc (0:180:0.4cm) [thick];
\filldraw  [black] (.1,.25) circle (2.5pt);
 \node at (-.1,.05) {\tiny $\l$};
 %% CAP
 \draw[thick, <-] (-.75,-.75)  .. controls ++(0,.8) and ++(0,.8) ..   (-1.75,-.75);
 %% LABEL
 \node at (1.1,.55) {\tiny $\lambda$};
\end{tikzpicture}}
\\
& \qquad
 % TRIPLE ARROW
\; \xy (-3.5,0)*{\;}="1"; (3.5,0)*{\;}="2";  {\ar@3{=}  "1"; "2"  }; \endxy \;
%%
%% 4th diagram
%%
\hackcenter{\begin{tikzpicture}[scale=0.8]
%% BUBBLE
\draw  (0,.5) arc (180:360:0.4cm) [thick];
 \draw[,<-](.8,.5) arc (0:180:0.4cm) [thick];
\filldraw  [black] (.1,.25) circle (2.5pt);
 \node at (-.1,.05) {\tiny $\l$};
 %% CAP
 \draw[thick,  ->] (0,-.75)  .. controls ++(0,.8) and ++(0,.8) ..   (1,-.75);
 %% LABEL
 \node at (1.1,.55) {\tiny $\lambda$};
\end{tikzpicture}}
 % TRIPLE ARROW
\; \xy (-3.5,0)*{\;}="1"; (3.5,0)*{\;}="2";  {\ar@3{=}  "1"; "2"  }; \endxy \;
%%
%% 5th diagram
%%
\hackcenter{\begin{tikzpicture}[scale=0.8]
%% BUBBLE
\draw  (0,.5) arc (180:360:0.4cm) [thick];
 \draw[,<-](.8,.5) arc (0:180:0.4cm) [thick];
\filldraw  [black] (.1,.25) circle (2.5pt);
 \node at (-.1,.05) {\tiny $\l$};
 %% CAP
 \draw[thick,  ->] (1,-.75)  .. controls ++(0,.8) and ++(0,.8) ..   (2,-.75);
 %% LABEL
 \node at (1.1,.55) {\tiny $\lambda$};
\end{tikzpicture}}
 % TRIPLE ARROW
\; \xy (-3.5,0)*{\;}="1"; (3.5,0)*{\;}="2";  {\ar@3{=}^{SInt}  "1"; "2"  }; \endxy \;
%%
%% 6th diagram
%%
\hackcenter{\begin{tikzpicture}[scale=0.8]
%% BUBBLE
\draw  (0,0) arc (180:360:0.4cm) [thick];
 \draw[,<-](.8,0) arc (0:180:0.4cm) [thick];
\filldraw  [black] (.1,-.25) circle (2.5pt);
 \node at (-.3,-.45) {\tiny $\l$};
 %% CAP
 \draw[thick, ->] (1.5,-.75) to (1.5,.75).. controls ++(0,.5) and ++(0,.5) .. (2.5,.75) to (2.5,-.75);
 %% LABEL
 \node at (1.1,.55) {\tiny $\lambda$};
\end{tikzpicture}}
\end{align*}

\begin{remark}
  No quasi-reduced monomial in $\mathcal{P}_2^s$ can be rewritten as a linear combination of other non-equivalent quasi-reduced monomials.
\end{remark}

% - - - - - - - - - - - - - - - -
\subsubsection{Weight functions and quasi-normal forms}
% - - - - - - - - - - - - - - - -
\label{subsec:weightfunctions}

\begin{definition}
  Let $C$ be a $2$-supercategory, then a \emph{weight function} on $C$ is a function $\tau \maps C_2 \to \N$ such that
  \begin{enumerate}
    \item $\tau(u \star_i v) = \tau(u) + \tau(v)$
    \item $\tau(u)=\text{max} \{ \tau(u_i) \mid u_i\in \text{Supp}(u) \}$
  \end{enumerate}
\end{definition}

When $C$ presented by $(3,2)$-superpolygraph $P$,  such a weight function is uniquely determined by its values on generating 2-cells $u$ of $P_2$.
This allows us to define a quasi-ordering $\gtrsim$ on $P_2^s$ by $u \gtrsim v$ if $\tau(u) \geq \tau(v)$.

\medskip
We define a weight function on $\mathbf{Osl(2)}_2^s$ by:
\begin{align*}
  \tau ( \cupr{}{} )=\tau(\cupl{}{})=\tau(\capr{}{})=\tau(\capl{}{})=0, \; & \tau(\udott{})=\tau(\ddott{})=0, & \tau(\crossup{}{})=\tau(\crossdn{}{})=3
\end{align*}

Then for all 3-cells $\alpha \in E_3 \backslash \{ dc \}$, we have $\tau(s_2(\alpha))=\tau (h)$ for all $h\in \text{Supp}(t_2(\alpha))$, so that all isotopy 3-cells but $dc$ preserve the weight function. In the procedure below, we  only use $dc$ from left to right, and stop the procedure whenever a $2$-cell $u$ is $0$.
Then starting with a monomial $u$ of $\mathbf{Osl(2)}_2^s$ that does not contain any negative degree bubble, and that is not $E$-equivalent to $0$:
\begin{itemize}
  \item While $u$ is not $0$ and can be rewritten with respect to $\ER$ into a 2-cell $u'$ such that $\tau(u) > \tau(u')$, then assign $u$ to $u'$.
  \item While $u$ is not $0$ and can be rewritten with respect to $\ER$ into a 2-cell $u'$ without any of the rewriting sequences in the definition of quasi-reduced monomial, namely $\Gamma_\lambda$, $on_1$, $on_2$ outside of indexed cycles, infinite Grassmannians, reduction of bubbles of degree $0$, bubble slide with a through strand, assign $u$ to $u'$.
\end{itemize}
This procedure terminates since $\gtrsim$ is well-founded, $R-\{s_\lambda^{\pm},r_\l^{\pm}, P_\lambda, P'_\lambda, Q_\lambda, Q'_\lambda \}$ is terminating by Lemma~\ref{lem:R'-terminates} and a bubble can only go through a finite number of through strands. It produces a linear combination of quasi-reduced monomials in $\mathbf{Osl(2)}_2^s$, on which one can only apply cycles generated by \eqref{eq:cyclePlambda1} and \eqref{eq:cyclePlambda2} and bubble slide through a cap or cup. Thus, $\ER$ is quasi-terminating. Moreover, we will fix a choice of preferred quasi-normal form with respect to these cycles by the following:
\begin{itemize}
\item slide the bubble outside of caps and cups, and slide them to the rightmost region of the diagram,
\item keep sideways crossings using their definition in terms of upward crossings \eqref{eq:crossl-gen-cyc}, use the cyclicity $3$-cell $P'_\lambda$ provided the number of leftward caps and cups is decreasing, and replace every downward crossing with its value in terms of upward crossings rightward caps and cups as in \eqref{eq:cyclic} using $Q'_\lambda$.
\end{itemize}

% ---------------------------------
\subsection{Confluence modulo}
%----------------------------------

In this section, we will prove that the $(3,2)$-superpolygraph modulo ${}_E R$ is confluent modulo $E$ by showing decreasing confluence of its critical branchings with respect to the quasi-normal form labelling for the quasi-normal forms fixed in Section \ref{subsec:weightfunctions}. We first start by enumerating many $3$-cells that can be derived from the generating $3$-cells of $\mathbf{Osl(2)}$, and that will be helpful for the proof of confluence of critical branchings and for the determination of the basis elements.

\subsubsection{Additional 3-cells}
\label{subsubsec:additional3cells}
From the definition of the (3,2)-superpolygraph $\mathbf{Osl(2)}$ we can derive the following 3-cells in $E^s$ or ${}_E R^s$.
We will often simplify summations involving bubbles by removing the terms involving negative degree bubbles by applying $b_\lambda^0$ or $c_\lambda^0$ to each term in a summation containing a negative bubble.  To make these types of 3-cells transparent in our notation we introduce a shorthand $b_\lambda'$ or $c_\lambda'$ to denote such application of $b_\lambda^0$ or $c_\lambda^0$.  For example,
\begin{align*}
  \sum\limits_{n=0}^\lambda
  \hackcenter{\begin{tikzpicture}[baseline=0]
	\draw[<-,thick,black] (0.3,-.45) to[out=90, in=0] (0,-0.05);
	\draw[-,thick,black] (0,-0.05) to[out = 180, in = 90] (-0.3,-.45);
    %--end of cap
  \draw[-,thick,black] (0,0.45) to[out=180,in=90] (-.2,0.25);
  \draw[<-,thick,black] (0.2,0.25) to[out=90,in=0] (0,.45);
 \draw[-,thick,black] (-.2,0.25) to[out=-90,in=180] (0,0.05);
  \draw[-,thick,black] (0,0.05) to[out=0,in=-90] (0.2,0.25);
   \node at (0.45,0) {$\scriptstyle{\lambda}$};
   \node at (-0.2,0.25) {$\color{black}\bullet$};
   \node at (-0.5,0.25) {$\color{black}\scriptstyle{n}$};
\end{tikzpicture}}
  \;\;\overset{b_\lambda'}{\Rrightarrow}\;\;
  \hackcenter{\begin{tikzpicture}[baseline=0]
	\draw[<-,thick,black] (0.3,-.45) to[out=90, in=0] (0,-0.05);
	\draw[-,thick,black] (0,-0.05) to[out = 180, in = 90] (-0.3,-.45);
    %--end of cap
  \draw[-,thick,black] (0,0.45) to[out=180,in=90] (-.2,0.25);
  \draw[<-,thick,black] (0.2,0.25) to[out=90,in=0] (0,.45);
 \draw[-,thick,black] (-.2,0.25) to[out=-90,in=180] (0,0.05);
  \draw[-,thick,black] (0,0.05) to[out=0,in=-90] (0.2,0.25);
   \node at (0.45,0) {$\scriptstyle{\lambda}$};
   \node at (-0.2,0.25) {$\color{black}\bullet$};
   \node at (-0.5,0.25) {$\color{black}\scriptstyle{\lambda -1}$};
\end{tikzpicture}}
\;\;+\;\;
\hackcenter{\begin{tikzpicture}[baseline=0]
	\draw[<-,thick,black] (0.3,-.45) to[out=90, in=0] (0,-0.05);
	\draw[-,thick,black] (0,-0.05) to[out = 180, in = 90] (-0.3,-.45);
    %--end of cap
  \draw[-,thick,black] (0,0.45) to[out=180,in=90] (-.2,0.25);
  \draw[<-,thick,black] (0.2,0.25) to[out=90,in=0] (0,.45);
 \draw[-,thick,black] (-.2,0.25) to[out=-90,in=180] (0,0.05);
  \draw[-,thick,black] (0,0.05) to[out=0,in=-90] (0.2,0.25);
   \node at (0.45,0) {$\scriptstyle{\lambda}$};
   \node at (-0.2,0.25) {$\color{black}\bullet$};
   \node at (-0.5,0.25) {$\color{black}\scriptstyle{\lambda}$};
\end{tikzpicture}}
\end{align*}
demonstrates how we will utilize this notation.

For $\lambda > 0$, define $A_{\lambda}'$ to be the 3-cell:
\[
  \tfishdrp{}{n} \overset{A_{\lambda}'}{\Rrightarrow} \left\{
    \begin{array}{ll}
       0 & \text{if $n<\lambda$} \\
        (-1)^{\lfloor\frac{\lambda+1}{2}\rfloor}\cupr{} & \text{if $n=\lambda$}
    \end{array}
    \right.
\]

We can use this to describe another 3-cell $A_\lambda ''$ for $\lambda >0$, defined by:
\[
  \tfishdrpMEinside{}{n} \overset{A_{\lambda}''}{\Rrightarrow} \left\{
    \begin{array}{ll}
       0 & \text{if $n<\lambda$} \\
        \cupr{} & \text{if $n=\lambda$}
    \end{array}
    \right.
\]

 For $\lambda > 0$, let $B_{\lambda}'$ be the 3-cell: 
$ \quad
  \tfishulpfME{}{n} \overset{B_{\lambda}'}{\Rrightarrow} \left\{
    \begin{array}{ll}
       0 & \text{if $n<\lambda$} \\
        \capl{} & \text{if $n=\lambda$}
    \end{array}
    \right.
$

For $\lambda < 0$, let $C_{\lambda}'$ be the 3-cell:
$ \quad
  \tfishurp{}{n} \overset{C_{\lambda}'}{\Rrightarrow} \left\{
    \begin{array}{ll}
       0 & \text{if $n<-\lambda$} \\
        \capr{} & \text{if $n=-\lambda$}
    \end{array}
    \right.
$

For $\lambda < 0$, let $D_{\lambda}'$ be the 3-cell:
$ \quad
  \tfishdlp{}{n} \overset{D_{\lambda}'}{\Rrightarrow} \left\{
    \begin{array}{ll}
       0 & \text{if $n<-\lambda$} \\
        \cupl{} & \text{if $n=-\lambda$}
    \end{array}
    \right.
$

As an illustration of how to derive these 3-cells, let us actually describe the process for creating $A_\lambda'$.
Given that $\lambda>0$, we slide the dots in the source of $A_\lambda'$ through all crossings possible.
\begin{align*}
  \hackcenter{\begin{tikzpicture}[baseline = 0,scale=0.7]
\draw[-,thick,black] (0.8,-0.5) to[out=-90, in=0] (0.5,-0.9);
	\draw[->,thick,black] (0.5,-0.9) to[out = 180, in = -90] (0.2,-0.5);
\draw[-,thick,black] (0.2,-.5) to (-0.3,.5);
	\draw[-,thick,black] (-0.2,-.2) to (0.2,.3);
        \draw[-,thick,black] (0.2,.3) to[out=50,in=180] (0.5,.5);
        \draw[-,thick,black] (0.5,.5) to[out=0,in=90] (0.8,-.5);
        \draw[-,thick,black] (-0.2,-.2) to[out=230,in=0] (-0.5,-.5);
        \draw[-,thick,black] (-0.5,-.5) to[out=180,in=-90] (-0.8,.5);
        \node at (1,0.5) {$\scriptstyle{\lambda}$};
         \node at (0.7,-0.8) {$\color{black}\bullet$};
   \node at (0.9,-1) {$\color{black}\scriptstyle{n}$};
        %\node at (-0.2,-0.7) {$\scriptstyle{#1}$};
\end{tikzpicture}}
\;\;\overset{SInt}{=}\;\;
(-1)^n
\hackcenter{\begin{tikzpicture}[baseline = 0,scale=0.7]
\draw[-,thick,black] (0.8,-0.5) to[out=-90, in=0] (0.5,-0.9);
	\draw[->,thick,black] (0.5,-0.9) to[out = 180, in = -90] (0.2,-0.5);
\draw[-,thick,black] (0.2,-.5) to (-0.3,.5);
	\draw[-,thick,black] (-0.2,-.2) to (0.2,.3);
        \draw[-,thick,black] (0.2,.3) to[out=50,in=180] (0.5,.5);
        \draw[-,thick,black] (0.5,.5) to[out=0,in=90] (0.8,-.5);
        \draw[-,thick,black] (-0.2,-.2) to[out=230,in=0] (-0.5,-.5);
        \draw[-,thick,black] (-0.5,-.5) to[out=180,in=-90] (-0.8,.5);
        \node at (1,0) {$\scriptstyle{\lambda}$};
         \node at (0.7,0.4) {$\color{black}\bullet$};
   \node at (0.9,.5) {$\color{black}\scriptstyle{n}$};
        %\node at (-0.2,-0.7) {$\scriptstyle{#1}$};
\end{tikzpicture}}
\;\;\overset{((i_\lambda^{2})^{\star_2 n})^-. on_2^n}{\Rrightarrow}\;\;
(-1)^{\lfloor \frac{n}{2} \rfloor}\left(
\hackcenter{\begin{tikzpicture}[baseline = 0,scale=0.7]
\draw[-,thick,black] (0.8,-0.5) to[out=-90, in=0] (0.5,-0.9);
	\draw[->,thick,black] (0.5,-0.9) to[out = 180, in = -90] (0.2,-0.5);
\draw[-,thick,black] (0.2,-.5) to (-0.3,.5);
	\draw[-,thick,black] (-0.2,-.2) to (0.2,.3);
        \draw[-,thick,black] (0.2,.3) to[out=50,in=180] (0.5,.5);
        \draw[-,thick,black] (0.5,.5) to[out=0,in=90] (0.8,-.5);
        \draw[-,thick,black] (-0.2,-.2) to[out=230,in=0] (-0.5,-.5);
        \draw[-,thick,black] (-0.5,-.5) to[out=180,in=-90] (-0.8,.5);
        \node at (1,0.5) {$\scriptstyle{\lambda}$};
         \node at (-0.2,-0.2) {$\color{black}\bullet$};
   \node at (-.4,-0.1) {$\color{black}\scriptstyle{n}$};
        %\node at (-0.2,-0.7) {$\scriptstyle{#1}$};
\end{tikzpicture}}
+
\sum\limits_{\overset{\scs r+s}{\scs =n-1}}
(-1)^{s+1}
\:\:
  \hackcenter{\begin{tikzpicture}[baseline=0,scale=1]
\begin{scope} [ x = 10pt, y = 10pt, join = round, cap = round, thick, scale=3]
  \draw[-,thick,black] (0.1,0.3) to[out=180,in=90] (-.1,0.1);
  \draw[<-,thick,black] (0.3,0.1) to[out=90,in=0] (0.1,.3);
 \draw[-,thick,black] (-.1,0.1) to[out=-90,in=180] (0.1,-0.1);
  \draw[-,thick,black] (0.1,-0.1) to[out=0,in=-90] (0.3,0.1);
  \node at (-.1,.1) {$\bullet$};
  \node at (-.3,0.1) {$ \scriptstyle{s}$};
  \node at (0.6,0.3) {$\scriptstyle{\lambda}$};
  %%%%%% of bubble part
 \draw[-,thick,black] (-.1,0.6) to[out=-90,in=180] (0.1,0.4);
  \draw[-,thick,black] (0.1,0.4) to[out=0,in=-90] (0.3,0.6);
  \draw[thick] (-.1,.6) -- (-.1,.8);
  \draw[thick,->] (.3,.6) -- (.3,.8);
\node at (.3,.6) {$\bullet$};
\node at (.5,.6) {$\scriptstyle{r}$};
 \end{scope}
\end{tikzpicture}} \right)
\end{align*}
The first term rewrites to 0 by the 3-cell
\begin{align*}
  (-1)^{\lfloor \frac{n}{2}\rfloor}
  \hackcenter{\begin{tikzpicture}[baseline = 0,scale=0.7]
\draw[-,thick,black] (0.8,-0.5) to[out=-90, in=0] (0.5,-0.9);
	\draw[->,thick,black] (0.5,-0.9) to[out = 180, in = -90] (0.2,-0.5);
\draw[-,thick,black] (0.2,-.5) to (-0.3,.5);
	\draw[-,thick,black] (-0.2,-.2) to (0.2,.3);
        \draw[-,thick,black] (0.2,.3) to[out=50,in=180] (0.5,.5);
        \draw[-,thick,black] (0.5,.5) to[out=0,in=90] (0.8,-.5);
        \draw[-,thick,black] (-0.2,-.2) to[out=230,in=0] (-0.5,-.5);
        \draw[-,thick,black] (-0.5,-.5) to[out=180,in=-90] (-0.8,.5);
        \node at (1,0.5) {$\scriptstyle{\lambda}$};
         \node at (-0.2,-0.2) {$\color{black}\bullet$};
   \node at (-.4,-0.1) {$\color{black}\scriptstyle{n}$};
        %\node at (-0.2,-0.7) {$\scriptstyle{#1}$};
\end{tikzpicture}}
\;\;\overset{((i_\l^{1})^{\star_2 n})^- \cdot A_\l}{\Rrightarrow}
0
\end{align*}
For $n<\lambda$, the second term rewrites to 0 by $b_\lambda'$
\begin{align*}
\sum\limits_{\overset{\scs r+s}{\scs =n-1}}
(-1)^{\lfloor \frac{n}{2}\rfloor+s+1}
\:\:
  \hackcenter{\begin{tikzpicture}[baseline=0,scale=1]
\begin{scope} [ x = 10pt, y = 10pt, join = round, cap = round, thick, scale=3]
  \draw[-,thick,black] (0.1,0.3) to[out=180,in=90] (-.1,0.1);
  \draw[<-,thick,black] (0.3,0.1) to[out=90,in=0] (0.1,.3);
 \draw[-,thick,black] (-.1,0.1) to[out=-90,in=180] (0.1,-0.1);
  \draw[-,thick,black] (0.1,-0.1) to[out=0,in=-90] (0.3,0.1);
  \node at (-.1,.1) {$\bullet$};
  \node at (-.3,0.1) {$ \scriptstyle{s}$};
  \node at (0.6,0.3) {$\scriptstyle{\lambda}$};
  %%%%%% of bubble part
 \draw[-,thick,black] (-.1,0.6) to[out=-90,in=180] (0.1,0.4);
  \draw[-,thick,black] (0.1,0.4) to[out=0,in=-90] (0.3,0.6);
  \draw[thick] (-.1,.6) -- (-.1,.8);
  \draw[thick,->] (.3,.6) -- (.3,.8);
\node at (.3,.6) {$\bullet$};
\node at (.5,.6) {$\scriptstyle{r}$};
 \end{scope}
\end{tikzpicture}}
\;\;\overset{b_\l'}{\Rrightarrow}\;\;
0
\end{align*}

\begin{comment}
\begin{align*}
  \sum\limits_{\overset{\scs r+s}{\scs =n-1}}
(-1)^{\lfloor \frac{n}{2}\rfloor+s+1}
\:\:
\hackcenter{
\begin{tikzpicture}[baseline=0, scale=0.7]
  \draw[thick,black,->] (0,0) to[out=90, in=180] (.5,.5);
  \draw[thick,black,-] (.5,.5) to[out=0, in=90] (1,0);
  \draw[thick,black,-] (1,0) to[out=-90, in=0] (.5,-0.5);
  \draw[thick,black,-] (.5,-.5) to[out=180, in=-90] (0,0);
  %
  \node at (.2,.35) {$\bullet$};
  \node at (.03,.45) {$\scriptstyle{s}$};
  %
  \draw[,->] (-.8,0) arc (180:360:.25)[thick];
  \draw[thick, black,-] (-.8,0) to (-.8,.5);
  \draw[thick,black,-] (-.3,0) to (-.3,.5);
  \draw[thick] (-.8,.5) -- (-.8,.75);
  \draw[thick] (-.3,.5) -- (-.3,.75);
  \node at (-.3,.15) {$\bullet$};
  \node at (-.5,.25) {$\scriptstyle{r}$};
  %
  \node at (1,.6) {$\scriptstyle{\lambda}$};
\end{tikzpicture}}
\;\;\overset{b_\l'}{\Rrightarrow}\;\;
0
\end{align*}
\end{comment}
%
For $n=\l$, only the $s=\l-1$ term remains non-zero after applying $b_\lambda'$ and we can apply $b_\lambda^1$ to this term to obtain
\begin{comment}
\begin{align*}
  (-1)^{\lfloor \frac{\l}{2}\rfloor+\l}
  \:\:\:
  \hackcenter{
\begin{tikzpicture}[baseline=0, scale=0.7]
  \draw[thick,black,->] (0,0) to[out=90, in=180] (.5,.5);
  \draw[thick,black,-] (.5,.5) to[out=0, in=90] (1,0);
  \draw[thick,black,-] (1,0) to[out=-90, in=0] (.5,-0.5);
  \draw[thick,black,-] (.5,-.5) to[out=180, in=-90] (0,0);
  %
  \node at (.23,.35) {$\bullet$};
  \node at (-.1,.5) {$\scriptstyle{\l-1}$};
  %
  \draw[,->] (-1,0) arc (180:360:.25)[thick];
  \draw[thick, black,-] (-1,0) to (-1,.5);
  \draw[thick,black,-] (-.5,0) to (-.5,.5);
  %
  \node at (1,.7) {$\scriptstyle{\lambda}$};
  %\node at (-.2,.15) {$\bullet$};
  %\node at (-.4,.25) {$\scriptstyle{}$};
\end{tikzpicture}}
\;\;\overset{b_\l^{1}}{\Rrightarrow}\;\;
(-1)^{\lfloor \frac{\l+1}{2}\rfloor}
\vcenter{\hbox{$\cupr{}$}}
\end{align*}
\end{comment}
\begin{align*}
(-1)^{\lfloor \frac{\l}{2}\rfloor+\l}
\:\:
  \hackcenter{\begin{tikzpicture}[baseline=0,scale=1]
\begin{scope} [ x = 10pt, y = 10pt, join = round, cap = round, thick, scale=3]
  \draw[-,thick,black] (0.1,0.3) to[out=180,in=90] (-.1,0.1);
  \draw[<-,thick,black] (0.3,0.1) to[out=90,in=0] (0.1,.3);
 \draw[-,thick,black] (-.1,0.1) to[out=-90,in=180] (0.1,-0.1);
  \draw[-,thick,black] (0.1,-0.1) to[out=0,in=-90] (0.3,0.1);
  \node at (-.1,.1) {$\bullet$};
  \node at (-.5, .1) {$ \scriptstyle{\lambda-1}$};
  \node at (0.6,0.4) {$\scriptstyle{\lambda}$};
  %%%%%% of bubble part
 \draw[-,thick,black] (-.1,0.6) to[out=-90,in=180] (0.1,0.4);
  \draw[-,thick,black] (0.1,0.4) to[out=0,in=-90] (0.3,0.6);
  \draw[thick] (-.1,.6) -- (-.1,.8);
  \draw[thick,->] (.3,.6) -- (.3,.8);
  %\node at (.3,.6) {$\bullet$};
  %\nodde at (.5,.6) {$\scriptstyle{r}$}
 \end{scope}
\end{tikzpicture}}
\;\;\overset{b_\l^{1}}{\Rrightarrow}\;\;
(-1)^{\lfloor \frac{\l+1}{2}\rfloor}
\vcenter{\hbox{$\cupr{}$}}
\end{align*}
Hence, for $\lambda>0$, we obtain a 3-cell $A_{\lambda}'$ given by
\begin{align*}
  \tfishdrp{}{n} \overset{A_{\lambda}'}{\Rrightarrow} \left\{
    \begin{array}{ll}
       0 & \text{if $n<\lambda$} \\
        (-1)^{\lfloor\frac{\lambda+1}{2}\rfloor}\cupr{} & \text{if $n=\lambda$}
    \end{array}
    \right.
\end{align*}

Using the bubble slide 3-cells of \ref{eq:ccbubslideLR}, we define a 3-cell $s_{\l,n}'$ that appears in some of the more complicated computations.
\begin{equation}
\sum\limits_{r\geq 0}(2r+1) \raisebox{-2mm}{$\negbubdfffsl{}{\substack{n-2r\\+\ast}}$} \identdotsu{}{2r} \: \overset{s_{\l,n}'}{\Rrightarrow} \: \identusl{} \negbubdff{}{n+ \ast}
\end{equation}

We have a $3$-cell in $E^s$ given by
\begin{equation}
 \raisebox{-7mm}{$\begin{tikzpicture}[scale=0.55,thick,black]
\draw[<-] (0.00,2.65)--(0.00,2.00) ;
\draw[<-] (1.00,2.65)--(1.00,2.00);
\draw[-] (2.00,2.25)--(2.00,1.75) ; \draw (0.00,2.00)--(1.00,1.50) (1.00,2.00)--(0.00,1.50) ;  \draw (0.00,1.50)--(0.00,1.00) (1.00,1.50)--(1.00,1.25) (2.00,1.75)--(2.00,1.25) ;  \draw (1.00,1.25)--(2.00,0.75) (2.00,1.25)--(1.00,0.75) ; \draw (0.00,1.00)--(0.00,0.50) (1.00,0.75)--(1.00,0.50) (2.00,0.75)--(2.00,0.25) ; \draw (0.00,0.50)--(1.00,0.00) (1.00,0.50)--(0.00,0.00) ;  \draw (0.00,0.00)--(0.00,-0.25) (1.00,0.00)--(1.00,-0.65) (2.00,0.25)--(2.00,-0.65) ;
%%% ROUND PART
	\draw[-,thick,black] (0,-0.25) to[out=-90, in=0] (-0.4,-0.65);
	\draw[-,thick,black] (-0.4,-0.65) to[out = 180, in = -90] (-0.8,-0.25);
	\draw[-,thick,black] (-0.8,-0.25) to (-0.8,2.65);
		\draw[<-,thick,black] (2.8,2.25) to[out=90, in=0] (2.4,2.65);
	\draw[-,thick,black] (2.4,2.65) to[out = 180, in = 90] (2,2.25);
		\draw[-,thick,black] (2.8,2.25) to (2.8,-0.65);
\end{tikzpicture}$}
\: \overset{yb}{\Rrightarrow} \:
\raisebox{-7mm}{$\begin{tikzpicture}[scale=0.55,thick,black]
\draw[<-] (0.00,2.65)--(0.00,1.50);
\draw[<-] (1.00,2.65)--(1.00,1.75) ;
\draw[-] (2.00,2.25)--(2.00,1.75) ;  \draw (1.00,1.75)--(2.00,1.25) (2.00,1.75)--(1.00,1.25) ; \draw (0.00,1.50)--(0.00,1.00) (1.00,1.25)--(1.00,1.00) (2.00,1.25)--(2.00,0.75) ; \draw (0.00,1.00)--(1.00,0.50) (1.00,1.00)--(0.00,0.50) ;  \draw (0.00,0.50)--(0.00,0.00) (1.00,0.50)--(1.00,0.25) (2.00,0.75)--(2.00,0.25) ;  \draw (1.00,0.25)--(2.00,-0.25) (2.00,0.25)--(1.00,-0.25) ; \draw (0.00,0.00)--(0.00,-0.25) (1.00,-0.25)--(1.00,-0.50) (2.00,-0.25)--(2.00,-0.50) ;
%%% ROUND PART
	\draw[-,thick,black] (0,-0.25) to[out=-90, in=0] (-0.4,-0.65);
	\draw[-,thick,black] (-0.4,-0.65) to[out = 180, in = -90] (-0.8,-0.25);
	\draw[-,thick,black] (-0.8,-0.25) to (-0.8,2.65);
		\draw[<-,thick,black] (2.8,2.25) to[out=90, in=0] (2.4,2.65);
	\draw[-,thick,black] (2.4,2.65) to[out = 180, in = 90] (2,2.25);
		\draw[-,thick,black] (2.8,2.25) to (2.8,-0.65);
\end{tikzpicture}$}
\end{equation}
which allows, up to isotopy and using sideways crossings as defined in \eqref{eq:crossl-gen-cyc}, to give an orientation for the Yang-Baxter relation for upward-upward-downward strands, corresponding to \cite[Equation (3.8)]{BE2}:
\[
\hackcenter{\begin{tikzpicture}[scale=0.8]
    \draw[thick, ->] (0,0) .. controls ++(0,1) and ++(0,-1) .. (1.5,2);
    \draw[thick, ] (.75,0) .. controls ++(0,.5) and ++(0,-.5) .. (1.5,1);
    \draw[thick, ->] (1.5,1) .. controls ++(0,.5) and ++(0,-.5) .. (0.75,2);
    \draw[thick, <-] (1.5,0) .. controls ++(0,1) and ++(0,-1) .. (0,2);
        \node at (1.75,.6) { $\lambda$};
%    \node at (-.2,.15) {\tiny $i$};
%    \node at (.95,.15) {\tiny $j$};
%    \node at (1.75,.15) {\tiny $k$};
\end{tikzpicture}}
\: \Rrightarrow
  \hackcenter{\begin{tikzpicture}[scale=0.8]
    \draw[thick, ->] (0,0) .. controls ++(0,1) and ++(0,-1) .. (1.5,2);
    \draw[thick, -] (.75,0) .. controls ++(0,.5) and ++(0,-.5) .. (0,1);
    \draw[thick, ->] (0,1) .. controls ++(0,.5) and ++(0,-.5) .. (0.75,2);
    \draw[thick, <-] (1.5,0) .. controls ++(0,1) and ++(0,-1) .. (0,2);
        \node at (1.75,.6) { $\lambda$};
%    \node at (-.2,.15) {\tiny $i$};
%    \node at (.95,.15) {\tiny $j$};
%    \node at (1.75,.15) {\tiny $k$};
\end{tikzpicture}}. \]
We actually can derive such $3$-cells either in $E^s$ using $yb$ or in ${}_E R^s$ using $\Gamma_\lambda$ to fix an orientation for all the possible configurations of Yang-Baxter $3$-cells.

Using the $3$-cell $Q'_\lambda$ to convert a downward crossing into an upward crossing with rightward caps and cups, and the odd isotopy 3-cells along with the 3-cells from superpolygraph $\mathbf{ONH}$, one can derive the following $3$-cells of ${}_E R^s$:
\[
\hspace{-0.5cm} \vcenter{\hbox{$\dbcrossingdown{}{}$}} \: \overset{dc^\lambda_-}{\Rrightarrow} 0,
\qquad
\vcenter{\hbox{$\ybrightdown{}{}{}$}} \: \overset{yb^\lambda_-}{\Rrightarrow} \: \vcenter{\hbox{$\ybleftdown{}{}{}$}},
\qquad
\begin{array}{l}
\hackcenter{\begin{tikzpicture}[scale=0.8]
    \draw[thick, <-] (0,0) .. controls ++(0,.55) and ++(0,-.5) .. (.75,1)
        node[pos=.25, shape=coordinate](DOT){};
    \draw[thick, <-] (.75,0) .. controls ++(0,.5) and ++(0,-.5) .. (0,1);
    \filldraw  (DOT) circle (2.5pt);
%    \node at (-.2,.15) {\tiny $i$};
%    \node at (.95,.15) {\tiny $j$};
\end{tikzpicture}}
\: \overset{on_{1,\lambda}^-}{\Rrightarrow} -
\hackcenter{\begin{tikzpicture}[scale=0.8]
    \draw[thick, <-] (0,0) .. controls ++(0,.55) and ++(0,-.5) .. (.75,1)
        node[pos=.75, shape=coordinate](DOT){};
    \draw[thick, <-] (.75,0) .. controls ++(0,.5) and ++(0,-.5) .. (0,1);
    \filldraw  (DOT) circle (2.5pt);
%    \node at (-.2,.15) {\tiny $i$};
%    \node at (.95,.15) {\tiny $j$};
\end{tikzpicture}}
 -\;\;
\hackcenter{\begin{tikzpicture}[scale=0.8]
    \draw[thick, <-] (0,0) to  (0,1);
    \draw[thick, <-] (.75,0)to (.75,1) ;
%    \node at (-.2,.15) {\tiny $i$};
%    \node at (.95,.15) {\tiny $i$};
\end{tikzpicture}},
\\
\hackcenter{\begin{tikzpicture}[scale=0.8]
    \draw[thick, <-] (0,0) .. controls ++(0,.55) and ++(0,-.5) .. (.75,1);
    \draw[thick, <-] (.75,0) .. controls ++(0,.5) and ++(0,-.5) .. (0,1) node[pos=.25, shape=coordinate](DOT){};
    \filldraw  (DOT) circle (2.5pt);
%    \node at (-.2,.15) {\tiny $i$};
%    \node at (.95,.15) {\tiny $j$};
\end{tikzpicture}}
    \: \overset{on_{2,\lambda}^-}{\Rrightarrow} -
\hackcenter{\begin{tikzpicture}[scale=0.8]
    \draw[thick, <-] (0,0) .. controls ++(0,.55) and ++(0,-.5) .. (.75,1);
    \draw[thick, <-] (.75,0) .. controls ++(0,.5) and ++(0,-.5) .. (0,1) node[pos=.75, shape=coordinate](DOT){};
    \filldraw  (DOT) circle (2.5pt);
%    \node at (-.2,.15) {\tiny $i$};
%    \node at (.95,.15) {\tiny $j$};
\end{tikzpicture}}
\; - \; \hackcenter{\begin{tikzpicture}[scale=0.8]
    \draw[thick, <-] (0,0) to  (0,1);
    \draw[thick, <-] (.75,0)to (.75,1) ;
%    \node at (-.2,.15) {\tiny $i$};
%    \node at (.95,.15) {\tiny $i$};
\end{tikzpicture}}
\end{array}
\]

% - - - - - - - - - - - - - - - -
\subsubsection{Critical branchings modulo of $\mathbf{Osl(2)}$}
% - - - - - - - - - - - - - - - -

We prove that ${}_E R$ is confluent modulo $E$ by showing that its critical branchings modulo are confluent and decreasing with respect to the quasi-normal form labelling for the fixed quasi-normal forms. All its critical branchings are proved confluent in Appendix \ref{appendix:criticalbranchingsfull}, and every rewriting step in these decrease the labelling to the quasi-normal form by~$1$. The classification of critical branchings modulo follows from \cite{DUP19bis}. Note that from the convergent presentation of the odd nilHecke $2$-supercategory given in Section~\ref{subsec:ConvergentPresentationONH}, all the critical branchings modulo involving two odd nilHecke $3$-cells are confluent. There is no critical branching implying the degree condition $3$-cells and infinite Grassmannians since these only reduce bubbles of positive degree by assumption, and branchings between degree condition $3$-cells and bubble slide $3$-cells are trivially confluent since the degree remains negative. There are critical branchings between infinite Grassmannians and bubble slide $3$-cells, that are proved confluent in Appendix \ref{appendix:criticalbranchingsfullOsl2}. Moreover, the critical branchings implied by $P_\lambda$ or $P'_\lambda$ with another $3$-cell given by modifying an upward crossing are trivially confluent, since there is a way to deform again the new crossing into the upward one, so that one gets back to the original $2$-cell and can apply the other $3$-cell of the branching to reach a confluence.

The remaining critical branchings are split into two families:
\begin{itemize}
\item branchings coming from the odd nilHecke $3$-cells, that is, those involving  a $3$-cell of $\mathbf{Osl(2)}$ and $on_1$ or $on_2$, and branchings that are given by applying two $3$-cells on terms that are equal modulo application of $yb$. These branchings are proved confluent in Appendix \ref{appendix:criticalbranchingsfullONH}.
\item branchings between the $3$-cells $A_\lambda, B_\lambda, C_\lambda, D_\lambda, E_\lambda, F_\lambda$ and $\Gamma_\lambda$. These ones are proved confluent modulo $E$ in Appendix \ref{appendix:criticalbranchingsfullOsl2}.
\end{itemize}

% #################################
\section{A basis theorem for odd categorified $sl(2)$} \label{sec:basis}
% #################################
 Split the (3,2)-superpolygraph $\mathbf{Osl(2)}$ into $E$ and $R$ as defined in section~\ref{sec:splittingofOsl2}. We have proved the following statement:
\begin{theorem}
The (3,2)-superpolygraph ${}_E  R$ is quasi-terminating and confluent modulo $E$.
\end{theorem}

The quasi-normal forms resulting from the (3,2)-superpolygraph %$\mathbf{Osl(2)}$
modulo  ${}_E  R$  can be described in a diagrammatic fashion.   The space 2-morphisms from $\cal{E}_{\und{\epsilon}}\1_{\l} = \cal{E}_{\epsilon_1}\dots \cal{E}_{\epsilon_k}\1_{\l}$ to $\cal{E}_{\und{\epsilon}'}\1_{\l}= \cal{E}_{\epsilon_1}\dots \cal{E}_{\epsilon_m}\1_{\l}$, when nonzero, consists of planar diagrams with $k$ points at the bottom equipped with upward/downward oriented collar neighborhoods for each $+\-$ sign $\epsilon_1, \dots , \epsilon_k$, and $m$ points at the top with collar neighborhoods determined by signs $\epsilon_1, \dots, \epsilon_m$.  These endpoints are connected by smoothly immersed directed strands whose endpoints connect the $(k+m)$ vertices compatibly with the orientation on the collar neighborhoods.
  Further,
\begin{itemize}
  \item we require that there are no triple intersections and no tangencies;
  \item no strand intersects itself, and intersects any other strand at most once;
  \item dots on a given strand appear only in a small interval near the negatively oriented endpoint of a strand connecting the vertices;
  \item all closed diagrams have been reduced to a product of non-nested dotted bubbles with a counterclockwise orientation.  Dots on bubbles are pushed to the rightmost edge of each bubble.
  \item If any three strands are such that each strand intersects the other two to create a triangle, then the triangle must be in the normal form with respect to the $(3,2)$-superpolygraph $\mathbf{SIso}$ given by one of the following:
\end{itemize}
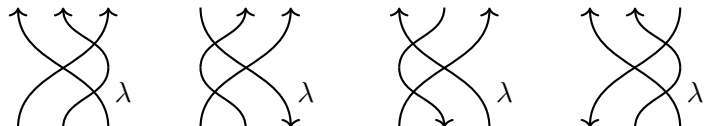
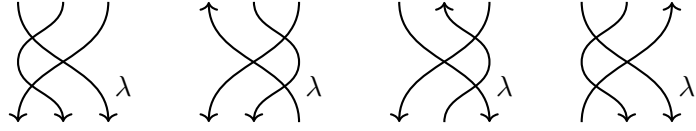
\begin{alignat*}{4}
& \hackcenter{\begin{tikzpicture}[scale=0.8]
    \draw[thick, ->] (0,0) .. controls ++(0,1) and ++(0,-1) .. (1.5,2);
    \draw[thick, ] (.75,0) .. controls ++(0,.5) and ++(0,-.5) .. (1.5,1);
    \draw[thick, ->] (1.5,1) .. controls ++(0,.5) and ++(0,-.5) .. (0.75,2);
    \draw[thick, ->] (1.5,0) .. controls ++(0,1) and ++(0,-1) .. (0,2);
        \node at (1.75,.6) { $\lambda$};
%    \node at (-.2,.15) {\tiny $i$};
%    \node at (.95,.15) {\tiny $j$};
%    \node at (1.75,.15) {\tiny $k$};
\end{tikzpicture}}
\qquad
 && \hackcenter{\begin{tikzpicture}[scale=0.8]
    \draw[thick, ->] (0,0) .. controls ++(0,1) and ++(0,-1) .. (1.5,2);
    \draw[thick, -] (.75,0) .. controls ++(0,.5) and ++(0,-.5) .. (0,1);
    \draw[thick, ->] (0,1) .. controls ++(0,.5) and ++(0,-.5) .. (0.75,2);
    \draw[thick, <-] (1.5,0) .. controls ++(0,1) and ++(0,-1) .. (0,2);
        \node at (1.75,.6) { $\lambda$};
%    \node at (-.2,.15) {\tiny $i$};
%    \node at (.95,.15) {\tiny $j$};
%    \node at (1.75,.15) {\tiny $k$};
\end{tikzpicture}}
\qquad
   &&
\hackcenter{\begin{tikzpicture}[scale=0.8]
    \draw[thick, ->] (0,0) .. controls ++(0,1) and ++(0,-1) .. (1.5,2);
    \draw[thick, <-] (.75,0) .. controls ++(0,.5) and ++(0,-.5) .. (0,1);
    \draw[thick, -] (0,1) .. controls ++(0,.5) and ++(0,-.5) .. (0.75,2);
    \draw[thick, ->] (1.5,0) .. controls ++(0,1) and ++(0,-1) .. (0,2);
        \node at (1.75,.6) { $\lambda$};
%    \node at (-.2,.15) {\tiny $i$};
%    \node at (.95,.15) {\tiny $j$};
%    \node at (1.75,.15) {\tiny $k$};
\end{tikzpicture}}
\qquad
 && \hackcenter{\begin{tikzpicture}[scale=0.8]
    \draw[thick, <-] (0,0) .. controls ++(0,1) and ++(0,-1) .. (1.5,2);
    \draw[thick, ] (.75,0) .. controls ++(0,.5) and ++(0,-.5) .. (1.5,1);
    \draw[thick, ->] (1.5,1) .. controls ++(0,.5) and ++(0,-.5) .. (0.75,2);
    \draw[thick, ->] (1.5,0) .. controls ++(0,1) and ++(0,-1) .. (0,2);
        \node at (1.75,.6) { $\lambda$};
%    \node at (-.2,.15) {\tiny $i$};
%    \node at (.95,.15) {\tiny $j$};
%    \node at (1.75,.15) {\tiny $k$};
\end{tikzpicture}}
\\
& \hackcenter{\begin{tikzpicture}[scale=0.8]
    \draw[thick, <-] (0,0) .. controls ++(0,1) and ++(0,-1) .. (1.5,2);
    \draw[thick, <-] (.75,0) .. controls ++(0,.5) and ++(0,-.5) .. (0,1);
    \draw[thick, -] (0,1) .. controls ++(0,.5) and ++(0,-.5) .. (0.75,2);
    \draw[thick, <-] (1.5,0) .. controls ++(0,1) and ++(0,-1) .. (0,2);
        \node at (1.75,.6) { $\lambda$};
%    \node at (-.2,.15) {\tiny $i$};
%    \node at (.95,.15) {\tiny $j$};
%    \node at (1.75,.15) {\tiny $k$};
\end{tikzpicture}}
\qquad
  && \hackcenter{\begin{tikzpicture}[scale=0.8]
    \draw[thick, <-] (0,0) .. controls ++(0,1) and ++(0,-1) .. (1.5,2);
    \draw[thick, <-] (.75,0) .. controls ++(0,.5) and ++(0,-.5) .. (1.5,1);
    \draw[thick, -] (1.5,1) .. controls ++(0,.5) and ++(0,-.5) .. (0.75,2);
    \draw[thick, ->] (1.5,0) .. controls ++(0,1) and ++(0,-1) .. (0,2);
        \node at (1.75,.6) { $\lambda$};
%    \node at (-.2,.15) {\tiny $i$};
%    \node at (.95,.15) {\tiny $j$};
%    \node at (1.75,.15) {\tiny $k$};
\end{tikzpicture}}
\qquad
  &&
\hackcenter{\begin{tikzpicture}[scale=0.8]
    \draw[thick, <-] (0,0) .. controls ++(0,1) and ++(0,-1) .. (1.5,2);
    \draw[thick, ] (.75,0) .. controls ++(0,.5) and ++(0,-.5) .. (1.5,1);
    \draw[thick, ->] (1.5,1) .. controls ++(0,.5) and ++(0,-.5) .. (0.75,2);
    \draw[thick, <-] (1.5,0) .. controls ++(0,1) and ++(0,-1) .. (0,2);
        \node at (1.75,.6) { $\lambda$};
%    \node at (-.2,.15) {\tiny $i$};
%    \node at (.95,.15) {\tiny $j$};
%    \node at (1.75,.15) {\tiny $k$};
\end{tikzpicture}}
\qquad
 && \hackcenter{\begin{tikzpicture}[scale=0.8]
    \draw[thick, ->] (0,0) .. controls ++(0,1) and ++(0,-1) .. (1.5,2);
    \draw[thick, <-] (.75,0) .. controls ++(0,.5) and ++(0,-.5) .. (0,1);
    \draw[thick, - ] (0,1) .. controls ++(0,.5) and ++(0,-.5) .. (0.75,2);
    \draw[thick, <-] (1.5,0) .. controls ++(0,1) and ++(0,-1) .. (0,2);
        \node at (1.75,.6) { $\lambda$};
%    \node at (-.2,.15) {\tiny $i$};
%    \node at (.95,.15) {\tiny $j$};
%    \node at (1.75,.15) {\tiny $k$};
\end{tikzpicture}}
\end{alignat*}
We can further reduce the ambiguity of our chosen basis by making a preferred choice of each super interchange class of diagram.  For example, choosing dots and crossings to decrease in height from right to left, with dots appearing above crossings when related by super interchange.

An example of the normal form of a 2-morphism from $\cal{E}_{-}\cal{E}_{+}\cal{E}_{+}\cal{E}_{-}\cal{E}_{+}\onel$ to
$\cal{E}_{+}\cal{E}_{+}\cal{E}_{-}\cal{E}_{-}\cal{E}_{+}\cal{E}_{+}\onel$ is given in the first diagram below, while the second would not be in normal form as it does not have the correct Yang-Baxter representative.
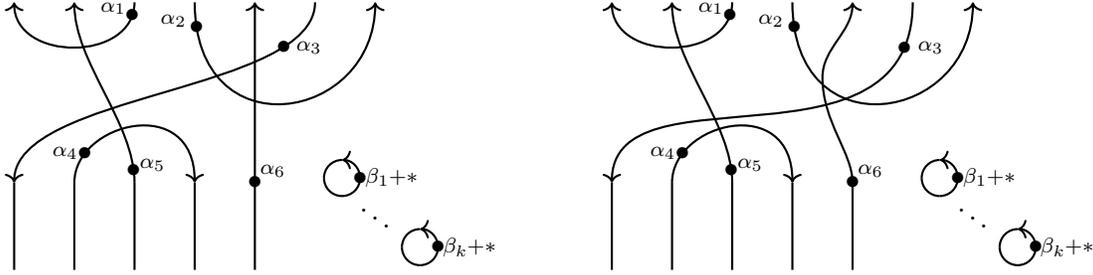
\[
\hackcenter{\begin{tikzpicture}[scale=0.8]
    \draw[thick,<-] (0,0) .. controls ++(0,1.5) and ++(0,-1.5) .. (5,3) node[pos=.8, shape=coordinate](DOT1){};
    \draw[thick,->] (1,0) .. controls ++(0,1) and ++(0,1.5) .. (3,0)  node[pos=.18, shape=coordinate](DOT2){};
    \draw[thick,<-] (0,3) .. controls ++(0,-1) and ++(0,-1) .. (2,3) node[pos=.92, shape=coordinate](DOT3){};
    \draw[thick,->] (2,0) .. controls ++(0,1) and ++(0,-1) .. (1,3) node[pos=.07, shape=coordinate](DOT4){};
   \draw[thick,->] (3,3) .. controls ++(0,-2.5) and ++(0,-2) .. (6,3) node[pos=.06, shape=coordinate](DOT5){};
    \draw[thick] (0,-1.45) -- (0,0)   (1,-1.45) -- (1,0)  (2,-1.45) -- (2,0)  (3,-1.45) -- (3,0) ;
   \draw[thick, ->] (4,-1.45) to (4,3);
            \node at (DOT1) {$\bullet$};\node at (DOT2) {$\bullet$}; \node at (DOT3) {$\bullet$};
            \node at (DOT4) {$\bullet$}; \node at (DOT5) {$\bullet$}; \node at (4,0) {$\bullet$};
    \node at (6,0) {$\negbubdfffsl{}{\beta_1+ \ast}$};
    \node at (6, -.45) {$\ddots$ };
    \node at (7.3,-1.15) {$\negbubdfffsl{}{\beta_k+ \ast}$};
    \node at (1.65, 2.9) {$\scs \alpha_1$ };
    \node at (2.65, 2.7) {$\scs \alpha_2$ };
    \node at (4.9, 2.25) {$\scs \alpha_3$ };
    \node at (.85, .5) {$\scs \alpha_4$ };
    \node at (2.3, .3) {$\scs \alpha_5$ };
   \node at (4.3, .2) {$\scs \alpha_6$ };
\end{tikzpicture}}
\qquad \quad
\hackcenter{\begin{tikzpicture}[scale=0.8]
    \draw[thick,<-] (0,0) .. controls ++(0,2) and ++(0,-3) .. (5,3) node[pos=.9, shape=coordinate](DOT1){};
    \draw[thick,->] (1,0) .. controls ++(0,1) and ++(0,1.5) .. (3,0)  node[pos=.18, shape=coordinate](DOT2){};
    \draw[thick,<-] (0,3) .. controls ++(0,-1) and ++(0,-1) .. (2,3) node[pos=.92, shape=coordinate](DOT3){};
    \draw[thick,->] (2,0) .. controls ++(0,1) and ++(0,-1) .. (1,3) node[pos=.07, shape=coordinate](DOT4){};
   \draw[thick,->] (3,3) .. controls ++(0,-2.5) and ++(0,-2) .. (6,3) node[pos=.06, shape=coordinate](DOT5){};
    \draw[thick] (0,-1.45) -- (0,0)   (1,-1.45) -- (1,0)  (2,-1.45) -- (2,0)  (3,-1.45) -- (3,0) ;
   \draw[thick, ->] (4,-1.45) to (4,0) .. controls ++(0,.5) and ++(0,-.75) .. (3.5,1.75)
                .. controls ++(0,.5) and ++(0,-.5) .. (4,3);
            \node at (DOT1) {$\bullet$};\node at (DOT2) {$\bullet$}; \node at (DOT3) {$\bullet$};
            \node at (DOT4) {$\bullet$}; \node at (DOT5) {$\bullet$}; \node at (4,0) {$\bullet$};
    \node at (6,0) {$\negbubdfffsl{}{\beta_1+ \ast}$};
    \node at (6, -.45) {$\ddots$ };
    \node at (7.3,-1.15) {$\negbubdfffsl{}{\beta_k+ \ast}$};
    \node at (1.65, 2.9) {$\scs \alpha_1$ };
    \node at (2.65, 2.7) {$\scs \alpha_2$ };
    \node at (5.3, 2.25) {$\scs \alpha_3$ };
    \node at (.85, .5) {$\scs \alpha_4$ };
    \node at (2.3, .3) {$\scs \alpha_5$ };
   \node at (4.3, .2) {$\scs \alpha_6$ };
\end{tikzpicture}}
\]
%the quasi-normal form produces a basis for the space of 2-morphisms of $\mf{U}$ proves
Hence, we have proven the non-degeneracy conjecture for the odd 2-category $\mf{U}$ from \cite[Section 8]{BE2}.

\begin{theorem}[Nondegeneracy Conjecture]
Fixing a choice of representative for each super interchange class of elements from the quasi normal form of the (3,2)-superpolygraph $\mathbf{Osl(2)}$ gives a basis for each Hom space $\Hom_{\mf{U}}(\cal{E}_{\und{\epsilon}}, \cal{E}_{\und{\epsilon'}})$.  In particular, $\Hom_{\mf{U}}(\cal{E}_{\und{\epsilon}}, \cal{E}_{\und{\epsilon'}})$ is a free right ${\sf Sym}[d]$-module with ${\sf Sym}[d]$ the bubble algebra defined in Remark~\ref{rem:Symd}.
\end{theorem}

\begin{corollary}
The conjectural classification of dg-structures on the super 2-category $\mf{U}(\mf{sl}_2)$ from \cite[Proposition 7.1]{Odd-diff} is a complete classification.
\end{corollary}

% ==============================================================================
% REFERENCES
%

\bibliographystyle{plain}
\bibliographystyle{amsalpha}
\bibliography{bib_bub}

%
% ==============================================================================

\clearpage

\csname @addtoreset\endcsname{section}{part}
\csname @addtoreset\endcsname{subsection}{part}

\appendix
\setcounter{secnumdepth}{2}

\section{Critical branchings for $\mathbf{SIso(\mathfrak{g})}$}
\label{appendix:criticalbranchingsSIso}

\subsection{Regular critical branchings}
Here we verify the critical branchings for the (3,2)-superpolygraph $\mathbf{SIso}(\mathfrak{g})$. For every $3$-cell other than $\alpha_{m,k}$ and $\beta_{m,k}$, every strand in both the source and target is labelled with $i$, so for branchings that don't use $\alpha_{m,k}$ and $\beta_{m,k}$, we often write $(-1)^{\l}$ instead of $(-1)^{\l_i}$. The classification of critical branchings is analogous to that of the 3-polygraph of pearls~\cite[Section 5.5]{GM09}, with one extra regular critical branching involving the $3$-cells $I_0$ and   $\alpha_{0,1_{1_{-2}}}$, and two extra indexed critical branchings involving the $3$-cells $\alpha_{m,k}$ and $\beta_{m,k}$, coming from the definition of the odd bubble.
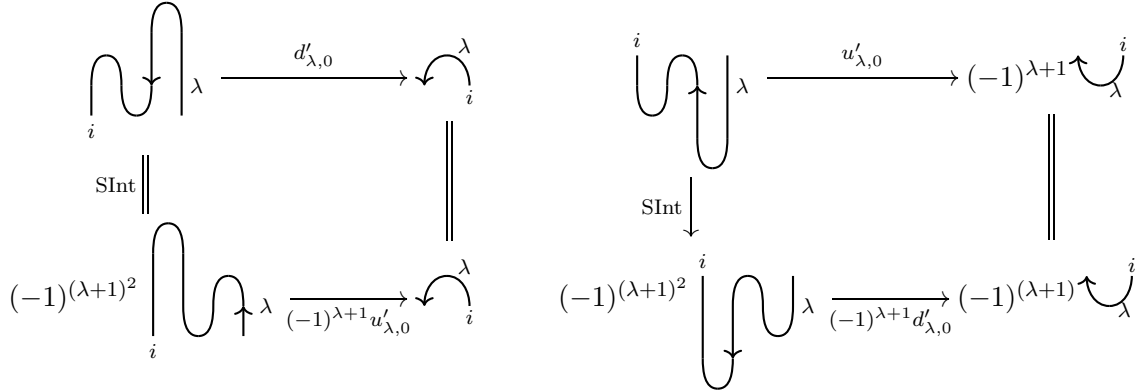
\[  \xymatrix@R=2em@C=4em{ \begin{tikzpicture}[baseline = 0]
\draw[-,thick,black] (0.7,0.7) to (0.7,-0.4);
	\draw[-,thick,black] (0.7,0.7) to[out=90, in=0] (0.5,1.1);
	\draw[-,thick,black] (0.5,1.1) to[out = 180, in =90] (0.3,0.7);
  \draw[<-,thick,black] (0.3,0) to (0.3,.7);
	\draw[-,thick,black] (0.3,0) to[out=-90, in=0] (0.1,-0.4);
	\draw[-,thick,black] (0.1,-0.4) to[out = 180, in = -90] (-0.1,0);
	\draw[-,thick,black] (-0.1,0) to[out=90, in=0] (-0.3,0.4);
	\draw[-,thick,black] (-0.3,0.4) to[out = 180, in =90] (-0.5,0);
  \draw[-,thick,black] (-0.5,0) to (-0.5,-.4);
   \node at (-0.5,-.6) {$\scriptstyle{i}$};
   \node at (0.9,0) {$\scriptstyle{\lambda}$};
\end{tikzpicture} \ar@{=}[d] _-{\text{SInt}} \ar [r] ^-{d'_{\lambda,0}} &  \capl{i}  \ar@{=} [d]
    \\
(-1)^{(\l+1)^2} \begin{tikzpicture}[baseline = 0]
    \draw[-,thick,black] (-0.1,0.7) to (-0.1,0);
	\draw[-,thick,black] (0.7,0) to[out=90, in=0] (0.5,0.4);
	\draw[-,thick,black] (0.5,0.4) to[out = 180, in =90] (0.3,0);
  \draw[<-,thick,black] (0.7,0) to (0.7,-0.4);
	\draw[-,thick,black] (0.3,0) to[out=-90, in=0] (0.1,-0.4);
	\draw[-,thick,black] (0.1,-0.4) to[out = 180, in = -90] (-0.1,0);
	\draw[-,thick,black] (-0.1,0.7) to[out=90, in=0] (-0.3,1.1);
	\draw[-,thick,black] (-0.3,1.1) to[out = 180, in =90] (-0.5,0.7);
  \draw[-,thick,black] (-0.5,0.7) to (-0.5,-.4);
   \node at (-0.5,-.6) {$\scriptstyle{i}$};
   \node at (1,0) {$\scriptstyle{\lambda}$};
\end{tikzpicture} \ar [r] _-{(-1)^{\l+1} u'_{\lambda,0}} & \capl{i} }
\qquad  \xymatrix@R=2em@C=4em{
\begin{tikzpicture}[baseline = 0]
    \draw[-,thick,black] (0.7,-0.7) to (0.7,0.4);
	\draw[-,thick,black] (0.7,-0.7) to[out=-90, in=0] (0.5,-1.1);
	\draw[-,thick,black] (0.5,-1.1) to[out = 180, in =-90] (0.3,-0.7);
  \draw[<-,thick,black] (0.3,0) to (0.3,-.7);
	\draw[-,thick,black] (0.3,0) to[out=90, in=0] (0.1,0.4);
	\draw[-,thick,black] (0.1,0.4) to[out = 180, in = 90] (-0.1,0);
	\draw[-,thick,black] (-0.1,0) to[out=-90, in=0] (-0.3,-0.4);
	\draw[-,thick,black] (-0.3,-0.4) to[out = 180, in =-90] (-0.5,0);
  \draw[-,thick,black] (-0.5,0) to (-0.5,.4);
   \node at (-0.5,.6) {$\scriptstyle{i}$};
   \node at (0.9,0) {$\scriptstyle{\lambda}$};
\end{tikzpicture}
\ar [d] _-{\text{SInt}} \ar [r] ^-{u'_{\lambda, 0}} & (-1)^{\l+1} \cupl{i} \ar@{=} [d] \\
(-1)^{(\l+1)^2} \begin{tikzpicture}[baseline = 0]
    \draw[-,thick,black] (0.7,0) to (0.7,0.4);
	\draw[-,thick,black] (0.7,0) to[out=-90, in=0] (0.5,-0.4);
	\draw[-,thick,black] (0.5,-0.4) to[out = 180, in =-90] (0.3,0);
  \draw[->,thick,black] (-0.1,0) to (-0.1,-.7);
	\draw[-,thick,black] (0.3,0) to[out=90, in=0] (0.1,0.4);
	\draw[-,thick,black] (0.1,0.4) to[out = 180, in = 90] (-0.1,0);
	\draw[-,thick,black] (-0.1,-0.7) to[out=-90, in=0] (-0.3,-1.1);
	\draw[-,thick,black] (-0.3,-1.1) to[out = 180, in =-90] (-0.5,-0.7);
  \draw[-,thick,black] (-0.5,-0.7) to (-0.5,.4);
   \node at (-0.5,.6) {$\scriptstyle{i}$};
   \node at (0.9,0) {$\scriptstyle{\lambda}$};
\end{tikzpicture} \ar [r] _-{(-1)^{\l+1} d'_{\lambda, 0}} & (-1)^{(\l+1)} \cupl{i} }
\]
since $(-1)^{(\l+1)^2} = (-1)^{\l^2 +1} = (-1)^{\l+1}$. Diagrams with reverse orientations give the same critical branchings as in the even case, since the use of superinterchange do not create any sign.
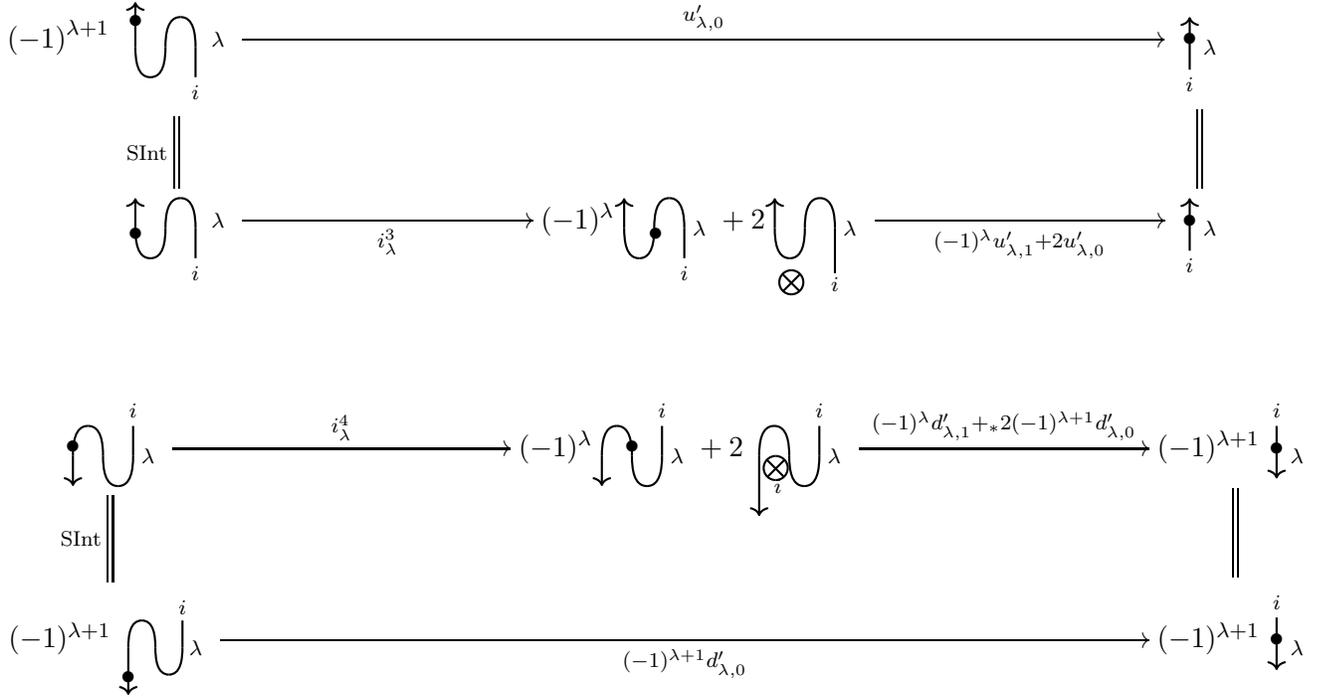
\[
(-1)^{\l+1} \xymatrix@R=2em@C=10em{ \begin{tikzpicture}[baseline = 0]
  \draw[-,thick,black] (0.3,0) to (0.3,-.4);
	\draw[-,thick,black] (0.3,0) to[out=90, in=0] (0.1,0.4);
	\draw[-,thick,black] (0.1,0.4) to[out = 180, in = 90] (-0.1,0);
	\draw[-,thick,black] (-0.1,0) to[out=-90, in=0] (-0.3,-0.4);
	\draw[-,thick,black] (-0.3,-0.4) to[out = 180, in =-90] (-0.5,0);
  \draw[->,thick,black] (-0.5,0) to (-0.5,.6);
  \node at (-0.5,0.35) {$\bullet$};
   \node at (0.3,-.6) {$\scriptstyle{i}$};
   \node at (0.6,0.1) {$\scriptstyle{\lambda}$};
\end{tikzpicture} \ar@{=} [d] _-{\text{SInt}}\ar[rr] ^-{u'_{\lambda,0}}  & &  \upd{i}  \ar@{=} [d] ^-{} \\
\begin{tikzpicture}[baseline = 0]
  \draw[-,thick,black] (0.3,0) to (0.3,-.4);
	\draw[-,thick,black] (0.3,0) to[out=90, in=0] (0.1,0.4);
	\draw[-,thick,black] (0.1,0.4) to[out = 180, in = 90] (-0.1,0);
	\draw[-,thick,black] (-0.1,0) to[out=-90, in=0] (-0.3,-0.4);
	\draw[-,thick,black] (-0.3,-0.4) to[out = 180, in =-90] (-0.5,0);
  \draw[->,thick,black] (-0.5,0) to (-0.5,.4);
  \node at (-0.5,-0.08) {$\bullet$};
   \node at (0.3,-.6) {$\scriptstyle{i}$};
   \node at (0.6,0.1) {$\scriptstyle{\lambda}$};
\end{tikzpicture} \ar [r]  _-{i_\lambda^3} & (-1)^\l \sdld{i} + 2  \begin{tikzpicture}[baseline = 0]
  \draw[-,thick,black] (0.3,0) to (0.3,-.6);
	\draw[-,thick,black] (0.3,0) to[out=90, in=0] (0.1,0.4);
	\draw[-,thick,black] (0.1,0.4) to[out = 180, in = 90] (-0.1,0);
	\draw[-,thick,black] (-0.1,0) to[out=-90, in=0] (-0.3,-0.4);
	\draw[-,thick,black] (-0.3,-0.4) to[out = 180, in =-90] (-0.5,0);
  \draw[->,thick,black] (-0.5,0) to (-0.5,.4);
   \node at (0.3,-.75) {$\scriptstyle{i}$};
   \node at (0.5,0) {$\scriptstyle{\lambda}$};
   %% Odd bubble
   \draw[-,thick,black,scale=0.8] (-0.35,-0.7) to[out=180,in=90] (-0.55,-0.9);
  \draw[-,thick,black,scale=0.8] (-0.15,-0.9) to[out=90,in=0] (-0.35,-0.7);
 \draw[-,thick,black,scale=0.8] (-0.55,-0.9) to[out=-90,in=180] (-0.35,-1.1);
  \draw[-,thick,black,scale=0.8] (-0.35,-1.1) to[out=0,in=-90] (-0.15,-0.9);
  \draw[-,thick,black,scale=0.8] (-0.21,-1.05) to (-0.49,-0.77);
  \draw[-,thick,black,scale=0.8] (-0.21,-0.77) to (-0.49,-1.05);
\end{tikzpicture} \ar [r] _-{(-1)^\l u'_{\lambda, 1} + 2 u'_{\lambda, 0}} & \upd{i}  }
\]

\[
\xymatrix@R=2em@C=10em{ \begin{tikzpicture}[baseline = 0]
  \draw[-,thick,black] (0.3,0) to (0.3,.4);
	\draw[-,thick,black] (0.3,0) to[out=-90, in=0] (0.1,-0.4);
	\draw[-,thick,black] (0.1,-0.4) to[out = 180, in = -90] (-0.1,0);
	\draw[-,thick,black] (-0.1,0) to[out=90, in=0] (-0.3,0.4);
	\draw[-,thick,black] (-0.3,0.4) to[out = 180, in =90] (-0.5,0);
  \draw[->,thick,black] (-0.5,0) to (-0.5,-.4);
  \node at (-0.5,0.12) {$\bullet$};
   \node at (0.3,.6) {$\scriptstyle{i}$};
   \node at (0.5,0) {$\scriptstyle{\lambda}$};
\end{tikzpicture} \ar@{=} [d] _-{\text{SInt}} \ar [r] ^-{i_\lambda^4}
& (-1)^\l \begin{tikzpicture}[baseline = 0]
  \draw[-,thick,black] (0.3,0) to (0.3,.4);
	\draw[-,thick,black] (0.3,0) to[out=-90, in=0] (0.1,-0.4);
	\draw[-,thick,black] (0.1,-0.4) to[out = 180, in = -90] (-0.1,0);
	\draw[-,thick,black] (-0.1,0) to[out=90, in=0] (-0.3,0.4);
	\draw[-,thick,black] (-0.3,0.4) to[out = 180, in =90] (-0.5,0);
  \draw[->,thick,black] (-0.5,0) to (-0.5,-.4);
  \node at (-0.1,0.12) {$\bullet$};
   \node at (0.3,.6) {$\scriptstyle{i}$};
   \node at (0.5,0) {$\scriptstyle{\lambda}$};
\end{tikzpicture}  + 2 \:  \begin{tikzpicture}[baseline = 0]
  \draw[-,thick,black] (0.3,0) to (0.3,.4);
	\draw[-,thick,black] (0.3,0) to[out=-90, in=0] (0.1,-0.4);
	\draw[-,thick,black] (0.1,-0.4) to[out = 180, in = -90] (-0.1,0);
	\draw[-,thick,black] (-0.1,0) to[out=90, in=0] (-0.3,0.4);
	\draw[-,thick,black] (-0.3,0.4) to[out = 180, in =90] (-0.5,0);
  \draw[->,thick,black] (-0.5,0) to (-0.5,-.8);
   \node at (0.3,.6) {$\scriptstyle{i}$};
   \node at (0.5,0) {$\scriptstyle{\lambda}$};
   % Odd bubble
     \draw[-,thick,black,scale=0.8] (-0.35,0) to[out=180,in=90] (-.55,-0.2);
  \draw[-,thick,black,scale=0.8] (-0.15,-0.2) to[out=90,in=0] (-0.35,0);
 \draw[-,thick,black,scale=0.8] (-.55,-0.2) to[out=-90,in=180] (-0.35,-0.4);
  \draw[-,thick,black,scale=0.8] (-0.35,-0.4) to[out=0,in=-90] (-0.15,-0.2);
  \draw[-,thick,black,scale=0.8] (-0.21,-0.35) to (-0.49,-0.07);
  \draw[-,thick,black,scale=0.8] (-0.21,-0.07) to (-0.49,-0.35);
 \node at (-0.25,-0.4) {$\scriptstyle{i}$};
\end{tikzpicture}
\ar [r] ^-{(-1)^\l d'_{\lambda,1} + {}_\ast 2(-1)^{\l+1} d'_{\lambda,0}} & (-1)^{\l+1} \dpd{i}  \ar@{=} [d] ^-{}
 \\
(-1)^{\l+1} \begin{tikzpicture}[baseline = 0,scale=0.9]
  \draw[-,thick,black] (0.3,0) to (0.3,.4);
	\draw[-,thick,black] (0.3,0) to[out=-90, in=0] (0.1,-0.4);
	\draw[-,thick,black] (0.1,-0.4) to[out = 180, in = -90] (-0.1,0);
	\draw[-,thick,black] (-0.1,0) to[out=90, in=0] (-0.3,0.4);
	\draw[-,thick,black] (-0.3,0.4) to[out = 180, in =90] (-0.5,0);
  \draw[->,thick,black] (-0.5,0) to (-0.5,-.7);
  \node at (-0.5,-0.45) {$\bullet$};
   \node at (0.3,.6) {$\scriptstyle{i}$};
   \node at (0.5,0) {$\scriptstyle{\lambda}$};
\end{tikzpicture} \ar[rr] _-{(-1)^{\l+1} d'_{\lambda,0}}  & & (-1)^{\l+1} \dpd{i}  } \]
where the $\ast$ symbol before the rewriting step $d'_{\lambda,0}$ means that we used a superinterchange relation, between the odd bubble and the leftward cup before applying the rewriting step, creating the sign $(-1)^{\l+1}$. Moreover, the critical branching involving $I_0$ and $\alpha_{0,1_{1_{-2}}}$ is proved confluent as follows:
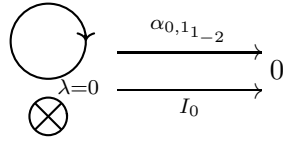
\[
\xymatrix{
\begin{tikzpicture}[baseline=0]
  \draw[thick, black,<-] (0,.5) arc (0:370:.5);
  \draw[thick, black] (-.25,-.5) arc (0:370:.25);
  \draw[thick, black] (-.32,-.32) to (-.68,-.68);
  \draw[thick, black] (-.32,-.68) to (-.68,-.32);
  \node at (-.1,-.1) {$\scriptstyle{\lambda=0}$};
\end{tikzpicture}
\ar@<-.6pc>[rr]_-{I_0} \ar@<.6pc>[rr]^-{\alpha_{0,1_{1_{-2}}}} & & 0
}
\]

% - - - - - - - - - - - - - - - -
\subsubsection{Shortened notation for critical branchings} \label{sec:shortnotation}
% - - - - - - - - - - - - - - - -
In order to avoid drawing all the critical branchings entirely, we introduced a shortened diagrammatic representation for these, encoding the minimal amount of data that we need in order to reconstruct the actual branching. We   only draw the diagrammatic source and the diagrammatic normal form (or chosen quasi-normal form) of the critical branchings, and indicate between brackets the two rewriting sequences that lead from the source to the common target. If one has to apply super-interchange relations at the source of the critical branching, we will indicate this by adding the element $\text{SInt}$ at the beginning of one of the rewriting paths. If one has to apply super-interchange relation in the middle of a rewriting path, we will indicate this by writing a symbol $\ast$ before applying the rewriting step with the correct sign brought by super interchange. Later, when rewriting modulo isotopy, we will indicate using a $3$-cell $e$ of the super-isotopy polygraph $E$ before applying a rewriting step $f$ of $R$ by $e \cdot f$.

 For example, the last critical branching above is depicted in our shorthand as follows:
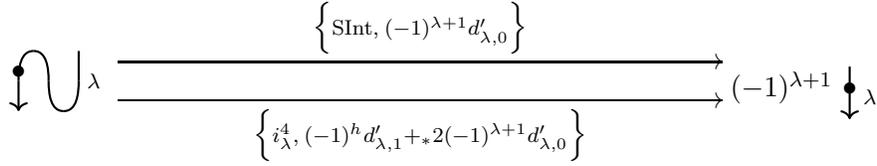
\[
\xy
  (-50,10)*+{ \begin{tikzpicture}[baseline = 0]
  \draw[-,thick,black] (0.3,0) to (0.3,.4);
	\draw[-,thick,black] (0.3,0) to[out=-90, in=0] (0.1,-0.4);
	\draw[-,thick,black] (0.1,-0.4) to[out = 180, in = -90] (-0.1,0);
	\draw[-,thick,black] (-0.1,0) to[out=90, in=0] (-0.3,0.4);
	\draw[-,thick,black] (-0.3,0.4) to[out = 180, in =90] (-0.5,0);
  \draw[->,thick,black] (-0.5,0) to (-0.5,-.4);
  \node at (-0.5,0.12) {$\bullet$};
   \node at (0.5,0) {$\scriptstyle{\lambda}$};
\end{tikzpicture}
       }="TL";
  (50,10)*+{  (-1)^{\l+1} \dpd{} }="TR";
     {\ar@<.6pc> ^-{
        \Big\{  \text{SInt}, \: (-1)^{\l+1} d'_{\lambda,0} \Big\}}  "TL";"TR"};
     {\ar@<-0.6pc>_-{
        \Big\{ i_\lambda^4, \: (-1)^h d'_{\lambda,1} + {}_\ast 2 (-1)^{\l+1} d'_{\lambda,0}   \Big\} }    "TL";"TR"};
\endxy
\]

We assume that if the two different reductions on a given diagram are applied at different heights, the upper branch of the critical branchings will represent the rewriting sequence corresponding to the application of the uppermost first rewriting step.
From now on, unless reconstructing the final result of a given critical branching is difficult for one branch of reductions, we will represent the critical branchings and critical branchings modulo using this notation.

\subsection{Indexed critical branchings}
The classification of indexed critical branchings follows from the indexed critical branchings for the $3$-polygraphs of pearls in \cite{GM09}, for all possible orientation of strands. Let us draw the ones that differ from the even case, labelled by some odd $i \in I$.
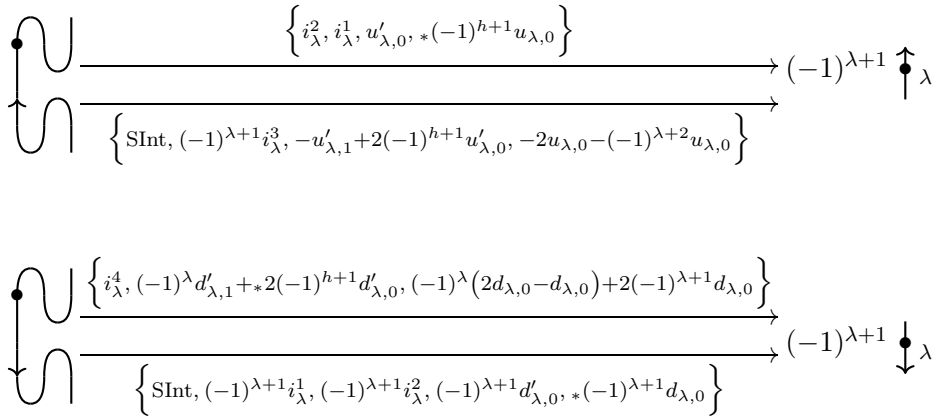
\[
\xy
  (-55,10)*+{ \begin{tikzpicture}[baseline = 0,scale=0.9]
  \draw[-,thick,black] (0.3,0) to (0.3,.4);
	\draw[-,thick,black] (0.3,0) to[out=-90, in=0] (0.1,-0.4);
	\draw[-,thick,black] (0.1,-0.4) to[out = 180, in = -90] (-0.1,0);
	\draw[-,thick,black] (-0.1,0) to[out=90, in=0] (-0.3,0.4);
	\draw[-,thick,black] (-0.3,0.4) to[out = 180, in =90] (-0.5,0);
  \draw[-,thick,black] (-0.5,0) to (-0.5,-.7);
  \draw[-,thick,black] (-0.5,-0.7) to (-0.5,-0.8);
 %%%%
  \draw[-,thick,black] (0.3,-1.2) to (0.3,-1.6);
	\draw[-,thick,black] (0.3,-1.2) to[out=90, in=0] (0.1,-0.8);
	\draw[-,thick,black] (0.1,-0.8) to[out = 180, in = 90] (-0.1,-1.2);
	\draw[-,thick,black] (-0.1,-1.2) to[out=-90, in=0] (-0.3,-1.6);
	\draw[-,thick,black] (-0.3,-1.6) to[out = 180, in =-90] (-0.5,-1.2);
  \draw[->,thick,black] (-0.5,-1.2) to (-0.5,-0.8);
  %%%
  \node at (-0.5,0) {$\bullet$};
\end{tikzpicture}
       }="TL";
  (55,10)*+{  (-1)^{\l+1} \upd{} }="TR";
     {\ar@<.6pc> ^-{
        \Big\{  i_{\lambda}^2, \: i_{\lambda}^1, \: u'_{\lambda, 0}, \: {}_\ast (-1)^{h+1} u_{\lambda,0} \Big\}}  "TL";"TR"};
     {\ar@<-0.6pc>_-{
        \Big\{ \text{SInt}, \: (-1)^{\l+1} i_\lambda^3, \: - u'_{\lambda,1} + 2 (-1)^{h+1} u'_{\lambda,0} , \: -2 u_{\lambda,0} - (-1)^{\l+2} u_{\lambda,0}  \Big\} }    "TL";"TR"};
\endxy
\]

\[
\xy
  (-55,10)*+{ \begin{tikzpicture}[baseline = 0,scale=0.9]
  \draw[-,thick,black] (0.3,0) to (0.3,.4);
	\draw[-,thick,black] (0.3,0) to[out=-90, in=0] (0.1,-0.4);
	\draw[-,thick,black] (0.1,-0.4) to[out = 180, in = -90] (-0.1,0);
	\draw[-,thick,black] (-0.1,0) to[out=90, in=0] (-0.3,0.4);
	\draw[-,thick,black] (-0.3,0.4) to[out = 180, in =90] (-0.5,0);
  \draw[-,thick,black] (-0.5,0) to (-0.5,-.7);
  \draw[-,thick,black] (-0.5,-0.7) to (-0.5,-0.8);
 %%%%
  \draw[-,thick,black] (0.3,-1.2) to (0.3,-1.6);
	\draw[-,thick,black] (0.3,-1.2) to[out=90, in=0] (0.1,-0.8);
	\draw[-,thick,black] (0.1,-0.8) to[out = 180, in = 90] (-0.1,-1.2);
	\draw[-,thick,black] (-0.1,-1.2) to[out=-90, in=0] (-0.3,-1.6);
	\draw[-,thick,black] (-0.3,-1.6) to[out = 180, in =-90] (-0.5,-1.2);
  \draw[<-,thick,black] (-0.5,-1.2) to (-0.5,-0.8);
  %%%
  \node at (-0.5,0) {$\bullet$};
\end{tikzpicture}
       }="TL";
  (55,10)*+{  (-1)^{\l+1} \dpd{} }="TR";
     {\ar@<.6pc> ^-{
        \Big\{  i_{\lambda}^4, \: (-1)^\l d'_{\lambda , 1} + {}_\ast 2 (-1)^{h+1} d'_{\lambda, 0}, \: (-1)^\l \left( 2 d_{\lambda, 0} - d_{\lambda, 0} \right) + 2(-1)^{\l+1} d_{\lambda,0}  \Big\}}  "TL";"TR"};
     {\ar@<-0.6pc>_-{
        \Big\{\text{SInt}, \: (-1)^{\l+1} i_\lambda^1, \: (-1)^{\l+1} i_\lambda^2, \:  (-1)^{\l+1} d'_{\lambda, 0}, \: {}_\ast (-1)^{\l+1} d_{\lambda,0}   \Big\} }    "TL";"TR"};
\endxy
\]

\[
\xy
  (-45,10)*+{
        \begin{tikzpicture}[baseline = 0,scale=0.8]
      \draw[-,thick,black] (0.3,0) to (0.3,-.4);
    	\draw[-,thick,black] (0.3,0) to[out=90, in=0] (0.1,0.4);
    	\draw[-,thick,black] (0.1,0.4) to[out = 180, in = 90] (-0.1,0);
    	\draw[-,thick,black] (-0.1,0) to[out=-90, in=0] (-0.3,-0.4);
    	\draw[-,thick,black] (-0.3,-0.4) to[out = 180, in =-90] (-0.5,0);
      \draw[<-,thick,black] (-0.5,0) to (-0.5,.4);
       \node at (0.9,0) {$\scriptstyle{\lambda}$};
       %%% BIG CAP
       	\draw[-,thick,black] (0.7,0.4) to[out=90, in=0] (0.1,0.9);
    	\draw[-,thick,black] (0.1,0.9) to[out = 180, in = 90] (-0.5,0.4);
    	\draw[-,thick,black] (0.7,0.4) to (0.7,-0.4);
    	\node at (-0.45,0.5) {$\bullet$};
        \end{tikzpicture}}="TL";
  (45,10)*+{ (-1)^\l \capldr{}{} + 2 \: \raisebox{-3mm}{$\scalebox{0.7}{\capoddbubblenest{}{}}$}}="TR";
     {\ar@<.6pc> ^-{
        \Big\{  i_{\lambda}^4,\:  (-1)^\l d_{\lambda+2 , 0} + 2 d_{\lambda+2 , 0} \Big\}}  "TL";"TR"};
     {\ar@<-0.6pc>_-{
        \Big\{ SInt, i_{\lambda}^1, \: i_{\lambda}^2, \: d_{\lambda+2 , 0}, \: i_{\lambda}^4 \Big\} }    "TL";"TR"};
\endxy
\]

\[
\xy
  (-45,10)*+{
  \begin{tikzpicture}[baseline = 0,scale=0.8]
  \draw[-,thick,black] (0.3,0) to (0.3,.4);
	\draw[-,thick,black] (0.3,0) to[out=-90, in=0] (0.1,-0.4);
	\draw[-,thick,black] (0.1,-0.4) to[out = 180, in = -90] (-0.1,0);
	\draw[-,thick,black] (-0.1,0) to[out=90, in=0] (-0.3,0.4);
	\draw[-,thick,black] (-0.3,0.4) to[out = 180, in =90] (-0.5,0);
  \draw[->,thick,black] (-0.5,0) to (-0.5,-.4);
   \node at (0.9,0) {$\scriptstyle{\lambda}$};
  \draw[-,thick,black] (-0.5,-0.4) to (-0.5,-0.7);
  	\draw[-,thick,black] (0.7,-0.7) to[out=-90, in=0] (0.1,-1.2);
	\draw[-,thick,black] (0.1,-1.2) to[out = 180, in = -90] (-0.5,-0.7);
	\draw (0.7,-0.7) to (0.7,0.4);
	\node at (-0.5,0.15) {$\bullet$};
\end{tikzpicture}
        }="TL";
  (45,10)*+{ (-1)^{\l+1} \cuprdr{i}{}}="TR";
     {\ar@<.6pc> ^-{
        \Big\{ i_\lambda^4, \: (-1)^\l d'_{\lambda+2,1} + {}_\ast 2 (-1)^{\l +1} d'_{\lambda+2,0}, \: (-1)^{\l+1} i_{\lambda}^1 \Big\}}  "TL";"TR"};
     {\ar@<-0.6pc>_-{
        \Big\{ SInt,  \: (-1)^{\l+1} i_\lambda^1, \: (-1)^{\l+1} d'_{\lambda+2,0}, \: (-1)^{\l+1} i_\lambda^1 \Big\} }    "TL";"TR"};
\endxy
\]

\[
\xy
  (-45,10)*+{
 \begin{tikzpicture}[baseline = 0,scale=0.8]
  \draw[-,thick,black] (0.3,0) to (0.3,.4);
	\draw[-,thick,black] (0.3,0) to[out=-90, in=0] (0.1,-0.4);
	\draw[-,thick,black] (0.1,-0.4) to[out = 180, in = -90] (-0.1,0);
	\draw[-,thick,black] (-0.1,0) to[out=90, in=0] (-0.3,0.4);
	\draw[-,thick,black] (-0.3,0.4) to[out = 180, in =90] (-0.5,0);
  \draw[<-,thick,black] (-0.5,0) to (-0.5,-.4);
   \node at (0.3,.6) {$\scriptstyle{i}$};
   \node at (0.9,0) {$\scriptstyle{\lambda}$};
  \draw[-,thick,black] (-0.5,-0.4) to (-0.5,-0.7);
  	\draw[-,thick,black] (0.7,-0.7) to[out=-90, in=0] (0.1,-1.2);
	\draw[-,thick,black] (0.1,-1.2) to[out = 180, in = -90] (-0.5,-0.7);
	\draw (0.7,-0.7) to (0.7,0.4);
	\node at (-0.5,0.15) {$\bullet$};
\end{tikzpicture}
        }="TL";
  (45,10)*+{ (-1)^\l  \cupldr{}{} + 2 \raisebox{-6mm}{$\cupoddbubble{i}{i}$} }="TR";
     {\ar@<.6pc> ^-{
        \Big\{ i_\lambda^2, \: i_\lambda^1, \: u_{\lambda-2,0}, \: i_{\lambda}^3 \Big\}}  "TL";"TR"};
     {\ar@<-0.6pc>_-{
        \Big\{ SInt,  \: i_\lambda^3, \: (-1)^\lambda u_{\lambda-2,0} + 2 u_{\lambda-2,0}  \Big\} }    "TL";"TR"};
\endxy
\]

\[
\xy
  (-45,10)*+{
\begin{tikzpicture}
\begin{scope} [ x = 10pt, y = 10pt, join = round, cap = round, thick, scale=3.2]
  \draw[->,thick,black] (0.7,0.2) to[out=90,in=0] (0,.9);
  \draw[-,thick,black] (0,0.9) to[out=180,in=90] (-.7,0.2);
\draw[-,thick,black] (-.7,0.2) to[out=-90,in=180] (0,-0.5);
  \draw[-,thick,black] (0,-0.5) to[out=0,in=-90] (0.7,0.2);
 \node at (0,-.75) {$\scriptstyle{i}$};
  \node at (0.7,1) {$\scriptstyle{\lambda}$};
 % End
 \node at (0.6,0.55) {$\bullet$};
 \node at (0.8,0.55) {$\scriptstyle{m}$};
 \node at (-0.45,0.7) {$\bullet$};
 %%%% Cell k
 \draw[-,thick,black] (-0.3,-0.35) to (-0.3,0) to (0.3,0) to (0.3,-0.35) to (-0.3,-0.35);
 \node at  (0,-0.2) {$\scriptstyle{k}$};
 \end{scope}
\end{tikzpicture}
        }="TL";
  (45,10)*+{
      (-1)^{m}\raisebox{-3mm}{$\targetbetamk{i}{i}{m+1}{k}$}} ="TR";
     {\ar@<.6pc> ^-{
        \Big\{  i_4^\lambda, \: {}_* (-1)^m \beta_{m,k} \Big\}}  "TL";"TR"};
     {\ar@<-0.6pc>_-{
        \Big\{  \text{SInt}, \: (-1)^{m+|k|} i_\lambda^1, \: \text{SInt}  \Big\} }    "TL";"TR"};
\endxy
\]

\[
\xy
  (-45,10)*+{
\begin{tikzpicture}
\begin{scope} [ x = 10pt, y = 10pt, join = round, cap = round, thick, scale=3.2]
  \draw[-<,thick,black] (0.7,0.2) to[out=90,in=0] (0,.9);
  \draw[-,thick,black] (0,0.9) to[out=180,in=90] (-.7,0.2);
\draw[-,thick,black] (-.7,0.2) to[out=-90,in=180] (0,-0.5);
  \draw[-,thick,black] (0,-0.5) to[out=0,in=-90] (0.7,0.2);
 \node at (0,-.75) {$\scriptstyle{i}$};
  \node at (0.7,1) {$\scriptstyle{\lambda}$};
 % End
 \node at (0.6,0.55) {$\bullet$};
 \node at (0.8,0.55) {$\scriptstyle{m}$};
 \node at (-0.45,0.7) {$\bullet$};
 %%%% Cell k
 \draw[-,thick,black] (-0.3,-0.35) to (-0.3,0) to (0.3,0) to (0.3,-0.35) to (-0.3,-0.35);
 \node at  (0,-0.2) {$\scriptstyle{k}$};
 \end{scope}
\end{tikzpicture}
        }="TL";
  (45,10)*+{
      \raisebox{-3mm}{$\targetalphamk{i}{i}{m+1}{k}$}} ="TR";
     {\ar@<.6pc> ^-{
        \Big\{  i_\lambda^2 \Big\}}  "TL";"TR"};
     {\ar@<-0.6pc>_-{
        \Big\{  \text{SInt}, \: (-1)^{m+|k|}i_\l^3, (-1)^{m+|k|}\alpha_{m,k}, \: \text{SInt}  \Big\} }    "TL";"TR"};
\endxy
\]

\[
\xy
  (-50,10)*+{
 \begin{tikzpicture}
\begin{scope} [ x = 10pt, y = 10pt, join = round, cap = round, thick, scale=3]
  \draw[-,thick,black] (0.1,0.3) to[out=180,in=90] (-.1,0.1);
  \draw[-,thick,black] (0.3,0.1) to[out=90,in=0] (0.1,.3);
 \draw[-,thick,black] (-.1,0.1) to[out=-90,in=180] (0.1,-0.1);
  \draw[-,thick,black] (0.1,-0.1) to[out=0,in=-90] (0.3,0.1);
  \draw[-,thick,black] (0.24,-0.05) to (-0.04,0.23);
  \draw[-,thick,black] (0.24,0.23) to (-0.04,-0.05);
 \node at (0.1,-.25) {$\scriptstyle{i}$};
  \node at (0.5,0.25) {$\scriptstyle{\lambda}$};
  %%%%%% End of odd bubble part
  \draw[<-,thick,black,scale=1.3] (0.1,1.2) to[out=180,in=90] (-.3,0.8);
  \draw[-,thick,black,scale=1.3] (0.5,0.8) to[out=90,in=0] (0.1,1.2);
 \draw[-,thick,black,scale=1.3] (-.3,0.8) to[out=-90,in=180] (0.1,0.4);
  \draw[-,thick,black,scale=1.3] (0.1,0.4) to[out=0,in=-90] (0.5,0.8);
 \node at (-0.45,0.8) {$\scriptstyle{i}$};
 \node at (0.58,1.3) {$\bullet$};
 \node at (0.8,1.3) {$\scriptstyle{m}$};
 \node at (-0.2,1.4) {$\bullet$};
 \draw[-,thick,black] (0.4,1.1) to (-0.1,1.1) to (-0.1,0.7) to (0.40,0.7) to (0.40,1.1);
 \node at (0.15,0.9) {$\scriptstyle{k}$};
 \end{scope}
\end{tikzpicture}
        }="TL";
  (50,10)*+{ \left\{
    \begin{array}{ll}
       (-1)^{m +1+ |k|} \raisebox{-3mm}{$\targetalphamk{i}{i}{m+2}{k}$} & \text{if $m+\l_i+2 $ is even} \\
        0 & \text{if $m+\l_i+2$ is odd}
    \end{array}
\right.}="TR";
     {\ar@<.6pc> ^-{
        \Big\{  i_\lambda^2, \: \alpha_{m+1,k} \Big\}}  "TL";"TR"};
     {\ar@<-0.6pc>_-{
        \Big\{  \text{SInt}, \: (-1)^{m+|k|} \: i_3^\lambda, {}_\ast \alpha_{m+1,k}  \Big\} }    "TL";"TR"};
\endxy
\]

\[
\xy
  (-50,10)*+{
\begin{tikzpicture}
\begin{scope} [ x = 10pt, y = 10pt, join = round, cap = round, thick, scale=3.2]
  \draw[->,thick,black] (0.7,0.2) to[out=90,in=0] (0,.9);
  \draw[-,thick,black] (0,0.9) to[out=180,in=90] (-.7,0.2);
\draw[-,thick,black] (-.7,0.2) to[out=-90,in=180] (0,-0.5);
  \draw[-,thick,black] (0,-0.5) to[out=0,in=-90] (0.7,0.2);
 \node at (0,-.75) {$\scriptstyle{i}$};
  \node at (0.7,1) {$\scriptstyle{\lambda}$};
  %% Odd bubble part
   \draw[-,thick,black,scale=0.8] (0,0.5) to[out=180,in=90] (-.2,0.3);
  \draw[-,thick,black,scale=0.8] (0.2,0.3) to[out=90,in=0] (0,0.5);
 \draw[-,thick,black,scale=0.8] (-.2,0.3) to[out=-90,in=180] (0,0.1);
  \draw[-,thick,black,scale=0.8] (0,0.1) to[out=0,in=-90] (0.2,0.3);
  \draw[-,thick,black,scale=0.8] (0.14,0.15) to (-0.14,0.43);
  \draw[-,thick,black,scale=0.8] (0.14,0.43) to (-0.14,0.15);
 \node at (-0.3,0.3) {$\scriptstyle{i}$};
 % End
 \node at (0.6,0.55) {$\bullet$};
 \node at (0.8,0.55) {$\scriptstyle{m}$};
 \node at (-0.45,0.7) {$\bullet$};
 %%%% Cell k
 \draw[-,thick,black] (-0.3,-0.35) to (-0.3,0) to (0.3,0) to (0.3,-0.35) to (-0.3,-0.35);
 \node at  (0,-0.2) {$\scriptstyle{k}$};
 \end{scope}
\end{tikzpicture}
        }="TL";
  (50,10)*+{\left\{
    \begin{array}{ll}
      (-1)^{m}\raisebox{-3mm}{$\targetbetamk{i}{i}{m+2}{k}$}  & \text{if $m+\l_i+2$ is even} \\
        0 & \text{if $m+\l_i+2$ is odd}
    \end{array}
\right.}="TR";
     {\ar@<.6pc> ^-{
        \Big\{  i_\lambda^4, \: (-1)^{\l_i} \beta_{m+1,k} \Big\}}  "TL";"TR"};
     {\ar@<-0.6pc>_-{
        \Big\{  \text{SInt}, \: (-1)^{m+1+|k|} i_\lambda^1, \: {}_\ast (-1)^m \beta_{m+1,k}  \Big\} }    "TL";"TR"};
\endxy
\]

\section{Critical branchings of $\mathbf{ONH}$}
\label{appendix:criticalbranchingsONH}

% - - - - - - - - - - - - - - - -
\subsection{Helpful 3-cells  of ONH}
% - - - - - - - - - - - - - - - -
We introduce some additional 3-cells in ONH that will be helpful in analyzing the critical branchings.
To simplify the description of these critical branchings we make use of the following 3-cells obtained by iterative application of $on_1$ and $on_2$, see also \cite[Lemma 3.1]{BE2}.
\begin{align}
&\hackcenter{
}
\end{align}
Using these 3-cells we also introduce 3-cells
\begin{align}
 \Theta_{x,y} &:= \left( \hackcenter{
%
}
\right)
\end{align}
and more generally we have
\begin{align}
&\Phi_{x,y}:= \left(\hackcenter{
%
} \right)
\nn
\end{align}
It will also be convenient to define a 3-cell $\Upsilon$ given by
\begin{align}
&\Upsilon := \left( \hackcenter{
%
 }
\end{align}
where in the first summation is zero unless $b>c$ and the second is zero unless $b<c$. We will show that the target of this 3-cell is zero.  Swapping the $b,c$ variables in the second summation the right hand side can be written
\begin{align} \label{eq:shouldBzero2}
\sum_{\overset{\scs a+b+c}{=n-2}}
  \sum_{j=c}^{b-1}
  (-1)^{j}\left[ (-1)^{jc} - (-1)^{jb}\right]
\hackcenter{
 }
\end{align}
where
\[
\left[ (-1)^{jc} - (-1)^{jb}\right]
=
\left\{
  \begin{array}{cl}
    2, & \hbox{if $b$ is even, $c$ is odd, and $j$ is odd,} \\
    -2, & \hbox{if $b$ is odd, $c$ is even, and $j$ is odd, } \\
    0 & \hbox{otherwise.}
  \end{array}
\right.
\]
Breaking the $b$ and $c$ summations in \eqref{eq:shouldBzero2} into a sum over even and odd terms, the only nonvanishing terms are
\begin{align}
& \sum_{\ell_1=0}^{\lfloor \frac{n-2}{2} \rfloor} \sum_{\ell_2=0}^{\lfloor \frac{2\ell_1-1}{2} \rfloor }  \sum_{j=2\ell_2+1}^{2\ell_1-1}
  (-1)^{j}\left[ (-1)^{j(2\ell_2+1)} - (-1)^{j(2\ell_1)}\right]
\hackcenter{
%
 }
\end{align}
Now observe that since $j$ is assumed to be odd , we can remove the $\ell_2=\ell_1$ term in the second summation since $2\ell_1 \leq j \leq 2\ell_1$ would imply $j$ was even. Similarly, since $j$ is odd the $j$ summation index in the second term can start at $2\ell_2+1$ and end at $2\ell_1-1$ so that the above terms cancel out and the target of the 3-cell from \eqref{eq:shouldBzero} is zero.
% - - - - - - - - - - - - - - - -
\subsection{Regular critical branchings of ONH}
% - - - - - - - - - - - - - - - -
In this section we study the critical branchings of the (3,2)-superpolygraph {\rm ONH}.  We begin with the relatively straightforward regular branchings using the shorthand notation introduced in Appendix~\ref{sec:shortnotation}.
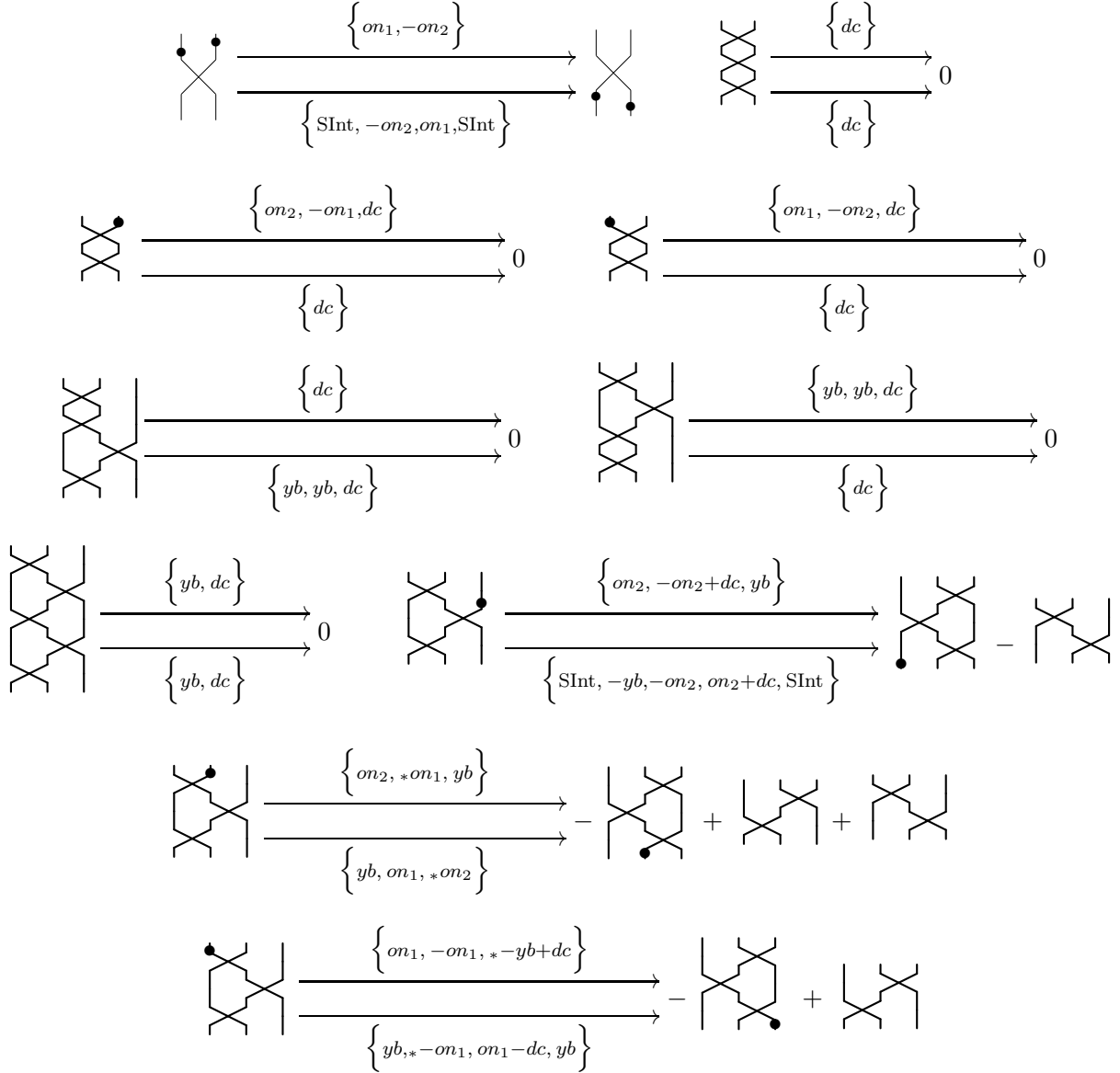
\[
\xy
  (-30,10)*+{
       \scalebox{0.5}{\raisebox{-1cm}{$\begin{tikzpicture}
\draw (0,0) to (1,1);
\draw (1,1) to (1,1.75);
\draw (0,1) to (0,1.75);
\draw (1,0) to (0,1);
\draw (0,0) to (0,-0.75);
\draw (1,0) to (1,-0.75);
\node[scale=1.75] at (0,1.2) {$\bullet$};
\node[scale=1.75] at (1,1.5) {$\bullet$};
\end{tikzpicture}$}}}="TL";
  (30,10)*+{ \scalebox{0.5}{\raisebox{-1cm}{$\begin{tikzpicture}
\draw (0,0) to (1,1);
\draw (1,1) to (1,1.75);
\draw (0,1) to (0,1.75);
\draw (1,0) to (0,1);
\draw (0,0) to (0,-0.75);
\draw (1,0) to (1,-0.75);
\node[scale=1.75] at (1,-0.5) {$\bullet$};
\node[scale=1.75] at (0,-0.2) {$\bullet$};
\end{tikzpicture}$}} }="TR";
     {\ar@<.6pc> ^-{
        \Big\{  on_1, - on_2  \Big\}}  "TL";"TR"};
     {\ar@<-0.6pc>_-{
        \Big\{ \text{SInt}, \: -on_2, on_1, \text{SInt} \Big\} }    "TL";"TR"};
\endxy
\qquad
\xy
  (-15,10)*+{
       \scalebox{0.911}{\raisebox{-5mm}{$\triplecrossing{}{}$}} }="TL";
  (15,10)*+{ 0 }="TR";
     {\ar@<.6pc> ^-{
        \Big\{  dc  \Big\}}  "TL";"TR"};
     {\ar@<-0.6pc>_-{
        \Big\{ dc \Big\} }    "TL";"TR"};
\endxy
\]
\[
\xy
  (-30,10)*+{
       \raisebox{-6mm}{$\tcrossur{}{}$} }="TL";
  (30,10)*+{ 0 }="TR";
     {\ar@<.6pc> ^-{
        \Big\{  on_2, \: - on_1, dc  \Big\}}  "TL";"TR"};
     {\ar@<-0.6pc>_-{
        \Big\{ dc \Big\} }    "TL";"TR"};
\endxy
\qquad
\xy
  (-30,10)*+{
       \raisebox{-6mm}{$\tcrossul{}{}$} }="TL";
  (30,10)*+{ 0 }="TR";
     {\ar@<.6pc> ^-{
        \Big\{  on_1, \: - on_2, \: dc  \Big\}}  "TL";"TR"};
     {\ar@<-0.6pc>_-{
        \Big\{ dc \Big\} }    "TL";"TR"};
\endxy
\]
\[
\xy
  (-30,10)*+{
       \raisebox{-10mm}{$\ybgcrossu$} }="TL";
  (30,10)*+{ 0 }="TR";
     {\ar@<.6pc> ^-{
        \Big\{  dc  \Big\}}  "TL";"TR"};
     {\ar@<-0.6pc>_-{
        \Big\{ yb, \: yb, \: dc \Big\} }    "TL";"TR"};
\endxy
\qquad
\xy
  (-30,10)*+{
       \raisebox{-10mm}{$\ybgcrossd{}{}{}$} }="TL";
  (30,10)*+{ 0 }="TR";
     {\ar@<.6pc> ^-{
        \Big\{  yb, \: yb, \: dc  \Big\}}  "TL";"TR"};
     {\ar@<-0.6pc>_-{
        \Big\{ dc  \Big\} }    "TL";"TR"};
\endxy
\]
\[
\xy
  (-20,10)*+{
       \raisebox{-7mm}{$\doubleybg{}{}{}$} }="TL";
  (20,10)*+{ 0 }="TR";
     {\ar@<.6pc> ^-{
        \Big\{ yb, \: dc  \Big\}}  "TL";"TR"};
     {\ar@<-0.6pc>_-{
        \Big\{ yb, \: dc \Big\} }    "TL";"TR"};
\endxy
\qquad
\xy
  (-40,10)*+{
       \raisebox{-5mm}{$\ybgur{}{}{}$} }="TL";
  (40,10)*+{ \raisebox{-5mm}{$\ybdddl{}{}{}$} \: - \: \raisebox{-5mm}{$\cgcd{}{}{}$} }="TR";
     {\ar@<.6pc> ^-{
        \Big\{  on_2, \: - on_2 + dc, \: yb \Big\}}  "TL";"TR"};
     {\ar@<-0.6pc>_-{
        \Big\{ \text{SInt}, \: -yb, -on_2, \: on_2 + dc, \: \text{SInt} \Big\} }    "TL";"TR"};
\endxy
\]
\[
\xy
  (-40,10)*+{
       \raisebox{-5mm}{$\ybgum{}{}{}$} }="TL";
  (40,10)*+{ - \: \raisebox{-7mm}{$\ybddm{}{}{}$} \: + \: \raisebox{-7mm}{$\cdcg{}{}{}$} + \: \raisebox{-5mm}{$\cgcd{}{}{}$} }="TR";
     {\ar@<.6pc> ^-{
        \Big\{  on_2, \: {}_\ast on_1, \: yb  \Big\}}  "TL";"TR"};
     {\ar@<-0.6pc>_-{
        \Big\{ yb, \: on_1, \: {}_\ast on_2 \Big\} }    "TL";"TR"};
\endxy
\]
\[
\xy
  (-40,10)*+{
       \raisebox{-5mm}{$\ybgul{}{}{}$} }="TL";
  (40,10)*+{ - \: \raisebox{-6mm}{$\ybddr{}{}{}$} \: + \: \raisebox{-8mm}{$\cdcg{}{}{}$} }="TR";
     {\ar@<.6pc> ^-{
        \Big\{  on_1, \: -on_1, \: {}_\ast -yb + dc \Big\}}  "TL";"TR"};
     {\ar@<-0.6pc>_-{
        \Big\{ yb, {}_\ast -on_1, \: on_1 - dc, \: yb  \Big\} }    "TL";"TR"};
\endxy
\]
% - - - - - - - - - - - - - - - -
\subsection{Indexed critical branchings of ONH}
% - - - - - - - - - - - - - - - -
We now verify the indexed critical branchings of $\mathbf{ONH}$, whose classification is the same as in \cite{DUP19bis}, by spelling out in greater detail the required steps as they are somewhat subtle and differ notably from the corresponding calculations in the even setting. The first indexed critical branching is obtained by reducing the diagram
\[
\hackcenter{
$}
\]
in two possible ways; the top branch is obtained by first applying the Yang-Baxter 3-cell to the bottom half, then sliding the $n$ dots to the bottom, while the bottom branch is obtained by first doing super interchange law, applying Yang-Baxter to the top half of the diagram, and sliding the $n$ dots to the bottom.  In detail, the top branch is given by the following.
\begin{align} \nn
&
\raisebox{-8mm}{$%
$}
\right] \\ \nn
& \overset{\{ dc, dc, yb, id \}}{\Rrightarrow}
\sum\limits_{a+b=n-1} (-1)^{b+n} \raisebox{-8mm}{$%
$} \label{eq:topbranchIC1}
\end{align}
The bottom branch is given by
\begin{align*}
\hspace{-1cm}
\raisebox{-8mm}{$%
$}
\end{align*}
where the first summand reduces by $\Phi_2$ into the right-hand side of \eqref{eq:topbranchIC1}
%
% \[ \sum\limits_{a+b+c=n-2} (-1)^{n-1-b} \raisebox{-8mm}{$%
$} \]
recovering the normal form obtained in the top part of the branching, and the second summand reduces to $0$ using the $3$-cell  in \eqref{eq:shouldBzero}.
\medskip

The final indexed critical branching for $\mathbf{ONH}$ is obtained from two branches obtained from the left-hand side below
\[
\xy
  (-30,10)*+{ \raisebox{-7mm}{$%
$} }="TR";
     {\ar@<.6pc> ^-{
         }  "TL";"TR"};
     {\ar@<-0.6pc>_-{
          }    "TL";"TR"};
\endxy
\]
obtained by applying the Yang-Baxter 3-cell to the top of the diagram, then sliding the $n$ dots to the bottom using 3-cells $on_2$.  The bottom branch is obtained by applying super interchange, applying the Yang-Baxter 3-cell to the bottom of the diagram, then sliding the $n$ dots to the bottom of the diagram using $on_2$.
In more detail, the top branching is
\begin{align*}
\raisebox{-8mm}{$%
$}
\end{align*}
Note that the last summand is a normal form, and the first three terms reduce respectively using the rewriting paths $\{ yb, \: SInt, \: yb, \: yb \}$, $\{ yb \}$ and $\{ SInt, \: -yb, \: -dc \}$ so that this branch of the branching gives
\begin{equation}
\label{E:TopPartLastIndexationONH}
 \raisebox{-9mm}{$%
$}
    \end{equation}

The bottom branch is given by applying the 3-cells below.

  \begin{align} \nn
&
     \raisebox{-8mm}{$%
$}  \label{eq:IC2-bottom}
        \end{align}

Now we relate the terms in \eqref{eq:IC2-bottom} to those appearing in the top branch~\eqref{E:TopPartLastIndexationONH}.
\begin{itemize}
  \item The first summand of \eqref{eq:IC2-bottom} reduces using the rewriting path $\{ SInt, \: -yb, \: SInt, \: yb, \: yb \}$ to the first summand of \eqref{E:TopPartLastIndexationONH}.
  \item The second summand of \eqref{eq:IC2-bottom} reduces using the rewriting path $\{ SInt, - \sum\limits_{a+b=n-1} \Phi_{1,b} \}$
  into
  \begin{equation} \label{eq:2nd-term}
 \sum\limits_{a+c+d = n-2} (-1)^{d+1}  \: \:  \raisebox{-8mm}{$
$} \end{equation}
\item The third summand of \eqref{eq:IC2-bottom} reduces as follows:
\begin{align*}
&
\sum\limits_{a+b = n-1} (-1)^{b+n+1} \raisebox{-8mm}{$%
$}
\end{align*}
so that the first sum gives the second sum of \eqref{E:TopPartLastIndexationONH}, and the second sum gives the third sum of \eqref{E:TopPartLastIndexationONH}. As a consequence, using the first and third summand on the bottom branch, we recover all the elements from the top branch.
\item The fourth summand of \eqref{eq:IC2-bottom} reduces using the $3$-cell $\Theta_{b,c}$ as follows:
\begin{align} \nn
-(-1)^{n+1}\sum\limits_{\overset{\scs a+b+c}{\scs =n-2}} (-1)^{b} \: \raisebox{-8mm}{$%
$}
\end{align}
where, after applying the 3-cell $yb$, the first term on the right-hand side cancels with \eqref{eq:2nd-term}, and the second summation above reduces to zero as in the computation of \eqref{eq:shouldBzero}.
The second sum
\begin{align}\label{needtobezeroII}
\sum\limits_{\overset{\scs a+b+c}{\scs =n-2}} \sum\limits_{j= \text{min}(b,c)}^{\text{max}(b,c)-1} (-1)^{c+b+j + jc + \delta_{\text{max}(b,c),b}} \: \raisebox{-7mm}{$
 %
$}
\end{align}
Swapping the $b,c$ variables in the first summation, this is equal to
\begin{align} \nn
  \sum\limits_{\overset{\scs a+b+c}{\scs =n-2}} \sum\limits_{j= c}^{b-1} (-1)^{b+j+c}[(-1)^{jb}-(-1)^{jc}] \: \raisebox{-7mm}{$
 $}
\end{align}
where
\[
\left[ (-1)^{jc} - (-1)^{jb}\right]
=
\left\{
  \begin{array}{cl}
    2, & \hbox{if $b$ is even, $c$ is odd, and $j$ is odd,} \\
    -2, & \hbox{if $b$ is odd, $c$ is even, and $j$ is odd, } \\
    0 & \hbox{otherwise.}
  \end{array}
\right.
\]
In particular, $b+c$ must be odd, so we write this summation as
\begin{align}\label{eq:lastindexonhzeroIII}
  \sum\limits_{b=1}^{n-2} \sum\limits_{c=1}^{b-1} \sum\limits_{j= c}^{b-1} (-1)^{j}[(-1)^{jc}-(-1)^{jb}] \: \raisebox{-7mm}{$
 %
$}
\end{align}
Breaking the $b$ and $c$ summations in \eqref{eq:lastindexonhzeroIII} into a sum over even and odd terms, the only non vanishing terms are
\begin{align}
& \sum_{\ell_1=0}^{\lfloor \frac{n-2}{2} \rfloor} \sum_{\ell_2=0}^{\lfloor \frac{2\ell_1-1}{2} \rfloor }  \sum_{j=2\ell_2+1}^{2\ell_1-1}
  (-1)^{j}\left[ (-1)^{j(2\ell_2+1)} - (-1)^{j(2\ell_1)}\right]
\hackcenter{
%
}
\end{align}
Now observe that since $j$ is assumed to be odd, we can remove the $\ell_2=\ell_1$ term in the second summation since $2\ell_1 \leq j \leq 2\ell_1$ would imply $j$ was even. Similarly, since $j$ is odd the $j$ summation index in the second term can start at $2\ell_2+1$ and end at $2\ell_1-1$ so that the above terms cancel out and the target of the 3-cell from \eqref{needtobezeroII} is zero.
  \end{itemize}

\section{Critical branchings modulo for the full $2$-category}
\label{appendix:criticalbranchingsfull}
In this Section, we prove that the critical branchings modulo for the $(3,2)$-superpolygraph $\mathbf{Osl(2)}$ are confluent modulo $E$.

\subsection{Critical branchings from $3$-cells of $\mathbf{ONH}$}
\label{appendix:criticalbranchingsfullONH}
We prove that the critical branchings implying a $3$-cell of $\mathbf{Osl(2)}$ with $on_1$, $on_2$ and two $3$-cells of $\mathbf{Osl(2)}$ on two terms that are equal up to $yb$ are confluent modulo $E$.

\medskip
\noindent \textbf{Critical branchings $(A_{\lambda}, on_{1,\lambda-2})$}
For any $\lambda$ in $\Z$ , the critical branchings $(A_{\lambda}, on_{1,\lambda-2})$ are confluent modulo super isotopies as follows:

\[
\xy
  (-65,10)*+{ \cbadot{}
       }="TL";
  (65,10)*+{   \sum\limits_{n=0}^{- \lambda} (-1)^n \scalebox{1.211}{\druleaf{}{n+1}{-n-1}} }="TR";
     {\ar@<.6pc> ^-{
        \Big\{  A_\lambda \Big\}}  "TL";"TR"};
     {\ar@<-0.6pc>_-{
        \Big\{ \text{SInt}, \: on_{1,\lambda-2}, \: (i_\lambda^3 \star_2 {}_\ast (i_\lambda^2)^- ) \cdot on_{2,\lambda-2}, \: (i_\lambda^1)^-  {}_\ast \cdot A_{\lambda}, \sum\limits_{n=0}^{-\lambda} (-1)^n \alpha_{-n-1,1_{1_\lambda}}, \: \sum\limits_{n=0}^{-\lambda} i_\lambda^1 \Big\} }    "TL";"TR"};
\endxy
\]

The last $3$-cell used in a bottom sequence is a $3$-cell of $E^s$, needed to close the confluence diagram modulo on the right. Note that, when applying the $3$-cell $(i_\lambda^1)^- \cdot {}_\ast A_{\lambda}$, we obtain the following $2$-cell:

\[
(-1)^\l \sum\limits_{n=0}^{- \lambda} (-1)^n \scalebox{1.011}{\druleafb{}{n}{-n-1}} - 2 \sum\limits_{n=0}^{- \lambda} (-1)^n \scalebox{1.011}{\druleafodd{}{n}{-n}} + (1+(-1)^\lambda) \scalebox{1.011}{\druleafempty{}}
\]
which reduces using the $3$-cell $\sum\limits_{n=0}^{-\lambda} (-1)^n \alpha_{-n-1,1_{1_\lambda}}$ on the second summand into

\[
(-1)^\lambda \left( \sum\limits_{n=0}^{-\lambda} \scalebox{1.011}{\druleafb{}{n}{-n-1}} - 2 \sum\limits_{\substack{n=0, \\ n - \lambda \: \text{even}}}^{-\lambda} \scalebox{1.011}{\druleafwodd{}{n}{-n}} \right) + (1+(-1)^\lambda) \scalebox{1.011}{\druleafempty{}}
\]
We can then use the isotopy $3$-cell $\sum\limits_{n=0}^{-\lambda} i_\lambda^1$ to move dots on the right of the cup of the first summand, and obtain the following term:
\[
(-1)^\lambda \left( \sum\limits_{n=0}^{-\lambda} \scalebox{1.011}{\druleaf{}{n+1}{-n-1}} - 2 \sum\limits_{\substack{n=0, \\ n - \lambda \: \text{even}}}^{-\lambda} \scalebox{1.011}{\druleafwodd{}{n}{-n}} \right) + (1+(-1)^\lambda) \scalebox{1.011}{\druleafempty{}}
\]
If $\lambda$ is even, this quantity is
\begin{align*}
\sum\limits_{n=0}^{-\lambda} \scalebox{1.011}{\druleaf{}{n+1}{-n-1}} - 2 \sum\limits_{\substack{n=0, \\ n  \: \text{even}}}^{-\lambda} \scalebox{1.011}{\druleafwodd{}{n}{-n}} + 2 \: \scalebox{1.011}{\druleafempty{}} \: & = \:
\sum\limits_{n=0}^{-\lambda} \scalebox{1.011}{\druleaf{}{n+1}{-n-1}} - 2 \sum\limits_{\substack{n=1, \\ n  \: \text{even}}}^{-\lambda} \scalebox{1.011}{\druleafwodd{}{n}{-n}} \\
& = \sum\limits_{n=0}^{-\lambda+1} \scalebox{1.011}{\druleafwodd{}{n}{-n}} - 2 \sum\limits_{\substack{n=1, \\ n  \: \text{even}}}^{-\lambda +1}  \scalebox{1.011}{\druleafwodd{}{n}{-n}} \\
& = \sum\limits_{n=0}^{-\lambda} (-1)^n \scalebox{1.011}{\druleaf{}{n+1}{-n-1}}
\end{align*}
If $\lambda$ is odd, this quantity is
\begin{align*}
- \sum\limits_{n=0}^{-\lambda} \scalebox{1.011}{\druleaf{}{n+1}{-n-1}} + 2 \sum\limits_{\substack{n=0, \\ n  \: \text{odd}}}^{-\lambda} \scalebox{1.011}{\druleafwodd{}{n}{-n}} \: & = \:
- \sum\limits_{n=1}^{-\lambda+1} \scalebox{1.011}{\druleafwodd{}{n}{-n}} + 2 \sum\limits_{\substack{n=1, \\ n  \: \text{odd}}}^{-\lambda+1} \scalebox{1.011}{\druleafwodd{}{n}{-n}} \\
& = \sum\limits_{n=0}^{-\lambda} (-1)^n \scalebox{1.011}{\druleaf{}{n+1}{-n-1}}
\end{align*}
which proves that the result is the same as in the top branch, so that this branching is confluent modulo $E$.

\medskip
\noindent \textbf{Critical branchings $(B_{\lambda}, i_4^2 \cdot on_{1,\lambda-2})$}

\[
\xy
  (-60,10)*+{ \cbbdot{}
       }="TL";
  (60,10)*+{  \sum\limits_{n=0}^{-\lambda} (-1)^n  \druleb{}{n+1}{-n-1} - (1 + (-1)^\lambda) \drulebempty{} }="TR";
     {\ar@<.6pc> ^-{
        \Big\{  B_\lambda , \: {}_\ast \sum\limits_{n=0}^{-\lambda} i_\lambda^4 \cdot  \sum\limits_{n=0}^{-\lambda} 2 s_{\lambda,1}, \: \sum\limits_{n=0}^{-\lambda} \alpha_{-n-1} \Big\}}  "TL";"TR"};
     {\ar@<-0.6pc>_-{
        \Big\{  \text{SInt}, \: i_\lambda^4 \cdot (on_{1,\lambda -2} + {}_\ast 2 s_{\lambda,1}), \: (i_\lambda^3 \star_2 {}_\ast (i_\lambda^2)^- ) \cdot on_{2,\lambda-2}, {}_\ast B_\lambda \Big\} }    "TL";"TR"};
\endxy
\] 
Note that after using the odd bubble slides in the top sequence, we use similar arguments as above to prove that the target of the rewriting step is equal to the expected result.

\medskip
\noindent \textbf{Critical branchings $(C_{\lambda}, (i_\lambda^2)^- \cdot on_{2,\lambda-2})$}

\[
\hspace{-1cm} \xy
  (-40,10)*+{ \cbcdot{}
       }="TL";
  (40,10)*+{  \sum\limits_{n=0}^{\lambda} (-1)^n  \begin{tikzpicture}[baseline=0]
	\draw[<-,thick,black] (0.3,-.65) to[out=90, in=0] (0,-0.05);
	\draw[-,thick,black] (0,-0.05) to[out = 180, in = 90] (-0.3,-.65);
   \node at (-.27,-.25) {$\bullet$};
   \node at (-.5,-.25) {$\scriptstyle{n}$};
   \node at (0.5,-.35) {$\color{black}\scriptstyle{1}$};
   \node at (0.30,-0.4) {$\bullet$};
  \draw[-,thick,black] (0,0.45) to[out=180,in=90] (-.2,0.25);
  \draw[->,thick,black] (0.2,0.25) to[out=90,in=0] (0,.45);
 \draw[-,thick,black] (-.2,0.25) to[out=-90,in=180] (0,0.05);
  \draw[-,thick,black] (0,0.05) to[out=0,in=-90] (0.2,0.25);
   \node at (-0.55,0) {$\scriptstyle{\lambda}$};
   \node at (0.2,0.25) {$\color{black}\bullet$};
   \node at (0.7,0.25) {$\color{black}\scriptstyle{-n-1}$};
\end{tikzpicture} }="TR";
     {\ar@<.6pc> ^-{
        \Big\{  C_\lambda   \Big\}}  "TL";"TR"};
     {\ar@<-0.6pc>_-{
        \Big\{  X \Big\} }    "TL";"TR"};
\endxy
\]
where $X$ is the rewriting sequence given by
\[ \{ \text{SInt}, -(i _\lambda^2)^- \cdot on_{2,\lambda-2}, \: - \gamma, \: {}_\ast (-1)^\lambda C_\lambda + {}_\ast 2 (-1)^\lambda C_\lambda, \: {}_\ast 2(-1)^\lambda \sum\limits_{n=0}^{\lambda} (-1)^n r_{\lambda,1}, \: 2 \sum\limits_{n=0}^{\lambda} \beta_{-n-1,1_{\lambda+2}} , \: \sum\limits_{n=0}^{\lambda}  i_\lambda^2 {}_\ast \}, \]
where the $3$-cell $\gamma$ is defined as follows:
\begin{align}
\hspace{-0.5cm} \begin{tikzpicture}[baseline = 0,scale=0.7]
  \draw[-,thick,black] (-0.2,0.5) to[out=90, in=0] (-0.5,0.9);
	\draw[->,thick,black] (-0.5,0.9) to[out = 180, in = 90] (-0.8,0.5);
    \node at (1,0.5) {$\scriptstyle{\lambda}$};
	\draw[-,thick,black] (0.3,-.5) to (-0.2,.5);
	\draw[-,thick,black] (-0.2,-.2) to (0.2,.3);
        \draw[-,thick,black] (0.2,.3) to[out=50,in=180] (0.5,.5);
        \draw[->,thick,black] (0.5,.5) to[out=0,in=90] (0.8,-.9);
        \draw[-,thick,black] (-0.2,-.2) to[out=230,in=0] (-0.5,-.5);
        \draw[-,thick,black] (-0.5,-.5) to[out=180,in=-90] (-0.8,.5);
        \node at (-0.2,-0.3) {$\bullet$};
\end{tikzpicture} &
\underset{(i_\lambda^1)^-}{\equiv}
\begin{tikzpicture}[baseline = 0,scale=0.7]
  \draw[-,thick,black] (-0.2,0.5) to[out=90, in=0] (-0.5,0.9);
	\draw[->,thick,black] (-0.5,0.9) to[out = 180, in = 90] (-0.8,0.5);
    \node at (1,0.5) {$\scriptstyle{\lambda}$};
	\draw[-,thick,black] (0.3,-.5) to (-0.2,.5);
	\draw[-,thick,black] (-0.2,-.2) to (0.2,.3);
        \draw[-,thick,black] (0.2,.3) to[out=50,in=180] (0.5,.5);
        \draw[->,thick,black] (0.5,.5) to[out=0,in=90] (0.8,-.9);
        \draw[-,thick,black] (-0.2,-.2) to[out=230,in=0] (-0.5,-.5);
        \draw[-,thick,black] (-0.5,-.5) to[out=180,in=-90] (-0.8,.5);
        \node at (-0.8,-0.3) {$\bullet$};
\end{tikzpicture} \underset{i_\lambda^4}{\equiv} (-1)^{\lambda} \begin{tikzpicture}[baseline = 0,scale=0.7]
  \draw[-,thick,black] (-0.2,0.5) to[out=90, in=0] (-0.5,0.9);
	\draw[->,thick,black] (-0.5,0.9) to[out = 180, in = 90] (-0.8,0.5);
    \node at (1,0.5) {$\scriptstyle{\lambda}$};
	\draw[-,thick,black] (0.3,-.5) to (-0.2,.5);
	\draw[-,thick,black] (-0.2,-.2) to (0.2,.3);
        \draw[-,thick,black] (0.2,.3) to[out=50,in=180] (0.5,.5);
        \draw[->,thick,black] (0.5,.5) to[out=0,in=90] (0.8,-.9);
        \draw[-,thick,black] (-0.2,-.2) to[out=230,in=0] (-0.5,-.5);
        \draw[-,thick,black] (-0.5,-.5) to[out=180,in=-90] (-0.8,.5);
        \node at (-0.2,0.5) {$\bullet$};
\end{tikzpicture} + 2
\begin{tikzpicture}[baseline = 0,scale=0.7]
  \draw[-,thick,black] (-0.2,0.5) to[out=90, in=0] (-0.5,0.9);
	\draw[->,thick,black] (-0.5,0.9) to[out = 180, in = 90] (-0.8,0.5);
    \node at (1,0.5) {$\scriptstyle{\lambda}$};
	\draw[-,thick,black] (0.3,-.5) to (-0.2,.5);
	\draw[-,thick,black] (-0.2,-.2) to (0.2,.3);
        \draw[-,thick,black] (0.2,.3) to[out=50,in=180] (0.5,.5);
        \draw[->,thick,black] (0.5,.5) to[out=0,in=90] (0.8,-.9);
        \draw[-,thick,black] (-0.2,-.2) to[out=230,in=0] (-0.5,-.5);
        \draw[-,thick,black] (-0.5,-.5) to[out=180,in=-90] (-0.8,.5);
        %%% Odd bubble part
 \draw[-,thick,black,scale=1] (-0.5,0.6) to[out=180,in=90] (-.7,0.4);
  \draw[-,thick,black,scale=1] (-0.3,0.4) to[out=90,in=0] (-0.5,0.6);
 \draw[-,thick,black,scale=1] (-.7,0.4) to[out=-90,in=180] (-0.5,0.2);
  \draw[-,thick,black,scale=1] (-0.5,0.2) to[out=0,in=-90] (-0.3,0.4);
  \draw[-,thick,black,scale=1] (-0.36,0.25) to (-0.64,0.53);
  \draw[-,thick,black,scale=1] (-0.36,0.53) to (-0.64,0.25);
\end{tikzpicture}
\nn \\ \nn
&
\underset{on_{1,\lambda-2} +2 s_{\lambda,1}}{\Rrightarrow}
(-1)^\lambda \begin{tikzpicture}[baseline = 0,scale=0.7]
  \draw[-,thick,black] (-0.2,0.5) to[out=90, in=0] (-0.5,0.9);
	\draw[->,thick,black] (-0.5,0.9) to[out = 180, in = 90] (-0.8,0.5);
    \node at (1,0.5) {$\scriptstyle{\lambda}$};
	\draw[-,thick,black] (0.3,-.5) to (-0.2,.5);
	\draw[-,thick,black] (-0.2,-.2) to (0.2,.3);
        \draw[-,thick,black] (0.2,.3) to[out=50,in=180] (0.5,.5);
        \draw[->,thick,black] (0.5,.5) to[out=0,in=90] (0.8,-.9);
        \draw[-,thick,black] (-0.2,-.2) to[out=230,in=0] (-0.5,-.5);
        \draw[-,thick,black] (-0.5,-.5) to[out=180,in=-90] (-0.8,.5);
        \node at (0.15,-0.2) {$\bullet$};
\end{tikzpicture} + (-1)^{\lambda+1} \begin{tikzpicture}[baseline=0]
	\draw[<-,thick,black] (0.3,-.65) to[out=90, in=0] (0,-0.05);
	\draw[-,thick,black] (0,-0.05) to[out = 180, in = 90] (-0.3,-.65);
  \draw[-,thick,black] (0,0.45) to[out=180,in=90] (-.2,0.25);
  \draw[->,thick,black] (0.2,0.25) to[out=90,in=0] (0,.45);
 \draw[-,thick,black] (-.2,0.25) to[out=-90,in=180] (0,0.05);
  \draw[-,thick,black] (0,0.05) to[out=0,in=-90] (0.2,0.25);
   \node at (-0.55,0) {$\scriptstyle{\lambda}$};
\end{tikzpicture} -2 \begin{tikzpicture}[baseline = 0,scale=0.7]
  \draw[-,thick,black] (-0.2,0.5) to[out=90, in=0] (-0.5,0.9);
	\draw[->,thick,black] (-0.5,0.9) to[out = 180, in = 90] (-0.8,0.5);
    \node at (1,0.5) {$\scriptstyle{\lambda}$};
	\draw[-,thick,black] (0.3,-.5) to (-0.2,.5);
	\draw[-,thick,black] (-0.2,-.2) to (0.2,.3);
        \draw[-,thick,black] (0.2,.3) to[out=50,in=180] (0.5,.5);
        \draw[->,thick,black] (0.5,.5) to[out=0,in=90] (0.8,-.9);
        \draw[-,thick,black] (-0.2,-.2) to[out=230,in=0] (-0.5,-.5);
        \draw[-,thick,black] (-0.5,-.5) to[out=180,in=-90] (-0.8,.5);
        %%% Odd bubble part
 \draw[-,thick,black,scale=1] (0.5,1) to[out=180,in=90] (0.3,0.8);
  \draw[-,thick,black,scale=1] (0.7,0.8) to[out=90,in=0] (0.5,1);
 \draw[-,thick,black,scale=1] (0.3,0.8) to[out=-90,in=180] (0.5,0.6);
  \draw[-,thick,black,scale=1] (0.5,0.6) to[out=0,in=-90] (0.7,0.8);
  \draw[-,thick,black,scale=1] (0.64,0.65) to (0.36,0.93);
  \draw[-,thick,black,scale=1] (0.64,0.93) to (0.36,0.65);
\end{tikzpicture}
\end{align}
and the proof that the final result of the bottom sequence is the same as the one obtained in the top branch is made similarly, using bubble slide through a downward strand to make the odd bubble go back into the regular bubble before applying the $3$-cell $2 \sum\limits_{n=0}^\lambda \beta_{-n-1,1_{1_{\lambda+2}}}$.

\medskip
\noindent \textbf{Critical branchings $(D_{\lambda}, on_2^\lambda)$}
\[
\hspace{-2cm} \xy
  (-40,10)*+{ \cbddot{}
       }="TL";
  (40,10)*+{ \sum\limits_{n=0}^{\lambda} (-1)^n\ruledbul{}{n}  }="TR";
     {\ar@<.6pc> ^-{
        \Big\{  D_\lambda \Big\}}  "TL";"TR"};
     {\ar@<-0.6pc>_-{
        \Big\{  X  \Big\} }    "TL";"TR"};
\endxy
\]
where $X$ is the rewriting path defined by
\[ \{ \text{SInt}, \: (-1)^{\lambda+1} on_{2,\lambda -2}, \: (-1)^\lambda \delta, \: i_\lambda^3 \cdot {}_\ast D_\lambda, \: (i_\lambda^3)^- \cdot -2(-1)^\lambda \sum\limits_{n=0}^{h} r_{\lambda,1}, \: -2(-1)^\lambda \sum\limits_{n=0}^\lambda \beta_{-n-1,1_{1_\lambda}}, \} \]
and the $3$-cell $\delta$ is defined as follows:
\begin{align*} \hspace{-1cm}

\end{align*}
and we then prove that we obtain the same result as in the top branch by using $i_\lambda^3$ on the first summand, creating an extra term that cancel the third summand above, and reducing the first summand with $D_\lambda$. After applying these $3$-cells, it remains
\[
\sum\limits_{n=0}^{\lambda} (-1)^n %
\]
and we prove after using $\sum\limits_{n=0}^\lambda (i_\lambda^3)^-$ to place all dots on the left that the final result is equal to the top result using similar arguments.

\medskip
\noindent  \textbf{Critical branchings $(\Gamma_\lambda, on_{1,\lambda})$}
This critical branching has source
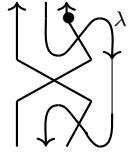
\[ \raisebox{0.6mm}{$%
$}
\]
One can use superinterchange and move the dots to the bottom of the diagram: this process gives minus a diagram on which we can apply $\Gamma$ with a dot at the bottom of the leftmost strand, and two extra terms on which we can apply the $3$-cell $E_\lambda$. By applying these $3$-cells, we obtain up to isotopy the following terms
\begin{align}
\label{E:CBGammaon1}
- \: \raisebox{0.8mm}{$%
 \nonumber
\end{align}
For the other branch, we first apply the $3$-cell $\Gamma$ and then move the dot to the bottom, creating two extra terms on which we can apply the $3$-cell $F_{\lambda+2}$, giving
\begin{align}
\label{E:CBGammaon2}
- \: \raisebox{0.6mm}{$%
 \nonumber
\end{align}

The first, fourth and fifth terms of \eqref{E:CBGammaon1} and \eqref{E:CBGammaon2} match. Moreover, one proves that extra terms in both \eqref{E:CBGammaon1} and \eqref{E:CBGammaon2} simplify to give
\[
\sum\limits_{r,s,t \geq 0} (-1)^{r+s}
%
 \]
Indeed, consider for instance the case $\lambda < 0$. In \eqref{E:CBGammaon2}, the third term reduces to $0$ using degree of bubble $3$-cells, and $6$th and $7$th terms are $0$ since sums are increasing.
In \eqref{E:CBGammaon1}, the third term also reduces to $0$. Moreover, changing variables to $s' = s+1$ and $t' = t-1$ gives a sum that is similar to the second element of \eqref{E:CBGammaon2}, up to extra terms given by $-$ the term for $s'=-1$ and $+$ the term for $t' = 0$, which cancel the $6$th and $7$th terms of \eqref{E:CBGammaon1}. We proceed similarly for $\lambda > 0$, where the second element of \eqref{E:CBGammaon1} reduces to $0$ using bubble $3$-cells.
Note that there is another critical branching implying $\Gamma$ and $on_1$, given by putting a dot on top of the other upward oriented strand, however this one would be proved confluent in a similar manner.

\medskip
\noindent \textbf{Critical branchings $(F_{\lambda}, on_{2,\lambda-2})$ and $(E_{\lambda}, on_{1,\lambda-2})$.} Let us denote by $on_{\lambda-2}$ the following composition of $3$-cells of $\ER^s$:
\[   \raisebox{3mm}{$\tdcrossrldtest{}{}$} \overset{on_{1,\lambda-2} \:}{\Rrightarrow} \: - \raisebox{3mm}{$\tdcrossrldm{}{}$} + \raisebox{3mm}{$\cuprb{}$} \overset{{}_\ast (-1)^h on_{2,\lambda-2}, \: \text{SInt} \:}{\Rrightarrow}  \: \raisebox{3mm}{$\tdcrossrldd{}{}$} + \raisebox{3mm}{$\cuprb{}$} - \raisebox{3mm}{$\acapl{}$}  \]

We then prove the critical branching $(F_{\lambda}, on_{2,\lambda-2})$ confluent modulo $E$ as follows:
\[   \xy
  (0,10)*+{ \tdcrosslrdd{}{}
       }="TL";
  (90,10)*+{  \sum\limits_{n=0}^{\lambda -1}  \sum\limits_{r = 0}^{\lambda-1} (-1)^{n+r} \stdbf{}{n+1}{-n-r-2}{r} - (-1)^{\lambda+1} \didupdowndl{}{}   }="TR";
     {\ar@<.6pc> ^-{
        \Big\{  F_\lambda, \: c'_\lambda \Big\}}  "TL";"TR"};
     {\ar@<-0.6pc>_-{
        \Big\{  \text{SInt}, \: on_{\lambda-2}, \: F_\lambda - A_\lambda + B_\lambda   \Big\} }    "TL";"TR"};
\endxy \]
where the $3$-cell $c'_\lambda$ is the $3$-cell defined in \ref{subsubsec:additional3cells}. Similarly, the critical branching $(E_{\lambda}, on_{1,\lambda-2})$ is proved confluent modulo $E$ as follows:
\[
  \xy
  (-10,10)*+{ \tdcrossrldtest{}{}
       }="TL";
  (110,10)*+{  \sum\limits_{n=0}^{ -\lambda -1} (-1)^{n+r} \sum\limits_{r = 0}^{-\lambda-1} \stda{}{n+1}{-n-r-2}{r} - (-1)^{\lambda+1} \diddownupdl{}{}   }="TR";
     {\ar@<.6pc> ^-{
        \Big\{  E_\lambda, \: b'_\lambda \Big\}}  "TL";"TR"};
     {\ar@<-0.6pc>_-{
        \Big\{  \text{SInt}, \: (-1)^{\lambda+1} on_{2,\lambda-2}, \: {}_\ast - on_{1,\lambda-2}, \: E_\lambda - D_\lambda + B_\lambda   \Big\} }    "TL";"TR"};
\endxy
\]
where the $3$-cell $b'_\lambda$ is the $3$-cell defined in \ref{subsubsec:additional3cells}.

\medskip

\noindent \textbf{Critical branchings $((u'_{\lambda,0} \star_2 u_{\lambda,0})^- \cdot F_{\lambda +2}, on_{2,\lambda}). $}
Starting from here, whenever we write $A\equiv B$ in the source of a branching between 3-cells $f$ and $e \cdot g$ we take $A$ to be the source and $B$ to be the result after applying $e$ to $A$.

\[
 \xy
  (-50,10)*++{ \isocbadot{}{} \equiv (-1)^{\lambda+1} \isocbbdotb{}{}
       }="TL";
  (50,10)*++{  \sum\limits_{n=0}^{\lambda+2} (-1)^n  \: \:\negbubdfsl{}{-n-1}  \identdotsu{}{n}   }="TR";
     {\ar@<.6pc> ^-{X }  "TL";"TR"};
     {\ar@<-0.6pc>_-{Y}    "TL";"TR"};
\endxy  \]
where
\begin{align*}
  X &:=  \Big\{  on_2, :\  \eta_\l,  \: -dc+\sum\limits_{n=0}^{\lambda + 2} u_{\lambda,0} \Big\}
\\
 Y &:=
        \Big\{ (u'_{\lambda,0} \star_2 u_{\lambda,0})^{-} \cdot {}_\ast F_{\lambda + 2},\: (u'_{\lambda,0} \star_2 u_{\lambda,0})_* \cdot s_{\l,\l+2}',\: c_{\l+2}'  \Big\}.
\end{align*}
%and the  $3$-cell $C'_\lambda$ in $\ER^s$ is defined as the following $\star_2$-composition of rewriting steps of $\ER$:
and the $3$-cell $\eta_\l:=(u_{\l,0}^-)_* \cdot C_{\l+2}$ in $\ER^s$ is defined as the first rewriting step in the following $\star_2$-composition of rewriting steps of $\ER$.
\begin{equation} \label{eq:defClambdap}
 \isocbadotbc{} \overset{(u_{\lambda,0}^-)_*}{\equiv}  \isocbadotbcb{}
 \: \overset{C_{\lambda+2}}{\Rrightarrow} \:
 \sum\limits_{n=0}^{\lambda+2} (-1)^n  \:
  \hackcenter{
\begin{tikzpicture} [scale =0.5, thick]
    \draw[<- ](0,1.5) .. controls ++(0,.6) and ++(0,.6) .. (.85,1.5);
    \draw (0,1.5) .. controls ++(0,-.6) and ++(0,-.6) .. (.85,1.5);
      \draw (0,-.5) to (0,0) .. controls ++(0,.5)  and ++(0,.5)  .. (1,0);
      \draw[->] (1,0) .. controls ++(0,-.75) and ++(0,-.75) .. (2,0) to (2,1.75);
      \node at (.85,1.5) {$\bullet$};
      \node at (1, .9) {$\scriptstyle -n-1$};
     \node at (2.5,1) {$\scs \l$};
      \node at (0,0) {$\bullet$};
      \node at (-0.4, 0) {$\scriptstyle n$};
%      \node at (2,.75) {$\bullet$};
%      \node at (2.35,.75) {$\scriptstyle b$};
\end{tikzpicture}}
 \: \overset{{}_* u_{\lambda,0}}{\equiv}
 \sum\limits_{n=0}^{\lambda+2} (-1)^n
%\raisebox{0mm}{$\isocbbcrcslfbis{}{n}$} \! \identdotsufdcb{}{n}
\hackcenter{
\begin{tikzpicture} [scale =0.5, thick]
    \draw[<- ](0,1.5) .. controls ++(0,.6) and ++(0,.6) .. (.85,1.5);
    \draw (0,1.5) .. controls ++(0,-.6) and ++(0,-.6) .. (.85,1.5);
      \draw[->] (2.5,0) to (2.5,2);
      %\draw[->] (1,0) .. controls ++(0,-.75) and ++(0,-.75) .. (2,0) to (2,1.75);
      \node at (.85,1.5) {$\bullet$};
      \node at (1.35, .95) {$\scriptstyle -n-1$};
    \node at (3.35,1.4) {$\scs \l$};
      \node at (2.5,.45) {$\bullet$};
      \node at (2.1, .35) {$\scriptstyle n$};
\end{tikzpicture}}
\end{equation}

Introduce the shorthand
\[
g(n) := \negbubdffffsl{}{\lambda+2-n+\ast}  \identdotsufdc{}{n}
\qquad
h(n)  := \identdotsufdc{}{n} \negbubdffffsl{}{\lambda+2-n+\ast}.
\]
Then the result of the top branch (and the critical branching) is $\sum\limits_{n= 0}^{\l+2}(-1)^n g(n)$.

In the bottom branch, the result after applying the steps up to and including $(u'_{\lambda,0} \star_2 u_{\lambda,0})_*$ is
\[
h(0)-\sum\limits_{n=0}^{\lambda+1}\sum\limits_{r\geq 0} (-1)^{nr} g(n+r+1)
\]

We can write this as follows:
\begin{align}\nn
  h(0)-\sum\limits_{n=0}^{\lambda+1}\sum\limits_{r\geq 0} (-1)^{nr} g(n+r+1)
  =
  h(0)-(\sum\limits_{z\geq 1} 2z g(2z)+\sum\limits_{z\geq 0}g(2z+1))
  \\
  =
  \qquad h(0)-\sum\limits_{z\geq 0}(2z+1)g(2z)+\sum\limits_{t\geq 0}(-1)^tg(t)
\end{align}

Then $h(0)-\sum\limits_{z\geq 0}(2z+1)g(2z) \overset{s_{\l,\l+2}'}{\Rrightarrow} 0$, so we have that
\begin{align}\nn
h(0)-\sum\limits_{n=0}^{\lambda+1}\sum\limits_{r\geq 0} (-1)^{nr} g(n+r+1)\overset{s_{\l,\l+2}'}{\Rrightarrow} \sum\limits_{t\geq 0}(-1)^t g(t)
 \overset{c_{\l+2}'}{\Rrightarrow} \sum\limits_{t=0}^{\l+2}(-1)^t g(t)=
 \sum\limits_{n=0}^{\lambda+2} (-1)^t
\hackcenter{
\begin{tikzpicture} [scale =0.5, thick]
    \draw[<- ](0,1.5) .. controls ++(0,.6) and ++(0,.6) .. (.85,1.5);
    \draw (0,1.5) .. controls ++(0,-.6) and ++(0,-.6) .. (.85,1.5);
      \draw[->] (2.5,0) to (2.5,2);
      %\draw[->] (1,0) .. controls ++(0,-.75) and ++(0,-.75) .. (2,0) to (2,1.75);
      \node at (.85,1.5) {$\bullet$};
      \node at (1.35, .95) {$\scriptstyle -t-1$};
    \node at (3.35,1.4) {$\scs \l$};
      \node at (2.5,.45) {$\bullet$};
      \node at (2.1, .35) {$\scriptstyle t$};
\end{tikzpicture}}
\end{align}

\medskip

\noindent \textbf{Critical branching $(F_\lambda, {}_* (u'_{\lambda+2,0} \star_2 u_{\lambda+2,0} \star_2 yb \star_2 u_{\l+2,0}^- ) \cdot {}_* C_{\lambda+4})$} Consider the critical branching
\begin{equation} \label{eq:branching0.1.8}
 \hspace{-2cm} \xy
  (-40,10)*++{ \isocbbcrsl{}{}
  \: \raisebox{-3mm}{$\equiv$} \:
   \begin{tikzpicture}[baseline = 0,scale=0.7]
		\draw[-,thick,black] (-0.2,.2) to[out=130,in=0] (-0.5,.5);
		\draw[->,thick,black] (-0.5,.5) to[out=180,in=-270] (-0.8,-.5);
		\draw[-,thick,black] (-0.5,-1.5) to[out=180,in=-90] (-0.8,-.5);
		\draw[-,thick,black] (-0.2,-1.2) to[out=230,in=0] (-0.5,-1.5);
			\draw[-,thick,black] (-0.2,0.2) to (-0.2,-0.2);
			\draw[-,thick,black] (-0.2,-1.2) to (-0.2,-0.8);
			\draw[-,thick,black] (-0.2,-0.2) to (0.3,-0.8) to (1.2,-1.5);
			\draw[-,thick,black] (-0.2,-0.8) to (0.2,-0.2);
			\draw[-,thick,black] (0.2,-0.2) to (0.2,-0.1);
			\draw[-,thick,black] (0.2,-0.8) to (0.2,-0.9);
			% Isotopies haut
			\draw[thick,-,black] (0.7,-0.1) to[out=90, in=0] (0.4,0.3);
			\draw[thick,-,black] (0.4,0.3) to[out = 180, in = 90] (0.2,-0.1);
			\draw[thick,-,black] (0.7,-0.1) to (0.7,-0.3);		
	\draw[-,thick,black] (1.2,-0.3) to[out=-90, in=0] (1,-0.6);
	\draw[-,thick,black] (1,-0.6) to[out = 180, in = -90] (0.7,-0.3);
	\draw[-,thick,black] (1.2,0.4) to (1.2,0.5);
	\draw[->,thick,black] (1.2,0.5) to (1.8,1.1);
	\draw[<-,thick,black] (1.2,1.1) to (1.8,0.5);
	\draw[-,thick,black] (1.2,-0.3) to (1.2,0.4);
	\draw[-,thick,black] (1.8,-0.7) to (1.8,0.5);
	      % Isotopies bas
	      	\draw[-,thick,black] (1.2,-1.5) to (1.8,-2.1);
	      	\draw[-,thick,black] (1.8,-1.5) to (1.2,-2.1);
	      	\draw[-,thick,black] (1.8,-1.5) to (1.8,-0.7);
		\node at (2.2,0) {$\scriptstyle{\lambda}$};
		\end{tikzpicture}
       }="TL";
  (40,10)*++{
      -\sum\limits_{r=0}^{\lambda+3}  \sum\limits_{s=0}^{r} (-1)^{r+s} \:
\hackcenter{
\begin{tikzpicture} [scale =0.5, thick]
    \draw[<- ](0,1.75) .. controls ++(0,.6) and ++(0,.6) .. (.85,1.75);
    \draw (0,1.75) .. controls ++(0,-.6) and ++(0,-.6) .. (.85,1.75);
      \draw (2.5,0) to (2.5,2);
      \draw (3.5,0) to (3.5,2);
\draw[->] (3.5,2) to (2.5,2.5) to (2.5,3);
\draw[->] (2.5,2) to (3.5,2.5) to (3.5,3);
      %\draw[->] (1,0) .. controls ++(0,-.75) and ++(0,-.75) .. (2,0) to (2,1.75);
      \node at (.85,1.75) {$\bullet$};
      \node at (1.45,1.3) {$\scriptstyle -r-2$};
    \node at (4.5,1.4) {$\scs \l$};
      \node at (2.5,.65) {$\bullet$};
      \node at (1.75, .55) {$\scriptstyle r-s$};
      \node at (3.5,.35) {$\bullet$};
      \node at (3.9,.45) {$\scriptstyle s$};
\end{tikzpicture}}
  }="TR";
     {\ar@<.6pc> ^-{ X}  "TL";"TR"};
     {\ar@<-0.6pc>_-{ Y }    "TL";"TR"};
\endxy
\end{equation}

with
\begin{align*}
  X&:= \Big\{ (u'_{\lambda+2,0} \star_2 u_{\lambda+2,0} \star_2 yb\star_2 u_{\l+2,0}^-) \cdot {}_* C_{\lambda+4}, \: {}_* u_{\l+2,0}\cdot \sum\limits_{n=0}^{\lambda+4} (-1)^{n} on_1^n,  \: \sum\limits_{n=0}^{\lambda+4} dc^\lambda  \Big\}
\\
  Y&:= \Big\{  \text{SInt}, \: - F_{\lambda + 4}, \: {}_\ast (u'_{\lambda+2,0} \star_2 u_{\lambda+2,0}) \cdot \sum\limits_{n=0}^{\lambda+3} \sum\limits_{r \geq 0} (-1)^{n+r+1+r\lambda} on_1^r, \text{SInt}, \: \gamma \Big\}
\end{align*}
where the $3$-cell $\eta_\l$ is defined in \eqref{eq:defClambdap} and the 3-cell $\gamma$ in the bottom branch will be defined in \eqref{E:ReductionLastCriticalBranching1} below.

Let us denote by $f(a,b)$ the monomial
\[ f(a,b)\;\;  := \;\;
\negbubdffffsl{}{\: \: \: \: \: \lambda+2 - (a+b) + \ast} \raisebox{-6mm}{$\isocbbcrcsldbis{}{}{a}{b}$}  \]
Then using $s_{\l+2,\l+2-a}'$, we get: 
\begin{equation} \label{eq:fsimp}
\sum\limits_{a=0}^{\lambda+2} \sum\limits_{b \geq 0} (-1)^a (2b+1)f(2b,a)
\: \overset{s_{\l+2,\l+2-a}'}{\Rrightarrow} \:
\sum\limits_{n=0}^{\lambda+2} (-1)^n %\raisebox{-8mm}{$\isocbbcrcslf{}{}{n}$} \! \identdotsufdcb{}{n}
\hackcenter{
\begin{tikzpicture} [scale =0.5, thick]
    \draw[<- ](0,1.5) .. controls ++(0,.6) and ++(0,.6) .. (.85,1.5);
    \draw (0,1.5) .. controls ++(0,-.6) and ++(0,-.6) .. (.85,1.5);
      \draw[->] (2.5,0) to (2.5,2);
    \draw[->] (-1,0) to (-1,2);
      %\draw[->] (1,0) .. controls ++(0,-.75) and ++(0,-.75) .. (2,0) to (2,1.75);
      \node at (.85,1.5) {$\bullet$};
      \node at (1.35, .95) {$\scriptstyle -n-1$};
    \node at (3.35,1.4) {$\scs \l$};
      \node at (2.5,.45) {$\bullet$};
      \node at (2.1, .35) {$\scriptstyle n$};
\end{tikzpicture}}
\end{equation}
and thus in particular we get that there is rewriting sequence $\gamma:=\{\eta_\l, u_{\l+2,0}\cdot s_{\l+2,\l+2-a}'\}$ of $\ER^s$ obtained from $\eta_\l$ and \eqref{eq:fsimp} as follows:
\begin{equation}
\label{E:ReductionLastCriticalBranching1}
\identusl{} \: \isocbadotbc{}
- \sum\limits_{a=0}^{\lambda+2} \sum\limits_{b \geq 0} (-1)^a (2b+1)f(2b,a) \overset{\gamma}{\Rrightarrow} 0
\end{equation}

Note that before applying the $3$-cell $\gamma$  in the bottom branch, we have obtained the polynomial
%\[
%\hspace{-1cm} (-1)^{h+1} \identusl{} \: \isocbadotbc{} - \sum\limits_{n=0}^{h+3} \sum\limits_{r \geq 0} (-1)^n \negbubdffffsl{}{\: \: \: \: \: h+3 - (n+r) + \ast} \raisebox{-6mm}{$\isocbbcrcsldcr{}{}{n}{r}$} + (-1)^h \sum\limits_{n=0}^{h+3} \sum\limits_{r \geq 1} \sum\limits_{s=0}^{r-1} (-1)^{s+sn+rn} \: \negbubdffffsl{}{\: \: \: \: \: \: h+3 - (n+r) + \ast}   \raisebox{-6mm}{$\isocbbcrcsldbisfl{}{}{n+r-1-s}{s}$}
%\]
\begin{equation}\label{eq:Bcancel}
 (-1)^{\lambda+1}
\hackcenter{
\begin{tikzpicture} [scale =0.5, thick]
\draw[->] (.5,-.25) to (.5,2.5);
\draw  (3.5,-.25) to (3.5,.75) .. controls ++(0,.3) and  ++(0,-.3) .. (2.5,1.5) ;
\draw[->]   (2.5,.75) .. controls ++(0,.3) and  ++(0,-.3) .. (3.5,1.5) to (3.5,2.5) ;
\draw[->] (2.5,1.5) .. controls ++(0,.65) and  ++(0,.65) ..(1.5,1.5) to (1.5,.75);
\draw (1.5,.75) .. controls ++(0,-.65) and  ++(0,-.65) ..(2.5,.75);
    \node at (4,1.4) {$\scs \l$};
\end{tikzpicture}}
\;
 -
\sum\limits_{n=0}^{\lambda+3} \sum\limits_{r \geq 0} (-1)^n
\hackcenter{
\begin{tikzpicture} [scale =0.5, thick]
    \draw[<- ](0,1.75) .. controls ++(0,.6) and ++(0,.6) .. (.85,1.75);
    \draw (0,1.75) .. controls ++(0,-.6) and ++(0,-.6) .. (.85,1.75);
      \draw (2.5,0) to (2.5,2);
      \draw (3.5,0) to (3.5,2);
\draw[->] (3.5,2) to (2.5,2.5) to (2.5,3);
\draw[->] (2.5,2) to (3.5,2.5) to (3.5,3);
      %\draw[->] (1,0) .. controls ++(0,-.75) and ++(0,-.75) .. (2,0) to (2,1.75);
      \node at (.85,1.75) {$\bullet$};
      \node at (.8, .8) {$\scriptstyle \overset{\scs \lambda+3-n}{\scs -r+\ast}$};
    \node at (4.5,1.4) {$\scs \l$};
      \node at (2.5,.65) {$\bullet$};
      \node at (2.2, .85) {$\scriptstyle n$};
      \node at (3.5,.35) {$\bullet$};
      \node at (3.9,.45) {$\scriptstyle r$};
\end{tikzpicture}}
+
(-1)^\lambda \sum\limits_{n=0}^{\lambda+3} \sum\limits_{r \geq 1} \sum\limits_{s=0}^{r-1} (-1)^{s+sn+rn}
\hackcenter{
\begin{tikzpicture} [scale =0.5, thick]
    \draw[<- ](0,2.5) .. controls ++(0,.6) and ++(0,.6) .. (.85,2.5);
    \draw (0,2.5) .. controls ++(0,-.6) and ++(0,-.6) .. (.85,2.5);
      \draw (2.5,0) to (2.5,2);
      \draw (3.5,0) to (3.5,2);
\draw[->] (2.5,2) to (2.5,2.5) to (2.5,3);
\draw[->] (3.5,2) to (3.5,2.5) to (3.5,3);
      %\draw[->] (1,0) .. controls ++(0,-.75) and ++(0,-.75) .. (2,0) to (2,1.75);
      \node at (.85,2.5) {$\bullet$};
      \node at (.8, 1.55) {$\scriptstyle \overset{\scs \lambda+3-n}{\scs -r+\ast}$};
    \node at (4.5,1.4) {$\scs \l$};
      \node at (2.5,.65) {$\bullet$};
      \node at (1.35, .5) {$\scriptstyle n+r-1$};
      \node at (3.5,.35) {$\bullet$};
      \node at (3.9,.45) {$\scriptstyle r$};
\end{tikzpicture}}
\end{equation}
We now show that  the first summand of \eqref{eq:Bcancel} cancels the third  using the 3-cell $\gamma$ from \eqref{E:ReductionLastCriticalBranching1}.
Using the $3$-cells $c_\lambda$ to remove the terms containing bubbles of negative degree, the last term reduces to
\begin{align*}
(-1)^\lambda
&\sum\limits_{n=0}^{\lambda+2} \sum\limits_{r \geq 0}^{\lambda+2-n} \sum\limits_{s=0}^{r} (-1)^{s+sn+rn+n} f(n+r-s,s)
 \\
&\quad =(-1)^\lambda
\sum\limits_{n=0}^{\lambda+2} \sum\limits_{r' =0}^{\lambda+2-n} \sum\limits_{a=0}^{\lambda+2-n-r'} (-1)^{a+an+n+(\lambda+2-n-r')n} f(\lambda+2-r'-a,a)
\\
& \quad =   (-1)^\lambda
 \sum\limits_{a=0}^{\lambda+2 } \sum\limits_{r' =0}^{\lambda+2-a} \sum\limits_{n=0}^{\lambda+2-r'-a}(-1)^{a+an+(\lambda+2-r')n} f(\lambda+2-r'-a,a)
\\
& \quad =
  (-1)^\lambda
 \sum\limits_{a=0}^{\lambda+2 } \sum\limits_{r' =0}^{\lambda+2-a} (-1)^{a}f(\lambda+2-r'-a,a) \left( \sum\limits_{n=0}^{\lambda+2-r'-a}(-1)^{(\lambda+2-r'-a)n}  \right)
\end{align*}
where we set $r'=\lambda+2-n-r$ and $s=a$ in the second equality and exchanged the summation order in the third.
Now let $b'=\lambda+2-a-r'$
\begin{align}
  &=  (-1)^\lambda
 \sum\limits_{a=0}^{\lambda+2 } \sum\limits_{b' =0}^{\lambda+2-a} (-1)^{a}f(b',a) \left( \sum\limits_{n=0}^{b'}(-1)^{b' n}  \right)
\end{align}
When $b'$ is odd the $n$ summation gives zero, but when $b'$ is even it gives a coefficient  $b+1$.  Keeping only the nonzero terms gives
\[
=  (-1)^\lambda
 \sum\limits_{a=0}^{\lambda+2 } \sum\limits_{b =0} (-1)^{a} (2b+1) f(2b,a)
\]
so that
\[
(-1)^\lambda
\sum\limits_{n=0}^{\lambda+2} \sum\limits_{r \geq 0}^{\lambda+2-n} \sum\limits_{s=0}^{r} (-1)^{s+sn+rn+n} f(n+r-s,s)
= (-1)^\lambda \sum\limits_{a=0}^{\lambda+2} \sum\limits_{b \geq 0} (-1)^a (2b+1)f(2b,a).
\]

Therefore, after applying the 3-cell $\gamma$ from  \eqref{E:ReductionLastCriticalBranching1}  to \eqref{eq:Bcancel}, only the second term remains
\begin{align*}
- \sum\limits_{n=0}^{\lambda+3} \sum\limits_{r \geq 0} (-1)^n
\:
\hackcenter{
\begin{tikzpicture} [scale =0.5, thick]
    \draw[<- ](0,1.75) .. controls ++(0,.6) and ++(0,.6) .. (.85,1.75);
    \draw (0,1.75) .. controls ++(0,-.6) and ++(0,-.6) .. (.85,1.75);
      \draw (2.5,0) to (2.5,2);
      \draw (3.5,0) to (3.5,2);
\draw[->] (3.5,2) to (2.5,2.5) to (2.5,3);
\draw[->] (2.5,2) to (3.5,2.5) to (3.5,3);
      %\draw[->] (1,0) .. controls ++(0,-.75) and ++(0,-.75) .. (2,0) to (2,1.75);
      \node at (.85,1.75) {$\bullet$};
      \node at (.8, .8) {$\scriptstyle \overset{\scs \lambda+3-n}{\scs -r+\ast}$};
    \node at (4.5,1.4) {$\scs \l$};
      \node at (2.5,.65) {$\bullet$};
      \node at (2.2, .85) {$\scriptstyle n$};
      \node at (3.5,.35) {$\bullet$};
      \node at (3.9,.45) {$\scriptstyle r$};
\end{tikzpicture}}
 &
%= - \sum\limits_{k=0}^{h+3} \sum\limits_{r \geq 0} (-1)^{k-r} \:
%\hackcenter{
%\begin{tikzpicture} [scale =0.5, thick]
%    \draw[<- ](-.25,1.5) .. controls ++(0,.6) and ++(0,.6) .. (.6,1.5);
%    \draw (-.25,1.5) .. controls ++(0,-.6) and ++(0,-.6) .. (.6,1.5);
%      \draw[->] (2.5,0) to (2.5,2);
%      \draw[->] (3.5,0) to (3.5,2);
%\draw (3.5,0) to (2.5,-.5) to (2.5,-.75);
%\draw (2.5,0) to (3.5,-.5) to (3.5,-.75);
%      \node at (.6,1.5) {$\bullet$};
%      \node at (.35, .6) {$\scriptstyle \overset{\scs h+3}{\scs -k+\ast}$};
%    \node at (4.5,1.4) {$\scs \l$};
%      \node at (2.5,.65) {$\bullet$};
%      \node at (1.8, .85) {$\scriptstyle k-r$};
%      \node at (3.5,.35) {$\bullet$};
%      \node at (3.9,.45) {$\scriptstyle r$};
%\end{tikzpicture}}
 = - \sum\limits_{r=0}^{\lambda+3} \sum\limits_{s \geq 0}^r (-1)^{r+s} \:
\hackcenter{
\begin{tikzpicture} [scale =0.5, thick]
    \draw[<- ](0,1.85) .. controls ++(0,.6) and ++(0,.6) .. (.85,1.85);
    \draw (0,1.85) .. controls ++(0,-.6) and ++(0,-.6) .. (.85,1.85);
      \draw (2.5,0) to (2.5,2);
      \draw (3.5,0) to (3.5,2);
\draw[->] (3.5,2) to (2.5,2.5) to (2.5,3);
\draw[->] (2.5,2) to (3.5,2.5) to (3.5,3);
      %\draw[->] (1,0) .. controls ++(0,-.75) and ++(0,-.75) .. (2,0) to (2,1.75);
      \node at (.85,1.85) {$\bullet$};
      \node at (.8,1.1) {$\scriptstyle \lambda+3-r+\ast$};
    \node at (4.5,1.4) {$\scs \l$};
      \node at (2.5,.65) {$\bullet$};
      \node at (1.8, .45) {$\scriptstyle r-s$};
      \node at (3.5,.35) {$\bullet$};
      \node at (3.9,.45) {$\scriptstyle s$};
\end{tikzpicture}}
\end{align*}
agreeing with the result in \eqref{eq:branching0.1.8} of the top branch,  establishing that this critical branching is confluent modulo $E$.

\medskip
\noindent \textbf{Critical branching $(E_\lambda,(u_{\l}\star_2 yb \star_2 (u_{\l,0}^-)_* \star_2 u_{\l+2,0}^-)_* \cdot \Gamma_\l)$}
Recalling the definition of sideways crossings from \ref{eq:crossl-gen-cyc}, we describe a critical branching between $E_\l$ and $(u_{\l}\star_2 yb \star_2 (u_{\l,0}^-)_* \star_2 u_{\l+2,0}^-)_* \cdot \Gamma_\l$
\[
  \xy
  (-5,10)*+{   \hackcenter{\begin{tikzpicture} [scale =0.5, thick]
      \draw (0,-.5) .. controls ++(0,.25)  and ++(0,-.25) ..    (1,.25)
          (1,-.5) .. controls ++(0,.25)  and ++(0,-.25) ..    (0,.25);
     \draw (0,.25) -- (0,1.75)  (1,.25) -- (1,1.75)  (3,-.5) -- (3,.5);
     \draw[<-] (2,-.5) -- (2,.25);
     \draw (4,.5) .. controls ++(0,-.5) and ++(0,-.5) .. (5,.5) -- (5,3);
     \draw (3,.5) .. controls ++(0,.25) and   ++(0,-.25) ..  (4,1.25)
            (4,.5).. controls ++(0,.25)  and ++(0,-.25) .. (3,1.25);
    \draw (4,1.25) -- (4,2.25) .. controls ++(0,.25)  and ++(0,-.25) .. (3,3) -- (3,3.25);
    \draw (2,.25) -- (2,1.25) .. controls ++(0,.5) and   ++(0,.5) ..  (3,1.25);
    \draw (0,1.75) --  (0,3.25)
          (1,1.75)   -- (1,3.5)
          (4,1.25) -- (4,2.25) .. controls ++(0,.25)  and ++(0,-.25) ..  (3,3) ;
    \draw (3,2.25) .. controls ++(0,.25)  and ++(0,-.25) ..  (4,3) ;
   \draw (2,2.25) .. controls ++(0,-.5) and ++(0,-.5) .. (3,2.25);
   \draw (2,2.25) -- (2,3.5)  .. controls ++(0,.5) and ++(0,.5) .. (1,3.5);
    \draw (4,3) .. controls ++(0,.5) and ++(0,.5) .. (5,3);
    \draw[->] (0,3.25) .. controls ++(0,.5) and ++(0,-.65) .. (2,5);
\draw[->] (3,3.25) .. controls ++(0,.5) and ++(0,-.65) .. (1,5);
  \end{tikzpicture}}
  \equiv
(-1)^{\l+1}
  \hackcenter{\begin{tikzpicture} [scale =0.55, thick]
    \draw (0,0) -- (0,1.75)  (1,0) -- (1,.75);
    \draw[<-] (2,0) .. controls ++(0,.45)  and ++(0,-.25) .. (3,.75);
    \draw (3,0) .. controls ++(0,.45)  and ++(0,-.25) .. (2,.75);
    \draw (2,.75)  .. controls ++(0,.25)  and ++(0,-.25) .. (1,1.5)
          (1,.75)  .. controls ++(0,.25)  and ++(0,-.25) .. (2,1.5)
          (0,.75) -- (0,1.5) (3,.75) -- (3,1.5);
    \draw (3,1.5)  .. controls ++(0,.25)  and ++(0,-.25) .. (2,2.25)
          (2,1.5)  .. controls ++(0,.25)  and ++(0,-.25) .. (3,2.25)
            (0,1.5) -- (0,2.25) (1,1.5) -- (1,2.25);
    \draw (2,2.25)  .. controls ++(0,.25)  and ++(0,-.25) .. (1,3)
          (1,2.25)  .. controls ++(0,.25)  and ++(0,-.25) .. (2,3)
          (0,2.25) -- (0,3) (3,2.25) -- (3,3);
    \draw (0,3)  .. controls ++(0,.5)  and ++(0,.5) .. (1,3);
    \draw[->] (2,3) -- (2,3.35);
     \draw[->] (3,3) -- (3,3.35);
  \end{tikzpicture}}
       }="TL";
  (110,10)*+{ (-1)^{\l}
 \hackcenter{
\begin{tikzpicture} [scale =0.5, thick]
    \draw (0,0)  .. controls ++(0,.25)  and ++(0,-.25) .. (1,.75)
          (1,0)  .. controls ++(0,.25)  and ++(0,-.25) .. (0,.75)
           (3,0) -- (3,.75);
    \draw[<-] (2,0) -- (2,.75);
    \draw (1,.75)  .. controls ++(0,.5)  and ++(0,.5) .. (2,.75)
          (0,.75)  .. controls ++(0,.25)  and ++(0,-.25) .. (1,1.5)
          (3,.75)  .. controls ++(0,.25)  and ++(0,-.25) .. (2,1.5);
    \draw[->] (2,1.5) .. controls ++(0,.25)  and ++(0,-.25) .. (1,2.25)--(1,2.5);
    \draw[->] (1,1.5) .. controls ++(0,.25)  and ++(0,-.25) .. (2,2.25)--(2,2.5);
  \end{tikzpicture} }
+
\sum\limits_{x,y,r \geq 0} (-1)^{x+r}
\hackcenter{
\begin{tikzpicture} [scale =0.5, thick]
     \node at (2.95,-.5) {$\bullet$};  \node at (3.3,-.5) {$\scs  \scs r$};
    \draw[<-] (2,-.75) .. controls ++(0,.65)  and ++(0,.65) .. (3,-.75) ;
    \draw (0,-.75) -- (0,1.75)  (1,-.75) -- (1,1.75);
    \draw (0,.75) -- (0,1.5)  (1,.75) -- (1,1.5);
    \node at (2,.75) {$\bullet$};  \node at (2.4,.2) {$\scs  \scs -x-y  -r-3 $};
    \draw[->] (2,.75) .. controls ++(0,.65)  and ++(0,.65) .. (3,.75) ;
    \draw (2,.75) .. controls ++(0,-.5)  and ++(0,-.5) .. (3,.75) ;
    \node at (0,1.25) {$\bullet$};  \node at (-.4,1.25) {$\scs  \scs x$};
    \node at (1,1.75) {$\bullet$};  \node at (1.3,1.75) {$\scs  \scs y$};
    \draw (1,1.75)  .. controls ++(0,.25)  and ++(0,-.25) .. (0,2.5)
          (0,1.75)  .. controls ++(0,.25)  and ++(0,-.25) .. (1,2.5);
    \draw[->] (0,2.5) -- (0,2.75);
    \draw[->] (1,2.5) -- (1,2.75);
  \end{tikzpicture} }  }="TR";
     {\ar@<.6pc> ^-{
        \Big\{  X \Big\}}  "TL";"TR"};
     {\ar@<-0.6pc>_-{
        \Big\{ Y  \Big\} }    "TL";"TR"};
\endxy
\]
with
\begin{align*}
  X&:= \Big\{ E_\l,\: {}_* u_{\l,0} \cdot \sum\limits_{n=0}^{-\l-1} \sum\limits_{r\geq 0} {}_* on_2^n, \:dc   \Big\}
\\
  Y&:= \Big\{  \text{SInt}, \: (u_{\l}\star_2 yb \star_2 (u_{\l,0}^-)_* \star_2 u_{\l+2,0}^-)_* \cdot \Gamma_\l, \: \Omega \Big\}
\end{align*}

Where $\Omega$ is the rewriting sequence described below.
The 3-cell
$\Omega$ has source
\[
\hackcenter{
\begin{tikzpicture} [scale =0.55, thick]
    \draw (0,0) -- (0,1.75)  (3,0) -- (3,.75);
    \draw (1,0) .. controls ++(0,.45)  and ++(0,-.25) .. (2,.75);
    \draw[<-] (2,0) .. controls ++(0,.45)  and ++(0,-.25) .. (1,.75);
    \draw (3,.75)  .. controls ++(0,.25)  and ++(0,-.25) .. (2,1.5)
          (2,.75)  .. controls ++(0,.25)  and ++(0,-.25) .. (3,1.5)
          (0,.75) -- (0,1.5) (1,.75) -- (1,1.5);
    \draw (2,1.5)  .. controls ++(0,.25)  and ++(0,-.25) .. (1,2.25)
          (1,1.5)  .. controls ++(0,.25)  and ++(0,-.25) .. (2,2.25)
          (0,1.5) -- (0,2.25) (3,1.5) -- (3,2.25);
    \draw (2,2.25)  .. controls ++(0,.25)  and ++(0,-.25) .. (1,3)
          (1,2.25)  .. controls ++(0,.25)  and ++(0,-.25) .. (2,3)
          (0,2.25) -- (0,3) (3,2.25) -- (3,3);
    \draw (0,3)  .. controls ++(0,.5)  and ++(0,.5) .. (1,3);
    \draw[->] (2,3) -- (2,3.35);
     \draw[->] (3,3) -- (3,3.35);
\end{tikzpicture}}
\quad - \quad
\sum_{r,s,t \geq 0} (-1)^{r+s}\;
\hackcenter{
\begin{tikzpicture} [scale =0.55, thick]
     \node at (1.05,-.5) {$\bullet$};  \node at (.7,-.5) {$\scs  \scs s$};
    \draw[->] (1,-.75) .. controls ++(0,.65)  and ++(0,.65) .. (2,-.75) ;
    \draw (0,-.75) -- (0,.75)  (3,-.75) -- (3,.75);
    \draw (0,.75) -- (0,1.5)  (3,.75) -- (3,1.5);
    \node at (2,.75) {$\bullet$};  \node at (1.4,.2) {$\scs  \scs -r-s  -t-3 $};
    \draw[<-] (1,.75) .. controls ++(0,.65)  and ++(0,.65) .. (2,.75) ;
    \draw (1,.75) .. controls ++(0,-.5)  and ++(0,-.5) .. (2,.75) ;
    \node at (1,2) {$\bullet$};  \node at (.6,2) {$\scs r$};
    \node at (3,2.25) {$\bullet$};  \node at (3.3,2.25) {$\scs t$};
    \draw   (1,2.25) -- (1,2)  .. controls ++(0,-.5)  and ++(0,-.5) .. (2,2) --(2,2.25)
          (0,1.5)  -- (0,2.25) (3,1.5)  -- (3,2.25);
    \draw (2,2.25)  .. controls ++(0,.25)  and ++(0,-.25) .. (1,3)
          (1,2.25)  .. controls ++(0,.25)  and ++(0,-.25) .. (2,3)
          (0,2.25) -- (0,3) (3,2.25) -- (3,3);
    \draw (0,3)  .. controls ++(0,.5)  and ++(0,.5) .. (1,3);
    \draw[->] (2,3) -- (2,3.35);
    \draw[->] (3,3) -- (3,3.35);
  \end{tikzpicture} }
\quad + \quad
\sum_{r,s,t \geq 0} (-1)^{r+s}\;
\hackcenter{
\begin{tikzpicture} [scale =0.55, thick]
     \node at (2.95,-.5) {$\bullet$};  \node at (3.3,-.5) {$\scs  \scs r$};
    \draw[<-] (2,-.75) .. controls ++(0,.65)  and ++(0,.65) .. (3,-.75) ;
    \draw (0,-.75) -- (0,.75)  (1,-.75) -- (1,.75);
    \draw (0,.75) -- (0,1.5)  (1,.75) -- (1,1.5);
    \node at (2,.75) {$\bullet$};  \node at (2.4,.2) {$\scs  \scs -r-s  -t-3 $};
    \draw[->] (2,.75) .. controls ++(0,.65)  and ++(0,.65) .. (3,.75) ;
    \draw (2,.75) .. controls ++(0,-.5)  and ++(0,-.5) .. (3,.75) ;
    \node at (1,2.25) {$\bullet$};  \node at (.6,2.25) {$\scs t$};
    \node at (3,1.9) {$\bullet$};  \node at (3.3,1.9) {$\scs s$};
    \draw   (2,2.25) -- (2,2)  .. controls ++(0,-.5)  and ++(0,-.5) .. (3,2) --(3,2.25)
          (0,1.5)  -- (0,2.25) (1,1.5)  -- (1,2.25);
    \draw (2,2.25)  .. controls ++(0,.25)  and ++(0,-.25) .. (1,3)
          (1,2.25)  .. controls ++(0,.25)  and ++(0,-.25) .. (2,3)
          (0,2.25) -- (0,3) (3,2.25) -- (3,3);
    \draw (0,3)  .. controls ++(0,.5)  and ++(0,.5) .. (1,3);
    \draw[->] (2,3) -- (2,3.35);
    \draw[->] (3,3) -- (3,3.35);
  \end{tikzpicture} }
\]
The second term rewrites to $0$ by either $A_{\l+2}''$ for $\l>-2$ or $c_{\l+2}^0$ for $\l \leq -2$.
The first term rewrites by $E_{\l+2}$ to
\[
(-1)^{\l}
 \hackcenter{
\begin{tikzpicture} [scale =0.55, thick]
    \draw (0,0)  .. controls ++(0,.25)  and ++(0,-.25) .. (1,.75)
          (1,0)  .. controls ++(0,.25)  and ++(0,-.25) .. (0,.75)
           (3,0) -- (3,.75);
    \draw[<-] (2,0) -- (2,.75);
    \draw (1,.75)  .. controls ++(0,.5)  and ++(0,.5) .. (2,.75)
          (0,.75)  .. controls ++(0,.25)  and ++(0,-.25) .. (1,1.5)
          (3,.75)  .. controls ++(0,.25)  and ++(0,-.25) .. (2,1.5);
    \draw[->] (2,1.5) .. controls ++(0,.25)  and ++(0,-.25) .. (1,2.25)--(1,2.5);
    \draw[->] (1,1.5) .. controls ++(0,.25)  and ++(0,-.25) .. (2,2.25)--(2,2.5);
  \end{tikzpicture} }
\]
plus an extra sum which is reduced to
\[
\sum\limits_{n=0}^{-\l-3}\sum\limits_{k\geq 0}
(-1)^{n+k}
\hackcenter{
\begin{tikzpicture}[scale=0.55,thick]
  \draw[black, ->] (-1,0) to (-1,3);
  \node at (-1,2.5) {$\bullet$};
  \node at (-.75,2.5) {$\scriptstyle{n}$};
  \draw[black,<-] (1.25,1.5) arc (0:370:.25);
  \node at (.75,1.5) {$\bullet$};
  \node at (.1,1.6) {$\color{black}\scriptstyle{\substack{-n-k \\ -2}}$};
  \draw[black, ->] (2,0) to (2,3);
  \draw[black,->] (3.25,0) arc (0:180:.5);
  \node at (3.2,.2) {$\bullet$};
  \node at (3.5,.2) {$\scriptstyle{k}$};
\end{tikzpicture}
}
\]
via the rewriting sequence
$\{{}_*u_{\l+2,0} \cdot on_1, (u'_{\l,0})^- \cdot {}_* F_{\l+2}+ C_{\l+2}' \}$.

Hence, first term rewrites to:
\begin{equation}\label{eq:EybI}
  (-1)^{\l}
 \hackcenter{
\begin{tikzpicture} [scale =0.55, thick]
    \draw (0,0)  .. controls ++(0,.25)  and ++(0,-.25) .. (1,.75)
          (1,0)  .. controls ++(0,.25)  and ++(0,-.25) .. (0,.75)
           (3,0) -- (3,.75);
    \draw[<-] (2,0) -- (2,.75);
    \draw (1,.75)  .. controls ++(0,.5)  and ++(0,.5) .. (2,.75)
          (0,.75)  .. controls ++(0,.25)  and ++(0,-.25) .. (1,1.5)
          (3,.75)  .. controls ++(0,.25)  and ++(0,-.25) .. (2,1.5);
    \draw[->] (2,1.5) .. controls ++(0,.25)  and ++(0,-.25) .. (1,2.25)--(1,2.5);
    \draw[->] (1,1.5) .. controls ++(0,.25)  and ++(0,-.25) .. (2,2.25)--(2,2.5);
  \end{tikzpicture} }
  \;\;+\;\;
\sum\limits_{n=0}^{-\l-3}\sum\limits_{k\geq 0}
(-1)^{n+k}
\hackcenter{
\begin{tikzpicture}[scale=0.55,thick]
  \draw[black, ->] (-1,0) to (-1,3);
  \node at (-1,2.5) {$\bullet$};
  \node at (-.75,2.5) {$\scriptstyle{n}$};
  \draw[black,<-] (1.25,1.5) arc (0:370:.25);
  \node at (.75,1.5) {$\bullet$};
  \node at (.1,1.6) {$\color{black}\scriptstyle{\substack{-n-k \\ -2}}$};
  \draw[black, ->] (2,0) to (2,3);
  \draw[black,->] (3.25,0) arc (0:180:.5);
  \node at (3.2,.2) {$\bullet$};
  \node at (3.5,.2) {$\scriptstyle{k}$};
\end{tikzpicture}
}
\end{equation}

For the third term, we use super isotopy and $on_2$ to move the $s$ dots through the sideways crossing and then move them below the $t$ dots to obtain
\[
\sum_{x,y,r \geq 0} (-1)^{x+r} \;
\hackcenter{
\begin{tikzpicture} [scale =0.55, thick]
     \node at (2.95,-.5) {$\bullet$};  \node at (3.3,-.5) {$\scs  \scs r$};
    \draw[<-] (2,-.75) .. controls ++(0,.65)  and ++(0,.65) .. (3,-.75) ;
    \draw (0,-.75) -- (0,1.75)  (1,-.75) -- (1,1.75);
    \draw (0,.75) -- (0,1.5)  (1,.75) -- (1,1.5);
    \node at (2,.75) {$\bullet$};  \node at (2.4,.2) {$\scs  \scs -x-y  -r-3 $};
    \draw[->] (2,.75) .. controls ++(0,.65)  and ++(0,.65) .. (3,.75) ;
    \draw (2,.75) .. controls ++(0,-.5)  and ++(0,-.5) .. (3,.75) ;
    \node at (0,1.25) {$\bullet$};  \node at (-.4,1.25) {$\scs  \scs x$};
    \node at (1,1.75) {$\bullet$};  \node at (1.3,1.75) {$\scs  \scs y$};
    \draw (1,1.75)  .. controls ++(0,.25)  and ++(0,-.25) .. (0,2.5)
          (0,1.75)  .. controls ++(0,.25)  and ++(0,-.25) .. (1,2.5);
    \draw[->] (0,2.5) -- (0,2.75);
    \draw[->] (1,2.5) -- (1,2.75);
\end{tikzpicture} }
\;\;+\;\;
\sum\limits_{a,b,r,t\geq 0}
(-1)^{r+a+bt-t}
\hackcenter{\begin{tikzpicture} [scale =0.55, thick]
     \node at (2.95,-.5) {$\bullet$};  \node at (3.3,-.5) {$\scs  \scs r$};
    \draw[<-] (2,-.75) .. controls ++(0,.65)  and ++(0,.65) .. (3,-.75) ;
    \draw (0,-.75) -- (0,1.75)  (1,-.75) -- (1,1.75);
    \draw (0,.75) -- (0,1.5)  (1,.75) -- (1,1.5);
    \node at (2,.75) {$\bullet$};  \node at (2.5,.2) {$\scs  \scs -a-b-t -r-4 $};
    \draw[->] (2,.75) .. controls ++(0,.65)  and ++(0,.65) .. (3,.75) ;
    \draw (2,.75) .. controls ++(0,-.5)  and ++(0,-.5) .. (3,.75) ;
    \node at (0,1.25) {$\bullet$};  \node at (-.4,1.25) {$\scs  \scs a$};
    \node at (1,1.75) {$\bullet$};  \node at (1.5,1.75) {$\scs  \scs b+t$};
    \draw [->] (1,1.75) to (1,2.75);
    \draw [->] (0,1.75) to (0,2.75);
\end{tikzpicture}}
\]
And one can check that the second term of this cancels with the second term of \ref{eq:EybI} once we apply bubble slides.

% - - - - - - - - - - - - - - - -
\subsection{Critical branchings from odd $sl(2$)-relations}
% - - - - - - - - - - - - - - - -
\label{appendix:criticalbranchingsfullOsl2}
We prove that the critical branching between two $3$-cells of the set
 $\{A_\l, B_\l, C_\l, D_\l, E_\l, F_\l, \Gamma_\l \}$ are confluent modulo $E$.

\medskip
\noindent \textbf{A and C }
\begin{enumerate}[{\bf i)}]

\item For $\lambda <0$:

\[
\xymatrix@R=1em@C=3em{\vcenter{\hbox{$\bcritacxy{}$}} \ar@<.6pc>[rr]^-{\{A_\lambda, c_\lambda', c_\lambda^1+b_\lambda^1\}} \ar@<-.6pc>[rr]_-{C_\lambda} & & 0}
\]

\item For $\lambda = 0$:

\[
\xymatrix{\bcritacxy{} \ar@<.6pc>[rr]^-{\{A_0, b_0^1\}} \ar@<-.6pc>[rr]_-{\{C_0, c_0^1,I_0\}} & & \raisebox{-3mm}{\oddbubble{}}}
\]

\item For $\lambda>0$, the calculation is similar to the case $\lambda<0$, except $A_\lambda$ takes it to 0 instead of $C_\lambda$.
\end{enumerate}

\noindent \textbf{A and F}
\begin{enumerate}[{\bf i)}]

\item For $\lambda <0$:

\[
\xymatrix{\vcenter{\hbox{$\bcritaf{}$}} \ar@<.6pc>[rr]^-{\{A_\lambda, \sum\limits_{n=0}^{-\lambda} D_{\lambda}',b_\lambda^1\}} \ar@<-.6pc>[rr]_-{F_\lambda} & & (-1)^{\lambda} \cupl{}}
\]

\item For $\lambda=0$:

\[
\xymatrix{\vcenter{\hbox{$\bcritaf{}$}} \ar@<.6pc>[rr]^-{\{A_0, b_0^1, D_0, c_0^1\}} \ar@<-.6pc>[rr]_-{F_0} & & \cupl{}}
\]

\item For $\lambda > 0$,

\[
\xymatrix@R=1em@C=4em{\vcenter{\hbox{$\bcritaf{}$}} \ar@<.6pc>[rr]^-{\{F_\lambda, b_\lambda',c_\lambda', b_\lambda^1,c_\lambda^1 \}} \ar@<-.6pc>[rr]_-{A_{\lambda}} & & 0}
\]

\end{enumerate}

\noindent \textbf{B and D}

\begin{enumerate}[{\bf i)}]

\item For $\lambda<0$,

\[
\xymatrix@R=1em@C=4em{\vcenter{\hbox{$\bcritbd{}$}} \ar@<.6pc>[rr]^-{\{B_\lambda, c_\lambda', c_\lambda^1+b_\lambda^1\}} \ar@<-.6pc>[rr]_-{D_\lambda} & & 0}
\]

\item For $\lambda=0$:

\[
\xymatrix{\bcritbd{} \ar@<.6pc>[rr]^-{\{B_0, b_0^1\}} \ar@<-.6pc>[rr]_-{\{D_0, c_0^1,I_0\}} & & \raisebox{-3mm}{\oddbubble{}}}
\]

\item For $\lambda >0$, we get a similar calculation as for $\lambda<0$ except $B_\lambda$ takes it to 0 instead of $D_\lambda$.

\end{enumerate}

\noindent \textbf{E and F}

\begin{enumerate}[{\bf i)}]

\item For $\lambda < 0$:

\[
\xymatrix{\vcenter{\hbox{$\bcritfe{}{}$}} \ar@<.6pc>[rr]^-{\{E_\lambda, \sum\limits_{n=0}^{-\lambda-1} \sum\limits_{r\geq 0} D_{\lambda}'\}} \ar@<-.6pc>[rr]_-{F_\lambda} & & (-1)^{\lambda} \tleftcrossw{}{}{} }
\]

\item For $\lambda=0$:
\[
\xymatrix{\vcenter{\hbox{$\bcritfe{}{}$}} \ar@<.6pc>[r]^-{E_0} \ar@<-.6pc>[r]_-{F_0}  & \tleftcrossw{}{}{} }
\]

\item For $\lambda > 0$:
\[
\xymatrix{\vcenter{\hbox{$\bcritfe{}{}$}} \ar@<.6pc>[rr]^-{\{F_\lambda, \sum\limits_{n=0}^{\lambda-1} \sum\limits_{r\geq 0} B_{\lambda}'\}} \ar@<-.6pc>[rr]_-{E_\lambda} & & (-1)^{\lambda} \tleftcrossw{}{}{} }
\]

\end{enumerate}

The other family of critical branchings with $F_\lambda$ and $E_\lambda$ would be proved to be confluent modulo $E$ in a similar manner.

\medskip
\noindent \textbf{B and F}
\begin{enumerate}[{\bf i)}]

\item For $\lambda > 0$:

\[
\xymatrix@R=1em@C=5em{ \vcenter{\hbox{$\bcritbf{}$}}  \ar@<-.6pc>[rr]_-{\{F_\lambda , b_\lambda' , c_\lambda', c_\lambda^1, b_\lambda^1\}} \ar@<.6pc>[rr]^-{B_\lambda}
& & 0}
\]

\item For $\lambda=0$:
\[
\xymatrix@R=1em@C=3em{ \vcenter{\hbox{$\bcritbf{}$}}  \ar@<-.6pc>[rr]_-{\{B_0 , b_0^1 , C_0 , c_0^1\}} \ar@<.6pc>[rr]^-{F_0} & & \capr{}}
\]

\item For $\lambda <0$:
\[
\xymatrix@R=1em@C=3em{ \vcenter{\hbox{$\bcritbf{}$}}  \ar@<-.6pc>[rr]_-{\{B_\lambda , \sum\limits_{n=0}^{-\lambda} C_{\lambda}', b_\lambda^1 \}} \ar@<.6pc>[rr]^-{F_\lambda} & & (-1)^{\lambda} \capr{}}
\]

\end{enumerate}

\noindent \textbf{E and D}

\begin{enumerate}[{\bf i)}]

\item For $\lambda \geq 0$:
\[
\xymatrix{\vcenter{\hbox{$\bcritde{}$}} \ar@<-.6pc>[rr]_-{\{D_\lambda ,\sum\limits_{n=0}^{\lambda} A_{\lambda} '', c_\lambda^1\}} \ar@<.6pc>[rr]^-{E_\lambda} & & (-1)^{\lambda}\cupr{}}
\]

\item For $\lambda<0$:

\[
\xymatrix@R=1em@C=5em{\vcenter{\hbox{$\bcritde{}$}} \ar@<-.6pc>[rr]_-{\{E_\lambda ,c_\lambda', b_\lambda', b_\lambda^1, c_\lambda^1 \}} \ar@<.6pc>[rr]^-{D_\lambda}
 & & 0}
\]

\end{enumerate}

\noindent \textbf{C and E}

\begin{enumerate}[{\bf i)}]

\item For $\lambda > 0$:

\[
\xymatrix{\vcenter{\hbox{$\bcritce{}$}}\ar@<-.6pc>[rr]_-{\{C_\lambda, \sum\limits_{n=0}^{\lambda} B_{\lambda}', c_\lambda^1\}} \ar@<.6pc>[rr]^-{E_0} & & (-1)^{\lambda}\capl{}}
\]

\item For $\lambda=0$:
\[
\xymatrix{\vcenter{\hbox{$\bcritce{}$}}\ar@<-.6pc>[rr]_-{\{C_0, c_0^1, B_0, b_0^1\}} \ar@<.6pc>[rr]^-{E_0} & & \capl{}}
\]

\item For $\lambda < 0$:

\[
\xymatrix@R=1em@C=5em{\vcenter{\hbox{$\bcritce{}$}} \ar@<-.6pc>[rr]_-{\{E_\lambda,  c_\lambda', b_\lambda', c_\lambda^1, b_\lambda^1 \}}
\ar@<.6pc>[rr]^-{C_\lambda}
 & & 0}
\]

\end{enumerate}

\medskip

\noindent \textbf{Critical branching $(\Gamma_\lambda, C_\lambda)$}
\[
\xy
  (-65,10)*+{ \begin{tikzpicture}[baseline = 0,scale=0.55]
	\draw[->,thick,black] (0.3,-.5) to (-0.3,.5);
	\draw[-,thick,black] (-0.2,-.2) to (0.2,.3);
        \draw[-,thick,black] (0.2,.3) to[out=50,in=180] (0.5,.5);
        \draw[->,thick,black] (0.5,.5) to[out=0,in=90] (0.8,-0.5);
        \draw[-,thick,black] (-0.2,-.2) to[out=230,in=0] (-0.5,-.4);
        \draw[-,thick,black] (-0.5,-.4) to[out=180,in=-90] (-0.8,.5);
        \draw[-,thick,black] (0.3,-1.5) to (-0.3,-2.6);
        \draw[-,thick,black] (0.5,-2.6) to[out=0,in=270] (0.8,-1.8);
        \draw[-,thick,black] (0.2,-2.3) to[out=130,in=180] (0.5,-2.6);
        \draw[-,thick,black] (-0.2,-1.8) to (0.2,-2.3);
        \draw[-,thick,black] (-0.2,-1.8) to[out=130,in=0] (-0.5,-1.6);
        \draw[->,thick,black] (-0.5,-1.6) to[out=180,in=-270] (-0.8,-2.6);
        \node at (1,0.5) {$\scriptstyle{\lambda}$};
        %% LEFT CROSSING PART
        \draw[-,thick,black] (0.3,-1.5) to (-1.5,-0.5);
        \draw[->,thick,black] (-1.5,-0.5) to (-1.5,0.8);
        \draw[-,thick,black]  (0.3,-0.5) to (-1.5,-1.5);
        \draw[-,thick,black] (-1.5,-1.5) to (-1.5,-2.6);
        \draw (0.8,-0.5) to (0.8,-2);
        %%% TOP CAP
     \draw[-,thick,black] (-0.2,0.4) to[out=90, in=0] (-0.5,0.8);
	\draw[->,thick,black] (-0.5,0.8) to[out = 180, in = 90] (-0.8,0.4);
\end{tikzpicture}
       }="TL";
  (65,10)*+{   \sum\limits_{n=0}^{ \lambda} (-1)^{n + \lambda}  \raisebox{-5mm}{$\begin{tikzpicture}
	\draw[<-,thick,black] (0.4,-0.1) to[out=90, in=0] (0.1,0.3);
	\draw[-,thick,black] (0.1,0.3) to[out = 180, in = 90] (-0.2,-0.1);
	\node at (-0.15,0.1) {$\bullet$};
	\node at (-0.35,0.1) {$\scriptstyle{n}$};
  \node at (1.4,0) {$\scriptstyle{\lambda}$};
  \draw[->,thick,black] (0.6,-0.1) to (0.6,1);
  %% Bubble part
    \draw[->,thick,black] (1.2,0.5) to[out=90,in=0] (1,.7);
  \draw[-,thick,black] (1,0.7) to[out=180,in=90] (0.8,0.5);
\draw[-,thick,black] (0.8,0.5) to[out=-90,in=180] (1,0.3);
  \draw[-,thick,black] (1,0.3) to[out=0,in=-90] (1.2,0.5);
   \node at (1.2,0.5) {$\color{black}\bullet$};
   \node at (1.7,0.5) {$\color{black}\scriptstyle{-n-1}$};
  \end{tikzpicture}$} }="TR";
     {\ar@<.6pc> ^-{
        \Big\{  C_\lambda, \: on_2^n, \: id + B'_\lambda, \: {}_\ast (u_{\lambda,0})^- \cdot E_\lambda \Big\}}  "TL";"TR"};
     {\ar@<-0.6pc>_-{
        \Big\{  \Gamma, \:  {}_\ast (u'_{\lambda,0})^- \cdot s_{\lambda, \lambda - s}' + \: {}_\ast (b'_\lambda \star_2 c'_\lambda),\: ({}_\ast (u'_{\lambda,0})^- \star_2 dc) \Big\} }    "TL";"TR"};
\endxy
\]
Note that there  is also a critical branching between $\Gamma_\l$ and $D_\l$ given by attaching to the source of $\Gamma_\l$ a rightward cup on bottom on the rightmost two strands. This one is proven confluent in a similar manner.

\smallskip

\noindent  \textbf{Critical branching $(\Gamma_\lambda, F_\lambda)$}
 \[
\xy
  (-65,10)*+{ \raisebox{0.6mm}{$\begin{tikzpicture}[baseline = 0,scale=0.45]
	\draw[-,thick,black] (0.3,-.5) to (-0.3,.5) to (-0.3,0.6);
	\draw[-,thick,black] (-0.2,-.2) to (0.2,.3);
        \draw[-,thick,black] (0.2,.3) to[out=50,in=180] (0.5,.5);
        \draw[->,thick,black] (0.5,.5) to[out=0,in=90] (0.8,-0.5);
        \draw[-,thick,black] (-0.2,-.2) to[out=230,in=0] (-0.5,-.4);
        \draw[-,thick,black] (-0.5,-.4) to[out=180,in=-90] (-0.8,.5);
        \draw[-,thick,black] (0.3,-1.5) to (-0.3,-2.6);
        \draw[-,thick,black] (0.5,-2.6) to[out=0,in=270] (0.8,-1.8);
        \draw[-,thick,black] (0.2,-2.3) to[out=130,in=180] (0.5,-2.6);
        \draw[-,thick,black] (-0.2,-1.8) to (0.2,-2.3);
        \draw[-,thick,black] (-0.2,-1.8) to[out=130,in=0] (-0.5,-1.6);
        \draw[->,thick,black] (-0.5,-1.6) to[out=180,in=-270] (-0.8,-2.6);
        \node at (1,0.5) {$\scriptstyle{\lambda}$};
        %% LEFT CROSSING PART
        \draw[-,thick,black] (0.3,-1.5) to (-1.5,-0.5);
        \draw[->,thick,black] (-1.5,-0.5) to (-1.5,1.6);
        \draw[-,thick,black]  (0.3,-0.5) to (-1.5,-1.5);
        \draw[-,thick,black] (-1.5,-1.5) to (-1.5,-2.6);
        \draw (0.8,-0.5) to (0.8,-2);
        %% TOP F PART
        	\draw[<-,thick,black] (0.3,1.6) to (-0.3,0.6);
	\draw[-,thick,black] (-0.2,1.3) to (0.2,0.8);
        \draw[-,thick,black] (0.2,0.8) to[out=130,in=180] (0.5,0.7);
        \draw[-,thick,black] (0.5,0.7) to[out=0,in=270] (0.8,1.6);
        \draw[-,thick,black] (-0.2,1.3) to[out=130,in=0] (-0.5,1.6);
        \draw[->,thick,black] (-0.5,1.6) to[out=180,in=-270] (-0.8,0.5);
\end{tikzpicture}$}
       }="TL";
  (65,10)*+{  (-1)^\lambda \: \raisebox{0.8mm}{$\begin{tikzpicture}[baseline = 0,scale=0.35]
        \draw[-,thick,black] (0.3,-1.5) to (-0.3,-2.6);
        \draw[-,thick,black] (0.5,-2.6) to[out=0,in=270] (0.8,-1.8);
        \draw[-,thick,black] (0.2,-2.3) to[out=130,in=180] (0.5,-2.6);
        \draw[-,thick,black] (-0.2,-1.8) to (0.2,-2.3);
        \draw[-,thick,black] (-0.2,-1.8) to[out=130,in=0] (-0.5,-1.6);
        \draw[->,thick,black] (-0.5,-1.6) to[out=180,in=-270] (-0.8,-2.6);
        \node at (1.2,0.5) {$\scriptstyle{\lambda}$};
        %% LEFT CROSSING PART
        \draw[-,thick,black] (0.3,-1.5) to (-1.4,-0.5);
        \draw[->,thick,black] (-1.4,-0.5) to (-1.4,0.5);
        \draw[-,thick,black]  (0.2,-0.5) to (-1.5,-1.5);
        \draw[-,thick,black] (-1.5,-1.5) to (-1.5,-2.6);
        \draw[-,thick,black] (0.8,-0.5) to (0.8,-2);
        \draw[->,thick,black] (0.2,-0.5) to (0.2,0.5);
        \draw[-,thick,black] (0.8,-1.8) to (0.8,0.5);
\end{tikzpicture}$}  \: + \: \sum\limits_{n=0}^{ \lambda-1}\sum\limits_{r\geq 0} (-1)^{n + r + \lambda}  \raisebox{-10mm}{$\begin{tikzpicture}[scale=0.85]
	\draw[<-,thick,black] (0.4,-0.1) to[out=90, in=0] (0.1,0.3);
	\draw[-,thick,black] (0.1,0.3) to[out = 180, in = 90] (-0.2,-0.1);
	\node at (-0.15,0.1) {$\bullet$};
	\node at (-0.35,0.1) {$\scriptstyle{n}$};
  \node at (1.4,0) {$\scriptstyle{\lambda}$};
  \draw[->,thick,black] (0.75,-0.1) to (0.75,1.5);
  %% Bubble part
    \draw[->,thick,black,scale=1.2] (1.4,0.5) to[out=90,in=0] (1.2,.7);
  \draw[-,thick,black,scale=1.2] (1.2,0.7) to[out=180,in=90] (1,0.5);
\draw[-,thick,black,scale=1.2] (1,0.5) to[out=-90,in=180] (1.2,0.3);
  \draw[-,thick,black,scale=1.2] (1.2,0.3) to[out=0,in=-90] (1.4,0.5);
   \node[scale=1.2] at (1.7,0.55) {$\color{black}\bullet$};
   \node[scale=1.2] at (2.2,0.55) {$\color{black}\scriptstyle{-n-1}$};
   %%% CUP
   	\node at (1.24,1.25) {$\color{black}\bullet$};
   	\node at (1.05,1.25) {$\scriptstyle{r}$};
	\draw[-,thick,black] (1.8,1.5) to[out=-90, in=0] (1.5,1);
	\draw[->,thick,black] (1.5,1) to[out = 180, in = -90] (1.2,1.5);
  \end{tikzpicture}$} }="TR";
     {\ar@<.6pc> ^-{
        \Big\{ F_\lambda, \: id + \sum\limits_{\substack{n=0, \\ r \geq 0}}^{\lambda-1} (-1)^{n+r} on_2^n, \: {}_\ast (u_{\lambda,0})^- \cdot E_{\lambda} + {}_\ast B'_\lambda   \Big\}}  "TL";"TR"};
     {\ar@<-0.6pc>_-{
        \Big\{  \Gamma_\lambda, \: \Omega \Big\} }    "TL";"TR"};
\endxy
\]
where $B'_\lambda$ is the $3$-cell that reduces the term
\[
\sum\limits_{n = 0}^{\lambda-1} \sum\limits_{r \geq 0} (-1)^{n+r} \sum\limits_{a+b=n-1} (-1)^a \tfishulpfME{}{b}
\]
into $0$, as defined in \ref{subsubsec:additional3cells}. The $3$-cell $\Omega$ is defined as the following composition of ${}_E R$-rewriting steps: when applying $\Gamma_\lambda$ we obtain the polynomial
\[
\raisebox{0.6mm}{$\begin{tikzpicture}[baseline = 0,scale=0.55]
	\draw[<-,thick,black] (0.3,.5) to (-0.3,-.5);
	\draw[-,thick,black] (-0.2,.2) to (0.2,-.3);
        \draw[-,thick,black] (0.2,-.3) to[out=130,in=180] (0.5,-.4);
        \draw[-,thick,black] (0.5,-.4) to[out=0,in=270] (0.8,.5);
        \draw[-,thick,black] (-0.2,.2) to[out=130,in=0] (-0.5,.5);
        \draw[-,thick,black] (-0.5,.5) to[out=180,in=-270] (-0.8,-.5);
      \draw[-,thick,black] (-0.3,-1.5) to (0.3,-2.5);
         \draw[-,thick,black] (-0.5,-2.5) to[out=180,in=-90] (-0.8,-1.5);
         \draw[-,thick,black] (-0.2,-2.2) to[out=230,in=0] (-0.5,-2.5);
         \draw[-,thick,black] (-0.2,-2.2) to (0.2,-1.7);
         \draw[-,thick,black] (0.2,-1.7) to[out=50,in=180] (0.5,-1.6);
          \draw[->,thick,black] (0.5,-1.6) to[out=0,in=90] (0.7,-2.5);
    \node at (1.9,0) {$\scriptstyle{\lambda}$};
    %% RIGHT CROSSING PART
        \draw[-,thick,black] (-0.3,-1.5) to (1.5,-0.5);
        \draw[->,thick,black] (1.5,-0.5) to (1.5,0.5);
        \draw[-,thick,black] (-0.3,-0.5) to (1.5,-1.5);
        \draw[-,thick,black] (1.5,-1.5) to (1.5,-2.5);
        \draw[-,thick,black]  (-0.8,-0.5) to (-0.8,-1.5);
            %% TOP F PART
        	\draw[<-,thick,black] (1.9,1.6) to (1.45,0.6) to (1.45,0.4);
	\draw[-,thick,black] (1.4,1.3) to (1.8,0.8);
        \draw[-,thick,black] (1.8,0.8) to[out=130,in=180] (2.3,0.7);
        \draw[-,thick,black] (2.1,0.7) to[out=0,in=270] (2.4,1.6);
        \draw[-,thick,black] (1.4,1.3) to[out=130,in=0] (1.1,1.6);
        \draw[->,thick,black] (1.1,1.6) to[out=180,in=-270] (0.8,0.5);
\end{tikzpicture}$}
\raisebox{-5mm}{$- \:
\sum\limits_{r,s,t \geq 0}
(-1)^{\lambda+r+s+1}
\begin{tikzpicture}[baseline = 0]
	\draw[-,thick,black] (0.3,0.7) to[out=-90, in=0] (0,0.3);
	\draw[-,thick,black] (0,0.3) to[out = 180, in = -90] (-0.3,0.7);
    \node at (1.4,-0.32) {$\scriptstyle{\lambda}$};
  \draw[->,thick,black] (0.2,0) to[out=90,in=0] (0,0.2);
  \draw[-,thick,black] (0,0.2) to[out=180,in=90] (-.2,0);
\draw[-,thick,black] (-.2,0) to[out=-90,in=180] (0,-0.2);
  \draw[-,thick,black] (0,-0.2) to[out=0,in=-90] (0.2,0);
   \node at (0.2,0) {$\color{black}\bullet$};
   \node at (.6,0) {$\color{black}\scriptstyle{\substack{-r-s\\-t-3}}$};
   \node at (-0.23,0.43) {$\color{black}\bullet$};
   \node at (-0.43,0.43) {$\color{black}\scriptstyle{r}$};
	\draw[<-,thick,black] (0.3,-.7) to[out=90, in=0] (0,-0.3);
	\draw[-,thick,black] (0,-0.3) to[out = 180, in = 90] (-0.3,-.7);
   \node at (-0.25,-0.5) {$\color{black}\bullet$};
   \node at (-.4,-.5) {$\color{black}\scriptstyle{s}$};
	\draw[-,thick,black] (.98,-0.7) to (.98,0.7);
   \node at (0.98,0.53) {$\color{black}\bullet$};
   \node at (1.15,0.53) {$\color{black}\scriptstyle{t}$};
   \draw[->,thick,black] (-0.3,0.7) to (-0.3,1.6);
          %% TOP F PART
        	\draw[<-,thick,black] (1.4,1.6) to (0.95,0.6);
	\draw[-,thick,black] (0.9,1.3) to (1.3,0.8);
        \draw[-,thick,black] (1.3,0.8) to[out=130,in=180] (1.6,0.7);
        \draw[-,thick,black] (1.6,0.7) to[out=0,in=270] (1.9,1.6);
        \draw[-,thick,black] (0.9,1.3) to[out=130,in=0] (0.6,1.6);
        \draw[->,thick,black] (0.6,1.6) to[out=180,in=-270] (0.3,0.5);
\end{tikzpicture} - \sum\limits_{r,s,t \geq 0}
(-1)^{\lambda+r+s+t}
\begin{tikzpicture}[baseline = 0]
	\draw[<-,thick,black] (0.3,0.7) to[out=-90, in=0] (0,0.3);
	\draw[-,thick,black] (0,0.3) to[out = 180, in = -90] (-0.3,0.7);
    \node at (.6,-0.32) {$\scriptstyle{\lambda}$};
  \draw[-,thick,black] (0.2,0) to[out=90,in=0] (0,0.2);
  \draw[<-,thick,black] (0,0.2) to[out=180,in=90] (-.2,0);
\draw[-,thick,black] (-.2,0) to[out=-90,in=180] (0,-0.2);
  \draw[-,thick,black] (0,-0.2) to[out=0,in=-90] (0.2,0);
   \node at (-0.2,0) {$\color{black}\bullet$};
   \node at (-0.6,0) {$\color{black}\scriptstyle{\substack{-r-s\\-t-3}}$};
   \node at (0.23,0.43) {$\color{black}\bullet$};
   \node at (0.43,0.43) {$\color{black}\scriptstyle{s}$};
	\draw[-,thick,black] (0.3,-.7) to[out=90, in=0] (0,-0.3);
	\draw[->,thick,black] (0,-0.3) to[out = 180, in = 90] (-0.3,-.7);
   \node at (0.25,-0.5) {$\color{black}\bullet$};
   \node at (.4,-.5) {$\color{black}\scriptstyle{r}$};
	\draw[->,thick,black] (-.98,-0.7) to (-.98,1.6);
   \node at (-0.98,0.58) {$\color{black}\bullet$};
   \node at (-1.15,0.58) {$\color{black}\scriptstyle{t}$};
            %% TOP F PART
        	\draw[<-,thick,black] (0.8,1.6) to (0.3,0.6);
	\draw[-,thick,black] (0.3,1.3) to (0.7,0.8);
        \draw[-,thick,black] (0.7,0.8) to[out=130,in=180] (1,0.7);
        \draw[-,thick,black] (1,0.7) to[out=0,in=270] (1.3,1.6);
        \draw[-,thick,black] (0.3,1.3) to[out=130,in=0] (0,1.6);
        \draw[-,thick,black] (0,1.6) to[out=180,in=-270] (-0.3,0.5);
\end{tikzpicture}$} \]
The third term reduces to $0$ using the $3$-cell $D'_\lambda$ into $0$ since $s < \lambda$, the first term reduces using $\{({}_*(u'_{\lambda,0})\star_2 yb \star_2 (u'_{\lambda,0})^-) \cdot {}_\ast F_{\lambda +2} \}$ into
\[  \raisebox{-5mm}{$(-1)^\lambda$} \:
 \raisebox{0.6mm}{$\begin{tikzpicture}[baseline = 0,scale=0.45]
        \draw[-,thick,black] (0.3,-1.5) to (-0.3,-2.6);
        \draw[-,thick,black] (0.5,-2.6) to[out=0,in=270] (0.8,-1.8);
        \draw[-,thick,black] (0.2,-2.3) to[out=130,in=180] (0.5,-2.6);
        \draw[-,thick,black] (-0.2,-1.8) to (0.2,-2.3);
        \draw[-,thick,black] (-0.2,-1.8) to[out=130,in=0] (-0.5,-1.6);
        \draw[->,thick,black] (-0.5,-1.6) to[out=180,in=-270] (-0.8,-2.6);
        \node at (1.2,0.5) {$\scriptstyle{\lambda}$};
        %% LEFT CROSSING PART
        \draw[-,thick,black] (0.3,-1.5) to (-1.4,-0.5);
        \draw[->,thick,black] (-1.4,-0.5) to (-1.4,0.5);
        \draw[-,thick,black]  (0.2,-0.5) to (-1.5,-1.5);
        \draw[-,thick,black] (-1.5,-1.5) to (-1.5,-2.6);
        \draw[-,thick,black] (0.8,-0.5) to (0.8,-2);
        \draw[->,thick,black] (0.2,-0.5) to (0.2,0.5);
        \draw[-,thick,black] (0.8,-1.8) to (0.8,0.5);
\end{tikzpicture}$} \]
plus an extra term that one might check is cancelled by the term obtained from the second summand when using super isotopies and making the $r$ dots move to the bottom of the crossing, so that it only remains the terms where the dots break the crossing, giving the summand
\[
\sum\limits_{a,b,s,t \geq 0} (-1)^{s+\lambda+b + (\lambda+a+b)(a+t)}
\raisebox{-10mm}{$\begin{tikzpicture}[scale=1]
	\draw[<-,thick,black] (0.4,-0.1) to[out=90, in=0] (0.1,0.3);
	\draw[-,thick,black] (0.1,0.3) to[out = 180, in = 90] (-0.2,-0.1);
	\node at (-0.15,0.1) {$\bullet$};
	\node at (-0.35,0.1) {$\scriptstyle{a}$};
  \node at (1.7,0) {$\scriptstyle{\lambda}$};
  \draw[->,thick,black] (1.5,-0.1) to (1.5,1.5);
  \node at (1.5,0.9) {$\bullet$};
  \node at (1.85,0.9) {$\scriptstyle{t+a}$};
  %% Bubble part
    \draw[->,thick,black,scale=1.2] (0.2,0.6) to[out=90,in=0] (0,.8);
  \draw[-,thick,black,scale=1.2] (0,0.8) to[out=180,in=90] (-0.2,0.6);
\draw[-,thick,black,scale=1.2] (-0.2,0.6) to[out=-90,in=180] (0,0.4);
  \draw[-,thick,black,scale=1.2] (0,0.4) to[out=0,in=-90] (0.2,0.6);
   \node at (0.25,0.65) {$\color{black}\bullet$};
   \node at (0.8,0.65) {$\color{black}\scriptstyle{\substack{-a-b \\ -s-t-4}}$};
   %%% CUP
   	\node at (1.94,1.35) {$\color{black}\bullet$};
   	\node at (1.75,1.35) {$\scriptstyle{b}$};
	\draw[-,thick,black] (2.5,1.6) to[out=-90, in=0] (2.2,1.1);
	\draw[->,thick,black] (2.2,1.1) to[out = 180, in = -90] (1.9,1.6);
  \end{tikzpicture}$}
\]
and one checks that this reduces using bubble slide $3$-cells $s_{\l,n}'$ into the second term of the final result. Note that there also is a critical branching between $\Gamma_\l$ and $E_\l$ given by attaching to the source of $\Gamma_\l$ a rightward crossing on bottom on the rightmost two strands. This one would be proved confluent in a similar manner.

\medskip
\noindent  \textbf{Critical branching $(ig_{2n}, s_{2n,\l}^+)$}
\[
\xymatrix{
\raisebox{-2mm}{$\posbubdfffsl{}{2n+ \ast}$} \identu{} \ar@<.6pc>[rrrr]^-{\{s_{\l,2n}^+, \sum\limits_{r=0}^{n-1}ig_{2n-2r,\l}\}} \ar@<-.6pc>[rrrr]_-{\{ig_{2n,\l+2}, s_{\l,2n-2\ell}^+, s_{\l,2\ell}^-, ig_{2n-2r,\l}\}} & & & & (2n+1)
\hackcenter{
\begin{tikzpicture}[scale=0.8]
    \draw [thick,black, ->] (0,-1) to (0,1);
    \node at (0,-.5) {$\bullet$};
    \node at (-.25,-.5) {$\scriptstyle{2n}$};
\end{tikzpicture}}
-
\sum\limits_{r=0}^{n-1} \sum\limits_{\ell=1}^{n-r}
(2r+1)
\hackcenter{
\begin{tikzpicture}[scale=0.8]
    \draw [thick,black, ->] (0,-1) to (0,1);
    \node at (0,-.5) {$\bullet$};
    \node at (-.25,-.5) {$\scriptstyle{2r}$};
\end{tikzpicture}}
\raisebox{2.5mm}{\hbox{$\posbubdfffsl{}{2n-2r-2\ell+\ast \qquad \;}$}}
\raisebox{-2.5mm}{\hbox{$\negbubdfffsl{}{\; 2\ell+ \ast }$}}
}
\]
where the last $ig$ $3$-cell in the bottom branch is only applied to terms without a counter-clockwise bubble of positive degree. Note that there is a similar branching between $ig_{2n}$ and $r_{2n,\l}^{-}$ given by changing the upward strand to the right of the bubble in the source of the last branching to a downward strand. This would be proved confluent in a similar manner.

\medskip

\end{document}